\documentclass{article}

\usepackage{latexsym}
\usepackage{authblk}
\usepackage{macro-phd}  
\usepackage[babel]{csquotes}

\usepackage{indentfirst}
\usepackage[english]{babel}
\usepackage{cancel}

\usepackage{lipsum}
\usepackage{array}
\usepackage{color} 
\usepackage{tabularx}
\usepackage{graphicx} 
\usepackage{amsfonts}	
\usepackage{amsmath}
\usepackage{amssymb}
\usepackage{moreverb}
\usepackage{dsfont}
\usepackage{tipa}
\usepackage{upgreek}
\usepackage{bm}
\usepackage{enumitem} 
\usepackage{multirow}
\usepackage{soul}
\usepackage{ textcomp }
\usepackage[svgnames]{xcolor}
\usepackage{mathtools}
\usepackage{setspace}
\usepackage{algorithm}
\usepackage{algorithmic}
\usepackage{caption}
\usepackage{subcaption}
\usepackage{booktabs}
\usepackage{relsize} 
\allowdisplaybreaks 
\newcommand\BibTeX{{\rmfamily B\kern-.05em \textsc{i\kern-.025em b}\kern-.08em
T\kern-.1667em\lower.7ex\hbox{E}\kern-.125emX}}

\usepackage{hyperref} 
\hypersetup{
    colorlinks=true,                          
    linkcolor=blue,
    citecolor=DarkRed,
    urlcolor=blue  }

\usepackage[backend=biber,giveninits=true,sorting=nyt,url=false,doi=true,eprint=false,isbn=false,
backref,backrefstyle=none,maxbibnames=99]{biblatex}
\DefineBibliographyStrings{english}{%
	backrefpage = {Cited on p\adddot},%
	backrefpages = {Cited on pp\adddot}%
}
\bibliography{references}

\renewbibmacro{in:}{%
	\ifboolexpr{%
		test {\ifentrytype{article}}%
	}{}{\printtext{\bibstring{in}\intitlepunct}}%
}

\newcommand{\pd}{\mathcal{\partial}}
\newcommand{\xx}{{\boldsymbol{x}}}
\newcommand{\XX}{{\boldsymbol{X}}}

\newcommand{\bfV}{\mathbf{V}}
\renewcommand{\S}{\mathbf{S}}

\newcommand{\F}{\mathbf{F}}

\newcommand{\B}{\mathbf{B}}

\newcommand{\rmd}{\mathrm{d}}  

\newcommand{\A}{\AAA}

\newcommand{\halb}{\frac{1}{2}}

\newcommand{\bdm}{\begin{displaymath}}
\newcommand{\edm}{\end{displaymath}}

\newcommand{\bea}{\begin{eqnarray} }
\newcommand{\eea}{\end{eqnarray} }

\newcommand{\AAA}{{\boldsymbol{A}}}

\newcommand{\GG}{{\boldsymbol{G}}}

\newcommand{\vv}{{\boldsymbol{v}}}
\newcommand{\uu}{{\boldsymbol{u}}}

\newcommand{\FF}{{\boldsymbol{F}}}

\newcommand{\QQ}{{\mathbf{Q}}}

\newcommand{\PP}{{\mathbf{P}}}

\newcommand{\Id}{{\boldsymbol{I}}}

\newcommand{\crho}{\varrho}
\newcommand{\UU}{{\boldsymbol{U}}}
\newcommand{\VV}{{\boldsymbol{V}}}
\newcommand{\ww}{{\boldsymbol{w}}}

\newcommand{\J}{{\boldsymbol{J}}}
\newcommand{\rhoE}{\mathcal{E}}
\newcommand{\ce}{{\varepsilon}}
\newcommand{\ps}{p}

\newcommand\rhoo{\rho \hspace{-0.015em}\mathrm{\scriptsize{o}}}
\newcommand\po{p \hspace{-0.015em}\mathrm{\scriptsize{o}}}

\newcommand\Cs{\mathrm{Cs}}
\newcommand\Ch{\mathrm{Ch}}
\newcommand\Cv{\mathrm{Cv}}
\newcommand\Co{\mathrm{Co}}

\newcommand{\Nph}{\mathrm{N}}

\let\svthefootnote\thefootnote
\newcommand\freefootnote[1]{%
  \let\thefootnote\relax%
  \footnotetext{#1}%
  \let\thefootnote\svthefootnote%
}

\newfont{\numerikEleven}{ecrm1000}
\newfont{\numerikTen}{cmss10}
\newfont{\numerikNine}{cmss9}
\newfont{\numerikEight}{cmss8}

\title{\vspace{-2.0cm}\bf A unified SHTC multiphase model of continuum mechanics}
\author[1]{Davide Ferrari}

\author[1]{Ilya Peshkov}

\author[2]{Evgeniy Romenski}

\author[1]{Michael Dumbser}

\affil[1]{\small Laboratory of Applied Mathematics, DICAM, University of Trento, via Mesiano 77, 38123 Trento, Italy}
\affil[2]{Sobolev Institute of Mathematics, 4 Acad. Koptyug avenue, Novosibirsk, Russia} 

\begin{document} 

\maketitle


\hfill%
\begin{minipage}{10cm}
\small The paper is dedicated to the memory of Prof. S.K. Godunov, 
who made an invaluable contribution to the development of many
areas of applied and computational mathematics. The content of
this article is largely based on his ideas, which are and will
be the basis of research for future generations of mathematicians.
\end{minipage}
\vspace{0.1cm}

\begin{abstract} \noindent In this paper, we present a unified
nonequilibrium model of continuum mechanics for compressible multiphase flows.
The model, which is formulated within the framework of Symmetric Hyperbolic
Thermodynamically Compatible (SHTC) equations, can describe arbitrary number of
phases that can be heat conducting \textit{inviscid} and \textit{viscous
fluids}, as well as \textit{elastoplastic solids}. The phases are allowed to
have different velocities, pressures, temperatures, and shear stresses, while
the material interfaces are treated as diffuse interfaces with the volume
fraction playing the role of the interface field. To relate our model to other
multiphase approaches, we reformulate the SHTC governing equations in terms of
the phase state parameters and put them in the  form of Baer-Nunziato-type models.
It is the Baer-Nunziato form of the SHTC equations which is then solved
numerically using a robust second-order path-conservative MUSCL-Hancock finite
volume method on Cartesian meshes. Due to the fact that the obtained governing
equations are very challenging we restrict our numerical examples to a
simplified version of the model, focusing on the isentropic limit for
three-phase mixtures. To address the stiffness properties of the relaxation source
terms present in the model, the implemented scheme incorporates a semi-analytical time integration
method specifically designed for the non-linear stiff source terms governing the
strain relaxation. The validation process involves a wide range of benchmarks
and several applications to compressible multiphase problems. Notably, results are presented for
multiphase flows in all the relaxation limit cases of the model, including 
inviscid and viscous Newtonian fluids, as well as non-linear hyperelastic and
elastoplastic solids. In all cases, the numerical results demonstrate good
agreement with established models, including the Euler or Navier-Stokes
equations for fluids and the classical hypo-elastic model with plasticity for solids. 
Importantly, however, this approach achieves these results within a unified multiphase framework of continuum mechanics.  

\end{abstract}

\section{Introduction} 
\label{sec.introduction}

\freefootnote{Submitted to: \textit{Jounral of Computational Physics}}
\freefootnote{
    \href{mailto:davide.ferrari@unitn.it}{davide.ferrari@unitn.it}, 
    \href{mailto:ilya.peshkov@unitn.it}{ilya.peshkov@unitn.it},
    \href{mailto:evrom@math.nsc.ru}{evrom@math.nsc.ru},
    \href{mailto:michael.dumbser@unitn.it}{michael.dumbser@unitn.it}
}

In this paper, we further develop the idea of a unified approach to continuum
mechanics, which is capable of describing the behavior of all states of matter
in a single system of governing equations. Here, we deal with multiphase
problems that may include an arbitrary number of phases, such as gases, liquids,
and solids and in arbitrary combinations. The approach we chose to tackle this
problem relies on the continuum mixture theory, which means that all the phases
are present in every material volume at all the times \cite{Truesdell1}. Despite
being popular in the fluid mechanics community for many decades now, the mixture
approach imposes a big challenge if one tries to generalize it to problems
involving the interaction of solids and fluids. The main difficulty is that
these two states of matter are generally studied by two distinct communities,
using different approaches. And it remained unclear, at least until recently,
how to combine these approaches (mathematical models) into a single unified and
thremodynamically compatible continuous mixture model. 

For example, the Eulerian description of the continuum is the most commonly used
to describe the motion of fluids, as they usually undergo large deformations
that are not suitable for a Lagrangian description. Instead, in solid mechanics,
the Lagrangian description is preferred, which allows very accurate description
of the solid body boundaries. Moreover, the fluid and solid mechanics models
have different number of degrees of freedom (state variables), e.g. the
classical fluid mechanics models do not need to evolve a deformation measure to
compute the stress tensor, while the solid mechanics models do. Such
incompatibilites in mathematical description have been a major obstacle to the
development of a general unified continuum mechanics model for multiphase
problems involving fluids and solids. Yet, in some limited cases the mixture
theory has been successfully applied to the description of the interaction
between fluids and solids, e.g. see
\cite{FavrieGavrilyukSaurel,NdanouFavrieGavrilyuk,Nikodemou2021,Nikodemou2023} and references therein. We
also remark that a conventional approach to the description of multiphase flows
involving solid-fluid interactions, such as in the Fluid-Structure-Interaction
(FSI) problems, is to use a combination of different models, one for each of the
phases involved, with additional coupling rules governing the interaction
between the different media, and even possibly different numerical schemes and
computational meshes for each phase. 

In recent years, we have been developing a unified Eulerian model of continuum
mechanics \cite{HPR2016,DUMBSER2016824,frontiers,HTC2022,ABGRALL2023127629} that allows
to eliminate the mentioned differences between the classical descriptions of
fluids and solids. According to this model, a fluid is treated as a special case
of an inelastic solid with a severe shear stress relaxation. Therefore, it
represents a very solid foundation for building a unified multiphase continuum
theory. On the other hand, the unified model belongs to a wider class of
equations that we call the Symmetric Hyperbolic Thermodynamically Compatible
(SHTC) equations
\cite{God1961,God1972MHD,GodRom1996b,Rom1998,Rom2001,GodRom2003,SHTC-GENERIC-CMAT}.
In particular, the SHTC class includes a multi-fluid model developed by Romenski
and co-authors in \cite{Romenski2004,Romenski2007,RomDrikToro2010} (the first
idea was published in \cite{GodRom1998,Rom2001}). Therefore, the main motivation
for this paper is to combine the two approaches and to develop a unified model
of continuum mechanics for multiphase problems with arbitrary composition of the
phases that may include inviscid and viscous fluids as well as elastic and
inelastic solids. Likewise, to demonstrate the practicality and advantages of
such a unified approach for the numerical simulation of fluid-fluid and
fluid-structure interaction problems.

It should be noted that despite the long history, the continuous mixture theory
still does not have a universally accepted mathematical formulation. It is a
hopeless task to try to review all the contributions to the development of the
mixture theory, however it worth mentioning the main directions that were taken
in order to obtain the balance equations for mixtures to position our approach
with respect to the others. One may roughly distinguish three main approaches:
the \textit{rational mechanics} approach, \textit{averaging} techniques, and
\textit{first-principle}-type approaches. However, let us immediately note that
in all these theories, the basic view on the mixture is the same, i.e. the
mixture is treated as a system of interpenetrating and interacting continua each
of which is described by their own state variables, equations of state, and the
evolution of constituents is governed by the phase balance equations similar to
that of a single-phase continua but having extra terms describing the
interactions between the phases. The main difference between these approaches is
the way the balance equations are derived, or more precisely, the way the
interaction terms are defined.

The \textit{rational mechanics} approach is a general approach based on some
principles mainly put forward by Truesdell, Coleman, and Noll
\cite{Truesdell1,TruesdellNoll2004,Coleman1963}, e.g. second law of
thermodynamics, objectivity, etc., e.g. see \cite{Muller1968}. However, the
generality of the approach does not provide specific instructions on how to
obtain the balance equations for mixtures in a closed form, i.e. the
interaction terms are in general not defined
\cite{Muller1968,Gurtin1971,Muller1984,NunzWalsh,StewartWendroff1984,Passman1984,Nigmatulin,MullerRuggeri1998}. This
explains why closed models can only be found in the literature for
application-specific two-phase flows, see e.g. deflagration-to-detonation
transition (DDT) in gas-permeable, reactive granular materials
\cite{baernunziato,Embid,Bdzil,POWERS1990264} and spray modeling \cite{Crowe}. 

The second approach for the derivation of balance laws for mixtures involves
overcoming the discrete nature of the phases that constitute a multiphase
mixture through the use of \textit{averaging} techniques. Within {averaged}
approaches, see e.g. \cite{Drew1971AveragedEF, Nigmatulin, Ishii}, the balance
equations of each constituent are derived by applying an averaging space and
time operation to the equations of motion for different phases separated by
interfaces across which the densities, velocities, etc., may jump. The main
difficulty of this way of deriving the balance equations is that it is almost
impossible to apply in a genuinely multiphysics setting, e.g. mixtures of more
than two constituents, different rheology of phases, phase transition,
electromagnetic forces (e.g. multispecies plasma), etc.

The third class of mixture theories consists of models whose governing equations
can be derived from a first-principle-type approach, e.g. a variational
principle (Hamilton's stationary action principle). Indeed, the success of
theoretical physics in the twentieth century taught us an important lesson that the most
successful theories should admit a variational formulation. However, since we
are interested in a macroscopic description of mixtures, such a description is
inevitably dissipative, and one cannot expect to derive the full mixture model
from a variational principle, but only its non-dissipative, or reversible, part.
Therefore, it is crucially important to formulate continuum models in such a way
that the reversible and irreversible parts can be clearly separated. This is one
of the main ideas behind the theory of Symmetric Hyperbolic Thermodynamically
Compatible (SHTC) equations \cite{SHTC-GENERIC-CMAT} that is used in this paper
to obtain the governing equations for multiphase and multilateral flows. In such
a theory, the dissipative processes can be formulated only via algebraic (no
partial derivatives) relaxation type terms, while the reversible part of SHTC
models is presented by the differential part of the time evolution equations.
And the latter can be derived from a variational principle. We believe that it
is important to have a variational formulation for such complex systems as
mixtures, because it allows to couple (define interactions) all parts of the
system in a consistent way. It is even more important to have such a
formulation in the multiphysics context, where the interactions between the
phases are not only mechanical, but also thermal, chemical, electromagnetic,
gravitational, etc. It is therefore the main objective of this paper to present
a unified model of continuum mechanics for multiphase problems and to test its practicality in numerical simulations of multiphase flows.

Despite being so fundamental in theoretical physics and other branches of
continuum mechanics, the variational principle is rarely used to derive
governing equations for mixtures. To the best of our knowledge, the few existing variational
formulations for mixtures of two fluids were proposed in
\cite{Gavrilyuk1999,GouinRugerri2003,SHTC-GENERIC-CMAT}.

We note that there are other first-principle type approaches such as the
Hamiltonian formulation of continuum mechanics and non-equilibrium
thermodynamics known as GENERIC (General Equations for Reversible-Irreversible
Coupling) and put forward by Grmela and \"Ottinger in \cite{GrmelaOttingerI,GrmelaOttingerII}, see also \cite{PKGBook2018,SHTC-GENERIC-CMAT}.
Despite the SHTC multiphase models can be also cast into the GENERIC
formulation, the GENERIC approach also allows to derive slightly different
governing equations for mixtures. A comparison of the two approaches was the
subject of \cite{Sykora2024}.

From the numerical view point, compressible multiphase problems also pose a great challenge
due to the complexity of the governing equations, various time and length
scales, and the presence of interfaces which require special treatment in many
cases. One could divide mesh-based methods into three different families for
dealing with multimaterial and multiphase problems: (i) Lagrangian and
Arbitrary-Lagrangian-Eulerian (ALE) methods on moving meshes, where free
surfaces of the fluid and fluid-fluid or fluid-solid interfaces are accurately
resolved by the moving computational grid, see e.g.
\cite{MaireCyl2,ShashkovMultiMat3,Kucharik2014,boscheri2015direct,ReALE2015,arepo};
(ii) Eulerian methods on fixed meshes with explicit interface
reconstruction, such as the volume of fluid method (VOF) originally developed by
Hirt and Nichols \cite{HirtNichols}, also see the works of Popinet \cite{Popinet2018} and Menshov et. al \cite{menshov1,menshov2,menshov3} which is also generalized to arbitrary number
of immiscible compressible fluids; (iii) Eulerian diffuse interface methods on
fixed grids, where the presence of each material is represented by a scalar
function, see e.g.
\cite{Saurel2018,DIM2D,GPR11,NdanouFavrieGavrilyuk,FavrieGavrilyukSaurel,gaburro2018diffuse,kemm2020simple}
and references therein. By its construction, the diffuse interface method is the
most natural method to represent interfaces within a continuous mixture model,
and it is the method we use in this paper.

In the compressible multiphase community, the Baer-Nunziato (BN) model
\cite{baernunziato} is one of the most popular mathematical models for
describing two-phase flows, and there are many works that solve this model
numerically. However, only a very limited number of publications exist on the
mathematical and computational issues of BN models for multiphase flows
describing \textit{more} than two phases, see e.g. \cite{Herard1,Herard2,Pelanti3phase}. It is therefore very attractive to
demonstrate how the BN-type models can be extended using the SHTC theory to
describe nonequilibrium compressible multiphase flows with more than two phases, that could be also solids,
and to show how such models can be solved numerically.

Therefore, in this paper we put the SHTC multiphase formulation into a BN-type
form. The resulting BN-type model yields a large system of nonlinear Partial Differential
Equation (PDE), which includes \textit{highly nonlinear stiff algebraic
source terms} as well as \textit{non-conservative} products. The different
complexities presented by the PDE system are addressed numerically following a
splitting approach, when the homogeneous part of the PDE system is discretized
with the aid of a robust second-order explicit MUSCL-Hancock finite volume
method on Cartesian meshes\cite{torobook}, combined with a path-conservative
technique of Castro and Pares \cite{Pares2006, Pimentel} for the treatment of
non-conservative products. Because the time scales associated with the
relaxation processes are much shorter than those given by the stability
condition of the explicit scheme, two different implicit methods are employed to
properly treat the relaxation terms. Namely, for the stiff but linear sources
related to velocity relaxation, a time integrator based on backward Euler is
employed, while the semi-analytical time integration method of Chiocchetti,
introduced in \cite{ExpInt1} for fracture modeling, and further developed in
\cite{ExpInt2,Chiocchetti2023,BOSCHERI2022110852}, is adopted for the nonlinear
stiff source terms governing the relaxation of shear stresses in viscous fluids
and plastic solids.

The paper is structured as follows: in Sec.\,\ref{intro.SHTC} we briefly remind
the main features of the SHTC framework, in Sec.\,\ref{sec.SHTC.model}, we
discuss the SHTC formulation for multiphase flows, while the variational
derivation is presented in \ref{sec.variation}. Furthermore, in
Sec.\,\ref{BN.SHTC.form} a BN-form of the SHTC model is presented. Sec.
\ref{SHTC.BN.num.meth} is devoted to the development of a robust numerical
scheme capable of addressing the various difficulties inherent in a unified
theory of compressible multiphase fluid and solids mechanics. Sec.
\ref{sec.num.res} provides and discusses an extensive collection of numerical
experiments, limited to three-phase flows in this paper, with the aim of
validating the numerical methods developed in this work, as well as providing
some rather unique results related to the behavior of multiphase flows for more
than two phases, described through the unified SHTC multiphase model of
continuum mechanics. Finally, in Sec. \ref{sec.conclusion} the main
achievements of this work are listed and we discuss future research directions regarding
numerical algorithms and modelling perspectives within the unified
SHTC multiphase model of continuum mechanics.

\section{SHTC class of equations}\label{intro.SHTC}

Perhaps, the main challenge in formulating continuum mechanics models for
complex phenomena in general and for multiphase flows in particular is
associated with the formulation of a closed model that satisfies \textit{a
priori} important physical and mathematical properties such as the principles of
invariance (Galilean or Lorentz invariance), conservation principle, the
principle of causality, the laws of thermodynamics and the possibility of having a
well-posed initial value problem (IVP). Methods for constructing equations
belonging to the SHTC class allow to build continuum mechanics models that
satisfy all these properties. Another general formulation of continuum mechanics
that is pursuing the same goals is the Rational Extended Thermodynamics by
M\"uller and Ruggeri \cite{MullerRuggeri1998}.

The origin of the SHTC formulation of continuum mechanics can be attributed to
the work of Godunov \cite{God1961}, who considered \enquote{an interesting
class} of nonlinear \textit{overdetermined} conservation laws. In this seminal
work, Godunov demonstrated an intricate connection between the
\textit{well-posedness} of the IVP for a nonlinear system of conservation laws
and principles of thermodynamics. The well-posedness was shown via putting the
system into a symmetric hyperbolic quasi-linear form,
which was a generalization of Friedrichs' linear symmetric hyperbolic systems
\cite{Friedrichs1958} to the nonlinear case. That is why we use the name Symmetric
Hyperbolic Thermodynamically Compatible systems. However, a very limited number
of continuum mechanics models admit a fully conservative formulation (e.g. Euler
equations, shallow water equations), which forced Godunov to extend his
observations to non-conservative systems such as the ideal magnetohydrodynamics
equations \cite{God1972MHD}. Later, in a series of works 
\cite{GodRom1995,GodRom1996b,Rom1998,Rom2001}, Godunov and Romenski extended
this idea further and showed that many continuum mechanics models (nonlinear
elasticity, binary mixture equations, electrodynamics equations, superfluidity
equations, etc.) can be put into an SHTC form. Later, it was also demonstrated
that all the SHTC equations admit a variational formulation  and can be also
viewed as a Hamiltonian formulation of continuum mechanics within the GENERIC
framework \cite{SHTC-GENERIC-CMAT}.

As discussed in \cite{GodRom1996b,SHTC-GENERIC-CMAT}, the SHTC equations have a
very peculiar structure, and recently, various numerical schemes were developed
to mimic the SHTC structure at the discrete level, see e.g. for a curl-free
discretization in \cite{SIGPR,Chiocchetti2023} and new class of
thermodynamically compatible schemes in
\cite{Busto2021b,HTC2022,Busto2022a,ABGRALL2023127629,Thomann2023}.

In the context of this paper, let us also note that it is possible to formulate
a compressible multiphase flow model within the SHTC framework. The seminal ideas
about an SHTC theory of mixtures was proposed by Romenski in \cite{Rom1998,Rom2001}
for the case of two fluids. It was then further developed in a series of works
\cite{Romenski2004,Romenski2007,RomDrikToro2010} and it was generalized to the
case of arbitrary number of constituents in \cite{Romenski2016}. The SHTC
mixture model was already discretized with various numerical schemes, e.g. a
MUSCL-Hancock scheme \cite{RomDrikToro2010,Romenski2016}, semi-implicit all-Mach
number schemes were developed in \cite{Thomann2022,Thomann2024}, high-order
discontinuous Galerkin and finite-volume schemes in the ADER framework were
developed in \cite{Rio-Martin2023} and a thermodynamically compatible
discretization was proposed in \cite{Thomann2023}. An exact Riemann solver for the barotropic two-fluid SHTC equations was introduced in \cite{Thein2022}.

An important step in demonstrating the potential of the SHTC theory was done in
\cite{HPR2016,DUMBSER2016824}, where the unified model of continuum fluid and
solid mechanics within the SHTC framework was proposed, see also its extensions
to electrodynamics of moving media \cite{DUMBSER2017298}, non-Newtonian fluids
\cite{GPR7,nonNewtonian2021}, poroelasticity \cite{Poroelast2020,Poroelasticity2}, relativistic continuum mechanics \cite{PTRSA2020},
flows with surface tension \cite{Chiocchetti2023,SHTC_surf2023}. The
mathematical model and its many variants have been referred to differently in
these various contexts, from Hyperbolic-Peshkov-Romenski (HPR) or
Godunov-Peshkov-Romenski (GPR) in \cite{DUMBSER2016824,
DUMBSER2017298,Jackson2017} to unified model of continuum mechanics (UMCM) in
\cite{BOSCHERI2022110852}. In this work, we adopt the generic terminology
Unified Continuum Mechanics Model, or when appropriate for brevity UMCM or GPR.
\section{Mixture characteristics}
\subsection{Composition characteristics}\label{sec.notaionts.mix}
We consider a mixture of $ \mathrm{N} $ constituents which are labeled by Latin indices $ 
a,b,\ldots=1,\ldots,\Nph $. Denoting by $M$ and $V$ the total mass 
and volume of the infinitesimal element of the mixture, respectively, we can write
\begin{equation}\label{eqn.MV}
	M = \sum_{a=1}^{\Nph} m_a,
	\qquad
	V = \sum_{a=1}^{\Nph} \nu_a,
\end{equation}
where $ m_a $ is the mass and $ \nu_a $ is the volume of the $a $-th constituent in the mixture 
control volume $V$. The mixture mass density is then defined as
\begin{equation}\label{eqn.rho.mix}
	\rho = \frac{M}{V} =  \frac{m_{1}+m_{2}+\ldots+m_{\Nph}}{V} =  \sum_{a=1}^{\Nph} \crho_a,
\end{equation}
where
\begin{equation}\label{eqn.truerho}
	\crho_a := \frac{m_a}{V}
\end{equation} 
denotes the density of the $a$-th phase inside the control volume $ V $.

 To characterise the volume and mass content of the $a$-th constituent inside the mixture control volume 
$V$, it is also convenient to introduce two non-dimensional scalars: the volume fraction
\begin{equation}\label{eqn.vf}
	\alpha_a:= \frac{\nu_a}{V},
	\qquad\qquad
	\sum_{a=1}^{\Nph}\alpha_a = 1,
\end{equation}
and the mass fraction 
\begin{equation}\label{eqn.mf}
	c_a:= \frac{m_a}{M}= \frac{\crho_a}{\rho},
	\qquad
	\sum_{a=1}^{\Nph}c_a = 1.
\end{equation}
Although $\crho_a$ represents the true mass density of the $a$-th constituent inside the control 
volume $V$, the equations of state of the constituents are usually given in the single-phase 
context, i.e. as if the phase $a$ would occupy the entire volume $V$. 
Therefore, to use the standard single-phase equations of state, we shall also need the mass density of the $a$-th phase not of the 
entire mixture control volume $V$, but of the partial volume $\nu_a$, i.e.
\begin{equation}
	\rho_a = \frac{m_a}{\nu_a} = \frac{m_a V}{\nu_a V} = \frac{\crho_a}{\alpha_a}.
\end{equation}
In other words, for phase $a$, its mass density $\rho_a$ w.r.t. the partial volume $\nu_a$ and its 
mass density $\crho_a$ w.r.t. the full control volume $V$ are related by
\begin{equation}
	\crho_a = \alpha_a \rho_a.
\end{equation}
The mixture entropy density $\eta = \rho S $ is defined as
\begin{equation}\label{eq.entropymix}
	\eta := \sum_{a=1}^{\Nph} \eta_a = \sum_{a=1}^{\Nph} \crho_a s_a
\end{equation}
where $s_a$ is the specific entropy of the $a$-th phase. Hence, the specific mixture entropy can be 
computed as
\begin{equation}
	S = \frac{\eta}{\rho} = c_1 s_1 + c_2 s_2 + \ldots + c_\Nph s_\Nph.
\end{equation}
\subsection{The kinematic quantities of mixtures}
Due to the conservation principle, the total momentum of a control volume is defined as the sum of the momenta of its parts. Thus the linear momentum $ \UU=\{U_k\}$ of the mixture control volume $V$, where $k$ denotes the component in space, is defined as the 
sum of the linear momenta $\uu_a = \{u_{a,k}\}:= \crho_a \vv_a $ of the constituents
\begin{equation}\label{eqn.mom}
	\UU := \uu_1 + \uu_2 + \ldots + 
	\uu_\Nph = \crho_1 \vv_1 + \crho_2 \vv_2 + \ldots + \crho_\Nph \vv_\Nph ,
\end{equation} 
where  $ \vv_a = \{v_{a,k}\} $ is the velocity of the $a$-th phase.
The velocity $ \VV =\{V_k\} := \UU/\rho $ of the mixture control volume is therefore defined as the 
center of mass velocity
\begin{equation}\label{eqn.moma}
	\VV := \frac{\UU}{\rho} = \frac{ \crho_1 \vv_1 + \crho_2 \vv_2 + \ldots + \crho_\Nph \vv_\Nph 
	}{\rho} = c_1 \vv_1 + c_2 \vv_2 + \ldots + c_\Nph \vv_\Nph.
\end{equation}
For the SHTC formulation of the mixture equations, in addition to the mixture momentum $\UU$, one 
also needs the relative velocity $\ww_a=\{w_{a,k}\}$ fields
\begin{equation}\label{eqn.w}
	\ww_a = \vv_a - \vv_\Nph,
	\qquad
	w_{a,k} = v_{a,k} - v_{\Nph,k},
	\qquad
	k = 1, \ldots, 3.
\end{equation}
which are defined with respect to the $ \Nph $-th constituent that can be chosen arbitrarily.
Whereas, in order to derive a BN-type formulation it is useful to define the relative velocity with respect to the mixture velocity
\begin{equation}\label{eqn.vbar}
	\bar{\ww}_{a} := \vv_a - \VV,
	\qquad
	\bar{w}_{a,k} = v_{a,k} - V_{k},
	\qquad
	k = 1, \ldots, 3.
\end{equation}
In the SHTC theory, the relative velocity $\ww_a$ is the preferred choice because it is 
dictated by the variational formulation as well as by the systematization of the governing 
equations.\\ 

\subsection{Deformation characteristics}
In order to describe the elastic and inelastic deformations of a single material in the SHTC framework, one needs to introduce the concept of the \textit{distortion field} $\A$, by means of which the 
evolution of elastic and elastoplastic solids and the dynamics of Newtonian and non-Newtonian 
fluids can be formulated in the SHTC formalism \cite{GodRom2003,HPR2016,Hyper-Hypo2018,GPR7,nonNewtonian2021}. 
In the classical formulation of ideal elastic solids the distortion matrix is interpreted as the 
inverse of the \textit{deformation gradient tensor}, commonly denoted in the literature as $\FF 
= \pd \xx/\pd\XX$, or, in index notations $F_{iK} = \pd x_i/\pd X_K$, 
and hence, $\A = 
\pd\XX/\pd\xx$ or $A_{Ki} = \pd X_K/\pd x_i$, for the dynamics of pure 
elastic solids. Here, as usual, we denote the coordinates of the reference configuration by $\mathbf{X}$ and the coordinates in the current configuration by $\mathbf{x}$. In the case of inelastic deformations (viscous flows, plastic deformations), the 
distortion field can be interpreted as the inverse of the elastic part $\FF^\mathrm{e}$ of the 
multiplicative decomposition $\FF = \FF^\mathrm{e} 
\FF^\mathrm{i}$ of the deformation gradient into elastic and inelastic part, e.g. see 
\cite{HPR2016, Hyper-Hypo2018}. Note that in the notation of the distortion matrix entries, we distinguish between the Eulerian (lowercase) index $i,j,k$ and the Lagrangian (uppercase) index $I,J,K$ so that the distortion matrix can be seen as a triad of three basis vectors $\A=\{\A_1,\A_2,\A_3\}$, where for each $K=1,2,3$, $\A_K$ is a 3-vector $\A_K=(A_{K1},A_{K2},A_{K3})$.

In the general setting, a multiphase medium in the SHTC framework should have multiple distortion 
fields $\A_a$, $a=1,2,\ldots,\Nph$. However, a rigorous derivation of such a model from the 
variational principle remains beyond our reach, even for the case of a mixture of ideal elastic 
solids. Therefore, following \cite{Poroelast2020,Poroelasticity2}, we present a simplified 
multiphase model with a single distortion field $\A_1=\A_2 = \ldots = \A_\Nph = 
\A$. For example, a multiphase single-distortion model was successfully used in 
\cite{Nikodemou2021,Nikodemou2023} 
for transient shock-dominated problems in 1D as 
well in multiple space dimensions.
Yet, 
from our computational experience with the model, especially for the problems with moderate Mach 
numbers, it has 
appeared that it is beneficial to evolve individual distortions 
for every phase. While a rigorous theoretical derivation of the multi-distortion model is still 
absent, we shall use it in the numerical experiments in Sec.\ref{sec.num.res}.

\subsection{Heat conduction characteristics}

To account for the heat conduction within the SHTC theory, it is necessary to 
consider another vector field, usually referred to as the thermal impulse. Thus, for each 
phase, we introduce the vector fields $\J_1,\J_2,\ldots, \J_\Nph$ that characterize the direction 
and intensity of the heat transfer in each phase.
\subsection{SHTC state variables for multiphase flows} 
Here, we list the set of SHTC state variables for mixtures which are partly different from the 
conventionally used state variables, for example, in the BN-type formulations, e.g. mixture 
momentum and relative velocities $\{\UU,\ww_a\}$, $a=1,\ldots,\Nph-1$ in the SHTC formulation 
versus phase momenta 
$\uu_a$, $a=1,\ldots,\Nph$ in the 
BN-type formulations. The SHTC choice is conditioned by the variational nature of the equations and 
their symmetrization procedure.
Thus, the vector of sought conservative SHTC variables is
\begin{equation}\label{eq.state.var}
	\QQ = 
	(\UU,\A,\rho, \crho_1,\ldots,\crho_{\Nph-1},\ww_1,\ldots,\ww_{\Nph-1},
	\eta_1,\ldots,\eta_\Nph, \J_1,\ldots, \J_\Nph,
	\phi_1,\ldots,\phi_{\Nph-1} )^{\mathsmaller{T}},
\end{equation}
which is related to the vector of primitive SHTC variables
\begin{equation}\label{eq.state.var.prim}
	\PP = 
	(\VV,\A,\rho, c_1,\ldots,c_{\Nph-1},\ww_1,\ldots,\ww_{\Nph-1},
	s_1,\ldots,s_\Nph, \J_1,\ldots, \J_\Nph,
	\alpha_1,\ldots,\alpha_{\Nph-1} )^{\mathsmaller{T}},
\end{equation}
as
\begin{equation}
	\phi_a = \rho \alpha_a, 
	\qquad
	\UU = \rho \VV,
	\qquad
	\crho_a = \rho c_a
	\qquad
	\eta_a = \rho c_a s_a.
\end{equation}
One should pay attention to that $\phi_{\Nph}$, $\crho_{\Nph}$, and $\ww_{\Nph}$ are excluded from 
the set of state 
variables 
because they can be expressed as
\begin{equation}
	\phi_{\Nph} = \rho - \phi_1 - \ldots - \phi_{\Nph-1},
	\qquad
	\crho_{\Nph} = \rho - \crho_1 - \ldots - \crho_{\Nph-1},	
	\qquad
	\ww_{\Nph} = 0,
\end{equation}
likewise their primitive counterparts
\begin{equation}
	\alpha_{\Nph} = 1 - \alpha_1 - \ldots - \alpha_{\Nph-1},
	\qquad
	c_{\Nph} = 1 - c_1 - \ldots - c_{\Nph-1},	
	\qquad
	\ww_{\Nph} = 0.
\end{equation}

\section{SHTC governing equations for multiphase flows}\label{sec.SHTC.model}

The SHTC equations for two-fluid mixtures were proposed by Romenski in \cite{Rom1998,Rom2001} (and 
later developed in \cite{Romenski2007,Romenski2007,RomDrikToro2010}) based on the special structure 
of the SHTC class of equations proposed in \cite{GodRom1995,GodRom1996b}. A variational 
formulation of the SHTC two-fluid model was discussed in \cite{SHTC-GENERIC-CMAT}, and its 
extension to the case of arbitrary number of fluid phases was presented in \cite{Romenski2016}, and 
extension to the case of mixtures of solid-fluid mixtures in the context of flows in porous media 
was given in \cite{Poroelast2020,Poroelasticity2}. 

In this section, we recall the single-distortion two-phase model
\cite{Romenski2016,Poroelast2020,Poroelasticity2} and extend it to the setting of an arbitrary 
number of fluid and 
solid constituents. We also give its variational formulation in \ref{sec.variation}. 
We then generalize the single-distortion model to a multi-distortion multiphase model whose 
rigorous theoretical derivation from a variational principle is still lacking and will be subject of future research.

Moreover, in this work, the SHTC multiphase model is extended with the additional physical effects 
such 
as heat conduction, viscosity of the fluid constituents, and plasticity of the solid ones. 
Specifically, in this generalization of the SHTC mixture theory, the viscosity and plasticity are  
incorporated via the unified hyperbolic model of 
continuum mechanics \cite{GodRom2003,HPR2016,DUMBSER2016824}, and the heat conduction is modeled via 
the 
hyperbolic heat equations \cite{Rom1986heat,SIGPR,PTRSA2020}, and thus retaining the first-order 
hyperbolic form of the equations. 

\subsection{Non-dissipative multiphase single-distortion model}

Ignoring for the moment all the dissipative processes, the single-distortion multiphase SHTC model 
reads (summation over repeated spatial $i,k,\ldots=1,2,3$ and phase $a,b=1,2,\ldots,\Nph$ indices 
is implied)
\begin{subequations}\label{eqn.SHTC.system2}
	\begin{align}	
		&\frac{\pd U_i}{\pd t}  + \frac{\pd \left( U_i V_k + P \delta_{ik} + A_{Jk} \rhoE_{A_{Ji}}  
		+ w_{a,i} \rhoE_{w_{a,k}}  + J_{a,i} \rhoE_{J_{a,k}} \right) }{\pd x_k} = 0, 
		\label{eqn.SHTC.U}\\[2.5mm]
		&\frac{\pd A_{Jk}}{\pd t}  + \frac{\pd A_{Jl} V_l }{\pd x_k} + V_i \left( \frac{\pd 
		A_{Jk}}{\pd x_i} - \frac{\pd A_{Ji}}{\pd x_k} \right)  = 0,
		\label{eqn.SHTC.A}\\[2.5mm]
		&\frac{\pd \crho_a}{\pd t}  + \frac{\pd \left(\crho_a V_k + \rhoE_{w_{a,k}}\right) }{\pd 
		x_k} = 0 ,  &a{=}1,\mathellipsis,\mathrm{N}{-}1,
		\label{eqn.SHTC.crho}\\[2mm]
		&\frac{\pd w_{a,k}}{\pd t}  + \frac{\pd \left(w_{a,l} V_l+ \rhoE_{\crho_a}\right) }{\pd 
		x_k}  + V_i \left( \frac{\pd w_{a,k}}{\pd x_i} - \frac{\pd w_{a,i}}{\pd x_k} \right)= 0, 
		\quad  &a{=}1,\mathellipsis,\mathrm{N}{-}1,
		\label{eqn.SHTC.w}\\[2mm]
		&\frac{\pd \eta_a}{\pd t}  + \frac{\pd \left(\eta_a V_k+ \rhoE_{J_{a,k}}\right) }{\pd 
		x_k} = 0 ,  &a{=}1,\mathellipsis,\mathrm{N},
		\label{eqn.SHTC.eta}\\[2mm]
		&\frac{\pd J_{a,k}}{\pd t}  + \frac{\pd \left(J_{a,l}V_l + \rhoE_{\eta_a}\right) }{\pd 
		x_k}  + V_i \left( \frac{\pd J_{a,k}}{\pd x_i} - \frac{\pd J_{a,i}}{\pd x_k} \right) = 0, 
		\quad  &a{=}1,\mathellipsis,\mathrm{N},
		\label{eqn.SHTC.J}\\[2mm]
		&\frac{\pd \rho}{\pd t}  + \frac{\pd \left(\rho V_k\right) }{\pd x_k} = 0, 
		\label{eqn.SHTC.rho}\\[2mm]
		&\frac{\pd \phi_a}{\pd t}  + \frac{\pd \left(\phi_a V_k\right) }{\pd x_k} = 0, \quad  
		&a{=}1,\mathellipsis,\mathrm{N}{-}1 ,
		\label{eqn.SHTC.phi}
		\end{align}
\end{subequations}
where
\begin{equation}\label{eqn.P}
	P(\QQ) := U_i \rhoE_{U_i}   + \rho 
		\rhoE_{\rho} + \crho_a \rhoE_{\crho_a} + \eta_a \rhoE_{\eta_a}   + \phi_a \rhoE_{\phi_a}  - 
		\rhoE 
\end{equation}
is the total thermodynamic mixture pressure. We note that, as usual, the thermodynamic pressure $P$ accounts 
only for contributions from the internal energies of the constituents (see \eqref{eqn.Psum0}--\eqref{eqn.Psum}), and, in general, is different from the total mechanical pressure.

As one can see, the fluxes 
\eqref{eqn.SHTC.U}, \eqref{eqn.SHTC.crho}, and \eqref{eqn.SHTC.eta} are defined 
in terms of the derivatives of the energy potential $\rhoE(\QQ)$ with respect 
to the state vector $\QQ$. Hence, to complete the model formulation, one needs 
to specify the energy and compute 
all the derivatives $\pd \rhoE/\pd \QQ = (\rhoE_{U_i},\rhoE_{A_{Ji}}, \rhoE_{\rho}, \rhoE_{{\crho_a}}, \rhoE_{w_{ak}}, \rhoE_{\eta_a}, \rhoE_{\J_a},  \rhoE_{\phi_a})$, which we do in the following section. 

\subsection{Closure relations}\label{sub.sec.closure.relations}
In this 
section, we summarize all the formulas for partial derivatives of the energy potential 
required in the formulation of the SHTC multiphase model. 

According to the principle of energy conservation, the total energy density $ \rhoE $ of 
the mixture, in the control volume $V$, can be defined as the sum of the  energy densities $\rhoE_a 
=\ce^{\mathrm{i}}_a + \ce^\mathrm{e}_a + \ce^\mathrm{t}_a + \ce^\mathrm{k}_{a}$ of its constituents 
\begin{equation}\label{eqn.rhoE0}
	\rhoE(\QQ)  = \sum_{a=1}^{\Nph} \rhoE_a, 
\end{equation}
where, for each phase $a=1,2,\ldots,\Nph$, $\ce^\mathrm{i}_a(\rho,\crho_a,\phi_a,\eta_a)$ is the internal energy, 
$\ce^\mathrm{e}_a(\crho_a,\A)$ is the elastic 
energy, $\ce^\mathrm{t}_a(\crho_a,\J_a)$ is the energy associated to the thermal impulse, and 
$\ce^\mathrm{k}_{a}(\crho_a,\uu_a)$ is 
the kinetic energy.
\subsubsection{Internal energy, $\ce^\mathrm{i}$}\label{eos.i}
The SHTC state variables \eqref{eq.state.var} are dictated by the variational 
formulation of the governing equations. However, these might be not the optimal 
choice of variables when it comes to expressing the fluxes in terms of the 
conventional fluid characteristics such as pressure, temperature, etc. 
Therefore, it is useful to express the total energy density of the mixture 
through two parametrizations. The first one is in terms of the state vector 
$\QQ$ in \eqref{eq.state.var}, and the second one in terms of the individual 
phase state parameters $\rho_a$ and $s_a$. The latter, only concerns the internal energies 
$\ce^\mathrm{i}_a$. Thus, we shall use the following notations
\begin{align}
	&\ce^\mathrm{i}_a(\rho, \crho_a,\phi_a,\eta_a) = \hat{\ce}^\mathrm{i}_a(\rho_a,s_a) = \crho_a 
	\hat{ e}^\mathrm{i}_a(\rho_a,s_a) = \crho_a \hat{ e}^\mathrm{i}_a\left (\frac{\crho_a 
	\rho}{\phi_a},\frac{\eta_a}{\crho_a}\right ), \qquad a=1,2,\ldots,\Nph-1\label{pot.hyd} 
\end{align}
and for $a=\Nph$
\begin{align}
	\hat{\ce}^\mathrm{i}_\Nph(\rho_\Nph,s_\Nph) = \bigg(\rho - \sum_{a=1}^{\Nph-1} \crho_a\bigg)  
	\hat{ e}^\mathrm{i}_\Nph\left( \frac{ \rho  \big(\rho - \sum_{a=1}^{\Nph-1} \crho_a\big) 
	}{\rho-\sum_{a=1}^{\Nph-1}\phi_a},\frac{\eta_\Nph}{\rho - \sum_{a=1}^{\Nph-1} \crho_a}\right).
\end{align}
With this parametrization of the internal energies, the phase pressures and temperatures are defined as
\begin{equation}\label{eqn.pa}
	p_a := \rho_a^2 \frac{\pd \hat{e}^\mathrm{i}_a}{\pd\rho_a},
	\qquad
	T_a := \frac{\pd \hat{e}^\mathrm{i}_a}{\pd s_a}.
\end{equation}


In this work, several test problems for multiphase flows of interacting gases, liquids, and solids will be presented. Each of these states of matter has its own equation of state which are described below. We remark that it is not the goal of the paper to provide a comprehensive list of equations of state for all possible materials, but rather to illustrate the flexibility of the SHTC model in handling different types of materials. 

\begin{itemize}
	\item For the gas phases, the equation of state of perfect gases is used in the form
\begin{equation}
	\hat{e}^\mathrm{i}_a(\rho_a,s_a) = \frac{ \Co_a ^2}{\gamma_a (\gamma_a-1)} \left( \frac{\rho_a}{\rhoo_{a}} 	\right)^{\gamma_a -1} e^{s_a/\Cv_a },
\end{equation}
where $\rhoo_{a}$ is the reference density, $\gamma_a$ is the adiabatic exponent, $\Co_{a}$ is the velocity of sound at normal atmospheric conditions, $\Cv_a $ is the specific heat capacity at constant volume. Then, according to \eqref{eqn.pa}, the pressure and temperature are computed as
\begin{align}
&p_a = \rho_a^2 \frac{\pd \hat{e}_a}{\pd\rho_a} = \frac{\rhoo_a  \Co_a  ^2  }{\gamma_a} \left( \frac{\rho_a}{\rhoo_a } 	\right)^{\gamma_a} e^{s_a/\Cv_a },\\
&T_a = \frac{\pd \hat{e}_a}{\pd s_a} = \frac{ \Co_a  ^2}{\Cv_a \gamma_a (\gamma_a-1)} \left( \frac{\rho_a}{\rhoo_a } 	\right)^{\gamma_a -1} e^{s_a/\Cv_a },
\end{align}
and the phase velocity of sound $\up{C}_a$ can be computed as 
\begin{equation}
	\up{C}_a^2 := \frac{\pd p_a}{\pd \rho_a} =  \Co_a  ^2 \left( \frac{\rho_a}{\rhoo_a } 	\right)^{\gamma_a -1} e^{s_a/\Cv_a }.
\end{equation}

\item For the liquid and solid phases, the stiffened gas equation of state will be used in the form
\begin{equation}\label{eqn.stiffened}
	\hat{e}^\mathrm{i}_a(\rho_a,s_a) = \frac{ \Co_a  ^2}{\gamma_a (\gamma_a-1)} \left( \frac{\rho_a}{\rhoo_a } 	\right)^{\gamma_a -1} e^{s_a/\Cv_a } + \frac{ \rhoo_a  \Co_a  ^2 - \gamma_a \po_{a} }{\gamma_a \rho_a},
\end{equation}
denoting with $\po_{a}$ the reference (atmospheric) pressure. In this case, the pressure and temperature are given by
\begin{align}
&p_a = \rho_a^2 \frac{\pd \hat{e}_a}{\pd\rho_a} = \frac{\rhoo_a \Co_a  ^2  }{\gamma_a} \left( \frac{\rho_a}{\rhoo_a } 	\right)^{\gamma_a} e^{s_a/\Cv_a } - \frac{ \rhoo_a  \Co_a  ^2 - \gamma_a p_{0a} }{\gamma_a },\\
	&T_a = \frac{\pd \hat{e}_a}{\pd s_a} = \frac{ \Co_a  ^2}{\Cv_a \gamma_a (\gamma_a-1)} \left( \frac{\rho_a}{\rhoo_a } 	\right)^{\gamma_a -1} e^{s_a/\Cv_a }
\end{align}
and the phase adiabatic sound speed $C_a$ results in 
\begin{equation}
	\up{C}_a^2 := \frac{\pd p_a}{\pd \rho_a} =  \Co_a  ^2 \left( \frac{\rho_a}{\rhoo_a } 	\right)^{\gamma_a -1} e^{s_a/\Cv_a }.
\end{equation}

\end{itemize}

\subsubsection{Elastic energy, $\ce^\mathrm{e}$}\label{eos.e}
We recall that according to the unified model of continuum mechanics 
\cite{HPR2016}, the Navier-Stokes equations can be considered as a the stiff relaxation 
limit of the SHTC viscoelastic model \cite{DUMBSER2016824}, and thus, like in elastic solids, their response 
to shear deformations is characterized by the elastic energy. 

In this work, the part of the energy density associated with the elastic-shear 
stress, $\ce^\mathrm{e} $, 
is assumed to be proportional to the second invariant of the deviator $\mathrm
{dev}{G}_{a,ij}={G}_{a,ij}-( {G}_{a,kk}/3)\delta_{ij}$ of the metric tensor of elastic deformations ${G}_{a,ij} =  {A}_{a,Ji}{A}_{a,Jj}$, and reads 
\begin{equation}\label{eqn.ce.s}
	\ce^\mathrm{e}_a(\crho_a,\A_a)  = \frac{1}{4}\crho_a  \Cs_a^2\
	\left(\mathrm{dev}{G}_{a,ik}  \mathrm{dev} {G}_{a,ki} \right),
\end{equation}
where $\Cs_a$ is a parameter representing the propagation speed of 
small-amplitude shear waves in $a$-th phase, and here it is referred to as 
\textit{shear sound velocity}. 

At the moment, an SHTC formulation for multiphase flows with different 
distortion fields $\A_a$ is unknown, and in the theoretical part, we assume a 
single-distortion approximation $\A=\A_1=\A_2 = \ldots =\A_\Nph$. In 
particular, it is unclear with which velocity the individual phase distortion fields $\A_a$ should 
be transported --- the mixture or phase velocity. Therefore, our initial 
intention was to use the single-distortion approximation. 
However, in numerous numerical experiments, we found out that the 
single-distortion approximation is not sufficient for obtaining good results in 
low-Mach problems. Instead, the best results were obtained for the case when 
the phase distortions $\A_a$ are assumed to be different and evolved according 
to different equations with phase velocities being the transport velocities, see Sec.\ref{sec.multiA}.


%
\subsubsection{Thermal energy, $\ce^\mathrm{t}$}\label{eos.t}
The energy associated to the heat conduction can be taken as 
\begin{equation}\label{eqn.ce.j}
	\ce^\mathrm{t}_a(\crho_a,\J_a)  = \crho_a\frac{1}{2} \Ch_a^2\ J_{a,k} J_{a,k},
\end{equation}
where $\Ch$ is a parameter representing the propagation of small-amplitude 
thermal perturbations. Note that one could use other forms for $\ce^\mathrm{t}_a$, see e.g. \cite{Dhaouadi2024}

\subsubsection{Kinetic energy, $\ce^\mathrm{k}$}\label{eos.k}
The phase kinetic energy 
\begin{align}
\ce^\mathrm{k}_{a} = \frac{1}{2 \crho_a} \Vert\uu_a \Vert^2,
\end{align}
is defined in terms of the phase momenta $\uu_a$. However, to compute the partial derivatives $\pd\rhoE/\pd\QQ$, after some algebra, it can be also expressed in terms of the SHTC variables $\QQ$ as
\begin{equation}\label{eqn.rhoE}
	\sum_{a=1}^{\Nph} \ce^\mathrm{k}_{a}  = \frac{1}{2\rho} \sum_{k=1}^{3}U_k^2 + 
	W(\rho,\crho_1,\ldots,\crho_{\Nph-1} ,\ww_1, \ldots, \ww_{\Nph-1}) 
\end{equation}
where the kinetic energy of relative motion $W$ is defined as
\begin{equation}\label{eqn.W}
	W(\rho,\crho_1,\ldots,\crho_{\Nph-1} ,\ww_a, \ldots, \ww_{\Nph-1}) :=  \frac{1}{2} 
	\sum_{k=1}^{3}\sum_{a=1}^{\Nph-1}\crho_a w_{a,k}^2 - \frac{1}{2\rho} 
	\sum_{k=1}^{3}\left (\sum_{a=1}^{\Nph-1}\crho_a w_{a,k} \right )^2.
\end{equation}

\subsubsection{Thermodynamic forces}

Keping in mind that 
\begin{align}
	\crho_\Nph = \rho - \sum_{a=1}^{\Nph-1} \crho_a \quad \mathrm{and} \quad \phi_\Nph = 
	\rho-\sum_{a=1}^{\Nph-1}\phi_a,
	\end{align}
the partial derivatives of the energy potential $\rhoE$ with respect to the state vector $\QQ$ are given by
\begin{subequations}\label{eqn.dEdQ0}
	\begin{align}
		\frac{\pd \rhoE}{\pd U_i} &= \frac{1}{\rho} U_i = V_i,
		\\[2mm]
		\frac{\pd \rhoE}{\pd A_{Jk}} &= \sum_{a=1}^{\Nph} \frac{\pd \ce^\mathrm{e}_a}{\pd A_{Jk}} = \sum_{a=1}^{\Nph} \crho_a \Cs_a^2 A_{Ji} 
		\mathrm{dev}G_{ik},
		\\[2mm]
		\frac{\pd \rhoE}{\pd \rho} &= \sum_{a=1}^{\Nph-1} \frac{\pd \hat{\ce}^\mathrm{i}_a}{\pd 
		\rho_a}\frac{\crho_a}{\phi_a} 
		+
		\frac{\pd \hat{\ce}^\mathrm{i}_\Nph}{\pd\rho_\Nph} \left(\frac{\rho  
		\phi_\Nph - \rho\crho_\Nph + \crho_\Nph \phi_\Nph}{\phi _\Nph^2}\right)
		-
		\frac{\pd \hat{\ce}^\mathrm{i}_\Nph}{\pd s_\Nph} \frac{\eta_\Nph}{\crho_\Nph^2}
		+ \frac{\ce^\mathrm{e}_\Nph}{\crho_\Nph} + \frac{\ce^\mathrm{t}_\Nph}{\crho_\Nph} + \\[2mm]
		&\frac{1}{2\rho^2} \sum_{k=1}^{3}\sum_{a=1}^{\Nph} (\crho_a w_{a,k})^2
		-
		\frac{1}{2\rho^2} \sum_{k=1}^{3}U_k^2,
		\\[2mm]
		\frac{\pd \rhoE}{\pd \phi_a} &= -\frac{\crho_a\rho}{\phi_a^2}\frac{\pd\hat{\ce}^\mathrm{i}_a}{\pd \rho_a}
		+
		\frac{\crho_\Nph\rho}{\phi_{\Nph}^2}\frac{\pd\hat{\ce}^\mathrm{i}_\Nph}{\pd \Nph}
		=
		-\frac{1}{\rho}\left(\rho_a^2 \frac{\pd \hat{e}_a}{\pd 
		\rho_a} - \rho_{\Nph}^2 \frac{\pd 
			\hat{e}_{\Nph}}{\pd \rho_{\Nph}}\right) =-\frac{p_a - p_{\Nph}}{\rho}
		,\qquad
		a=1,\ldots,N-1,
		\label{eqn.dEdphi}
		\\[2mm]
		\frac{\pd \rhoE}{\pd \crho_a} &=
		\frac{\rho}{\phi_a}\frac{\pd \hat{\ce}^\mathrm{i}_a}{\pd \rho_a} - \frac{\eta_a}{\crho_a^2}\frac{\pd 
		\hat{\ce}^\mathrm{i}_a}{\pd s_a}
		-
		\frac{\rho}{\phi_\Nph}\frac{\pd \hat{\ce}^\mathrm{i}_\Nph}{\pd \rho_\Nph} + 
		\frac{\eta_\Nph}{\crho_{\Nph}^2}\frac{\pd \hat{\ce}^\mathrm{i}_\Nph}{\pd s_\Nph}
		+ \frac{\ce^\mathrm{e}_a}{\crho_a} - \frac{\ce^\mathrm{e}_\Nph}{\crho_\Nph} + \frac{\ce^\mathrm{t}_a}{\crho_a} - \frac{\ce^\mathrm{t}_\Nph}{\crho_\Nph} + 
		\nonumber\\[2mm]
		& \hspace{1cm} \frac{1}{2}\sum_{k=1}^{3}w_{a,k}^2
		-
		\frac{1}{\rho}\sum_{k=1}^{3}\sum_{b=1}^{\Nph-1}\crho_bw_{b,k}w_{a,k},
		\quad
		a=1,\ldots,N-1,
		\label{eqn.dEdcrho}
		\\[2mm]
		\frac{\pd \rhoE}{\pd w_{a,k}} &=\crho_a w_{a,k} - 
		\frac{\crho_a}{\rho}\sum_{b=1}^{\Nph-1}\crho_b w_{b,k} = \crho_a  ( v_{a,k} - V_{k}) = \crho_a \bar{w}_{a,k},
		\qquad
		a=1,\ldots,N-1,
		\label{eqn.dEdw}
		\\[2mm]
		\frac{\pd \rhoE}{\pd \eta_a} &= \frac{1}{\crho_a}\frac{\pd \hat{\ce}^\mathrm{i}_a}{\pd s_a} = \frac{\pd \hat{e}^\mathrm{i}_a}{\pd s_a} =  T_a,\\[2mm]
		\frac{\pd \rhoE}{\pd J_{a,k}} &= \frac{\pd \ce^\mathrm{t}}{\pd J_{a,k}} = \crho_a \Ch_a^2 J_{a,k} .
	\end{align}
\end{subequations}

In particular, with formulas \eqref{eqn.dEdQ0}, one can show that the mixture 
pressure \eqref{eqn.P} can be computed as the following sum
\begin{equation}\label{eqn.Psum0}
	P =  \sum_{a=1}^{\Nph} \left( \rho_a \frac{\pd \hat{\ce}^\mathrm{i}_a}{\pd \rho_a} - \hat{\ce}^\mathrm{i}_a 
	\right)
\end{equation}
which suggests that the quantities
\begin{equation}\label{eqn.hat.pa}
	P_a := \rho_a \frac{\pd \hat{\ce}^\mathrm{i}_a}{\pd \rho_a} - \hat{\ce}^\mathrm{i}_a
\end{equation}
can be called the partial phase pressures. Moreover, due to the fact $ \rho_a \frac{\pd 
\hat{\ce}^\mathrm{i}_a}{\pd \rho_a} - \hat{\ce}^\mathrm{i}_a = \alpha_a \rho_a^2 
\frac{\pd 
	\hat{e}^\mathrm{i}_a}{\pd \rho_a}$, the partial phase pressures $P_a$ and 
the single phase pressures \eqref{eqn.pa} are related by
\begin{equation}
	P_a = \alpha_a p_a.
\end{equation}
In other words, according to the SHTC formulation for multiphase flows, the mixture pressure can be computed as 
\begin{equation}\label{eqn.Psum}
	P = P_1 + \ldots + P_{\Nph} = \alpha_1 p_1 + \ldots + \alpha_{\Nph} p_{\Nph},
\end{equation}
which is also known as Dalton's law of partial pressures in mixtures.

\begin{equation}\label{eqn.dEcrho2}
		\frac{\pd \rhoE}{\pd \crho_{a}} 
		= 
		\mu_{a} - \mu_{\Nph} 
		+
		e^\mathrm{e}_{a} - e^\mathrm{e}_{\Nph}
		+
		e^\mathrm{t}_{a} - e^\mathrm{t}_{\Nph}
		+
		\sum_{k=1}^{3} w_{a,k} \left( \bar{w}_{a,k} - \frac{1}{2}w_{a,k} \right),
\end{equation} 
where 
\begin{equation}
	\mu_{a} := e^\mathrm{i}_{a} + \frac{p_{a}}{\rho_{a}} - s_{a} T_{a} = \frac{\pd \ce^\mathrm{i}_a}{\pd \crho_a}
\end{equation}
is the chemical potential of the $a$-th constituent, and $e^\mathrm{e}_a = \ce^\mathrm{e}_a/\crho_a$, $e^\mathrm{t}_a = \ce^\mathrm{t}_a/\crho_a$.
\subsection{Irreversible dynamics, dissipative processes} \label{sub.sec.irr.dissipative}
In the SHTC theory, a dissipative process is associated with the irreversible 
part 
of the time evolution equations that increases the entropy of the system and 
that is modeled via algebraic relaxation source terms \cite{SHTC-GENERIC-CMAT}. 
They are defined in terms of the gradients of the energy $\rhoE_{\QQ}
$, (i.e. in terms of the conjugate state variables), 
thus the irreversible part of the SHTC equations can be called the gradient dynamics \cite{SHTC-GENERIC-CMAT, PKGBook2018}.

The following, at the moment, arbitrary functions of the conjugate state 
variables
\begin{equation}\label{eqn.source.dissipative}
	Z_{Jk}(\rhoE_{\QQ}), \quad \Lambda_{a,k}(\rhoE_{\QQ}), \quad \chi_a(\rhoE_{\QQ}), \quad \Gamma_{a,k}(\rhoE_{\QQ}), \quad \pi_a(\rhoE_{\QQ}),\quad \Phi_a(\rhoE_{\QQ}),
\end{equation}
can be added in the right-hand side of the system equations \eqref{eqn.SHTC.system2}
\begin{subequations}\label{eqn.SHTC.system3}
	\begin{align}	
		&\frac{\pd U_i}{\pd t}  + \frac{\pd \left( U_i V_k + P \delta_{ik} + A_{Ji} \rhoE_{A_{Jk}}  + w_{a,i} \rhoE_{w_{a,k}}  + J_{b,i} \rhoE_{J_{b,k}} \right) }{\pd x_k} = 0, 
		\label{eqn.SHTC.U3}\\[2.5mm]
		&\frac{\pd A_{Jk}}{\pd t}  + \frac{\pd A_{Jl} V_l }{\pd x_k} + V_i \left( \frac{\pd A_{Jk}}{\pd x_i} - \frac{\pd A_{Ji}}{\pd x_k} \right)  = Z_{Jk},
		\label{eqn.SHTC.A3}\\[2.5mm]
		&\frac{\pd \crho_a}{\pd t}  + \frac{\pd \left(\crho_a V_k + \rhoE_{w_{a,k}}\right) }{\pd 
		x_k} = \chi_a ,  &a{=}1,\mathellipsis,\mathrm{N}{-}1,
		\label{eqn.SHTC.crho3}\\[2mm]
		&\frac{\pd w_{a,k}}{\pd t}  + \frac{\pd \left(w_{a,l} V_l+ \rhoE_{\crho_a}\right) }{\pd 
		x_k}  + V_i \left( \frac{\pd w_{a,k}}{\pd x_i} - \frac{\pd w_{a,i}}{\pd x_k} \right)= \Lambda_{a,k}, \quad  &a{=}1,\mathellipsis,\mathrm{N}{-}1,
		\label{eqn.SHTC.w3}\\[2mm]
		&\frac{\pd \eta_a}{\pd t}  + \frac{\pd \left(\eta_a V_k+ \rhoE_{J_{a,k}}\right) }{\pd 
		x_k} = \pi_a + \Pi_a ,  &a{=}1,\mathellipsis,\mathrm{N},
		\label{eqn.SHTC.eta3}\\[2mm]
		&\frac{\pd J_{a,k}}{\pd t}  + \frac{\pd \left(J_{a,l}V_l + \rhoE_{\eta_a}\right) }{\pd 
		x_k}  + V_i \left( \frac{\pd J_{a,k}}{\pd x_i} - \frac{\pd J_{a,i}}{\pd x_k} \right) = \Gamma_{a,k}, \quad  &a{=}1,\mathellipsis,\mathrm{N},
		\label{eqn.SHTC.J3}\\[2mm]
		&\frac{\pd \rho}{\pd t}  + \frac{\pd \left(\rho V_k\right) }{\pd x_k} = 0, 
		\label{eqn.SHTC.rho3}\\[2mm]
		&\frac{\pd \phi_a}{\pd t}  + \frac{\pd \left(\phi_a V_k\right) }{\pd x_k} = \Phi_a , \quad  &a{=}1,\mathellipsis,\mathrm{N}{-}1 ,
		\label{eqn.SHTC.phi3}
	\end{align}
\end{subequations}

In what follows, we specify the form of the source terms in \eqref{eqn.source.dissipative} for the SHTC mixture model, and explain their physical meaning.

\subsubsection{Strain relaxation, $Z_{Jk}$}

The strain relaxation source term $Z_{Jk}$
\begin{equation}\label{eqn.Z}
	 Z_{Jk}:= - \frac{1}{\rho} \Upsilon  \rhoE_{A_{Jk}},
\end{equation} 
where $\Upsilon$ is a positive relaxation scaling function. For the elastic energy chosen in this paper, $\rhoE_{A_{Jk}}$ takes the form
\begin{equation}\label{eqn.rhoE.A}
	\rhoE_{A_{Jk}} = \sum_{a=1}^{\Nph}\frac{\pd \ce_a^\mathrm{e} }{\pd  A_{Jk}}  =   A_{Ji} 
	\mathrm{dev} G_{ik} \sum_{a=1}^{\Nph} \crho_a \Cs_a^2 = \rho A_{Ji} 
	\mathrm{dev} G_{ik} \sum_{a=1}^{\Nph} c_a \Cs_a^2,
	\quad
	\Upsilon = \frac{3}{\tau^\mathrm{e}} \left(\sum_{a=1}^{\Nph} c_a \Cs_a^2\right)^{-1}  \mathrm{det} (\A)^{5/3}.
\end{equation} 
The scaling parameter $\Upsilon$ is chosen in such a form so that to respect Newton's law of viscosity with a constant shear viscosity in the limit $\tau^\mathrm{e}\to 0$ \cite{HPR2016,DUMBSER2016824}, where $\tau^\mathrm{e}$ is the so-called \textit{strain relaxation time}, which governs the rate at which the strain, in a given control volume, is dissipated.
When $\tau^\mathrm{e} = 0$, the thermodynamic force $\rhoE_{A_{Jk}}$ vanishes 
instantaneously ($\rhoE_{A_{Jk}}=0$) and hence, the shear stress as well $A_{Jk}
\rhoE_{A_{Ji}}=0$, and therefore, inviscid flow is retrieved. On the other 
hand, for $\tau^\mathrm{e} = \infty$ there is no relaxation,
and the behavior of a pure elastic solid is retrieved. For intermediate values 
$0 < \tau^\mathrm{e} <\infty$, the medium is neither ideal fluid nor elastic 
solid, but a viscoelastic material, and with a proper choice of the function 
$\tau^\mathrm{e}(\QQ)$ various material responses can be recovered, e.g. 
Newtonian \cite{HPR2016,DUMBSER2016824} for $\tau^\mathrm{e}=const$ with the kinematic viscosity (for the elastic energy in the form \eqref{eqn.ce.s})
\begin{equation}\label{eqn.viscosity}
	\nu = \frac{1}{6}\Cs^2 \tau^\mathrm{e},
\end{equation}
non-Newtonian fluids \cite{GPR7}, 
viscoplastic Bingham-type fluids \cite{nonNewtonian2021}, elastoplastic solids 
\cite{GodRom2003,Hyper-Hypo2018}.

\subsubsection{Interfacial friction, $\Lambda_{a,k}$, and thermal impulse dissipation, $\Gamma_{a,k}$}
The second dissipative process relevant to multiphase flows is the interfacial 
friction.  Accounting for this process results in relaxation of the phase 
velocities $\vv_a$ towards a common value, which in turn leads to the 
relaxation of the relative velocities $\ww_a$ towards zero. 

The relaxation source terms in the thermal impulse equations are introduced to 
model the heat conduction process within a given phase and between the 
constituents. It has appeared that to get the heat conducting BN-type model as a relaxation limit 
of the SHTC multiphase equations, one needs to 
couple the relaxation of the phase velocities with the relaxation of the 
thermal impulse. Thus, following \cite{Romenski2007,Thomann2023}, the relaxation source terms $\Lambda_{a,k}$ and $\Gamma_{a,k}$ are introduced in the phase velocity and thermal impulse equations, respectively, as
\begin{subequations}\label{eqn.vel.relax1}
	\begin{align}
		\Lambda_{a,k}&:= -\frac{1}{\rho}\sum_{b=1}^{\Nph-1}\lambda_{ab,k} \rhoE_{w_{b,k}}
		- \frac{1}{\rho}\sum_{b=1}^{\Nph}\zeta_{ab,k} \rhoE_{J_{b,k}},
		\qquad
		\rhoE_{w_{b,k}} = \frac{\pd \rhoE}{\pd w_{b,k}} = \crho_b (v_{b,k} - V_{k}),
		\\
		\Gamma_{a,k}&:= - \frac{1}{\rho}\sum_{b=1}^{\Nph-1}\zeta_{ba,k} \rhoE_{w_{b,k}}
		- \frac{1}{\rho} \gamma_{a,k} \rhoE_{J_{a,k}},
		\qquad
		\qquad
		\rhoE_{J_{bk}} = \crho_b \frac{\pd  {e}^\mathrm{t}_{b} }{\pd J_{{b},k}} = \crho_b \Ch_b^2  J_{{b},k} ,
	\end{align}
\end{subequations}
where the kinetic coefficients $\lambda_{ab,k} = \lambda_{ab,k}(\QQ)$ are the entries of three 
($k=1,2,3$) symmetric positive semi-definite matrices, and 
\begin{equation}
	 \gamma_{a,k} = \frac{1}{\tau^\mathrm{t}_{a,k} \Ch_a^2} \geq 0
\end{equation} 
where $\tau^\mathrm{t}_{a,k}$ are the relaxation time scales governing the 
heat conduction process within the $a$-th constituent. Moreover, in accordance 
with the Onsager principle, we chose the matrices of dissipative kinetic 
coefficients to be symmetric, for consistency with the first and second laws of 
thermodynamics, we must require that the following three matrices ($k=1,2,3$)
\begin{equation}\label{eqn.matrices.wJ}
	\boldsymbol{\Psi}_k = \left(
		\begin{array}{cc}
		 \boldsymbol{\lambda}_k & \boldsymbol{\zeta}_k   \\
		 \boldsymbol{\zeta}_k^T & \boldsymbol{\gamma}_{k}	 \\
		\end{array}
		\right) \geq 0,
\end{equation} 
are positive semi-definite, where $\boldsymbol{\lambda}_k$ and $\boldsymbol{\zeta}_k$ are the $(\Nph-1)\times(\Nph-1)$ and $(\Nph-1)\times\Nph$ matrices with the entries $\lambda_{ab,k}$ and $\zeta_{ab,k}$, respectively, and $\boldsymbol{\gamma}_{k}$ are the diagonal matrices with the entries $\gamma_{a,k}$ on the diagonal.

\subsubsection{Kinetics of phase transformation, $\chi_a$}
For the sake of completeness, we also mention how the kinetics of phase transformation can be 
introduced into the SHTC mixture equations. In order to achieve this, it is 
necessary to introduce the sources $\chi_a$ into the phase mass balance 
equations, which are defined as
\begin{equation}
	\chi_a := -\rho\sum_{b=1}^{\Nph-1} \chi_{ab} \rhoE_{\crho_b},
	\qquad
	\sum_{a=1}^{\Nph}\chi_a = 0,
\end{equation}
where $\rhoE_{\crho_a}$ is given in \eqref{eqn.dEdcrho}, or \eqref{eqn.dEcrho2}, and the kinetic coefficients $\chi_{ab} = \chi_{ab}(\QQ)$ form a symmetric positive semi-definite matrix. \\

\subsubsection{Temperature relaxation, $\pi_a$}

The source terms $\pi_a$ in the {phase entropy} equations are defined as
\begin{equation}\label{eqn.pi_a}
	\pi_a := -\crho_a\frac{\rhoE_{\eta_a} - \sum\limits_{b=1}^{\Nph}\frac{\crho_b}{\rho} 
	\rhoE_{\eta_b}}{\tau \rhoE_{\eta_a}}
	=
	\crho_a\frac{T_a - \sum\limits_{b=1}^{\Nph}c_b T_b}{\tau T_a},
\end{equation}
and they model the phase temperature relaxation towards the common temperature 
\begin{equation}\label{eqn.T}
	T := \sum\limits_{a=1}^{\Nph}c_a T_a
\end{equation}
which can be called the temperature of the mixture control volume. Here, $\tau<0$ is the relaxation 
parameter that characterizes the rate at which the temperature equilibrium $T_1=\ldots =T_{\Nph} = 
T$ is approached by the 
system. \\

\subsubsection{Pressure relaxation, $\Phi_a$}

Finally, the dissipative process related to the {pressure relaxation} towards a 
common pressure are introduced as source terms in the volume fraction balance 
laws of the phases by the functions
\begin{equation}
	\Phi_a := - \rho \sum_{b=1}^{\Nph-1} \phi_{ab} \rhoE_{\phi_b},
	\qquad
	\rhoE_{\phi_a} = \frac{\pd \rhoE}{\pd \phi_a} 
	=
	-\frac{p_a - p_{\Nph}}{\rho},
\end{equation}
where $\phi_{ab} = \phi_{ab}(\QQ)$ are the kinetic coefficients which again are the entries of a symmetric positive semi-definite matrix.

\subsubsection{Entropy production terms, $\Pi_a$}

The remaining undefined dissipative terms, the entropy production terms $\Pi_a$, $a=1,\ldots,\mathrm{N}$, 
serve the goal of 
making the system compatible with the two laws of thermodynamics. 
Therefore, to fulfill the first and second 
law of thermodynamics (see \eqref{eqn.E.sum} and \eqref{eqn.entropy.prod} for details), $\Pi_a$ must be defined as
\begin{subequations}\label{eqn.Pi_a}
\begin{equation}
\Pi_a := \frac{1}{\rhoE_{\eta_a}} \left(\hat{\Pi}_a + c_a \tilde{\Pi}\right)
\end{equation}
where 
\begin{equation}
	\hat{\Pi}_a := \sum_{k=1}^{3} \frac{1}{\gamma_{a,k}} 
	\left[
		\left( 
			\gamma_{a,k} \rhoE_{J_{a,k}} 
			+ \sum_{b=1}^{\Nph-1} \zeta_{ba,k} \rhoE_{w_{b,k}}
		\right)^2 
		- \left( \sum_{b=1}^{\Nph-1} \zeta_{ba,k} \rhoE_{w_{b,k}}\right)^2 
		+ \sum_{b=1}^{\Nph-1} \zeta_{ba,k}^2
\rhoE_{w_{b,k}}^2\right],
\end{equation}
\begin{equation}
	\tilde{\Pi} := {\frac{1}{\rho} \rhoE_{A_{Ji}}\Upsilon  
	\rhoE_{A_{Jk}}  } + \sum_{b=1}^{\Nph-1} \sum_{c=1}^{\Nph-1} \left( 
		\sum_{k=1}^{3} \bigg( \rhoE_{w_{c,k}} \lambda_{cb,k} \rhoE_{w_{b,k}} - \frac{1}{\gamma_{c,k}} \rhoE_{w_{b,k}}^2 \zeta_{bc,k}^2 \bigg) 
	+
	\rho\rhoE_{\crho_b} \chi_{bc} \rhoE_{\crho_c} 	
	+
	\rho \rhoE_{\phi_b} \phi_{bc} \rhoE_{\phi_c} 
	\right) .
\end{equation}
\end{subequations}
Remark that the positive sign of each $\Pi_a\geq 0$ is guaranteed for the case 
when the coupling coefficients $\zeta_{ab,k}$ in \eqref{eqn.vel.relax1} vanish. 
In this case, \eqref{eqn.Pi_a} is the sum of quadratic forms with positive 
semi-definite matrices of coefficients, $\lambda_{ab,k}$, $\chi_{ab}$, $\phi_
{bc}$, and hence $\Pi_a\geq 0$. Apparently, for sufficiently small $\zeta_{ab,k}
$ we also have $\Pi_a\geq 0$ (because $\hat{\Pi}_a\geq 0$ and $\tilde{\Pi}\geq 0$), but such inequalities cannot be guarantied for 
arbitrary $\zeta_{ab,k}$. Yet, even if the production terms $\Pi_a\geq0$, the 
phase entropies may 
decrease due to the presence of the temperature relaxation terms $\pi_a$ that 
makes the sign of 
$\Pi_a-\pi_a$, in general, indefinite. This of course doesn't contradict the 
second law because the 
mixture constituents are not isolated systems. However, as discussed in the 
next section, this 
choice 
of $\Pi_a$ guarantees the fulfillment of the second law of thermodynamics for 
the entire mixture.

\subsubsection{Thermodynamic equilibrium}

The dissipative source terms $\Phi_a$, $\Lambda_{a,k}$, $\Gamma_{a,k}$ and $\chi_a$ are defined in such a way that 
they diminish the thermodynamic forces $\rhoE_{\phi_a}$, $\rhoE_{w_{a,k}}$, $\rhoE_{J_{a,k}}$, $\rhoE_{A_{Jk}}$, and $\rhoE_{\chi_a}$, 
i.e they lead the mixture towards a thermodynamic equilibrium state at which these forces must vanish 
$\rhoE_{\phi_a}=0$, $\rhoE_{w_{a,k}}=0$, $\rhoE_{J_{a,k}}= 0$, $\rhoE_{A_{Jk}}=0$, and $\rhoE_{\chi_a}=0$, while the temperature relaxation terms
$\pi_a $ tend to make the phase temperatures equal. 
\subsection{Consistency with the first and second laws of thermodynamics}\label{sec.consistency}
One may notice that the total energy conservation law (first law of thermodynamics) is 
not listed within the set of 
equations \eqref{eqn.SHTC.system3}. In fact, one of the main features of all the SHTC models 
\cite{GodRom1995,GodRom1996b,Rom1998,Rom2001,SHTC-GENERIC-CMAT} is that the energy conservation law
\begin{equation}\label{eq.Econs}
	{\frac{\pd \rhoE}{\pd t} + \frac{\pd}{\pd x_k} \left( V_k \rhoE  +V_k ( P \delta_{ik} + w_{a,i} 
	\rhoE_{w_{a,k}}  + A_{Ji} \rhoE_{A_{Jk}} ) + \rhoE_{\crho_a} \rhoE_{w_{a,k}} + \rhoE_{\eta_a} \rhoE_{J_{a,k}}  \right) = 0}
\end{equation}
 is automatically fulfilled for smooth solutions of system \eqref{eqn.SHTC.system3}. 
 In other words, the energy conservation law is a consequence of the governing equations and can be obtained as a linear combination of the governing equations 
\eqref{eqn.SHTC.system3} 
multiplied with the corresponding factors (the thermodynamic conjugate variables or main-field variables \cite{MullerRuggeri1998}). 
Thus, Eq.\eqref{eq.Econs} can be obtained as the following linear combination of equations \eqref{eqn.SHTC.system3} multiplied by the corresponding factors
\begin{align}\label{eqn.E.sum}
\begin{split}
	&\eqref{eq.Econs} \equiv 
					 	 \rhoE_{U_i}\cdot\eqref{eqn.SHTC.U3}
					 + \rhoE_{A_{Jk}}\cdot\eqref{eqn.SHTC.A3}
					 + \rhoE_{\crho_a}\cdot\eqref{eqn.SHTC.crho3}
					 + \rhoE_{w_{a,k}}\cdot\eqref{eqn.SHTC.w3}
					 + \rhoE_{\eta_a}\cdot\eqref{eqn.SHTC.eta3}
					 + \rhoE_{J_{a,k}}\cdot\eqref{eqn.SHTC.J3} + \\
					 &\phantom{\eqref{eq.Econs} \equiv 
					 	 \rhoE_{U_i}\cdot\eqref{eqn.SHTC.U3}
					 + \rhoE_{A_{Jk}}\cdot\eqref{eqn.SHTC.A3} }
					 + \rhoE_{\rho}\cdot\eqref{eqn.SHTC.rho3}
					 + \rhoE_{\phi_a}\cdot\eqref{eqn.SHTC.phi3}
.
\end{split}
\end{align}
However, we can obtain zero in the right hand-side of \eqref{eq.Econs} by these means only if we 
define the phase entropy production terms $\Pi_a$ as in \eqref{eqn.Pi_a}. Note that the temperature 
relaxation terms $\pi_a$ are defined in such a way that they sum up to zero in \eqref{eqn.E.sum}:
\begin{equation}
	\sum_{a=1}^{\Nph}\rhoE_{\eta_a} \pi_a = \sum_{a=1}^{\Nph}T_a \pi_a = 0.
\end{equation}

As we have already mentioned, if the relaxation processes discussed in the previous section are taken into account, our choice 
\eqref{eqn.Pi_a} of the phase entropy production terms $\Pi_a$ cannot guarantee positive sign of 
$\Pi_a-\pi_a $. However, the mixture itself (in the absence of exchange with the exterior) is an isolated 
system and the second law must hold. Thus, our choice of the phase entropy production terms 
$\Pi_a$ not only guarantees the energy conservation law for the entire mixture (the first law of 
thermodynamics) but also the second law. Indeed, 
the mixture entropy density is defined as
\begin{equation}
	\eta = \eta_1 + \ldots + \eta_{\Nph}
\end{equation}
and fulfills the entropy balance law
\begin{subequations} \label{eqn.entropy.prod}
	\begin{equation}
		\frac{\pd \eta}{\pd t} + \frac{\pd (\eta V_a)}{\pd x_k} = \Pi + \pi \geq 0,
	\end{equation}	
which is the sum of the phase entropy balance laws, and where the mixture entropy production $\Pi$ and $\pi$ are defined as
	\begin{equation}\label{eqn.Pi}
		\qquad \Pi := \sum_{a=1}^{\Nph}\Pi_a
		=
		\sum_{a=1}^{\Nph}\hat{\Pi}_a + 
		\left( \frac{c_1}{T_1} + \ldots + \frac{c_{\Nph}}{T_{\Nph}}\right) \tilde{\Pi} \geq 0,
	\end{equation}	
	\begin{equation}
		\pi = -\sum_{a=1}^{\Nph} \pi_a = \frac{1}{2} \text{tr} \left( \boldsymbol{M}^{\text T} \boldsymbol{M} \right) \geq 0,
	\end{equation}
\end{subequations}
with $\boldsymbol{M}$ being a symmetric matrix with the entries 
\begin{equation}
	M_{ab} = \sqrt{ \frac{\crho_a \crho_b}{\rho \tau} \frac{(T_a - T_b)^2}{T_a T_b}}, \qquad 
	a,b=1,\ldots,\mathrm{N}.
\end{equation}
The sign on the right hand side of \eqref{eqn.Pi}, and hence 
the signe of the overall mixture entropy production in \eqref{eqn.entropy.prod}, is guaranteed only for 
sufficiently small coupling coefficients $\zeta_{ab,k}$ in \eqref{eqn.vel.relax1}. 
\section{Baer-Nunziato form of the multiphase SHTC equations}\label{BN.SHTC.form}
As anticipated in the introduction, it is of particular interest to compare the 
structure of the proposed SHTC multiphase model \eqref{eqn.SHTC.system3} with 
the other approaches. However, the SHTC state variables dictated by its 
variational nature are different from those traditionally used in the 
popular Baer-Nunziato (BN) model \cite{baernunziato}, or in the
conventional formulation of balance laws for multiphase flows \cite{Muller1968,Gurtin1971,StewartWendroff1984,Passman1984,Nigmatulin}. Therefore, to 
make the comparison possible, we first need to put the SHTC equations in a form 
that is similar to the conventional way of writing multiphase models, i.e. in 
terms of the phase mass, momentum, and energy balance laws.

It is of particular interest to relate the SHTC formulation to the well-known 
Baer-Nunziato (BN) model \cite{baernunziato}, which was studied and used for example in \cite{SaurelAbgrall,AbgrallSaurel,AndrianovWarnecke,Schwendeman,DeledicquePapalexandris,SaurelCavitation}. However, to the best 
of our knowledge, there is no extension of the original two-fluid 
non-equilibrium (two-pressure, two-velocity) BN model 
to the case of an arbitrary number of phases, and therefore, a direct comparison 
is impossible. The authors are only aware of non-equilibrium BN-type models for three phases, see \cite{Herard1,Herard2}. An interesting alternative compressible multiphase model for a general number of $\mathrm{N}$ phases was recently presented and discussed in \cite{CompressibleNPhase}, but it does not take a BN-type form and thus is not directly comparable to the formulation considered here.   
Nevertheless, it is interesting to put the SHTC model in a form 
that is similar to the BN model to make such a comparison possible in the 
future. Moreover, the BN form of the SHTC model and not the original SHTC equations is discretized in this paper as 
it was our original intention to compare the SHTC model with the BN model in 
the context of numerical simulations. 

Let us first remind the structure of the BN model. In its original formulation 
\cite{baernunziato} it is a two-phase model specifically designed for 
applications 
describing the deflagration-to-detonation (DDT) transition in reactive, 
gas-permeable granular materials. However, nowadays this model is usually 
referred to by considering only the homogeneous part of the original system, 
i.e., without the algebraic phase interaction terms, see e.g. 
\cite{AndrianovWarnecke,SaurelAbgrall,Schwendeman,USFORCE2,
ReAbgrall2022}. 
The complete two-fluid seven-equation Baer-Nunziato model, without algebraic source terms, is a first-order system of nonlinear PDEs,  which in our notations, see Section\,\ref{sec.notaionts.mix},  reads 
\begin{subequations}\label{BN.system}
	\begin{align}
		&\frac{\pd \crho_a}{\pd t} + \frac{\pd u_{a,k}}{\pd x_k} = 0, \\ 
		&\frac{\pd u_{a,i}}{\pd t} + \frac{\pd}{\pd x_k}\left(u_{a,i} v_{a,k}  + P_a \delta_{ki} \right) =  \sum_{b=1}^{\mathrm{N}} p_{\mathrm{I},ab} \frac{\pd \alpha_b}{\pd x_i}, \\
		&\frac{\pd \rhoE_a}{\pd t} + \frac{\pd}{\pd x_k}\left(v_{a,k}  \rhoE_a  + v_{a,k} P_a \right) = -\sum_{b=1}^{\mathrm{N}} p_{\mathrm{I},ab} \frac{\pd \alpha_b}{\pd t} ,\\
		& \frac{\pd\alpha_a}{\pd t} + v_{\mathrm{I},a,k} \frac{\pd \alpha_a}{\pd x_k} = 0
	\end{align} 	
\end{subequations}
where the subscript $a=1,2$ lables the phase. As it is well known, the BN model 
requires a proper choice of the interface velocity $\vv_{\mathrm{I,a}}=\{v_{\mathrm{I,a},k}\}$ and the interface pressure $p_\mathrm{I,ab}$ \cite{SaurelStiff2007}. Such a structure of the governing equations for multiphase flows, in the 
form of mass, momentum, and energy balance laws for each phase, is referred to as the \textit{BN-type form}.

%
%
\subsection{ A BN-type form of the SHTC multiphase model}
To compare the structure of the multiphase SHTC model \eqref{eqn.SHTC.system3} 
with the BN-type structure, the system of PDEs in \eqref{eqn.SHTC.system3} is 
rewritten in terms of the phase mass and momentum. We deliberately keep all 
the dissipative source terms of the SHTC model in the BN-type form, as we 
believe that this is will be useful in the future once the multiphase 
extension of the BN model will be obtained.
In the following, we collect the individual phase related terms on the left-hand side, while the \textit
{interphase exchange terms} and the terms arising from the dissipative 
processes considered in \ref{sub.sec.irr.dissipative} are collected on the 
right-hand side.


\subsubsection{The phase mass balance equations} 
The evolution equations in \eqref{eqn.SHTC.crho3} can be immediately 
rewritten in a more traditional form using the definition of the mixture velocity $V_k$ and the expression of $\rhoE_{w_{a,k}}$, given in
\eqref{eqn.dEdw}: 
\begin{equation}\label{eqn.crho.BN}
	\frac{\pd \crho_a}{\pd t} + \frac{\pd \left( \crho_a v_{a,k} \right) }{\pd x_k} = \chi_a.
\end{equation}

\subsubsection{The phase volume fraction equations}
 The phase volume fraction equations can be retrieved from the equations \eqref{eqn.SHTC.phi3}, using the conservation of total mass and balance laws of the phase densities \eqref{eqn.crho.BN}, and read
\begin{equation}\label{eqn.alpha.BN}
	\frac{\pd \alpha_a}{\pd t} + V_{k}\frac{\pd \alpha_a }{\pd x_k} =  \frac{1}{\rho}\Phi_a.
\end{equation}
Thus, comparing with the BN model, one can conclude that the interface velocity $v_{\mathrm{I},k}$ 
of the BN model is replaced by the mixture velocity $V_k$ in the SHTC model, as was already noticed 
in \cite{Romenski2004,Romenski2007}.

\subsubsection{The phase momentum equations} 
The balance laws for the phase momenta can be obtained from the mixture 
momentum conservation and relative velocity equations in the following way. Let 
$\mathcal{D}^U_i$ represents the mixture momentum 
equation \eqref{eqn.SHTC.U3}, $\mathcal{D}^w_{a,i}$ represents the relative velocity equations 
\eqref{eqn.SHTC.w3}, and $\mathcal{D}^u_{a,i}$ represents the sought phase momentum equations. Then, the phase momentum balance equations can be obtained as
\begin{equation}\label{eqn.ua.idea}
	\mathcal{D}^u_{a,i} = \frac{\crho_a}{\rho} \mathcal{D}^U_i - \frac{\crho_a}{\rho} 
	\sum_{b=1}^{\Nph} \crho_b \mathcal{D}^w_{b,i} + \crho_a \mathcal{D}^w_{a,i}. 
\end{equation}
After lengthy but rather straightforward manipulations of the terms in \eqref{eqn.ua.idea}, the 
individual phase momentum balance equations can be written as
\begin{subequations}\label{eqn.BN.mom}
	\begin{align}
		\frac{\pd u_{a,i}}{\pd t} {+} \frac{\pd
		( u_{a,i} v_{a,k} + P_a \delta_{ki} +  \sigma^\mathrm{e}_{a,ki} + \sigma^\mathrm{t}_{a,ki} ) }{\pd x_k}= 
			&{-} c_a \sum_{b=1}^{\Nph}\ps_b \frac{\pd \alpha_b}{\pd x_i} + \ps_a \frac{\pd 
			\alpha_a}{\pd 
			x_i}\label{eqn.BN.p}   \\[2mm]
		& {-} c_a \sum_{b=1}^{\Nph} \crho_b \bar{w}_{b,k} \omega_{b,k,i} + \crho_a \bar{w}_{a,k} 
		\omega_{a,k,i} \label{eqn.BN.lift}\\[2mm]
		&- c_a \sum_{b=1}^{\Nph}\crho_b s_b \frac{\pd T_b}{\pd x_i} + \crho_a s_a \frac{\pd 
		T_a}{\pd x_i} \label{eqn.BN.T}\\[2mm]
		&- c_a \sum_{b=1}^{\Nph}\frac{\pd \left(\sigma^\mathrm{e}_{b,ki} + \sigma^\mathrm{t}_{b,ki}\right)}{\pd x_k} + \frac{\pd \left(\sigma^\mathrm{e}_{a,ki} + \sigma^\mathrm{t}_{a,ki} \right)}{\pd x_k}  \label{eqn.BN.sigma}\\[2mm]
		&+ c_a \sum_{b=1}^{\Nph}\crho_b \frac{\pd (e^\mathrm{e}_b + e^\mathrm{t}_b)}{\pd x_i} - \crho_a \frac{\pd (e^\mathrm{e}_a + e^\mathrm{t}_a)}{\pd x_i} \label{eqn.BN.ee}\\[2mm]
		& - c_a \sum_{b=1}^{\Nph}\crho_b \Lambda_{b,i} + \crho_a \Lambda_{a,i}  \label{eqn.BN.Vel.rel}\\[2mm]
		& - c_a \sum_{b=1}^{\Nph}v_{b,i} \chi_b + v_{a,i} \chi_a,\label{eqn.BN.p.trans}
	\end{align}
\end{subequations}
where
\begin{equation}
\sigma^\mathrm{e}_{a,ki} := A_{Ji} \frac{\pd \rhoE_{a}}{\pd A_{Jk}},
\qquad
\sigma^\mathrm{t}_{a,ki} := J_{a,i} \frac{\pd\rhoE_a}{\pd J_{a,k}},
\qquad
\omega_{a,k,i} := \frac{\pd w_{a,i}}{\pd x_k} - \frac{\pd w_{a,k}}{\pd x_i}
\end{equation}
and $e^\mathrm{e}_a = \ce^\mathrm{e}_a/\crho_a$, $e^\mathrm{t}_a = \ce^\mathrm{t}_a/\crho_a$,  and $\bar{w}_{a,k}$ was previously defined in \eqref{eqn.vbar}.
Here, the presence of the phase transformation terms $\chi_a$ is due to the 
appearance of $\pd \crho_a/\pd t$ and the need to replace them by their 
expressions from \eqref{eqn.crho.BN}. 

The phase momentum equations \eqref{eqn.BN.mom} derived from the SHTC mixture 
equations can be compared with those of the original BN model in order to 
understand some differences between the two approaches. Of course, the main differences are related 
to the phase interaction terms that are collected on the right hand side of \eqref{eqn.BN.mom}. 
Thus, the main difference concerns 
the fact that in the momentum equations \eqref{eqn.BN.mom} derived from the 
SHTC mixture theory, no closure problem arises for the interface quantities $p_{\mathrm{I},ab}$ and $v_{\mathrm{I},a,k}$ as in the BN model, as already noticed in \cite{Romenski2004,Romenski2007}. 
The second difference concerns the presence of the gradients of the relative 
velocities (lift forces) \eqref{eqn.BN.lift} and temperature gradients 
\eqref{eqn.BN.T} that are not present in the BN model. 

Overall, one can notice that the phase interaction terms are rather complicated and take into 
account all physical processes occurring in the mixture, including elastic and thermal stresses 
\eqref{eqn.BN.sigma}, interfacial friction \eqref{eqn.BN.Vel.rel}, and the phase transformation 
\eqref{eqn.BN.p.trans}.

\subsubsection{The phase energy equations} 
The phase energy balance laws can be obtained from the conservation equations of mixture momentum, 
relative velocity and entropy, similarly to the phase momentum equations, after a lengthy 
manipulation. However, this procedure is not illustrated in this work, since the energy equations 
depend directly on the momenta PDEs, and the latter already fulfill the comparative purpose of this 
section by making the interphase terms clear.

\subsection{Multi-distortion extension of the BN formulation of SHTC equations}\label{sec.multiA}

According to the unified model of continuum mechanics \cite{HPR2016} and SHTC 
formulation for multiphase flows \cite{Romenski2007,Romenski2016}, a true 
non-equilibrium multiphase model should have, in general, different pressures, 
temperatures, velocities, distortion fields, and etc. for each phase. However, 
at the moment, we only know how to derive the SHTC multiphase model for the 
case when the phase distortions are equal, $\A_1=\A_2=\ldots=\A_\Nph=\A$. 
Therefore, in this section we discuss a heuristic extension of the SHTC  
multiphase model in its BN form \eqref{eqn.BN.mom} to the case of different phase distortions $\A_a$, and this formulation was used to obtain all the numerical results presented in Section\,\ref{sec.num.res}. 

Let us note that the
single-distortion formulation still might be a reasonable approximation in many situations, for example, in 
high-energy transient problems like in \cite{Nikodemou2021,Nikodemou2023}. 
However, according to our experience, the single-distortion formulation is not 
suitable for low-Mach long-time fluid-structure interaction problems considered in Section \ref{sec.num.res}. Note that a multiphase 
single-velocity multi-distortion model was already proposed in 
\cite{NdanouFavrieGavrilyuk} in which the distortion fields $\A_a$ were transported by the mixture velocity, which is a natural choice for a single-velocity approximation. In our multi-velocity model we opt for the advection of the distortion fields by the phase velocities.

Thus, to extend \eqref{eqn.BN.mom} to the multi-distortion formulation, we employ a heuristic approach and simply introduce the phase distortion fields $\A_a=\{A_{a,Jk}\}$ for each phase $a=1,\ldots,\Nph$ and assume that each $\A_a$ is advected by the phase velocity $\vv_{a}$, i.e. we add to system \eqref{eqn.BN.mom} the following equations for the phase distortion fields:
\begin{equation}\label{eqn.PDE.Aa}
	\frac{\pd A_{a,Jk}}{\pd t} + \frac{\pd \left( A_{a,Jl} v_{a,l} \right) }{\pd x_k} + 
	v_{a,i}\left(\frac{\pd  A_{a,Jk} }{\pd x_{i}} - \frac{\pd  A_{a,Ji} }{\pd x_k} \right)=  
	Z_{a,Jk}.
\end{equation}

Further modifications concern the elastic stress $\sigma^\mathrm{e}_{a,ki}$ and thermodynamically conjugate variables of the phase distortion fields, which (for the elastic energies $\ce^\mathrm{e}_a$ used in this paper) are computed as
\begin{equation}\label{eqn.conj.multiA}
	\sigma^\mathrm{e}_{a,ki} = A_{a,Ji}\rhoE_{A_{a,Jk}} = A_{a,Ji}\frac{\pd\ce^\mathrm{e}_a}{\pd A_{a,Jk}},
	\qquad
	\frac{\pd \ce^\mathrm{e}_a}{\pd A_{a,Jk}}  =  \crho_a \Cs_a^2 A_{a,Ji} \mathrm{dev} G_{a,ik}.
\end{equation} 

Finally, the strain dissipation source term $Z_{a,Jk}$ in \eqref{eqn.PDE.Aa} is defined as 
\begin{equation}\label{eqn.rhoE.Aa}
	Z_{a,Jk}:= - \frac{1}{\crho_a} \Upsilon_a  \rhoE_{A_{a,Jk}},
	\qquad
	\Upsilon_a = \frac{3}{\tau^\mathrm{e}_a} \Cs_a^{-2}  \mathrm{det} (\A_a)^{5/3},
\end{equation} 
with $\tau^\mathrm{e}_a$ being the strain relaxation time of phase $a$, and in general is assumed being different for each phase.
This timescale defines the stiff nature of the strain relaxation source term towards an equilibrium state of material deformation, e.g.
\begin{enumerate}
\item for $\tau^\mathrm{e}_a = 0$, the so-called stiff relaxation limit is achieved instantaneously and therefore inviscid flow is retrieved;
\item for sufficiently small $\tau^\mathrm{e}_a$ with respect to the flow timescale, the model reproduces the Navier-Stokes equations of viscous fluids, for the chosen shear energy, it can be computed to fit the kinematic viscosity $\nu_a$ of a fluid as $\tau^\mathrm{e}_a = 6\nu_a/\Cs_a^2$ ;
\item  for $\tau^\mathrm{e}_a = \infty$ the relaxation is absent and the behavior of a pure elastic 
solid is retrieved.
\end{enumerate}
The relaxation time $\tau^\mathrm{e}_a$, in general is a function of the state variables. In the multiphase context, it is useful to define it as a function of the volume fraction $\alpha_a$ for each phase, by means of a smooth logarithmic interpolation, which can be computed as follows
\begin{align}
\label{eqn.relaxation.time}
	 \tau^\mathrm{e}_a = (\tau^\mathrm{e}_a)^\xi \ \tau_\up{o}^{1-\xi},
\end{align}
where $\xi$ can be evaluated by
\begin{align}
	\delta \alpha = 	\frac{ \alpha_a - \alpha_\up{m}}{\alpha_\up{M} - \alpha_\up{m}}, \quad \delta \alpha = \up{max} \big(  0,  \up{min}(  1, \delta \alpha )\big), \quad \xi = \delta \alpha^2 (3- 2 \delta \alpha),
\end{align}
which results in a smooth transition from $\tau^\mathrm{e}_a$ to $\tau_\up{o}$, where $\tau_\up{o}$ is usually assumed to be a small constant like $10^{-14}$. Then, $\alpha_\up{M}$ and $\alpha_\up{m}$ represent the extrema at which this operator makes the transition. In this way, where a phase is not present, the \textit{strain is dissipated instantaneously}, as for a prefect fluid, and no stresses are generated in the respective momentum conservation equation. This rescaling of the relaxation time is well suited to multimaterial problems, and in the rest of the paper is referred to as \textit{vanishing ghost solid relaxation time}.

\subsection{Three phase SHTC model in the BN-type form}
\label{sub.sec.reduced.model}

As the BN model is one of the most popular mathematical models for describing two-phase flows, there 
are many works in which BN-type equations are solved numerically. However, only a very limited number of 
publications exist on the mathematical and computational issues of BN models for flows with \textit{more} than two phases \cite{Herard1,Herard2}. It is therefore very attractive to numerically address 
the BN-type SHTC multiphase model presented in the previous section \ref{BN.SHTC.form}, which has 
been also generalized to fluid and solid mechanics in this work.

The BN-type form of the SHTC equations \eqref{eqn.BN.mom}, \eqref{eqn.PDE.Aa} admits an arbitrary number of phases, but in this paper, the numerical tests are restricted to three-phase problems in one and two dimensions. Moreover, we further simplify equations \eqref{eqn.BN.mom} by assuming isothermal conditions, absence of heat conduction and phase transitions, and isotropic phase velocity relaxation.

Such an approximation is reasonable for low-Mach number flows, i.e. shock waves 
are either totally absent, or very weak. Thus, instead of considering the phase 
energy balance laws, we assume the conservation of phase entropies on shocks, and consider phase entropy balance equations, which have a much simpler structure than 
the phase energy balance laws. Therefore, the \textit{phase entropy equations} 
can be retrieved from the SHTC equations 
\eqref{eqn.SHTC.eta3}, using the conservation of total mass \eqref{eqn.SHTC.rho3}, and read
\begin{equation}
	\frac{\pd s_a}{\pd t} + V_k \frac{\pd s_a}{\pd x_k}
	=  \Pi_a - \pi_a.
\end{equation}
Then, neglecting the phase pressure relaxation towards a common pressure for all the phases 
($\Phi_a= 0$), assuming the absence of phase transformation ($\chi_a=0$), temperature relaxation 
($\pi_a=0$), and assuming the isotropy of the phase velocity relaxation process 
($\lambda_{ab,k} $ becomes $\lambda_{a}$, $a=1,2,3$), 
the mathematical model 
that is solved numerically in this paper reads
\begin{subequations}\label{eqn.BN.NumMeth}
	\begin{align}
		&\frac{\pd \crho_a}{\pd t} + \frac{\pd u_{a,k} }{\pd x_k} =  0 
		\label{eqn.BN.NumMeth.rho.a}
		\\[1mm]
		&\frac{\pd u_{a,i}}{\pd t} + \frac{\pd}{\pd x_k}
		\big( u_{a,i} v_{a,k} + P_a \delta_{i,k} - \sigma^\mathrm{e}_{a,ik} \big) = - c_a 
		\sum_{b=1}^{\Nph}\ps_b \frac{\pd \alpha_b}{\pd x_i}+ \ps_a \frac{\pd \alpha_a}{\pd 
			x_i}  \label{eqn.BN.NumMeth.p}\\[1mm]
		&\phantom{\frac{\pd u_{a,i}}{\pd t} + \frac{\pd}{\pd x_k}
		\big( u_{a,i} v_{a,k} }  - c_a \sum_{b=1}^{\Nph} \crho_b \bar{w}_{b,k} \omega_{b,k,i} + 
		\crho_a \bar{w}_{a,k} 
		\omega_{a,k,i} 
		\label{eqn.BN.NumMeth.lift}\\[1mm]
		&\phantom{\frac{\pd u_{a,i}}{\pd t} + \frac{\pd}{\pd x_k}
		\big( u_{a,i} v_{a,k} } + c_a \sum_{b=1}^{\Nph}\lambda_b \crho_b  \left( v_{b,i} - V_i 
		\right) - \lambda_a 
		\crho_a 
		\left( v_{a,i} - V_i \right)  ,
		\\[1mm]
		&\frac{\pd s_a}{\pd t} + V_{k}\frac{\pd s_{a} }{\pd x_k} =  \frac{\lambda_a}{T_a} c_a 
		\left( v_{a,k} - V_{k} \right)^2, \label{eqn.BN.NumMeth.s.a}
		\\[1mm]
		&\frac{\pd \alpha_a}{\pd t} + V_{k}\frac{\pd \alpha_a }{\pd x_k} =  \ 0,\\[1mm]
		&\frac{\pd A_{a,Jk}}{\pd t} + \frac{\pd \left( A_{a,Jl} v_{a,l} \right) }{\pd x_k} + v_{a,i}\left(\frac{\pd  A_{a,Jk} }{\pd x_{i}} - \frac{\pd  A_{a,Ji} }{\pd x_k} \right)=  Z_{a,Jk}.\label{eqn.BN.NumMeth.A}		
	\end{align}
\end{subequations}

In order to simplify notation for discussing the numerical method, we introduce 
a compact matrix-vector notation so that system \eqref{eqn.BN.NumMeth} can be written as
\begin{equation}
\label{eqn.Matrix.Vector1}
	\pd_t \QQ + \nabla \cdot \F(\QQ) + \B(\QQ)\cdot \nabla\QQ = \S(\QQ)
\end{equation}
with $\QQ = \{\QQ_1,\QQ_2,\QQ_3\}$ being the vector of conservative state variables, and  $\QQ_a$ being the conservative variables for each phase $a = 1, 2, 3$:
\begin{equation}\label{eqn.cons.var.vect}
	\QQ_a = (\crho_a, \uu_a,  s_a, \alpha_a, \A_{a,1}, \A_{a,2}, \A_{a,3} ),
\end{equation}
where we use the fact that each distortion matrix $\A_a$ is in fact a triad of three basis vectors, i.e. $\A_a=\{\A_{a,1},\A_{a,2},\A_{a,3}\}$ and for each $J=1,2,3$,  $\A_{a,J} = \{A_{a,J1}, A_{a,J2}, A_{a,J3}\}$ is a 3-vector. 

The flux tensor $\F(\QQ)$ in \eqref{eqn.Matrix.Vector1} is decoupled with
respect to the phases, i.e. $\F(\QQ) = \{\F_1(\QQ_1),\F_2(\QQ_2),\F_3(\QQ_3)\}$,
and can be written as the sum of several contributions as follows
\begin{equation}
	 \F_a(\QQ_a) = \F^\mathrm{c}_a(\QQ_a) + \F^\mathrm{p}_a(\QQ_a) + \F^\mathrm{s}_a(\QQ_a) + \mathbf{F}^\mathrm{d}_a(\QQ_a) ,
\end{equation}
where each term depends only on $\QQ_a$ and is defined as
\begin{align}
\label{eqn.cons.fluxes0}
\F^\mathrm{c}_a(\QQ_a)  = \begin{pmatrix} \uu_a \\ \uu_a \otimes \vv_a \\ \bm{0}_{3\times1} \\ \bm{0}_{3\times1} \\ \bm{0}_{3\times1} \\ \bm{0}_{3\times1} \\ \bm{0}_{3\times1}\end{pmatrix}, 
\quad
\F^\mathrm{p}_a(\QQ_a)  = \begin{pmatrix} \bm{0}_{3\times1} \\ P_a \Id \\ \bm{0}_{3\times1} \\ \bm{0}_{3\times1} \\ \bm{0}_{3\times1} \\ \bm{0}_{3\times1} \\ \bm{0}_{3\times1} \end{pmatrix}, 
\quad
\F^\mathrm{s}_a(\QQ_a)  = \begin{pmatrix}  \bm{0}_{3\times1} \\ \bm{\sigma}^\mathrm{e}_a \\  \bm{0}_{3\times1} \\  \bm{0}_{3\times1} \\  \bm{0}_{3\times1} \\ \bm{0}_{3\times1} \\ \bm{0}_{3\times1} \end{pmatrix},
\quad
\F^\mathrm{d}_a(\QQ_a)  = \begin{pmatrix}  \bm{0}_{3\times1} \\  \bm{0}_{3\times3} \\  \bm{0}_{3\times1} \\  \bm{0}_{3\times1} \\ \A_{a,1} \vv_a \Id \\ \A_{a,2} \vv_a \Id \\ \A_{a,3} \vv_a \Id \end{pmatrix},
\end{align}
Thus, tensor $\F^\mathrm{c}_a(\QQ_a)$ contains the convective terms for the mass and momentum balance equations, $\F^\mathrm{p}_A(\QQ_a)$ is the phase related pressure flux tensor, $\F^\mathrm{s}_a(\QQ_a)$  is the flux tensor containing contribution due to shear viscous and elastic stresses, while $\F^\mathrm{d}_a(\QQ_a)$ contains advective terms from the distortion PDE. 

The so-called non-conservative matrix-vector product in
\eqref{eqn.Matrix.Vector1} contains the phase coupling terms and can be
presented as a sum of the following contributions
\begin{equation}
	\B(\QQ) \cdot \nabla\QQ = \left[ \B^\mathrm{c}(\QQ) + \B^\mathrm{p}(\QQ) + \B^\mathrm{w}(\QQ) + \B^\mathrm{d}(\QQ) \right] \cdot \nabla\QQ,
\end{equation}
where the convective part is given by (components restricted to phase $a$)
\begin{equation}
\left(\B^\mathrm{c}(\QQ)\cdot \nabla\QQ\right)_a  = ( 0, \bm{0}, \VV\cdot\nabla s_a, \VV\cdot\nabla \alpha_a, \bm{0}, \bm{0}, \bm{0}  ),
\end{equation}
while the non-conservative products related to the multiphase and multi-material nature of the model read
\begin{equation}
\left(\B^\mathrm{p}(\QQ) \cdot \nabla\QQ\right)_a  = \left( 0 ,  c_a \sum_{b=1}^{\Nph}\ps_b \nabla 
\alpha_b - \ps_a \nabla \alpha_a , 0 , 0 , \bm{0} , \bm{0} , \bm{0} \right), 
\end{equation}
\begin{equation}
	\left(\B^\mathrm{w}(\QQ) \cdot \nabla\QQ\right)_a =  \left( 0,  c_a \sum_{b=1}^{\Nph} \crho_b (\nabla \vv_b 
	{-}\nabla \vv_b^\transpose) ( \vv_b - \VV) - \crho_a (\nabla \vv_a - \nabla \vv_a^\transpose) ( \vv_a {-}\VV)  , 0 , 0 , \bm{0} , \bm{0}, \bm{0}\right),
\end{equation}
\begin{equation}
	\left(\B^\mathrm{d}(\QQ) \cdot \nabla\QQ\right)_a  = \left( 0 ,  \bm{0} , 0 , 0 , \left( \nabla \A_{a,1} - \nabla \A_{a,1}^\transpose \right) \vv_a ,
	\left( \nabla \A_{a,2} - \nabla \A_{a,2}^\transpose \right) \vv_a ,
	\left( \nabla \A_{a,3} - \nabla \A_{a,3}^\transpose \right) \vv_a \right)  .
\end{equation}

The source term vector $ \S(\QQ)$ can be written as the sum of two different vectors that will be discretized by two different approaches, 
\begin{align}
 \S(\QQ) = \S^\mathrm{w}(\QQ)  + \S^\mathrm{s}(\QQ),
\end{align}
where $\S^\mathrm{v}(\QQ)$ is stiff but linear in $\QQ$ (relative velocity
relaxation), while $\S^\mathrm{s}(\QQ)$ is related to the strain relaxation
source terms of the distortion matrix, which is nonlinear and can be extremely
stiff. These source term vectors read
\begin{align}
\label{eqn.Source.vect1}
\S^\mathrm{w}(\QQ_a) = \begin{pmatrix} 0 \\  c_a \sum_{b=1}^{\Nph}\lambda_b \crho_b  \left( 
\vv_{b} {-} \VV \right) {-} \lambda_a 
		\crho_a \left( \vv_{a} {-} \VV \right) \\ 0 \\  \frac{\lambda_a}{T_a} c_a 
		\left( \vv_{a} {-} \VV \right)^2 \\ \bm{0} \\ \bm{0} \\ \bm{0} \end{pmatrix}, \quad
\S^\mathrm{s}(\QQ_a) = \begin{pmatrix} 0 \\ 0 \\ 0 \\ 0 \\ \mathbf{Z}_{a,1} \\ \mathbf{Z}_{a,2} \\ \mathbf{Z}_{a,3}\end{pmatrix},
\end{align}
where $\bm{Z}_{a} = \{\bm{Z}_{a,1},\bm{Z}_{a,2},\bm{Z}_{a,3}\}$ is the phase strain relaxation matrix defined in \eqref{eqn.rhoE.Aa}.

Furthermore, to better address some of the specific issues usually encountered
in the numerical solution of multiphase flow models, it is convenient to
introduce the so-called primitive variables. For example, in the multiphase
context and high-order FV schemes, the use of only conservative variables, such
as $\crho_a$ and $\uu_a$, may result in non-physical discontinuities in the
reconstructed velocity and density fields, as well as violations of the
positivity constraint in the reconstructed mass fraction values
\cite{SHTCSurfaceTension}. Whereas, a reconstruction in the primitive variable
space, for a second-order MUSCL-Hancock TVD scheme, significantly mitigates these
problems, see e.g.  \cite{MunzPrimRec,Zanotti2016,SHTCSurfaceTension} and references therein. 

Therefore, alongside with the vector of conservative variables $\QQ = (\QQ_a,\QQ_2,\QQ_3)$ we consider the vector of primitive variables $\bfV = (\bfV_1,\bfV_2,\bfV_3)$, where for each phase $a = 1,2,3$, 
\begin{equation}\label{eqn.prim.vect0}
	\QQ_a = (\crho_a, \uu_a,  s_a, \alpha_a, \A_a ),
	\qquad
	\bfV_a = (\rho_a, \vv_{a},  p_a, \alpha_a, \A_a ),
\end{equation}

The \textit{primitive-to-conservative} transformation operator will be denoted
by $\mathcal{C}$ and its complementary \textit{conservative-to-primi\-ti\-ve} by
$\mathcal{P}$, i.e.
\begin{equation}
	\bfV_{a}(x,y) = \mathcal{P}[\QQ_{a}(x,y)], \qquad \mathrm{and} \qquad {\QQ}_{a}(x,y) = \mathcal{C}[\bfV_{a}(x,y)].
\end{equation}
In the mixture context, these operators must be defined with care to avoid
division by zero when a phase vanishes. In the following, we illustrate how
these conversion operators are defined in our numerical method to address this
issue and to satisfy the unit sum constraints on the volume fractions. 

First, a sum of the volume
fractions over the phases is evaluated
\begin{equation}
	\alpha_{tot} = \sum^\Nph_{a=1} \alpha_a, 
\end{equation}
from which a preliminary phase volume fraction is computed as 
\begin{equation}
	\alpha_a^* = \max( \epsilon, \mathrm{min} \left(1, \alpha_a /\alpha_{tot}^*) \right),
\end{equation}
where $\epsilon = 10^{-14}$ is a small constant introduced to avoid division by
zero in the following formula.  
Then the conservative-to-primitive $\mathcal{P}$ operator reads
\begin{equation}
	\mathcal{P}[\QQ_{a}] =  
	\begin{pmatrix} 
		\crho_a/\alpha_a^* \\ 
		\uu_a \crho_a/ ( \crho_a^2 + \epsilon^2 ) \\ 
		P[\crho_a/\alpha_a^*, s_a] \\ 
		\max( \epsilon, \mathrm{min} \left(1, \alpha_a /\alpha_{tot}^*) \right) \\ 
		\A_a 
	\end{pmatrix},
\end{equation}
where $P[\crho_a/\alpha_a^*, s_a]$ is the pressure function that can be defined
according to the EOS chosen for the phase.\\

\section{Explicit FV scheme for compressible multiphase fluid and solid mechanics} \label{SHTC.BN.num.meth}

In this section, we describe the way we adapt the well-known Finite-Volume (FV)
MUSCL-Hancock method \cite{leer5,torobook} for addressing the challenges encountered while solving the
BN form of the SHTC multifluid model \ref{sub.sec.reduced.model}. The multiphase
system is restricted to at most three phases. However, they can be freely chosen
as gaseous, viscous or inviscid liquid, or elastoplastic solid, e.g. the
multiphase system can be a gas-liquid-solid system or gas-solid-solid system,
etc. The phase interfaces are treated in the spirit of the \textit{diffuse
interface approach}, thanks to the suitable mathematical description of the
interfaces through the volume fractions $\alpha_a$. 

Due to the presence of source terms in \eqref{eqn.BN.NumMeth}, the MUSCL-Hancock
scheme will be implemented in an \textit{operator splitting} manner. Moreover,
to address the stiff character of the relaxation source terms, a specially
designed implicit discretization of the sources is incorporated at both the
predictor and corrector stages of the scheme. To address the presence of the
non-conservative product terms in \eqref{eqn.BN.NumMeth}, a path-conservative
variant of the MUSCL-Hancock scheme is employed.


%
%
\subsection{Data representation and reconstruction, slope limiting}
\label{sec.reconstruction.limiting}
The computational domain $\Omega \subset \mathbb{R}^2$ is partitioned in Cartesian elements, denoted by
\begin{equation}
\Omega_{ij} = \left[x_i - \frac{\Delta x_i}{2}, x_i + \frac{\Delta x_i}{2} \right] \times \left[y_j - \frac{\Delta y_j}{2}, y_j + \frac{\Delta y_j}{2} \right],
\end{equation}
where the indices $i$ and $j$ go from 1 to the total number of elements in each
direction. From now on, to avoid possible confusion of spatial and
discretization indices, we shall use $x$ for the direction $x_1$, and $y$ for
the direction $x_2$.

The discrete solution of the conservative and primitive state variables for a generic element $\Omega_{ij}$ at time $t^n$ is denoted by $\QQ_{i,j}^n$ and $\bfV_{ij}^n$, and which are defined as volume (area) averaged values, i.e.
\begin{equation}\label{eqn.average}
	\QQ_{i,j}^n = \frac{1}{|\Omega_{ij}|}\int_{\Omega_{ij} }\QQ(x,y,t^n)d\Omega,
	\qquad
	\bfV_{i,j}^n = \frac{1}{|\Omega_{ij}|}\int_{\Omega_{ij} }\bfV(x,y,t^n)d\Omega.
\end{equation}
 In order to achieve a second-order accuracy, it is necessary to perform a data
 reconstruction that, for each cell, yields a first-degree polynomial
 representation of the state variables, named $\QQ^r_{ij}(x, y, t)$ and $\bfV^r_{ij}(x, y, t)$.

A cell-local primitive variable polynomial reconstruction $\bfV^r_{ij}(x,
y, t)$ is now reconstructed, for each element
$\Omega_{ij}$, from the discrete primitive state vector $\bfV_{ij}^n$ at
time $t^n$. For each Cartesian cell $\Omega_{ij}$, a jump in primitive variables
through each edge can be evaluated. These are then combined in a non-linear
fashion in order to obtain a slope in the $x$ and $y$-direction respectively and
to guarantee non-oscillatory properties of the resulting scheme. For instance,
in the $x$-direction, left and right jumps are evaluated as 
\begin{equation}
\Delta\bfV_L = \bfV_{i,j}^n - \bfV_{i-1,j}^n \qquad \mathrm{and} \qquad \Delta\bfV_R = \bfV_{i+1,j}^n - \bfV_{i,j}^n,
\end{equation}
respectively. These are then combined in a non-linear fashion to obtain a
preliminary slope $\widetilde{\Delta \bfV}_i$ by means of a slope limiter. In
our implmentation,  we use a limiter that is usually referred to as the
\textit{Generalised minmod slope limiter}, and is given by
\begin{equation}
	\widetilde{\Delta\bfV}_i = \frac{\Delta\bfV_R \ \max\left[\  0, \ \mathrm{min} \left(\ \beta\Delta\bfV_R ^2, \ \Delta\bfV_R\Delta\bfV_L \right)\right]} { 2 \ \Delta\bfV_R^2 + \epsilon^2 } \\
	+ \frac{\Delta\bfV_L \ \max\left[\  0, \ \mathrm{min} \left(\ \beta\Delta\bfV_L ^2, \ \Delta\bfV_L\Delta\bfV_R \right)\right]} { 2 \ \Delta\bfV_L^2 + \epsilon^2 },
\end{equation}
where $\epsilon$ is a small constant that avoids division by zero, e.g.
$\epsilon = 10^{-14}$, and $\beta$ defines a family of \textit{minmod} limiters.
For $\beta = 1$, the classic minmod slope limiter is obtained, whereas it
reduces to the MUSCL-Barth-Jespersen limiter for $\beta = 3$, and $\beta = 2$
represents a good compromise between robustness and dissipation, and this value
will be assumed for all the subsequent numerical tests.

In fact the slope $\widetilde{\Delta \bfV}_i$ is not the final slope but a
\textit{preliminary} one since it was found useful to adopt a slope
\textit{rescaling} approach presented in \cite{Chiocchetti2023}. Thus, after
$\widetilde{\Delta \bfV}_i$ is computed it is then corrected to impose an upper
or lower limit for certain variables; in this way, the positivity of the
reconstructed density values is guaranteed and the upper and lower bounds of the
volume fractions of the phases are respected. We list below the steps to achieve
this rescaling
\begin{equation}\label{eq.slope.limited.rescaled}
	\Delta\bfV_i =\widetilde{\Delta \bfV}_i \mathrm{min} \left(	1, \Phi^+_i, \Phi^-_i	\right),
\end{equation}
with
\begin{align}
\begin{split}
	 \Phi^+_i =  \frac{ \left[ \left( \mid \widetilde{\Delta \bfV}_i \mid {+}  \widetilde{\Delta \bfV}_i		\right) ( \bfV_{M} - \bfV_i) + \left( \mid \widetilde{\Delta \bfV}_i \mid -  \widetilde{\Delta \bfV}_i		\right) ( \bfV_{m} - \bfV_i) 
	 \right]   \widetilde{\Delta \bfV}_i }
	 {2 \mid \widetilde{\Delta \bfV}_i^3 \mid + \epsilon^3}, \\
	 \Phi^-_i =  \frac{ \left[ \left( \mid \widetilde{\Delta \bfV}_i \mid -  \widetilde{\Delta \bfV}_i		\right) ( \bfV_{i} - \bfV_{M}) + \left( \mid \widetilde{\Delta \bfV}_i \mid +  \widetilde{\Delta \bfV}_i		\right) ( \bfV_{i} -  \bfV_{m}) 
	 \right]   \widetilde{\Delta \bfV}_i }
	 {2 \mid \widetilde{\Delta \bfV}_i^3 \mid + \epsilon^3},
\end{split}
\end{align}
where, the vectors $\bfV_{m}$ and $\bfV_{M}$ represent the lower and upper bounds for each variable of the primitive state vector, and are set, for each phase $a=1,2,3$, as
\begin{equation}
	\bfV_{m,a} = (0, -\bm{h}, -H, 0, -\bm{H})_a, \quad \mathrm{and} \quad \bfV_{M,a} = (H,\bm{h}, H, 1, \bm{H})_a,
\end{equation}
where the values of $H$, $\bm{h}$ and $\bm{H}$ should be set to certain bound
values following from the physical meaning of the corresponding quantity or to
arbitrary large value to represent the absence of bounds. The analogous steps
are performed in $y$-direction to calculate the reconstruction of the slope
$\Delta\bfV_j$.

After the spatial reconstruction, at a given time instant $t^n$, the cell-local
space-time primitive variable polynomial reconstruction, in each cell
$\Omega_{ij}$, is written in terms of a space-time Taylor series expanded about
$x_i$, $y_i$ and $t^n$ as
\begin{equation}
\label{eqn.spacetime.predict}
	\bfV^r_{i,j}(x, y, t) =  \bfV^n_{i,j} + (x- x_{i,j}) \frac{\Delta\bfV_i}{\Delta x} + (y- y_{i,j}) \frac{\Delta\bfV_j}{\Delta y} + (t-t^n) \ \pd_t\bfV_{i,j}.
\end{equation}
The time derivative, in \eqref{eqn.spacetime.predict}, is computed in terms of
primitive variables in two steps, through the following straight-forward
application of the operator splitting approach. 

Thus, to determine $\pd_t\bfV_{i,j}$, we consider 
\begin{equation} \label{eq.predictor12}
	\pd_t \bfV_{i,j} = - \ (\mathbf{C} \cdot\hat{\boldsymbol{n}}_x) \ \pd_x\bfV  - \ (\mathbf{C} \cdot\hat{\boldsymbol{n}}_y) \ \pd_y\bfV + \S,
\end{equation}
where 
\begin{equation}\label{eq.C0}
	\mathbf{C} = \left( \frac{\pd \QQ}{ \pd \bfV}\right)^{-1} \left( \frac{\pd  \F(\QQ)}{ \pd \bfV} +   \B(\QQ)\frac{\pd \QQ}{ \pd \bfV} \right) \qquad 
	\textnormal{ and } \qquad 
	\tilde{\mathbf{B}} =  \B(\QQ) \frac{\pd \QQ}{ \pd \bfV}.
\end{equation}
To solve \eqref{eq.predictor12}, we split it into the homogeneous part 
\begin{equation} \label{eq.predictor123}
	\pd_t \bfV_{i,j} = - \ (\mathbf{C} \cdot\hat{\boldsymbol{n}}_x) \ \pd_x\bfV  - \ (\mathbf{C} \cdot\hat{\boldsymbol{n}}_y) \ \pd_y\bfV,
\end{equation}
and the source part  
 \begin{equation}
	\frac{\rmd  \bfV_{i,j} }{\rmd \hat{t}} = \S(\bfV_{i,j}),
	\qquad
	\hat{t} \in[t^n, t].
\end{equation}
In order to approximate the spatial derivatives of the primitive state variables
in \eqref{eq.predictor123}, we use a central finite difference with respect to
the cell center  by using the boundary primitive reconstructed values from
within the cell $\Omega_{i,j}$, as
\begin{equation} \label{eq.predictor1234}
	\pd_t \bfV_{i,j} \approx \left( \mathbf{C}(\bfV_{i,j}^n) \cdot\hat{\boldsymbol{n}}_x 		\right)  \frac{\Delta\bfV_i}{\Delta x}  + \left( \mathbf{C}(\bfV_{i,j}^n) \cdot\hat{\boldsymbol{n}}_y 		\right)  \frac{\Delta\bfV_j}{\Delta y}.
\end{equation}
Hence, using \eqref{eq.predictor1234}, one can compute an update for each cell, such that $\bfV^{\hspace{+0.05em}*}_{i,j}$ is the solution of \eqref{eq.predictor123} at time $t$ with an initial value $\bfV^n_{i,j}$, as follows
\begin{equation}
	\bfV^{\hspace{+0.05em}*}_{i,j} = \bfV^n_{i,j} + (t-t^n)\ \pd_t\bfV_{i,j}.
\end{equation}
In the second step, we consider the contribution of the stiff source terms via
solving the initial value problem
 \begin{equation}\label{ini.value.probl}
	\frac{d  \bfV_{i,j} }{d \hat{t}} = \S(\bfV_{i,j}), \qquad \hat{t} \in[t^n, t], \qquad \bfV_{i,j}(t^n) = \bfV^{\hspace{+0.05em}*}_{i,j},
\end{equation}
whose solution at time $t$ is denoted by $\bfV_{i,j}^{**}$. This initial value problem is solved with two different implicit methods discussed in detail in Sec.\,\ref{section.sources}.

Finally, by introducing the discrete solution $\bfV_{i,j}^{**}$ of $\pd_t\bfV_{i,j}$ in \eqref{eqn.spacetime.predict}, the cell-local primitive variable polynomial reconstruction reads
\begin{equation}
\label{eqn.local.space.time.pred}
	\bfV^r_{i,j}(t, x, y) =  \bfV_{i,j}^{**}(t) + (x- x_{i,j}) \frac{\Delta\bfV_i}{\Delta x} + (y- y_{i,j}) \frac{\Delta\bfV_j}{\Delta y},
\end{equation}
which we also refer to as the \textit{cell-local space-time predictor}.
\subsection{Explicit finite volume discretization of the homogeneous system}
\label{sec.FV.discretization} 
After obtaining the local space-time predictor \eqref{eqn.local.space.time.pred}, the final solution $\QQ_{i,j}^{n+1}$ at $t^{n+1}$  of the MUSCL-Hancock scheme is also computed using the splitting approach, in which we first compute the solution of the homogeneous PDE system 
\begin{equation}
	\label{eqn.Matrix.Vector.homo}
	\pd_t \QQ+ \nabla \cdot \F(\QQ) + \B(\QQ)\cdot \nabla\QQ = \bm{0},
\end{equation}
with the initial data obtained by extrapolating the reconstructed polynomials
towards the cell boundaries and by applying the standard explicit FV update
formula to \eqref{eqn.Matrix.Vector.homo}. The latter is obtained by integrating
\eqref{eqn.Matrix.Vector.homo} over the space-time element and applying Gauss’s
theorem for integrating the divergence of fluxes in space:
\begin{equation}
	\int^{t^{n+1}}_{t^{n}} \int_{\Omega_{ij}}   \pd_t \QQ  \rmd\xx  \rmd t + \int^{t^{n+1}}_{t^{n}} \int_{\pd\Omega_{ij}}  \F(\QQ)\cdot \hat{\mathbf{n}}  \rmd S  \rmd t
	+ \int^{t^{n+1}}_{t^{n}}     \int_{\Omega_{ij}}     \B(\QQ)\cdot \nabla\QQ \ \rmd\xx \rmd t = 0,
\end{equation}
where $\hat{\mathbf{n}} $ defines the outward unit normal vector on the element boundary, and $\rmd \xx= \rmd x \rmd y$.

Then, by using the reconstructed polynomials $\bfV^r_{i,j}(t, x, y)$ and
treating the non-conservative terms using the path-conservative approach by Castro and Parés \cite{Pares2006,Castro2006}, we get the usual
\textit{path-conservative} FV discretization
\begin{align}
	{\int^{t^{n+1}}_{t^{n}}}    \int_{\Omega_{ij}}   \pd_t \QQ d\xx \rmd t &+ \int^{t^{n+1}}_{t^{n}}    \int_{\pd\Omega_{ij}}    \bigg(  \F (\bfV^{r,-}_{i,j},\bfV^{r,+}_{i,j} ) + \mathbf{D}(\bfV^{r,-}_{i,j},\bfV^{r,+}_{i,j} ) \bigg) \cdot \hat{\mathbf{n}}  \rmd S  \rmd t \ +\nonumber\\ 
	&+ \int^{t^{n+1}}_{t^{n}}    \int_{\Omega_{ij}\textbackslash \pd\Omega_{ij}}    \tilde{\B}(\bfV^r_{i,j})\cdot \nabla\bfV^r_{i,j} \ d\xx \rmd t = 0,
\end{align}
where, within the framework of path-conservative schemes, the new term
$\mathbf{D}$ was introduced to take into account the jumps of the primitive
variables $\bfV$ across the space-time element boundaries $\pd\Omega_{ij}$,
while the last term is the integral over the smooth part of the non-conservative
terms.

Using notations \eqref{eqn.average}, the fully discrete one-step update formula
for the solution $\QQ_{i,j}^{(1)}$ of the homogeneous part of the system at time
$t^{n+1}$ reads
\begin{align}
\label{eqn.discrete.update.formula}
\begin{split}
	\QQ_{i,j}^{(1)} = \QQ_{i,j}^{n} &- \frac{\Delta t}{\Delta x} \left( \ \F^{\mathrm{RS}}_{i+1/2,j} - \F^{\mathrm{RS}}_{i-1/2,j} + 	\mathbf{D}_{i+1/2,j} + \mathbf{D}_{i-1/2,j}	\right) + \\
	&- \frac{\Delta t}{\Delta y} \left( \ \F^{\mathrm{RS}}_{i,j+1/2} - \F^{\mathrm{RS}}_{i,j-1/2} + 	\mathbf{D}_{i,j+1/2} + \mathbf{D}_{i,j-1/2}	\right) +\\
	&+ \frac{\Delta t}{\Delta x} \tilde{\B}_1 \left[  \bfV^r_{i,j} \left(t^{\ n + 1/2}, x_{i}, y_j\right) 		\right] \Delta\bfV_i + \\
	&+ \frac{\Delta t}{\Delta y} \tilde{\B}_2 \left[  \bfV^r_{i,j} \left(t^{\ n + 1/2}, x_{i}, y_j\right) 		\right] \Delta\bfV_j,
\end{split}
\end{align}
where $\F^{\mathrm{RS}}$ is the generic conservative numerical flux, which can be compute with different approximate Riemann solvers. 

In order to describe easily each term in \eqref{eqn.discrete.update.formula}, we introduce a compact notation for the boundary-extrapolated primitive states $\bfV_R$ and $\bfV_L$, which can be evaluated from the solution of the cell-local space-time predictor \eqref{eqn.local.space.time.pred}. In particular, the space-time midpoint values for each face, of generic index $i + \frac{1}{2}, j$ in the $x$-direction or $i, j + \frac{1}{2}$ in the $y$-direction, read
 \begin{align}
 \begin{split} 
&(\bfV_L)_{i+\halb,j} = \bfV^r_{i,j}(t^{n + \halb}, x_{i+\halb}, y_j), \quad  (\bfV_R)_{i+\halb,j} = \bfV^r_{i+1,j}(t^{n + \halb}, x_{i+\halb}, y_j), \\
&(\bfV_L)_{i,j+\halb} = \bfV^r_{i,j}(t^{n + \halb}, x_{i}, y_{j+\halb}),\quad 
(\bfV_R)_{i,j+\halb} = \bfV^r_{i,j+1}(t^{n + \halb}, x_{i}, y_{j+\halb}).
\end{split} 
\end{align}
Thus, using this simpler notation, we illustrate how conservative numerical
$\F^{\mathrm{RS}}$ flows are defined. In this paper, we employ the simple Rusanov
flux:
\begin{align}
\begin{split}
\label{eqn.num.flux}
	\F^{\mathrm{RS}}_{i+1/2,j} (\bfV_L, \bfV_R) = \halb \Big(\F_1(\bfV_L) + \F_1(\bfV_R) \Big) - \halb s_1^{\max} \Big( \mathcal{C}[\bfV_R] - \mathcal{C}[\bfV_L] \Big),		\\
	\F^{\mathrm{RS}}_{i,j+1/2} (\bfV_L, \bfV_R) = \halb \Big(\F_2(\bfV_L) + \F_2(\bfV_R) \Big) - \halb s_2^{\max} \Big( \mathcal{C}[\bfV_R] - \mathcal{C}[\bfV_L] \Big),
\end{split}
\end{align}
where $\F_1$ and $\F_2$ are the conservative fluxes in the first and in the second space direction.

The Rusanov numerical flux requires the knowledge of an estimate for the maximum
wave velocity $s^{\max}$ for each direction. In this paper, keeping in mind that
we are interested in problems with not high Mach numbers, the absolute value of
the maximum eigenvalue of the PDE system  linearized at the states $\bfV_L$ and
$\bfV_R$ can be a good estimate for $s^{\max}$. Therefore, the maximum wave
speed estimates read
\begin{align}
\begin{split}
	s_1^{\max} (\bfV_L, \bfV_R) = \max \left(\ \lambda_1^{\max}(\bfV_L), \  \lambda_1^{\max}(\bfV_R) \ \right), \\
	s_2^{\max} (\bfV_L, \bfV_R) = \max \left(\ \lambda_2^{\max}(\bfV_L), \  \lambda_2^{\max}(\bfV_R) \ \right).
\end{split}
\end{align}
The maximum eigenvalues of \eqref{eqn.BN.NumMeth} can be estimated as described in \ref{sec.eigenvalues}.

The nonconservative products appearing in the BN model are treated within the framework of path-conservative
schemes \cite{Castro2006,Pares2006,Castro2008,CastroPardoPares,USFORCE2}. Thus, at each cell interface the following path integrals must be
prescribed
\begin{equation}
	\mathbf{D}_\Psi (\bfV_L, \bfV_R) \cdot \hat{\mathbf{n}} 	= \halb \int_0^1 \tilde{\mathbf{B}}\left[ \Psi(\bfV_L, \bfV_R, s)		\right]\cdot \hat{\mathbf{n}} \frac{\pd \Psi}{\pd s}  \rmd s = \halb \left(\int_0^1 \tilde{\mathbf{B}}\left[ \Psi(\bfV_L, \bfV_R, s)		\right]\cdot \hat{\mathbf{n}} \rmd s \right)\big(\bfV_R - \bfV_L\big),
\end{equation}
in which $\Psi(\bfV_L, \bfV_R, s) =  \bfV_L + s (\bfV_R- \bfV_L) $ is a simple
segment path function connecting the left and right states in the primitive
state space. These path integrals, which are denoted by $\mathbf{D}_{i+1/2,j}$
and $\mathbf{D}_{i,j+1/2}$ in \eqref{eqn.discrete.update.formula}, are computed
with a three-point Gauss-Legendre quadrature rule with points $s_k \in [0,1]$ and weights $w_k$ as follows (see \cite{USFORCE2}) 
\begin{align}
\begin{split}
	&\mathbf{D}_{i+1/2,j} = \halb \sum_{k=1}^3 w_k \ \tilde{\mathbf{B}}_1 \Big[ \Psi(\bfV_L, \bfV_R, s_k) \Big] \ \big(\bfV_R - \bfV_L \big),	\\
	&\mathbf{D}_{i,j+1/2} = \halb \sum_{k=1}^3 w_k \ \tilde{\mathbf{B}}_2 \Big[ \Psi(\bfV_L, \bfV_R, s_k) \Big] \ \big(\bfV_R - \bfV_L \big)\	.
\end{split}
\end{align}
Therefore, by these means we compute the preliminary state vector $\QQ^{(1)}_{i,j} $, which is the updated solution of the left hand side of \eqref{eqn.Matrix.Vector1}. To get the final solution $\QQ^{n+1}_{i,j}$, it remains to compute the solution of the relaxation source terms, which is done in the next section.

Before describing in detail the implicit solver for the relaxation source terms, we note that in order to guarantee stability of the explicit FV time-stepping described above, the time-step size is restricted by
\begin{equation} \label{CFL.SHTC}
	\Delta t \leq k_{\scriptscriptstyle{\up{CFL}}} \frac{1}{	{\Delta x}/{\lambda^{\scriptscriptstyle{M}}_1} + {\Delta y}/{\lambda^{\scriptscriptstyle{M}}_2}	},
\end{equation}
where $\lambda^{\scriptscriptstyle{M}}_k$ is the maximum absolute value of all
eigenvalues found in the domain, in the $x_k$-direction, which, for system \eqref{eqn.BN.NumMeth}, can be estimated
as detailed in \ref{sec.eigenvalues}. Also, $k_{\scriptscriptstyle{\up{CFL}}}
\leq 1$ is a Courant-type number \cite{courant1928}, which is typically chosen
as $k_{\scriptscriptstyle{\up{CFL}}}  = 0.9$ for all the simulations presented
in this work.

\subsection{Integration of relaxation sources }\label{section.sources} 

As previously mentioned, to account for the algebraic relaxation source terms in
the numerical solution, we adopt a splitting approach. It is a simple but robust
strategy since it allows to separate the contribution of relaxation terms from
the reversible part of the time evolution equations. Here, we discuss the
details of the implicit method that is used to solve the ordinary
differential equations (ODE) 
\begin{equation}\label{eqn.ODE.source}
	\frac{\up{d}  \QQ_{i,j} }{\up{d} {t}} = \S(\QQ_{i,j}), \qquad {t} \in[t^n, t^{n+1}], \qquad \QQ_{i,j}(t^n) =\QQ^{(1)}_{i,j}.
\end{equation}

In Section \ref{sub.sec.reduced.model}, the source terms were separated into
$\S^\mathrm{w}(\QQ_a)$, which contains the velocity relaxation terms and are
stiff, but linear with respect to the relative velocity, and
$\S^\mathrm{s}(\QQ_a)$, which contains the strain relaxation terms of the
distortion matrix $\mathbf{ A}_a$, and which is non-linear and can be very
stiff. The integration of these two different source terms is carried
out with two different implicit approaches.
\subsubsection{Relative velocity relaxation}
In the following, we describe the system of ODEs arising from the source vector $\S^\mathrm{w}(\QQ)$ in \eqref{eqn.Source.vect1}, related to the relative velocity relaxation. 
Since there are zeros in $\S^\mathrm{w}(\QQ_a)$ corresponding to the
conservation equations of mass, volume fraction, and distortion matrix, they are
remained constant over time in \eqref{eqn.ODE.source}. Therefore, these
quantities can be considered as constant parameters and can be omitted from the
state vector $\QQ_{i,j}$ of the initial value problem \eqref{eqn.ODE.source}.

The integration of the remaining quantities, the phase momenta $\uu_a$ and entropies $s_a$, is carried out in terms of the primitive state variables. More precisely, instead of \eqref{eqn.ODE.source} we consider the following reduced ODE system for the vector$ \mathbf{\widetilde{V}} = \left( \vv_{1}, \vv_{2}, \vv_{3}, s_1, s_2, s_3 \right)$:
\begin{equation}\label{eqn.ODE.w.s}
	\frac{\rmd \mathbf{\widetilde{V}}_{ij}  }{\rmd t} = \S(\mathbf{\widetilde{V}}_{ij}), \qquad {t} \in[t^n, t^{n+1}],\qquad \mathbf{\widetilde{V}}_{ij}(t^n) = \mathbf{\widetilde{V}}_{ij}^{(1)},
\end{equation}
where $\mathbf{\widetilde{V}}^{(1)}$ is the primitive variable reduced state
vector obtained in \eqref{eqn.discrete.update.formula} as the solution to the
the homogeneous PDE system. The later ODE system can be easily integrated by
means of the backward Euler method to obtain an updated solution
$\mathbf{\widetilde{V}}^{(2)}$ at time $t^{n+1}$. 

For instance, it can be seen that for each spatial direction $x_k$, $k = 1,2$,
one can decouple the phase velocity equations from the phase entropy ones, and
this velocity subsystem reads
\begin{align}
\label{eqn.ODE.1}
\begin{split}
	&\frac{\rmd  v_{1,k}}{\rmd t} =  \lambda  \bigg( c_1(v_{1,k} {-} V_k) + c_2(v_{2,k} {-} V_k)  + c_3(v_{3,k} {-} V_k)  - (v_{1,k} {-} V_k) \bigg)   , \\
	&\frac{\rmd v_{2,k}}{\rmd t}  = \lambda  \bigg( c_1(v_{1,k} {-} V_k) + c_2(v_{2,k} {-} V_k)  + c_3(v_{3,k} {-} V_k)  - (v_{2,k} {-} V_k) \bigg)   , \\
	&\frac{\rmd v_{3,k}}{\rmd t}  = \lambda  \bigg( c_1(v_{1,k} {-} V_k) + c_2(v_{2,k} {-} V_k)  + c_3(v_{3,k} {-} V_k)  - (v_{3,k} {-} V_k) \bigg)   , \\
\end{split}
\end{align}
where the phase kinetic coefficients $\lambda_{ab,k}$, that define the time
scale for friction relaxation dissipative process in \eqref{eqn.vel.relax1}, are
assumed to be equal throughout all phases and directions in this paper, i.e.
$\lambda_{ab,k} = \lambda$. It is also usually assumed that this parameter is
larger than $1/\Delta t$, where $\Delta t$ is the time-step given by the
stability condition in \eqref{CFL.SHTC}, and hence, we can say that we have a
stiff ODE system, and an implicit discretization is needed. For such
a system, in which all the sources are linear, the following
discretization can be written for each cell $\Omega_{ij}$
\begin{align}\label{source.update.vel}
	\begin{pmatrix} v_{1,k} \\ v_{2,k} \\ v_{3,k}\end{pmatrix}^{(2)}_{i,j} = (\up{I}- \Delta t \mathcal{M} )^{-1}\begin{pmatrix} v_{1,k} \\ v_{2,k} \\ v_{3,k}\end{pmatrix}^{(1)}_{i,j} 
\end{align}
where the matrix inverse can be evaluated analytically, and after defining $\lambda^* = \Delta t  \lambda$, it reads
\begin{align}
(\up{I}-\Delta t \mathcal{M} )^{-1} =\frac{1}{1+\lambda^*}
	\begin{pmatrix} 1+\lambda^* c_1 	&& \lambda^* c_2 	&& \lambda^* c_3 \\
				\lambda^* c_1		&& 1+\lambda^* c_2	&& \lambda^* c_3 \\
				\lambda^* c_1		&& \lambda^* c_2	&& 1+\lambda^* c_3 \\
	\end{pmatrix}.
\end{align}

After that, it remains to solve the independent ODE subsystem for the phase
entropies that reads
\begin{align}
\label{eqn.ODE.2}
\begin{split}
	&\frac{\rmd s_{1}}{\rmd t} = \frac{\lambda  c_1}{T_1}  \bigg( \big(v_{1,1} - V_1\big)^2 + \big(v_{1,2} - V_1\big)^2+ \big(v_{1,3} - V_1\big)^2 \bigg)	, \\
	&\frac{\rmd s_{2}}{\rmd t} = \frac{\lambda  c_2}{T_2}  \bigg( \big(v_{2,1} - V_1\big)^2 + \big(v_{2,2} - V_1\big)^2+ \big(v_{2,3} - V_1\big)^2 \bigg)	, \\
	&\frac{\rmd s_{3}}{\rmd t} = \frac{\lambda  c_3}{T_3}  \bigg( \big(v_{3,1} - V_1\big)^2 + \big(v_{3,2} - V_1\big)^2+ \big(v_{3,3} - V_1\big)^2 \bigg) , 
\end{split}
\end{align}
where we can use the updated velocities evaluated in \eqref{source.update.vel}.
Finally, to discretize the three ODEs in \eqref{eqn.ODE.2}, a generic implicit
backward Euler time integrator based on Newton's method can be used in order to
deal with the nonlinearity inherent to the definitions of the phase temperatures
$T_a(\rho_a,s_a)$.
\subsubsection{Strain relaxation}
In contrast to the relative velocity relaxation subsystem \eqref{eqn.ODE.w.s},
an accurate integration of the non-linear stiff source $\S^\mathrm{s}(\QQ)$
governing the strain relaxation of the distortion matrix $\A_a$ is a more
challenging task, especially in the context of multiphase flows. Let us begin
with some remarks on the evolution equations of the phase distortion fields
$\A_a$, which we recall to be defined for each phase $a = 1,2,3$. 

In the multiphase context, the evolution of the three distortion fields $\A_a$,
$a=,1,2,3$ given by 
\begin{equation}
\label{eqn.distortionMatrix}
	\pd_t\A_a + \nabla \big( \A_a \cdot \vv_a \big) + \big(\nabla \A_a - \nabla \A_a^{\transpose} \big) \cdot \vv_a = - \frac{3}{\tau^\mathrm{e}_a} \left( \mathrm{det} \A_a\right)^{5/3} \A_a  \mathrm{dev}  \big( \A^{\transpose}_a  \A_a \big).
\end{equation}
may occur over a very wide range of time scales in a single
computational cell $\Omega_{ij}$. Namely, there might be infinitely slow strain
relaxation time scale ($\tau^\mathrm{e}_a = 10^{14}$) in an elastic solid phase
and extremely fast relaxation of shear stresses in the inviscid
($\tau^\mathrm{e}_a=10^{-14}$) and viscous fluid phase
($\tau^\mathrm{e}_a\sim10^{-6}-10^{-3}$). These different time scales are
quantified by means of the relaxation time $\tau^\mathrm{e}_a$ in the evolution
equation of the phase distortion field

The interpretation of the strain relaxation timescale $\tau^\mathrm{e}_A$ and
its definition in the multiphase context were described in
Section\,\ref{sec.multiA}. From that description, it is clear that one of the
major difficulties in solving the unified multiphase model of continuum
mechanics is conditioned by the presence of these stiff and very non-linear
strain relaxation source terms. Therefore, it is necessary to solve the
associated ODE systems with care using an appropriate implicit time integrator.

Following the ideas in \cite{GodRom2003}, an efficient and robust method for a
semi-analytical implicit integration of the strain relaxation ODE systems was
introduced by Chiocchetti and co-authors in \cite{ExpInt1} in the context of
strain relaxation in the damaged solids, and further developed in \cite{ExpInt2,Chiocchetti2023} for finite-rate pressure and strain relaxation in multiphase
flows. The \textit{key idea} of this time integrator is a reduction of the
problem by using the \textit{polar decomposition} for each phase distortion
matrix $\A_a$:
\begin{align}\label{ev.eq.A.source}
	\A_a = \boldsymbol{R}_a \ \GG_a^{1/2} \qquad \mathrm{with} \qquad \GG_a^{1/2} = \boldsymbol{E}_a \  \hat{\GG}^{1/2}_a \ \boldsymbol{E}_a^{-1},
\end{align}
where $\boldsymbol{R}_a$ is an orthogonal matrix with a positive unit
determinant, while the matrix square root $\GG_a^{1/2}$ can be defined by means
of eigen-decomposition of the symmetric positive definite matrix $\GG_a$, where $\boldsymbol{E}_a$ is the matrix whose
columns are eigenvectors and $\hat{\GG}^{1/2}_a$ is the diagonal matrix whose
diagonal elements are the roots of the eigenvalues.

Indeed, the distortion field $\A_a$ represents three local basis vectors
representing the volume, shape, and the orientation of the phase control volume.
Its 9 independent components (degrees of freedom), therefore, encodes two
different types of information. Six degrees of freedom are strictly related to
the definition of the stress tensor $\bm{\sigma}^\mathrm{e}_a = \crho_a \Cs_a^2
( \GG_a \mathrm{dev} \GG_a)$ via the six independent components of the metric
tensor $\GG_a$, and the three remaining degrees of freedom that define the
angular orientation of the control volume.

 Numerically, the matrix $\GG_a^{1/2}$ can be simply evaluated using the
 Denman-Beavers algorithm. Thus, for any given state $\A_a$, one can easily
 compute $\GG_a$, its square root $\GG_a^{1/2}$, and eventually the inverse
 $\GG_a^{-1/2}$. After that, the rotation matrix can be computed as
\begin{align}
	\boldsymbol{R}_a = \A_a \ \GG_a^{-1/2}.
\end{align}
Moreover, the invariance of the rotational component of the distortion matrix
under strain relaxation can be proven following the arguments in
\cite{GodRom2003,Chiocchetti2023}, which means that during the strain relaxation step, one can use the evolution PDE for the metric tensor  
\begin{align}
\label{eqn.Metric.Tensor}
\pd_t\GG_a +\big( \nabla  \GG_a\big) \ \vv_a  + \GG_a \nabla \vv_a - (\nabla \vv_a)^{\transpose} \GG_a = - \frac{6}{\tau^\mathrm{e}_a} \left( \det \GG_a\right)^{5/6} \GG_a \mathrm{dev} \GG_a\ ,
\end{align}
instead of the PDE for the full distortion matrix.

We now have all the ingredients to describe the following steps in order to
obtain the final solution for the distortion matrix. First, we
calculate the update $\A_a^{(1)}$ of the distortion matrix $\A_a$ obtained
from the left hand side of the evolution equation \eqref{eqn.distortionMatrix},
as presented in \eqref{eqn.discrete.update.formula}. From this solution, a rotation matrix can be calculated
independently of the non-linear source terms as 
\begin{align}
\label{eqn.rotations.homog}
	\boldsymbol{R}_a^{(1)} = \A_a^{(1)} \ \GG_a^{{(1)} \ -1/2} \quad \mathrm{with} \quad \GG_a^{(1)} = \left(\A_a^{(1)}\right)^\transpose \A_a^{(1)},
 \end{align}
 with $\GG_a^{{(1)} \ -1/2}$ computed by means of the Denman-Beavers algorithm.
 Then, the following nonlinear ODE system should be solved
\begin{align}
\label{eqn.alternative.ODE}
	\frac{\rmd  \GG_a}{\rmd t} = \boldsymbol{L}_a^{(1)}  - \frac{6}{\tau^\mathrm{e}_a} \left( \ \mathrm{det}\ \GG_a\right)^{5/6} \ \GG_a \ \mathrm{dev} \ \GG_a\ ,
\end{align}
where $ \boldsymbol{L}_a^{(1)}$ is a constant convective/productive forcing term
evaluated simply as
\begin{align}
	\boldsymbol{L}_a^{(1)} = \frac{\GG_a^{(1)}  - \GG_a^n }{\Delta t},  \quad \mathrm{with} \quad \GG_a^{(1)} = \left(\A_a^{(1)}\right)^\transpose \A_a^{{(1)}}, 
	\quad
	\GG_a^{n} =  \GG_a\left(t^n\right). 
 \end{align}
This term, which takes into account the left-hand side of
\eqref{eqn.Metric.Tensor}, is introduced to converge to the asymptotically
correct state in the stiff limit of the equations. This alternative ODE problem
\eqref{eqn.alternative.ODE} is then solved by computing the analytical solution
of a sequence of linearized problems that approximate the original non-linear
ODE, according to the procedure outlined in \cite{Chiocchetti2023}. Once the
source term applied to the metric tensor is integrated, and thus obtaining
$\GG_a^{(2)}$ at time $t^{n+1}$, the information can be mapped back to get the
updated distortion field as 
\begin{align}
\label{eqn.final.update}
	\A_a^{(2)} = \boldsymbol{R}^{(1)}_a \ \left(\GG_a^{(2)} \right)^{1/2}.
\end{align}
\subsection{Further remarks on the distortion field }\label{remarks.dist.field}
Before assembling the final solution of the entire PDE system, we must make some
important remarks concerning the challenges related to the numerical computation
of the distortion field $\A_a$. 

\subsubsection{Algebraic determinant constraint}
\label{det.constr}
In the numerical solution of the evolution equation \eqref{ev.eq.A.source} for the phase distortion field $\A_a$, particularly when describing liquid phases, one must be careful with preserve the nonlinear algebraic constraint 
\begin{align}\label{constrain.A}
	\rho_a = \rhoo_a  \mathrm{det} \A_a.
\end{align}
This constraint stems from the fact that the phase mass balance equation
\eqref{eqn.BN.NumMeth.rho.a} is the consequence of the time evolution
\eqref{ev.eq.A.source} for the distortion matrix $\A_a$, see e.g. \cite{God1978,
Pesh2014}. Yet, in an Eulerian scheme, it is important to explicitly discretize
the mass balance equation to ensure that it is fulfilled at the discrete level
(important for the Rankine–Hugoniot conditions). Thus, the phase density can be
computed by two means: from the mass balance equations and from
\eqref{constrain.A}. Therefore, to ensure that the solution is consistent with
\eqref{constrain.A}, this constraint must be actively enforced in the numerical
scheme.

A simple but effective approach is to manually impose the constraint at each time iteration. Specifically, the distortion field obtained after solving the homogeneous problem is enforced to satisfy the ratio
\begin{align}\label{constrain.A1}
	 \mathrm{det} \A_a^{(1)} = \frac{\rho_a^{(1)}}{\rhoo_a }
\end{align}
as detailed in \cite{Chiocchetti2023}.
\subsubsection{Linear combination of pure rotational fields} 

Numerical discretization, in all its parts from data representation to explicit
discretization with the FV update formula, is applied to the vector and
tensorial quantities in a component-wise manner. This is not a problem for the
velocity field, but the distortion field requires more attention due to its
rotational component. From our experience, an improper treatment of the
rotational matrix $\bm{R}_a$ can lead to artificial stresses and other numerical
artifacts.

It is clear that even in the case of a simple component-by-component linear
combination of two rotational matrices, such as that presented by an average
operator for example, the resulting matrix is a rotational matrix only in the
case of infinitesimal rotations. However, the additivity rule does not in
general hold for finite rotations. For example, we have observed in our
numerical experiments, in which the solution of the distortion field is not
particularly smooth, that simple averaging of the rotational matrices results in
artificial stresses. This happened in tests such as the lid-driven cavity test,
Sec.\,\ref{lid.driven.distortion} if no special treatment is applied. In this
test, the boundary conditions produce a velocity gradient singularity at the
corners of the cavity, which results in a locally discontinuous distortion
field. On the other hand, no issue arises in the double shear layer problem,
Sec.\,\ref{double.shear.layer}. For this reason, a different approach would be
required to describe the rotational component of the information encoded by
$\A_a$, by means of a auxiliary mathematical representation of these quantities
that allows a component-by-component treatment even for rotations of finite
amplitude. 

Thus, in this paper, we propose a simple but effective approach to address this
issue. It employs the efficiency of the Chiocchetti semi-analytical solver
\cite{Chiocchetti2023} in the infinitely stiff relaxation regime,
$\tau^\mathrm{e}_a \to 0$. In particular, in this limit, the strains encoded in
the metric tensor $\GG_a$ dissipate almost instantaneously, resulting in a
distortion field that is represented by a pure rotation matrix, i.e.
$\A_a=\bm{R}_a$ and $\GG_a=\Id$. 

Therefore, the idea of decoupling of the evolution of the two types of
information encoded in $\A_a$ (at least numerically) is straightforward. This
can be done by adding an auxiliary evolution equation
\eqref{eqn.distortionMatrix} for a new auxiliary distortion field
$\widetilde{\A}_a$ subject to a relaxation timescale
$\tilde{\tau}^\mathrm{e}_a\to 0$. This auxiliary distortion field
$\widetilde{\A}_a$ carries only the information about the rotational component
of the original distortion $\A_a$, i.e. $\widetilde{\A}_a = \widetilde{\bm{R}}_a$. Moreover, thanks to this almost instantaneous
relaxation (numerically we use $\tilde{\tau}^\mathrm{e}_a=10^{-14}$), the
artificial stresses that may arise from the combination of finite amplitude
rotations are dissipated instantaneously, both at the level of the predictor,
Sec.\,\ref{sec.reconstruction.limiting}, and at the level of the final solution
at each time-step, Sec.\,\ref{sec.FV.discretization} and
Sec.\,\ref{section.sources}.

On the other hand, due to this decoupling, the original distortion field $\A_a$
at each time step carries only the information about the strains, i.e. $\A_a =
\GG_a^{1/2}$ and $ \boldsymbol{R}_a = \Id$, both at the level of the
predictor and the solution updated at time $t^{n+1}$. Therefore, expression
\eqref{eqn.final.update}, in the source term integration process for the
original distortion field $\A_a$, should be rewritten as 
\begin{align}
	\A_a^{(2)} =\Id \left(\GG_a^{(2)} \right)^{1/2}.
\end{align}
\subsection{Final solution for the complete problem  }
The final solution $\vec{Q}_{i,j}^{n+1}$ of the complete problem
\eqref{eqn.Matrix.Vector1} can now be retrieved by considering the solution of
the homogeneous problem, the contribution of the source terms, and the remarks
on rotational matrices discussed above. Thus, including formally the auxiliary
phase distortion field $\widetilde{\A}_a$ to the set of state parameters, the
final solution of each phase $a=1,2,3$ reads
\begin{align}
	\vec{Q}_{a,i,j}^{n+1} = \left(  \crho_a^{(1)}, \uu_a^{(2)},  s_a^{(2)}, \alpha_a^{(1)}, \A_{a}^{(3)}, \widetilde{\A}_{a}^{(2)}\right)_{i,j},
\end{align}
where the complete phase distortion field  $\A_{a}^{(3)}$ is computed as
\begin{align}
	\A_{a}^{(3)} = \widetilde{\boldsymbol{R}}^{(1)}_{a} \ \left(\GG_{a}^{(2)} \right)^{1/2}.
\end{align}
 \section{Numerical results}\label{sec.num.res}
This section provides the results obtained with the numerical scheme presented
in Section\,\ref{SHTC.BN.num.meth} addressing the three-phase reduced BN-type
SHTC model presented in \ref{sub.sec.reduced.model}. The numerical scheme
developed in this work considers the entropy inequalities of each phase
\eqref{eqn.BN.NumMeth.s.a}, instead of the phase energy balance laws as in the
most existing finite-volume discretizations. This choice was made in order to
make the mathematical model less complex. Indeed, the phase entropy balance laws
have a much simpler structure than the energy balance laws. This, however, comes
with the price that such a numerical scheme can not be applied to problems with
high Mach number. In particular, the following numerical test problems also aim
to demonstrate that considering the phase entropy balance laws rather than phase
energy ones leads to negligible errors for problems with relatively \textit{weak
shocks}.  Also, it should be emphasized that our future intention is to develop
a numerical scheme more suited to the original SHTC formulation
\eqref{eqn.SHTC.system3}, i.e. the thermodynamically compatible (HTC) type
schemes where the fully-discrete energy conservation is obtained as a mere
consequence of the discretized PDEs, see e.g.
\cite{ABGRALL2023127629,Thomann2023, HTC2022}. 

Furthermore, the presented results consist of a wide range of benchmarks and
problems that may occur in real life problems involving several phases. Some
results demonstrate the solution of the multiphase model in a single-phase
limit, and in the relaxation limits of the GPR model, e.g. the inviscid and
viscous fluid, as well as in the limit of nonlinear elasticity and plasticity.
In all the cases, the numerical results are comparable with results obtained
from established standard models, i.e. the Euler or Navier-Stokes equations for
fluids, or the classical hypo-elastic model with plasticity, but, notably, in
our case the solution is obtained within the unified multiphase model of
continuum mechanics.

In all the tests, the time step $\Delta t$ is computed according to the $\up{CFL}$ condition expressed in \eqref{CFL.SHTC}, in order to guarantee the stability of the explicit FV time-stepping. Furthermore, the initial conditions for volume fractions are defined with respect to a minimum value $\alpha_{\mathrm{min}}={\epsilon} = 10^{-6}$; i.e. when $a$-th phase is not present in a given computational cell, the volume fraction of that phase is set to $\alpha_a = \alpha_{\mathrm{min}}$.
\subsection{Numerical convergence study}\label{sec.vortex}
A numerical convergence study is presented by solving the isentropic vortex problem proposed in \cite{BALSARA2000405, HU199997}, considering the  one-phase limit of the model: $\alpha_1 = 1- 2 \epsilon$, $  \alpha_2 = \epsilon$, $\alpha_3 = \epsilon$. For this problem, there is an exact analytical solution for the compressible Euler equations, i.e. in the stiff inviscid limit $\tau^\mathrm{e}_1 \rightarrow 0 $ of the SHTC BN-model considered in this work. The initial condition consists of a linear superposition of a homogeneous background field and some $\delta$ perturbations, which in terms of primitive variables for the first phase read
\begin{equation}
\vec{V}_1 = (1+\delta\rho_1, 1+\delta v_{1,1}, 1+\delta v_{1,2},1 + \delta p_1, 1- 2\epsilon, \Id ),
\end{equation}
where the phase distortion field is initially set equal to the identity, while the quantities for the absent phases are set in the same way except for the volume fractions. The computational domain is $\Omega= [0; 10] \times [0; 10]$ and periodic boundary conditions are applied everywhere. In this domain, the perturbations of velocities $\delta v_{1,k}$ and temperature $\delta T_{1}$  are given by 
\begin{align}\label{num.conv1}
\begin{pmatrix}
\delta v_{1,1} \\
\delta v_{1,2}
\end{pmatrix} = \frac{5}{2 \pi}  e^{ 0.5(1-r^2)}\begin{pmatrix}
5-y \\
x-5
\end{pmatrix}, \qquad 
\delta T_{1} = - \frac{(\gamma_1 - 1)5^2}{8 \gamma_1 \pi^2} e^{1-r^2},
\end{align}
where $r$ is the distance frm the center of the vortex.
Aditionally, because we are considering an isentropic vortex, it is assumed that the
perturbation of the entropy $\delta s_{1}$ is zero, hence the perturbations for
density and pressure result in
\begin{align}\label{num.conv2}
\delta \rho_{1} =  {( 1+ \delta T_{1})}^{\frac{1}{\gamma_1 -1}} -1, \qquad \delta p_{1} =  {( 1+ \delta T_{1})}^{\frac{\gamma_1}{\gamma_1 -1}} -1.
\end{align}
The exact analytical solution of the problem represented by these initial
conditions for the compressible Euler equations is represented simply by the
time-shifted initial conditions \eqref{num.conv1}, \eqref{num.conv2}, convected
following the mean velocity $\bar{\vv} =(1,1)$. The equation of state parameters
that remain to be defined are assumed to be $\gamma_1 = 1.4$, $\Cv_1 = 1$,
$\Cs_1 = 0.5$, $\tau^\mathrm{e}_1 = 10^{-14}$.

This test is performed up to a final time of $t = 1.0$ using a sequence of successively refined equidistant meshes composed of $\up{N}_x \times \up{N}_y$ control volumes. The $L^1$ and $L^2$ error norms at the final time for the density $\rho_1$, the velocity component $v_{1,1}$ and the phase entropy $s_1$ are shown in Tab. \ref{tab.convergence1} and Tab. \ref{tab.convergence2} together with the corresponding convergence rates. From the results shown in the tables, it can be seen that second-order accuracy is achieved for this inviscid problem, i.e. in the stiff limit of the governing PDE system.
\begin{table}[ht]
	\centering
	\renewcommand{\tablename}{\footnotesize{Tab.}}
	\begin{tabular}{@{}lllllll@{}}
		\toprule
		$\up{N}_x \times \up{N}_y$    	&	 $L^1_{\rho_1}$& $L^1_{v_{1,1}}$ 	&	 $L^1_{s_{1}}$ 	& $\mathcal{O}_{\rho_1}$  &  $\mathcal{O}_{v_{1,1}}$  & $\mathcal{O}_{s_1}$ 	\\ \midrule
		32 			&	2.5094E-1		& 	5.1290E-1	   	&	1.3009E-2	 	& 			&			&		\\
		64  			&	5.2676E-2		&	1.1826E-1 	&	5.4240E-3		&	2.25		&	2.11		&	1.26	\\
		128			&	1.0012E-2 	&	2.7041E-2		&	9.9400E-4		&	2.39		&	2.12		&	2.44	\\
		256			&	1.8412E-3		&	6.3160E-3		&	1.5781E-4		&	2.44		&	2.09		&	2.65	\\
		\midrule
	\end{tabular}
	\caption{\footnotesize {Mesh elements, $L^1$-error norms and their respective numerical convergence rates for the density $\rho_1$, the velocity component $v_{1,1}$ and the phase entropy $s_1$, applied to the isentropic vortex problem.}  }
	\label{tab.convergence1}
\end{table}
\begin{table}[ht]
	\centering
	\renewcommand{\tablename}{\footnotesize{Tab.}}
	\begin{tabular}{@{}lllllll@{}}
		\toprule
		$\up{N}_x \times \up{N}_y$    	&	 $L^1_{\rho_1}$& $L^1_{v_{1,1}}$ 	&	 $L^1_{s_{1}}$ 	& $\mathcal{O}_{\rho_1}$  &  $\mathcal{O}_{v_{1,1}}$  & $\mathcal{O}_{s_1}$ 	\\ \midrule
		32 			&	6.6187E-2		& 	1.3959E-1	   	&	4.4700E-3	 	& 			&			&		\\
		64  			&	1.4075E-2		&	3.4710E-2		&	2.3585E-3		&	2.23		&	2.01		&	0.93	\\
		128			&	2.6702E-3 	&	8.5657E-3		&	4.9292E-4		&	2.40		&	2.02		&	2.26	\\
		256			&	4.8569E-4		&	2.0754E-3		&	7.6455E-5		&	2.46		&	2.05		&	2.69	\\
		\midrule
	\end{tabular}
	\caption{\footnotesize {Mesh elements, $L^2$-error norms and their respective numerical convergence rates for the density $\rho_1$, the velocity component $v_{1,1}$ and the phase entropy $s_1$, applied to the isentropic vortex problem.}  }
	\label{tab.convergence2}
\end{table}
\subsection{ Shear motion in solids and fluids}
In the context of this work, this test has a twofold purpose of showing that the
unified model for fluid and solid mechanics and the developed numerical scheme
can indeed model the behavior of viscous fluids and elastic solids at once. We
consider a simple shear motion in solids and fluids in the single-phase limit of
the entire multiphase model: $\alpha_1 = 1- 2 \epsilon$, $  \alpha_2 =
\epsilon$, $\alpha_3 = \epsilon$. Similar to the previous section
\ref{sec.vortex}, the time evolution of an incompressible shear layer is one of
the few test problems for which the exact analytical solution of the
non-stationary Navier-Stokes equations is known, and for the velocity component $v_{1,1}$ is given by the following error function
\begin{equation}\label{exact.first.stokes}
v_{1,1}(x,y,t) = v_{1,1}(x,y,0)\, \textnormal{erf}\left(\halb \frac{x}{\sqrt{\nu_1 t}}\right).
\end{equation} 
However, because we are
discretizing compressible equations with an explicit scheme, the best we can do
is to simulate the problem at sufficiently low Mach number, e.g. $\up{M}_1 = 0.1$
was sufficient to obtain an almost incompressible behavior.

The computational domain is $\Omega=[-0.5 ; 0.5] \times [-0.0625; 0.0625]$, with
the opposite velocities imposed on the left and right halfs of the domain in the
$x$-direction, while we use periodic boundary conditions in the $y$-direction.
The initial conditions of the problem for the first phase, are given by
\begin{align}\label{exact.first.stokes}
\begin{split}
&\alpha_1 = 1- 2 \epsilon,  	 \quad   \rho_1 = 1,  \quad p_1 = \frac{1}{\gamma_1}, \quad \A_1 = \Id \\
&v_{1,1} = 0, \quad
v_{1,2}(x,y) = \begin{cases} +0.1  &\mathrm{if}  \ \  x > 0 , \\
	-0.1 \ \ & \mathrm{if}  \ \  x \leq 0.
\end{cases}
\end{split}
\end{align} 
with the physical parameters set to $\gamma_1 = 1.4$, $\Cv_1 = 1$, $\Cs_1 = 1$.
The strain relaxation time $\tau^\mathrm{e}_1 = 6 \nu_1/\Cs_1^2$ is chosen
according to \eqref{eqn.viscosity} for various values of the fluid kinematic
viscosity $\nu_1$, while for the elastic solid limit is set to $\tau^\mathrm{e}_1 =
10^{14}$. 

For the elastic solid limit, this initial condition leads to two shear waves
travelling to the left and right with the shear sound speed. In this case, a
reference solution for the solid limit was obtained for the single-material GPR
model using a classical second-order MUSCL-Hancock scheme \cite{torobook} on a
fine mesh of 32000 cells, as done in \cite{HTC2022}.

Simulations are carried out on a grid composed of $256\times32$ control volumes up to a final time of $t = 0.4$. 
\begin{figure}[!htbp]
	\renewcommand{\figurename}{\footnotesize{Fig.}}
	{\includegraphics[width=.46\textwidth]{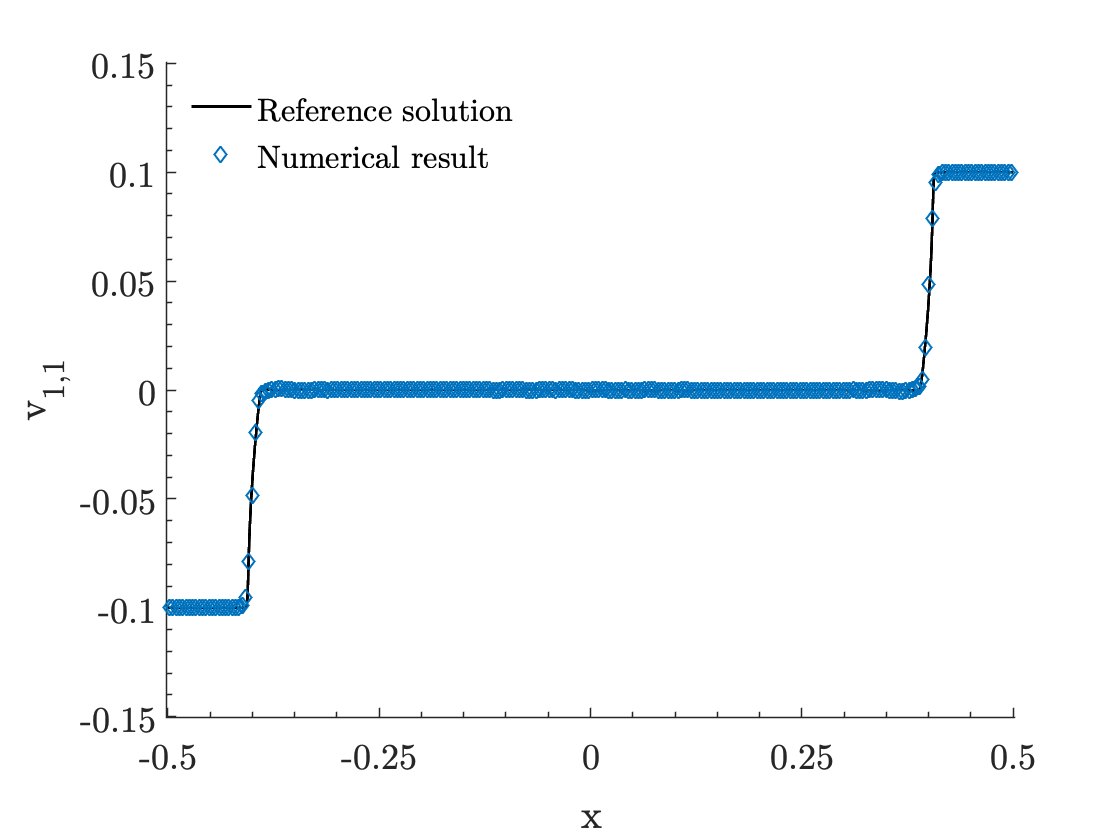}}\hfill
	{\includegraphics[width=.46\textwidth]{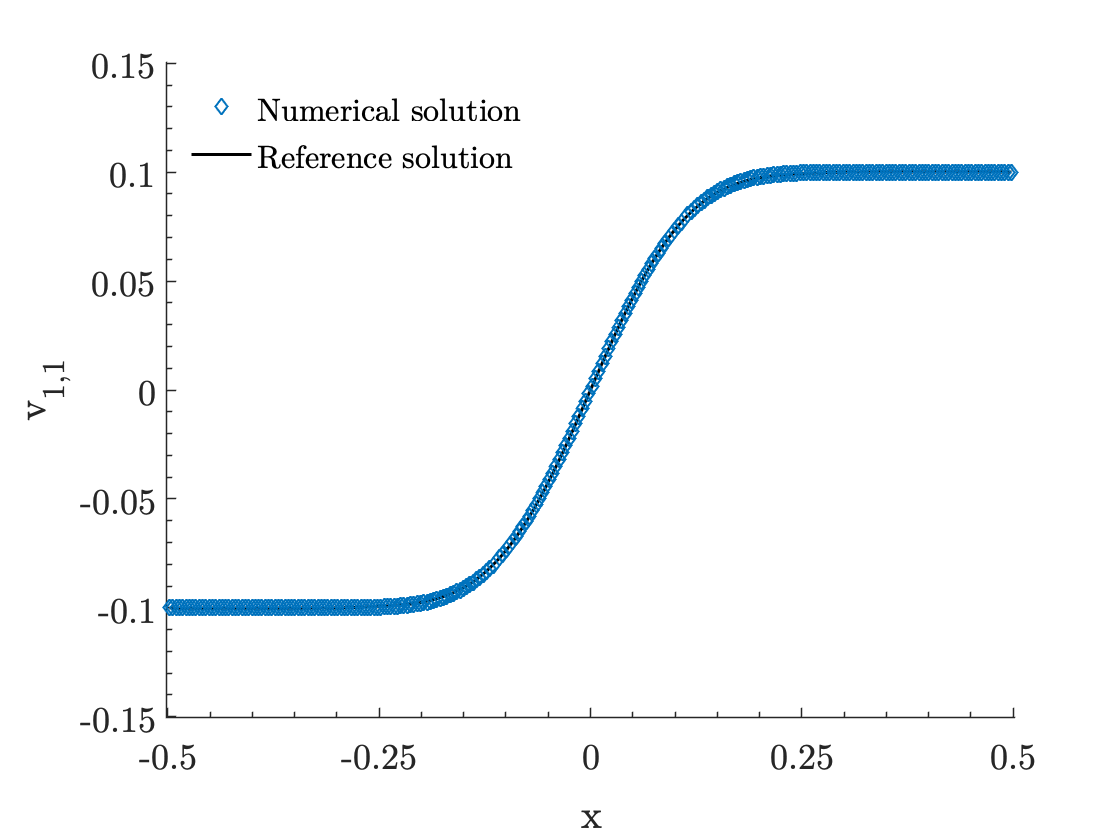}}\hfill 	
	{\includegraphics[width=.46\textwidth]{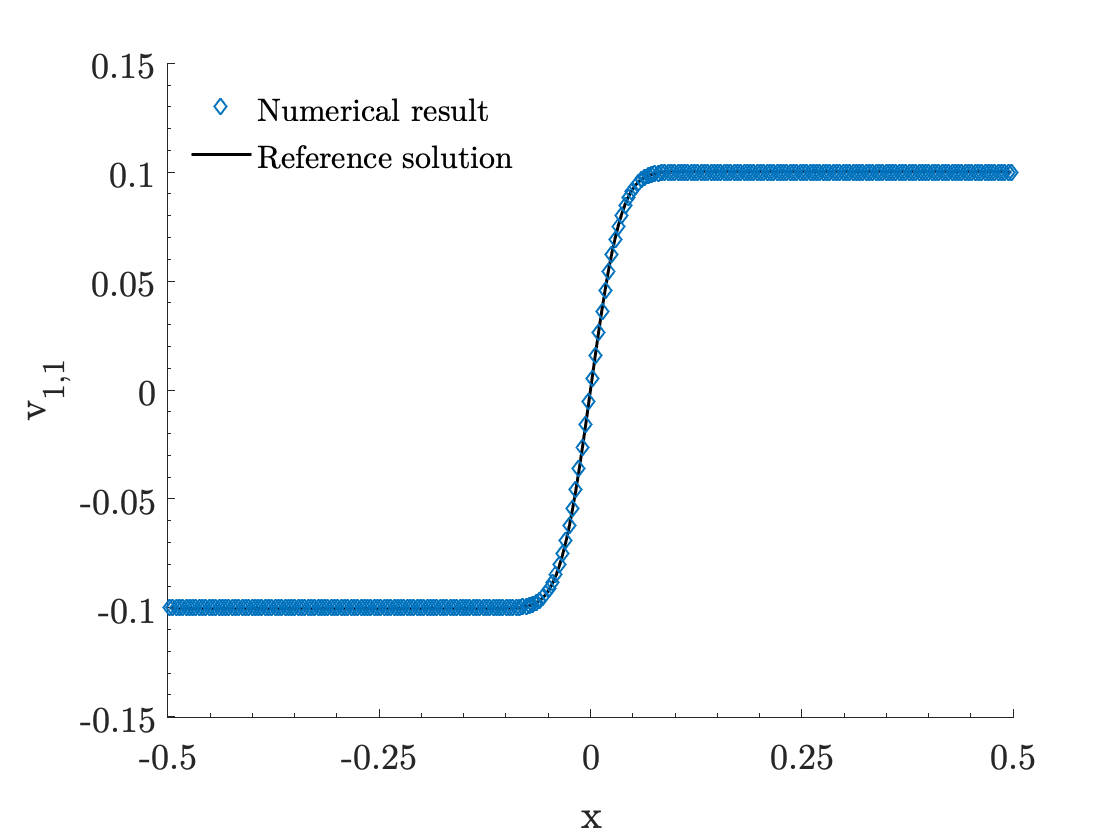}}\hfill
	{\includegraphics[width=.46\textwidth]{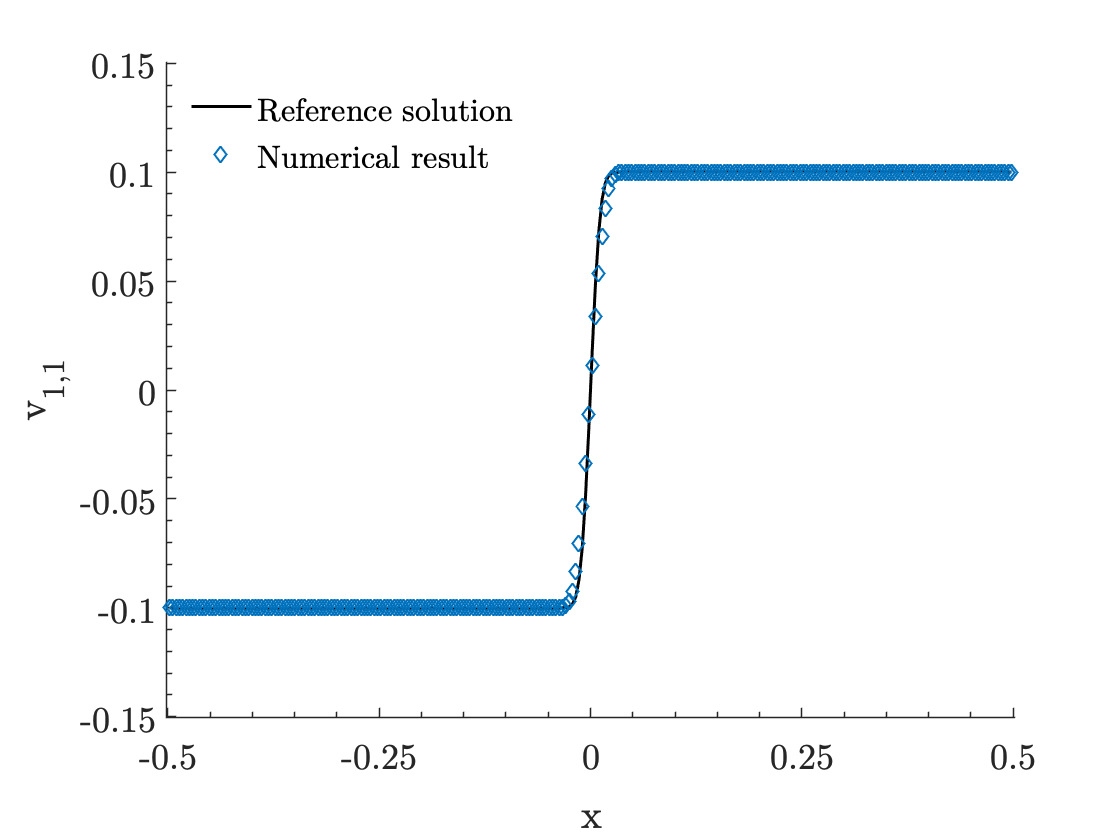}}\hfill 
	\caption{\footnotesize Numerical solution at time $t = 0.4$ obtained with the explicit FV scheme for compressible multiphase fluid and solid mechanics applied to a simple shear flow in fluids and in an elastic solid. Results for the solid limit (top left) and for fluids with different viscosities $\nu_1 = 10^{-2}$ (top right), $\nu_1 = 10^{-2}$ (bottom left) and $\nu_1 = 10^{-2}$ (bottom right). For fluids, the analytical solution of the first problem of Stokes is used as the reference solution.}
	\label{fig:shear.solid.liquid}
	\vspace{-1em}
\end{figure}
The comparison between the numerical results and the previously mentioned
reference solutions is presented in Fig. \ref{fig:shear.solid.liquid}, where an
excellent agreement between the two solutions can be observed for both solid and fluid
behaviour.
\subsection{Riemann problems}
We continue the validation of our numerical scheme with a set of Riemann
problems to quantify the error encountered by considering the balance equations
of phase entropies rather than of the phase energies, and to see that the
correct wave structure can still be reproduced for problems with relatively weak
shocks. In this section, we solve a series of Riemann problems with initial data
according to Tab. \ref{tab.RP}, for the Euler equations of compressible gas
dynamics, which can be retrieved in the stiff relaxation limit
$\tau^\mathrm{e}_1 \rightarrow 0$.
\begin{table}[ht]
	\centering
	\renewcommand{\tablename}{\footnotesize{Tab.}}
	\begin{tabular}{@{}lllllllll@{}}
		\toprule
		$\up{RP}$    		&	 ${\rho_1^L}$	& $ v^L_{1,1} $ 	&	 $ v^L_{1,2} $ 	& $ P^L_{1} $  	&	 ${\rho_1^R}$	& $ v^R_{1,1} $ 	&	 $ v^R_{1,2} $ 	& $ P^R_{1} $ 	\\ \midrule
		\up{RP1} 			&	1.0			& 	0.0	   		&	0.0			& 	1.0		&	0.125		&	0.0			&		0.0		&	0.1	\\
		\up{RP2}  			&	1.0			&	0.75			&	0.0			&	1.0		&	0.125		&	0.0			&		0.0		&	0.1	\\
		\up{RP3}			&	1.0	 		&	0.0			&	-0.2			&	1.0		&	0.5			&	0.0			&		0.2		&	0.5	\\
		\midrule
	\end{tabular}
	\caption{\footnotesize {Left initial state (L) and right initial state (R)
	for the quantities related to the first phase. In particular the density
	$\rho_1$,  velocity $\vv = (v_{1,1}, v_{1,2}, 0 )$ and  pressure $P$ are
	defined for three different Riemann problems. These Riemann problems
	(RP1), (RP2) and (RP3) can be referred to the solution of the Euler equations, i.e.
	$\tau^\mathrm{e}_1 = 10^{-14}$. } }
	\label{tab.RP}
\end{table}

The computational $\Omega=[-0.5 ; 0.5] \times [-0.0625; 0.0625]$ is partitioned
into two regions with constant states, left (L) and right (R), separated by a
discontinuity normal to the $x$-direction, located at $x_d$. The distortion
field is initially set equal to the identity $\A_1 = \Id$, while the equation of
state parameters are taken as $\gamma_1 = 1.4$, $\Cv_1 = 1.0$, $\Cs_1 = 1.0$ and
Riemann problems. Simulations are carried out on a grid composed of $512\times64$ control
volumes up to a final time of $t = 0.2$.

\begin{figure}[!htbp]
	\renewcommand{\figurename}{\footnotesize{Fig.}}
	{\includegraphics[width=.46\textwidth]{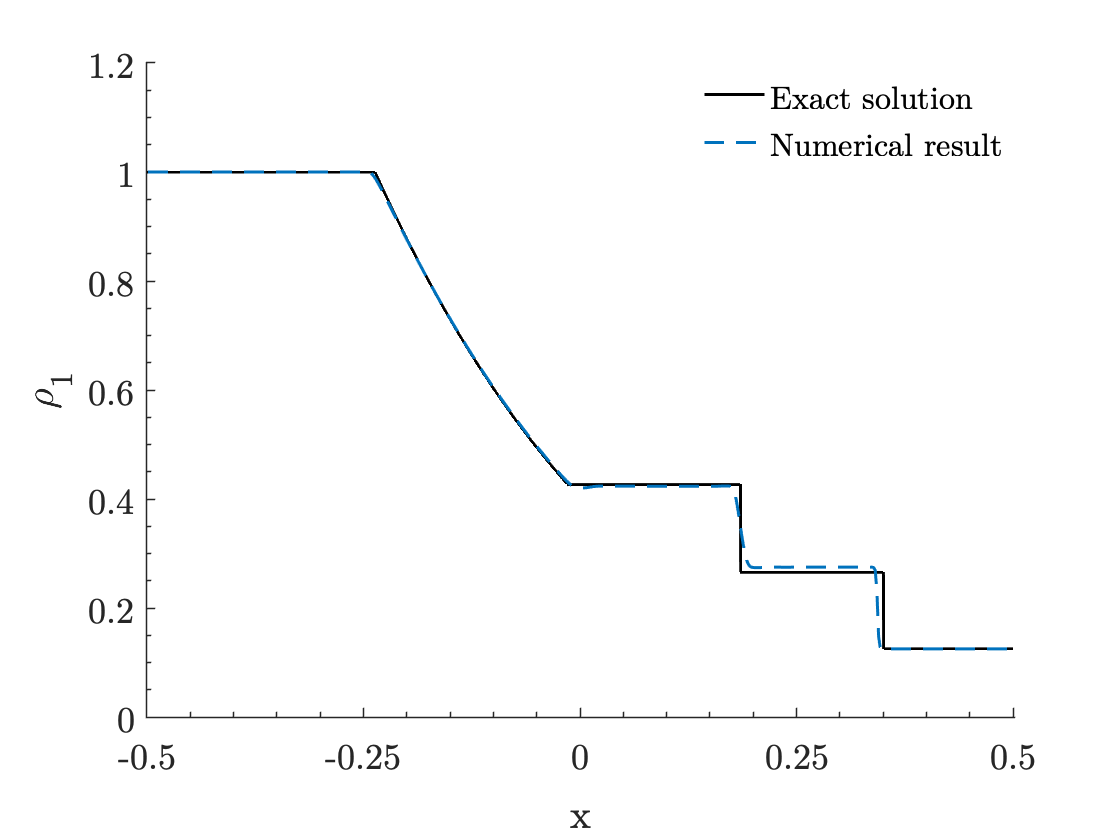}}\hfill
	{\includegraphics[width=.46\textwidth]{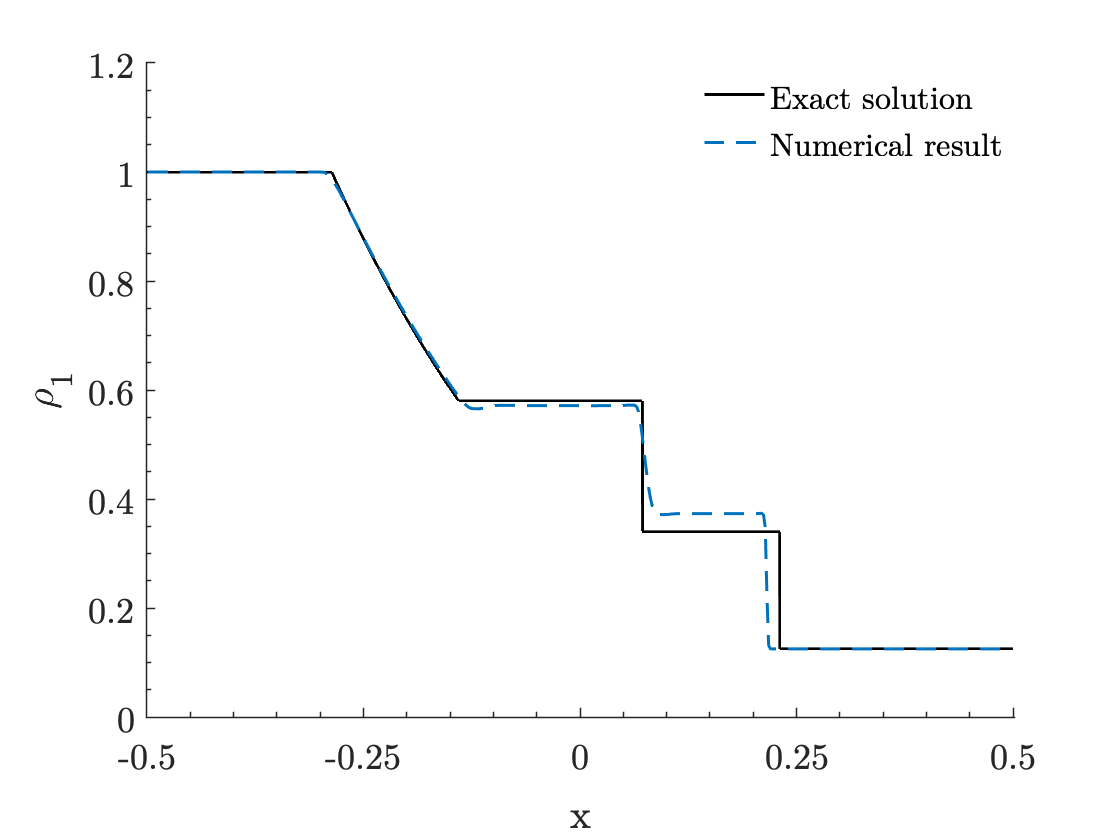}}\hfill 	
	{\includegraphics[width=.46\textwidth]{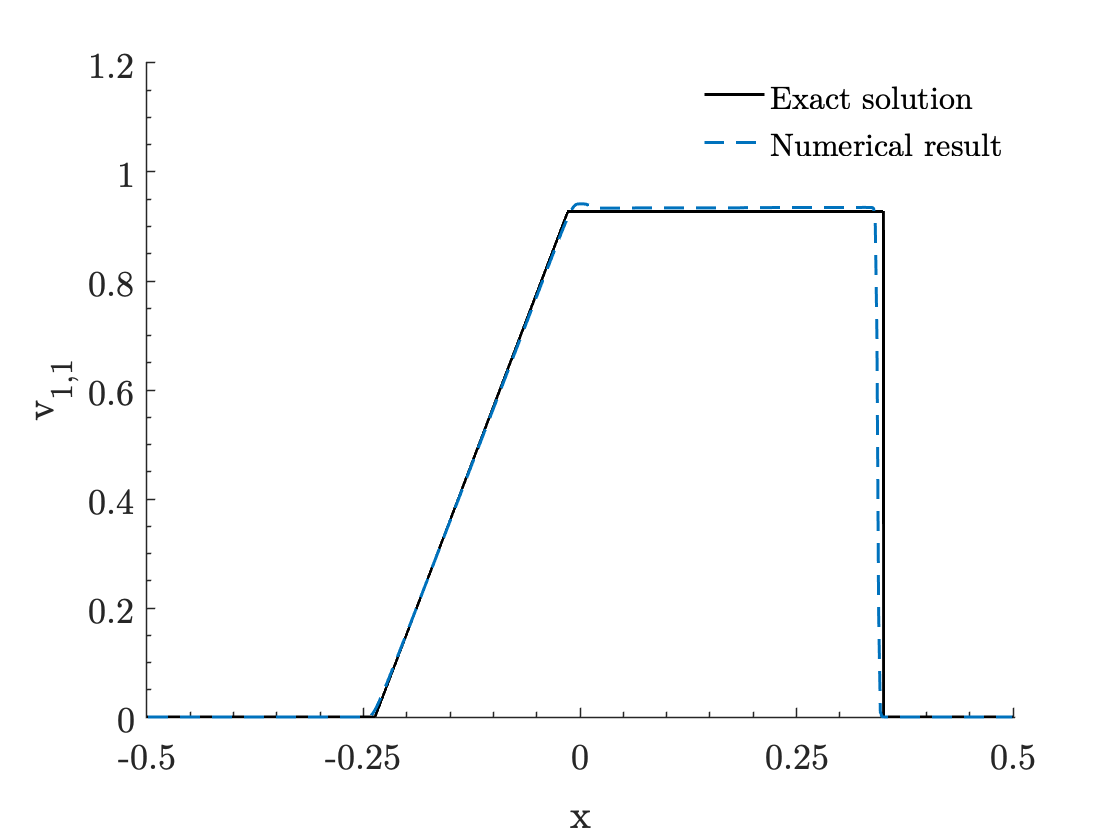}}\hfill
	{\includegraphics[width=.46\textwidth]{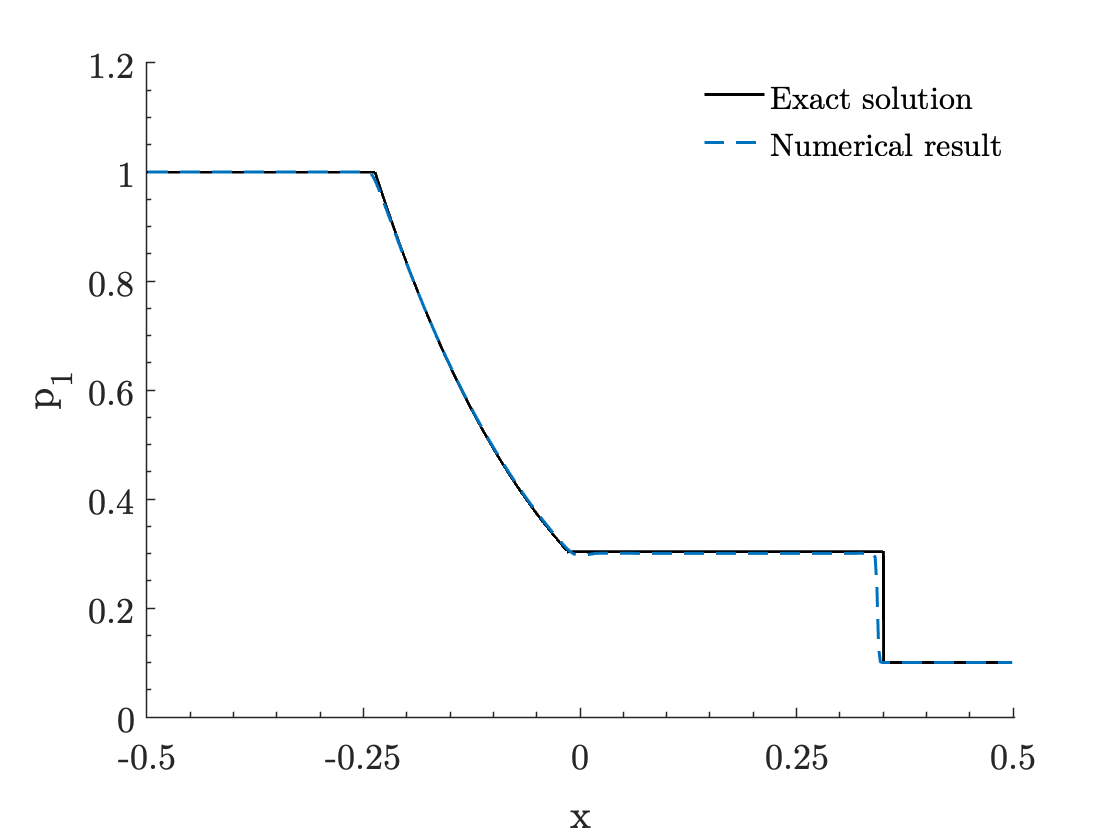}}\hfill 
	\caption{\footnotesize Numerical results (dashed line) for density $\rho_1$, velocity component $v_{1,1}$ and pressure $p_1$ in the inviscid limit $\tau_1 = 10^{-14}$, for the Riemann problem RP1 $(x_d = 0)$  (top left, bottom left and right), for the Riemann problem RP2 $(x_d = -0.2)$ (top right). The exact solution of the compressible Euler equations (black solid line).}
	\label{fig:RP1}
	\vspace{-1em}
\end{figure}

In Figure \ref{fig:RP1}, the one-dimensional profiles of the density $\rho_1$,
the $x$-component of the velocity field $v_{1,1}$ and pressure $p_1$ for the
Riemann problems RP1 $(x_d = 0)$ and RP2 $(x_d = -0.2)$ are shown. The results
are compared with the exact solution of the compressible Euler equations. From
the results it can be observed that the correct wave structure is overall
reproduced properly for the Riemann problem RP1, while as the shock wave becomes
stronger, as for RP2, the error introduced due to the use of phase entropy
balance laws increases. This is why we limit ourselves to low Mach number flows,
and RP2 clearly demonstrates the well-known fact that satisfying the
conservation of energy is essential to correctly solve problems involving shock
waves. RP2 was proposed by Toro in \cite{torobook} and includes a sonic
rarefaction, however this test cases is well resolved and does not present any
sonic glitches.

\begin{figure}[!htbp]
	\renewcommand{\figurename}{\footnotesize{Fig.}}
	{\includegraphics[width=.46\textwidth]{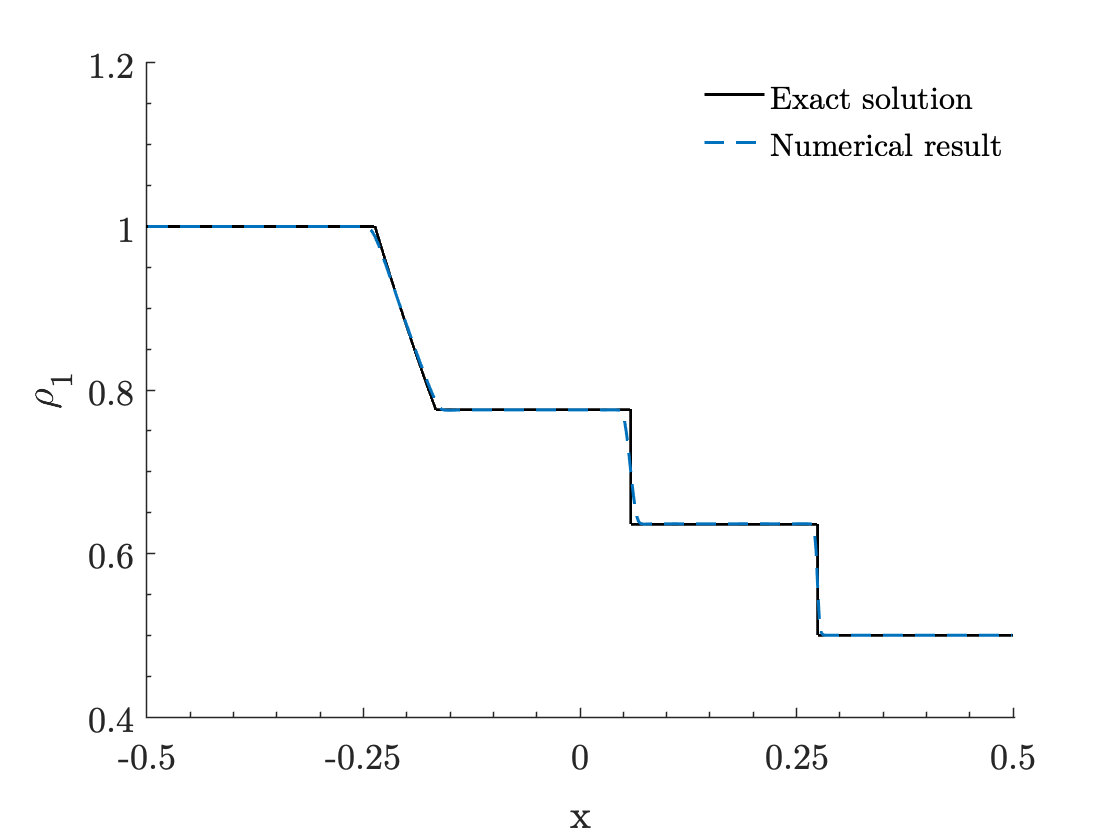}}\hfill
	{\includegraphics[width=.46\textwidth]{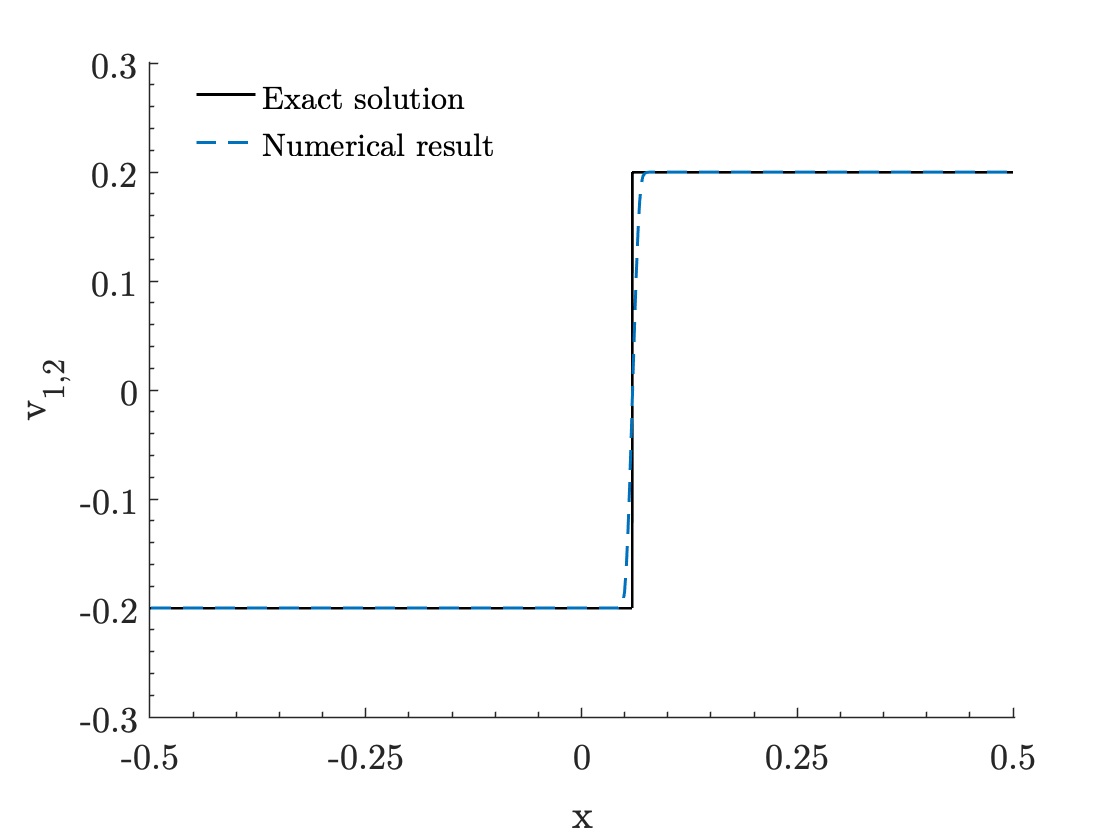}}\hfill 	
	\caption{\footnotesize Numerical results for density $\rho_1$ and velocity component $v_{1,2}$ in the inviscid limit $\tau_1 = 10^{-14}$, at time $t = 0.2$, for the Riemann problem RP3 $(x_d = 0)$ (dashed line). The exact solution of the compressible Euler equations (black solid line).}
	\label{fig:RP2}
\end{figure}

The numerical results obtained for the Riemann problem RP3 $(x_d = 0)$ is shown
in Fig. \ref{fig:RP2}. In this case, the shock present is even weaker, the
solution is very close to an isoentropic one, and therefore the numerical
solution is in very good agreement with the exact one.

\subsection{Double shear layer problem}\label{double.shear.layer} The numerical
scheme is now applied to solve the double shear layer test problem, see e.g.
\cite{Bell1989, DUMBSER2016824, Tavelli2015, HTC2022,Chiocchetti2023}. It is
another classical benchmark problem which is useful for the validation of the
model and numerical algorithm on viscous flows. Here, the model is still
considered in the single-phase limit with the volume fractions taken as
$\alpha_1 = 1- 2 \epsilon$, $ \alpha_2 = \epsilon$, $\alpha_3 = \epsilon$. For
this test, the 2D computational domain is $\Omega=[0 ; 1]^2$, with periodic
boundary conditions imposed everywhere. The initial conditions contain a steep
velocity gradient and are defined as follows
\begin{align}
\begin{split}
&\alpha_1 = 1- 2 \epsilon,  	 \quad   \rho_1 = 1,  \quad p_1 = \frac{100}{\gamma_1}, \quad \A_1 = \Id, \\[1.mm]
&v_{1,1}(x,y)  = \begin{cases} \mathrm{tanh}\big( (y-0.25)\bar{\rho} \big),  &\mathrm{if}  \ \  x \leq 0.5 , \\
	 \mathrm{tanh}\big((0.75-y)\bar{\rho} \big),\ \ & \mathrm{if}  \ \  x > 0.5,
\end{cases} \\[1.5mm]
&v_{1,2}(x,y) = \delta \mathrm{sin} (2 \pi x),
\end{split}
\end{align} 
where the parameters that determine the shape of the velocity field are set to $\delta=0.05$ and $\bar{\rho} = 30$. The other physical parameters are assumed to be $\gamma_1 = 1.4$, $\Cv_1 = 1$, $\Cs_1 = 8.0$ while two different viscosity coefficients were set in two separate runs of the test problem, $\nu_1= 2\times10^{-3} $ ($\reynolds \simeq 1000$) and $\nu_1= 2\times10^{-4} $ ($\reynolds \simeq 10000$) respectively, which result in $\tau^\mathrm{e}_1 = 1.875\cdot 10^{-4}$ and $\tau^\mathrm{e}_1 = 1.875\cdot 10^{-5}$, respectively. 

\begin{figure}[!htbp]
	\renewcommand{\figurename}{\footnotesize{Fig.}}
	\centering
	{\includegraphics[width=.40\textwidth]{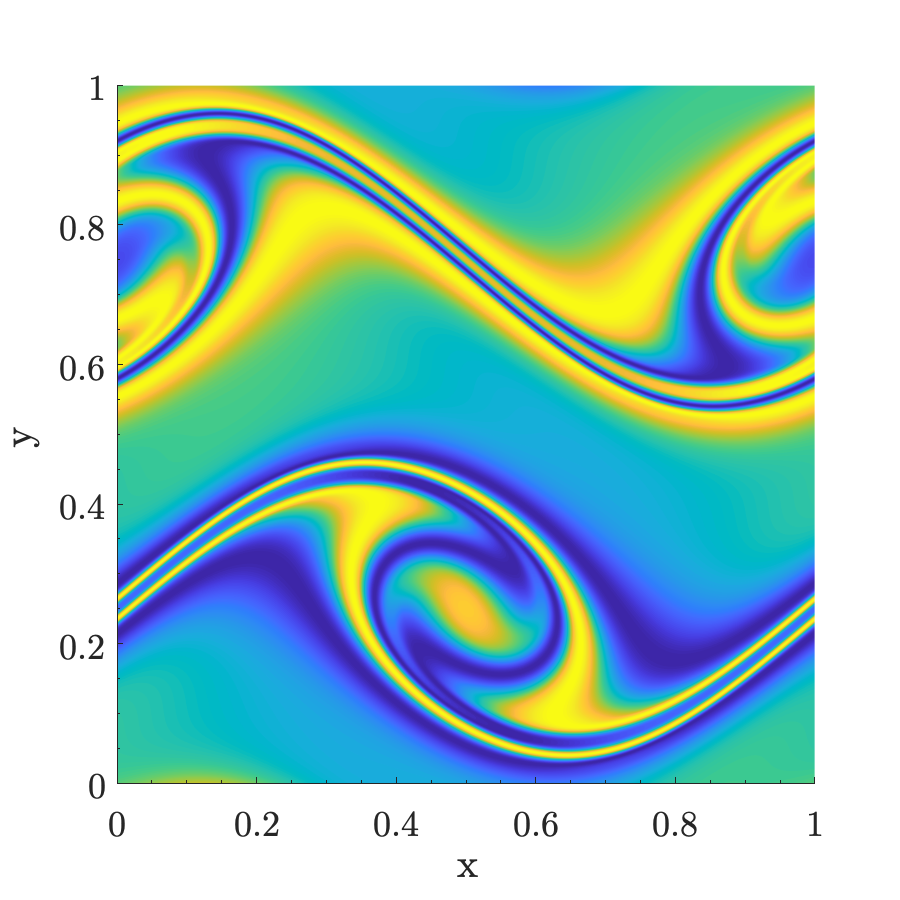}}
	{\includegraphics[width=.45\textwidth]{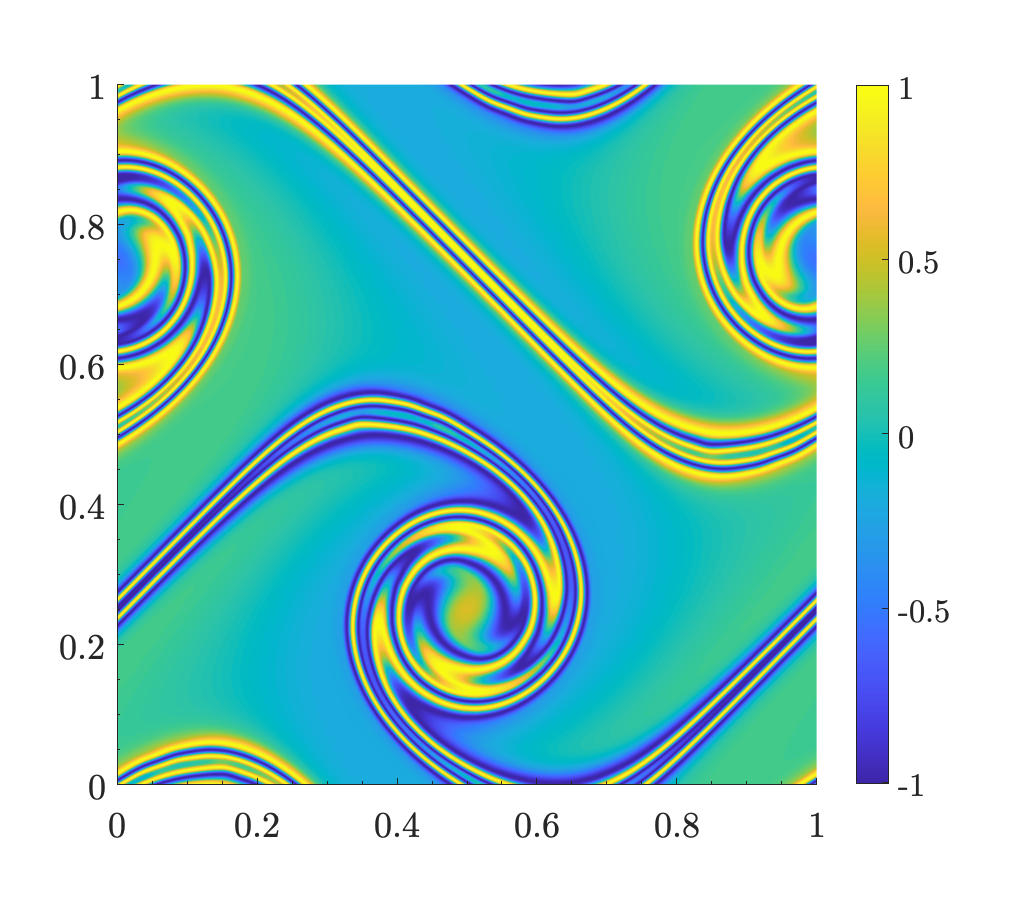}} 		
	{\includegraphics[width=.40\textwidth]{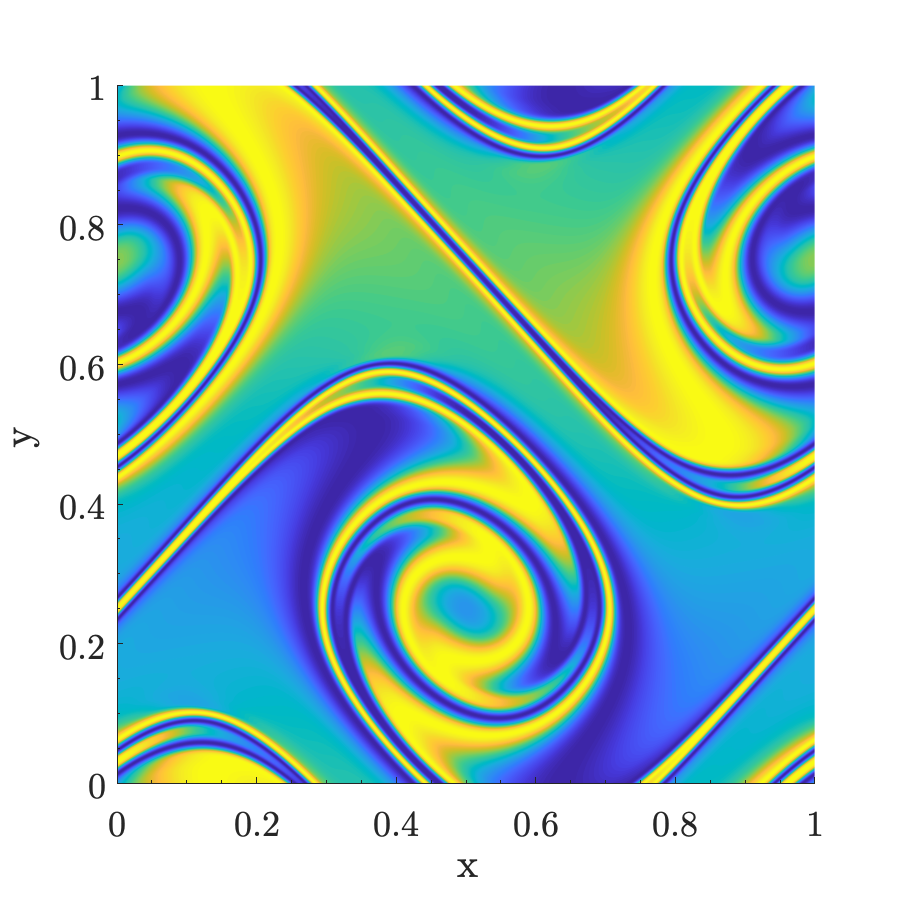}}
	{\includegraphics[width=.45\textwidth]{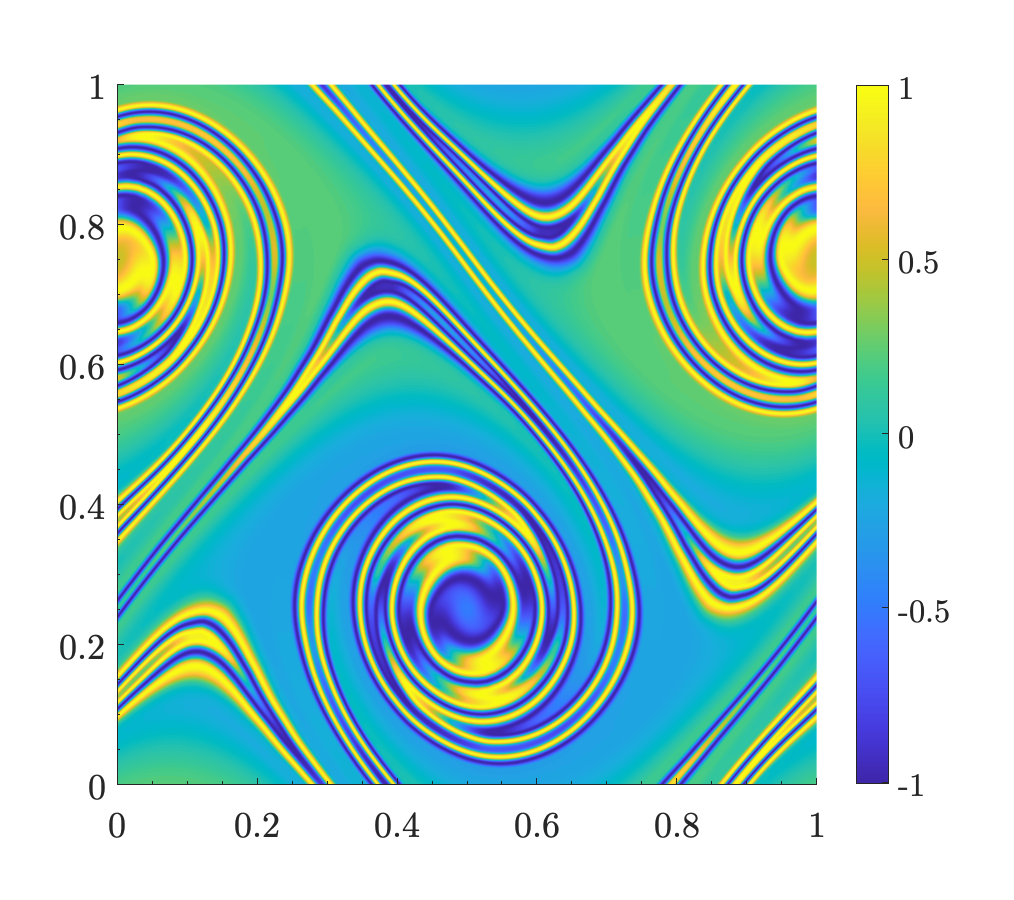}} 
	{\includegraphics[width=.40\textwidth]{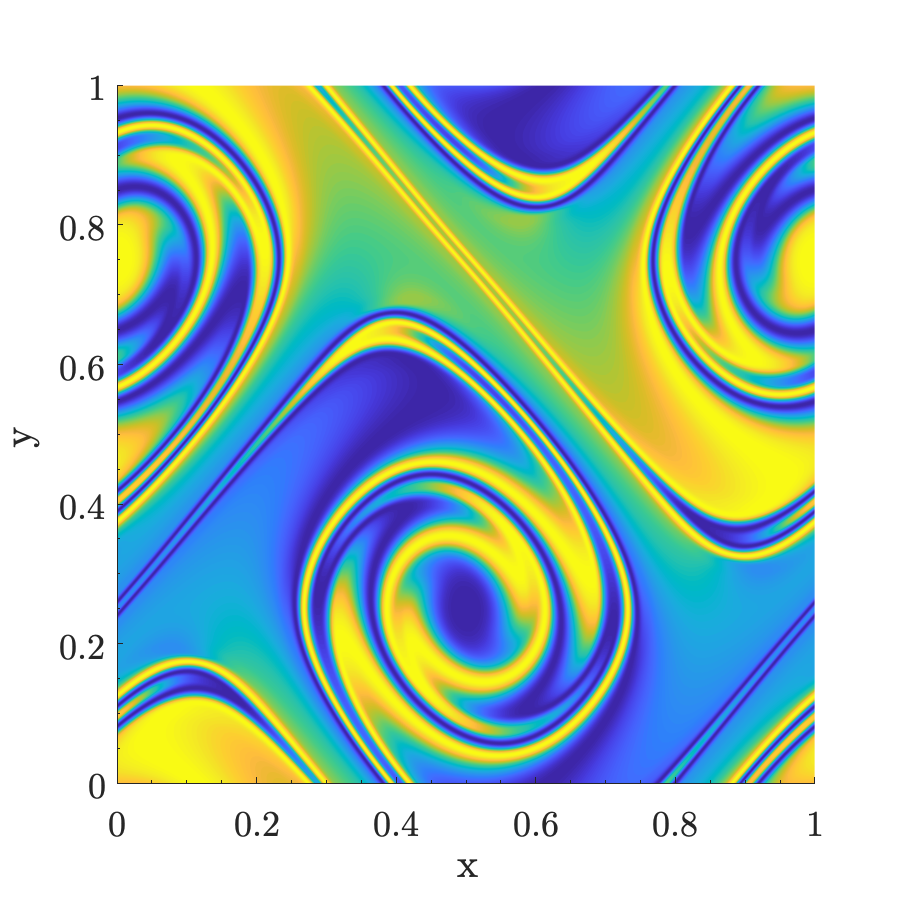}}
	{\includegraphics[width=.45\textwidth]{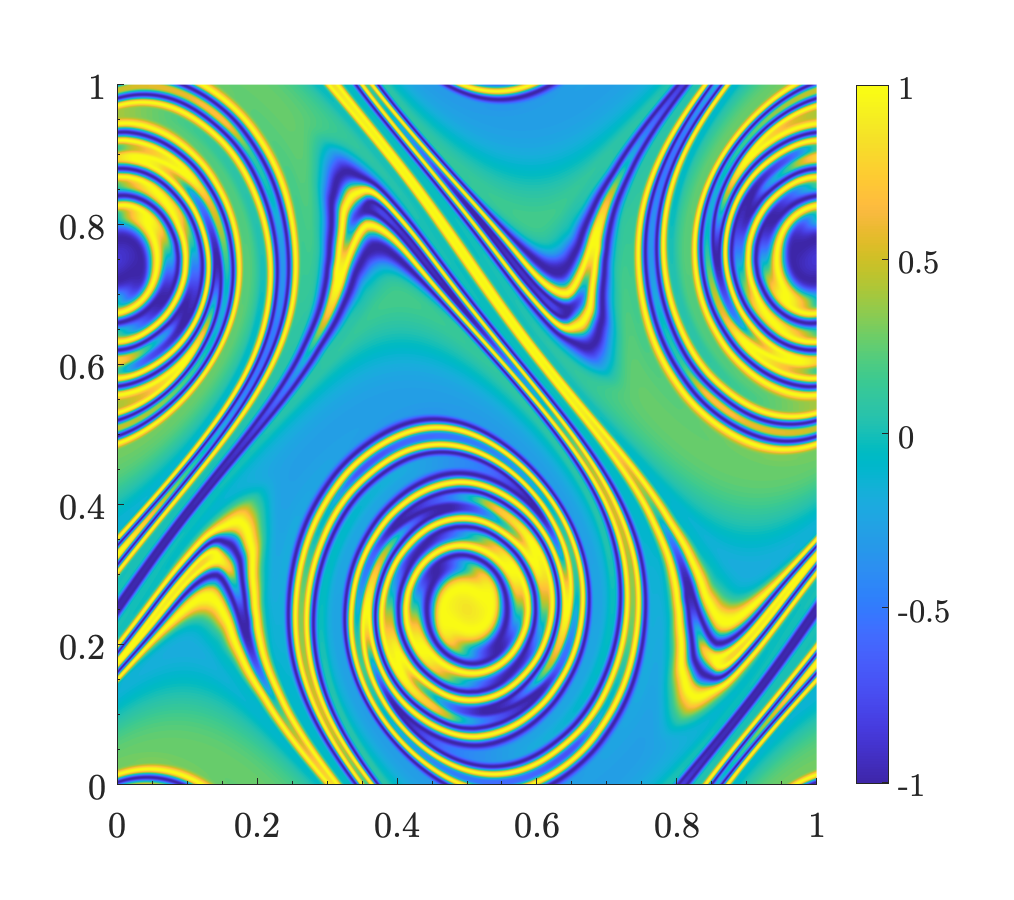}}
	\caption{\footnotesize Filled contours of one component of the distortion field $\A_1$, namely of the $A_{1,12}$ component, for the double shear layer problem at times $t=1.2$ (top), $t=1.6$ (center) and $t=1.8$ (bottom); for two values of kinematic viscosity $\nu_1= 2\times10^{-3} $ ($\reynolds \simeq 1000$) (left) and $\nu_1= 2\times10^{-4} $ ($\reynolds \simeq 10000$) (right).}
	\label{fig:DS}
\end{figure}

Simulations are carried out up to a final time of $t=1.8$ on a grid consisting of $1280\times1280$ control volumes. Figure \ref{fig:DS} shows the time evolution of the $A_{1,12}$ component of the distortion field at times $t=1.2$ (top), $t=1.6$ (center) and $t=1.8$ (bottom), for the two different viscosity coefficients considered (left) and (right), respectively. The dynamics of the flow, as already described in \cite{Bell1989, DUMBSER2016824, Tavelli2015, HTC2022,Chiocchetti2023}, is represented by the evolution of the initially perturbed shear layers into different vortices, which exhibit particularly complex flow structures.

The results in Fig. \ref{fig:DS}, highlight the incredible capability of the distortion field to describe the details of the flow structures, which in particular are encoded in the rotational component $\boldsymbol{R}_1$ of the distortion field $\A_1$. The results obtained are in excellent agreement with those obtained in \cite{HTC2022}, where a thermodynamically compatible scheme is used and with those in \cite{Chiocchetti2023} obtained through a semi-implicit structure-preserving scheme, despite the fact that in these works a four times finer grid was used.
\subsection{Lid-driven cavity}\label{lid.driven.distortion} As a last numerical
test considering a single-phase limit of the full model, we present the
lid-driven cavity problem, see \cite{Ghia1982}. It is a classical benchmark
problem for numerical methods applied to incompressible Navier-Stokes equations,
see \cite{Tavelli2014}, however it can be used to validate compressible flow
solvers in the low Mach number regime \cite{DumbserCasulli2016,
TavelliDumbser2017, BERMUDEZ12}. Moreover, it has already been successfully
solved with the GPR model in \cite{DUMBSER2016824, SIGPR} with a discontinuous Galerkin and a semi-implicit scehem, and with a
thermodynamically compatible scheme in \cite{HTC2022}. However, in these works,
high-order schemes or schemes that make use of a particular time or structure-preserving discretizations have
been used, e.g. on staggered grid or thermodynamically compatible
discretization. These could be the reasons why it would appear that the problem
associated with the discretization, i.e. the combination, of purely rotational
fields presented in Sec.\,\ref{remarks.dist.field} does not arise in these works,
whereas it does for the classical MUSCL-Hancock discretization, specifically in
this test where the boundary conditions produce a velocity gradient singularity
in the corners.

The computational domain is given by $\Omega=[-0.5 ; 0.5]\times[-0.5 ; 0.5]$ and the initial condition are simply 
\begin{align}
\alpha_1 = 1- 2 \epsilon,  	 \quad   \rho_1 = 1 , \quad \vv_1 = \boldsymbol{0},  \quad p_1 = \frac{100}{\gamma_1}, \quad \A_1 = \Id.
\end{align} 
The fluid flow inside the cavity is driven by the lid on the upper boundary,
whose velocity is set to $v_{1,1} = 1$, resulting in the
Mach number $\up{M}_1 = 0.1$.  On all the other
boundaries, a no-slip wall boundary condition with $\vv_1 = 0$ has been be
imposed. Furthermore, the parameters of the model are set to $\gamma_1 = 1.4$,
$\Cv_1 = 1$, $\Cs_1 = 8.0$ and the kinematic viscosity is chosen as $\nu_1=
10^{-2} $ so that the Reynolds number of the test problem is $\reynolds =
100$. 

Simulations are carried out up to a final time of $t=10$ on a grid consisting of
$512\times512$ control volumes. To correctly set the no-slip wall boundary
conditions, it is necessary to compute and prescribe a specific distortion field
$\boldsymbol{A}_{BC}$ using the values taken in the edge-adjacent cell. First,
the information encoded by the distortion field is expressed through
$\boldsymbol{R}_{BC\, (1)}$ and $\boldsymbol{G}_{BC\, (1)}$, by means of a polar
decomposition. Subsequently, the inverse of the rotational component can be
easily evaluated as $\boldsymbol{R}^{ -1}_{BC\, (1)}=\boldsymbol{R}^{
\transpose}_{BC\, (1)}$. At this point the information can be mapped back to
obtain the boundary condition for the distortion field as 
\begin{equation}
\boldsymbol{A}_{BC} = \boldsymbol{R}^{ \transpose}_{BC\, (1)} \ \GG_{BC\, (1)}^{1/2}.
\end{equation}
\begin{figure}[!htbp]
	\renewcommand{\figurename}{\footnotesize{Fig.}}
	\centering
	{\includegraphics[width=.45\textwidth]{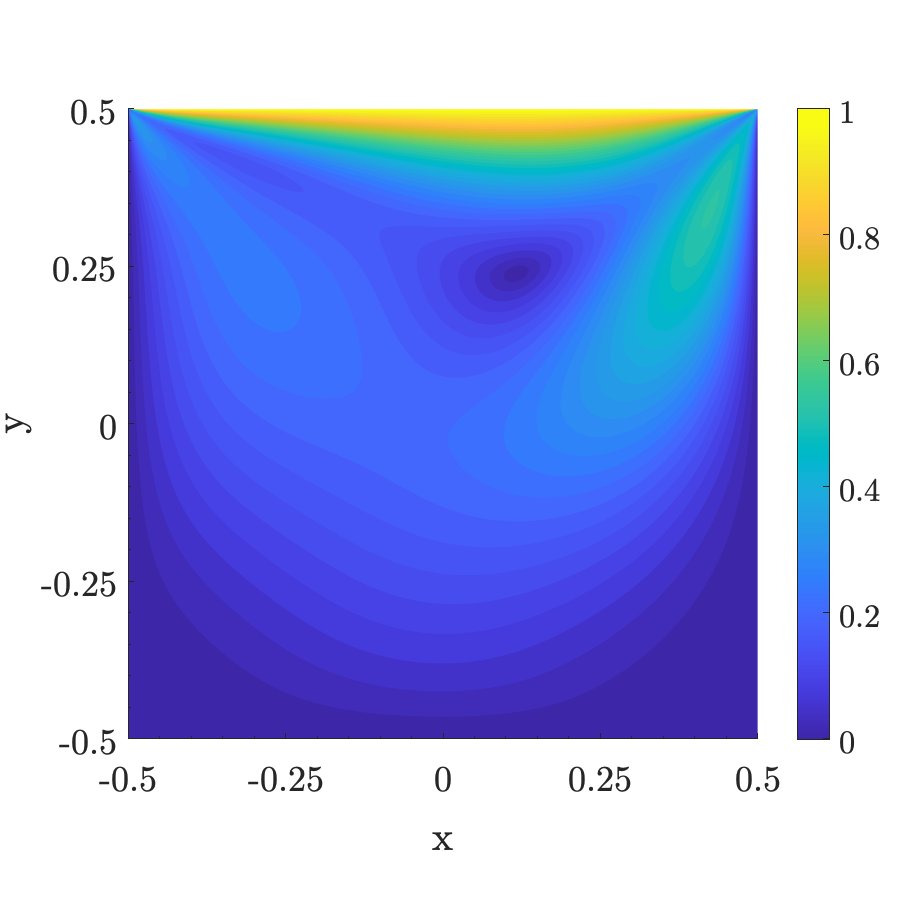}}
	{\includegraphics[width=.45\textwidth]{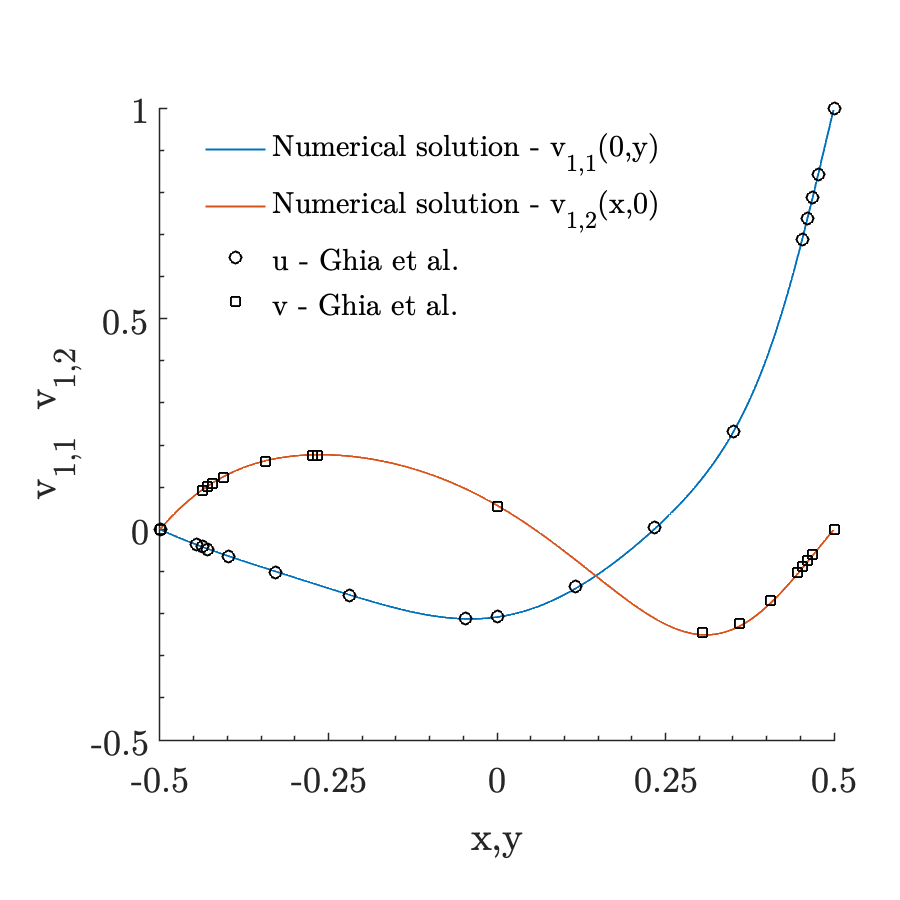}} 	
	\caption{\footnotesize Lid driven cavity at Reynolds number $\reynolds = 100$. Numerical results obtained at time $t=10.0$. Colour contours of the velocity module (left), and a comparison with the reference solution of Ghia et al. \cite{Ghia1982} of the velocity components $v_{1,1}$ and $v_{1,2}$ for 1D cuts along the $x$ and $y$ axis.}
	\label{fig:Lid1}
	\vspace{-1em}
\end{figure}

\begin{figure}[!htbp]
	\renewcommand{\figurename}{\footnotesize{Fig.}}
	\centering
	{\includegraphics[width=.40\textwidth]{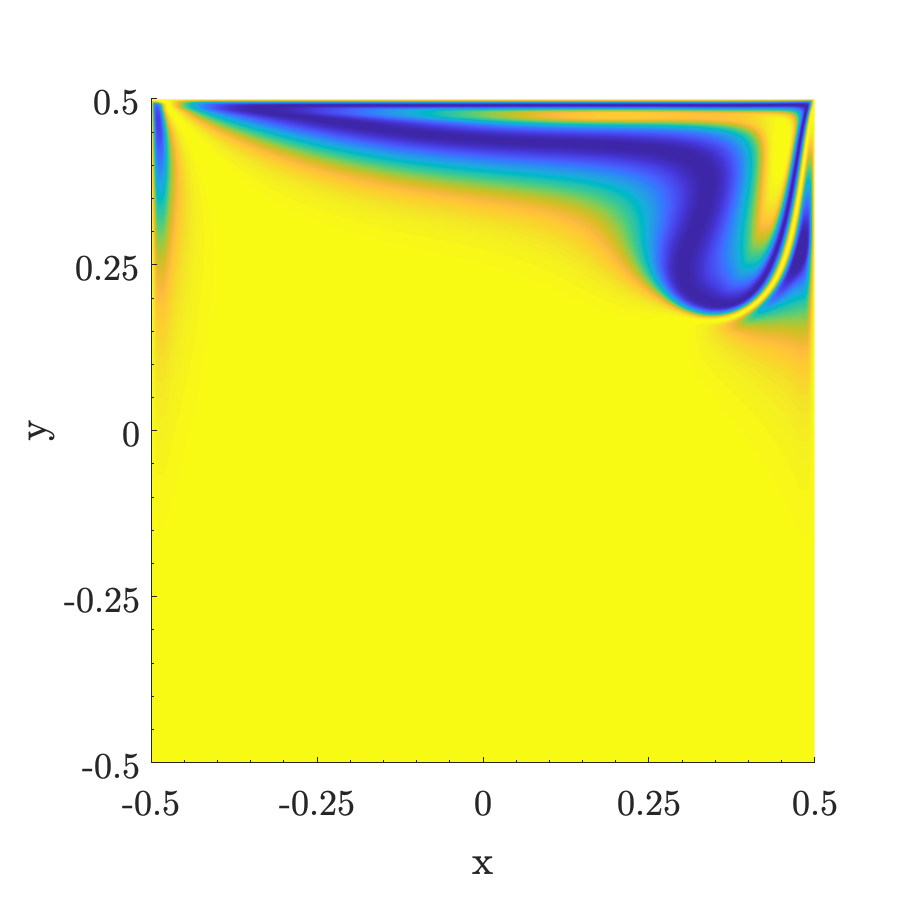}}
	{\includegraphics[width=.45\textwidth]{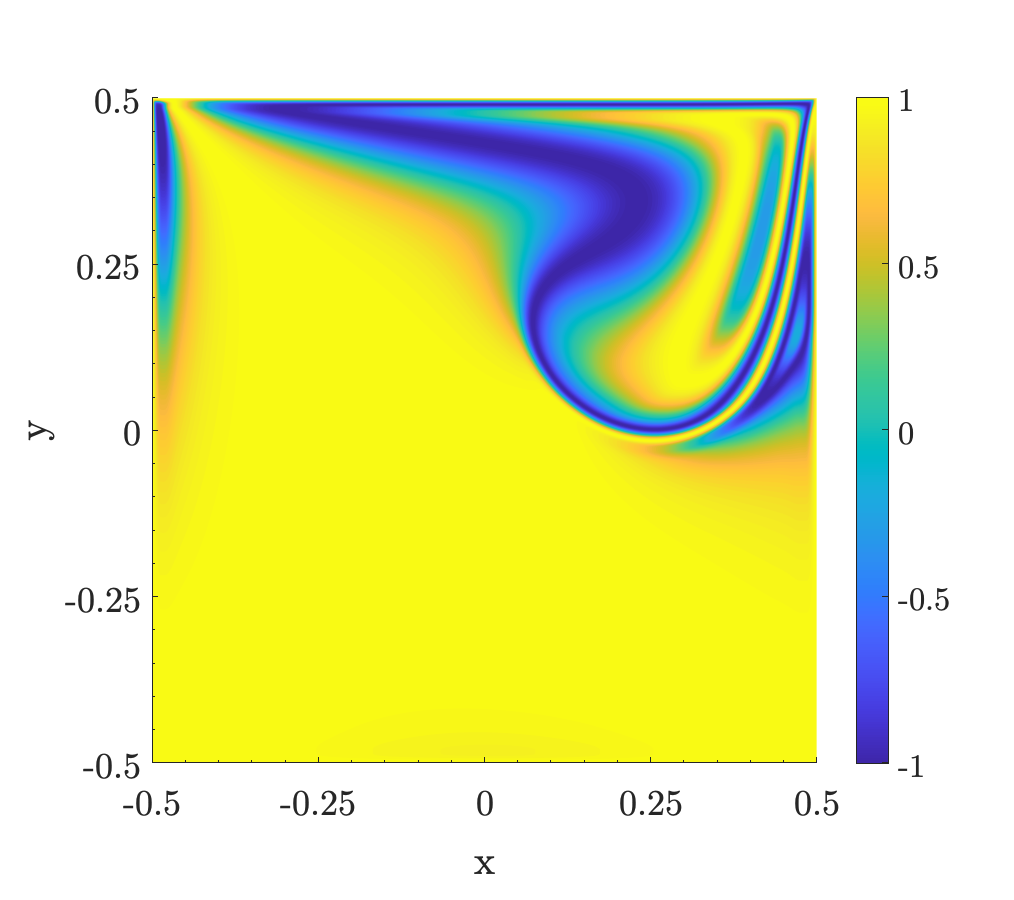}} 	
	{\includegraphics[width=.40\textwidth]{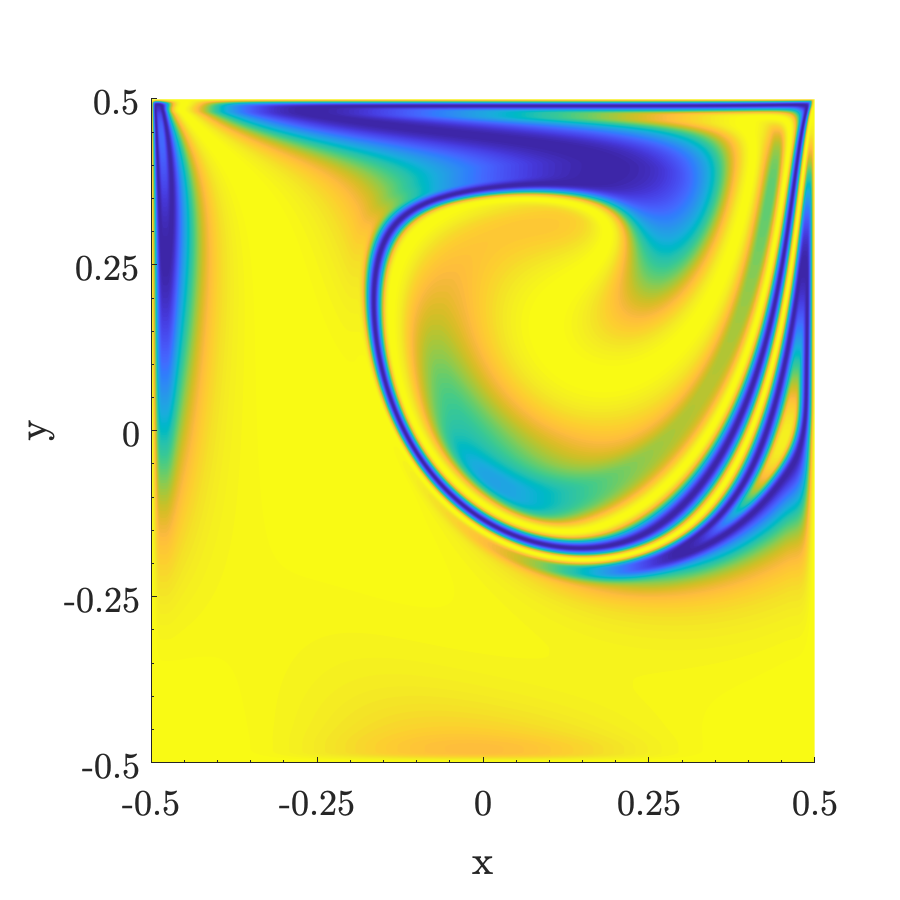}}
	{\includegraphics[width=.45\textwidth]{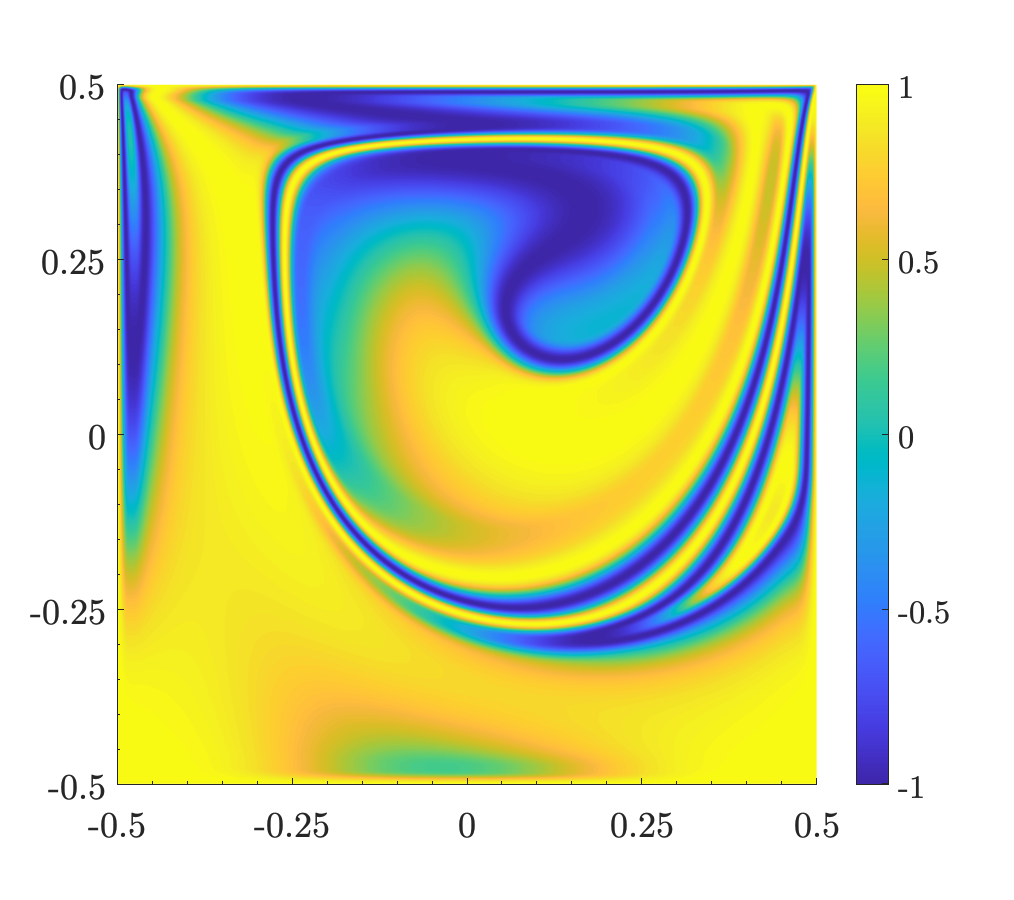}} 
	{\includegraphics[width=.40\textwidth]{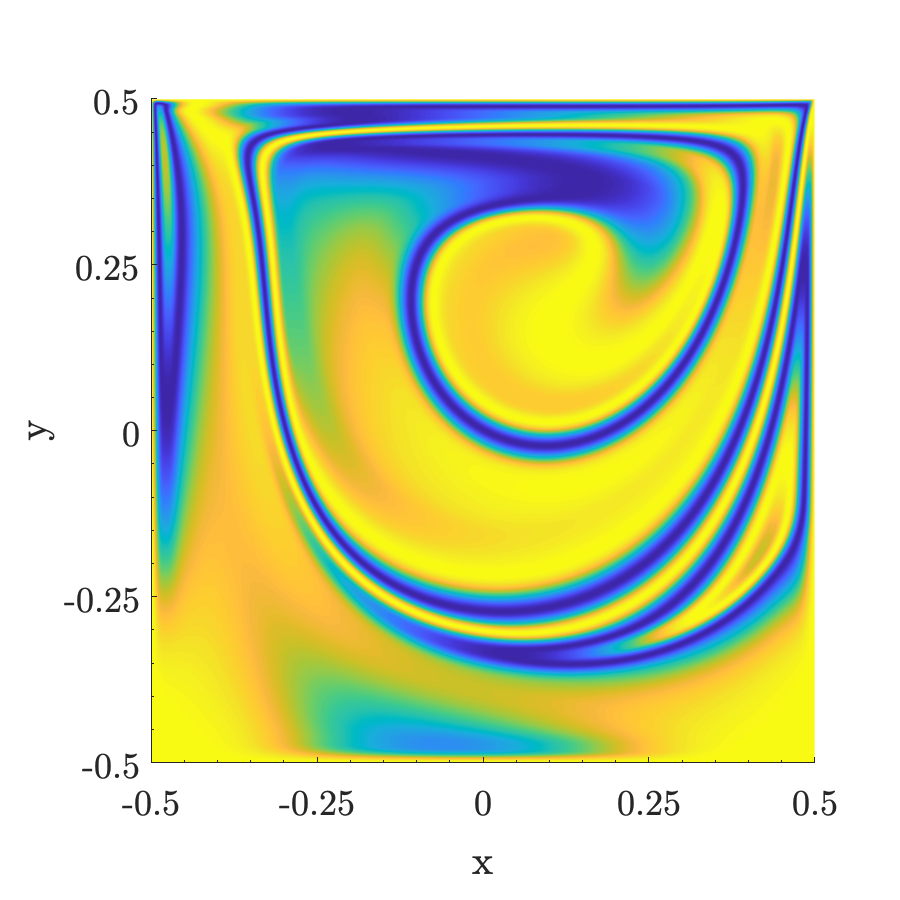}}
	{\includegraphics[width=.45\textwidth]{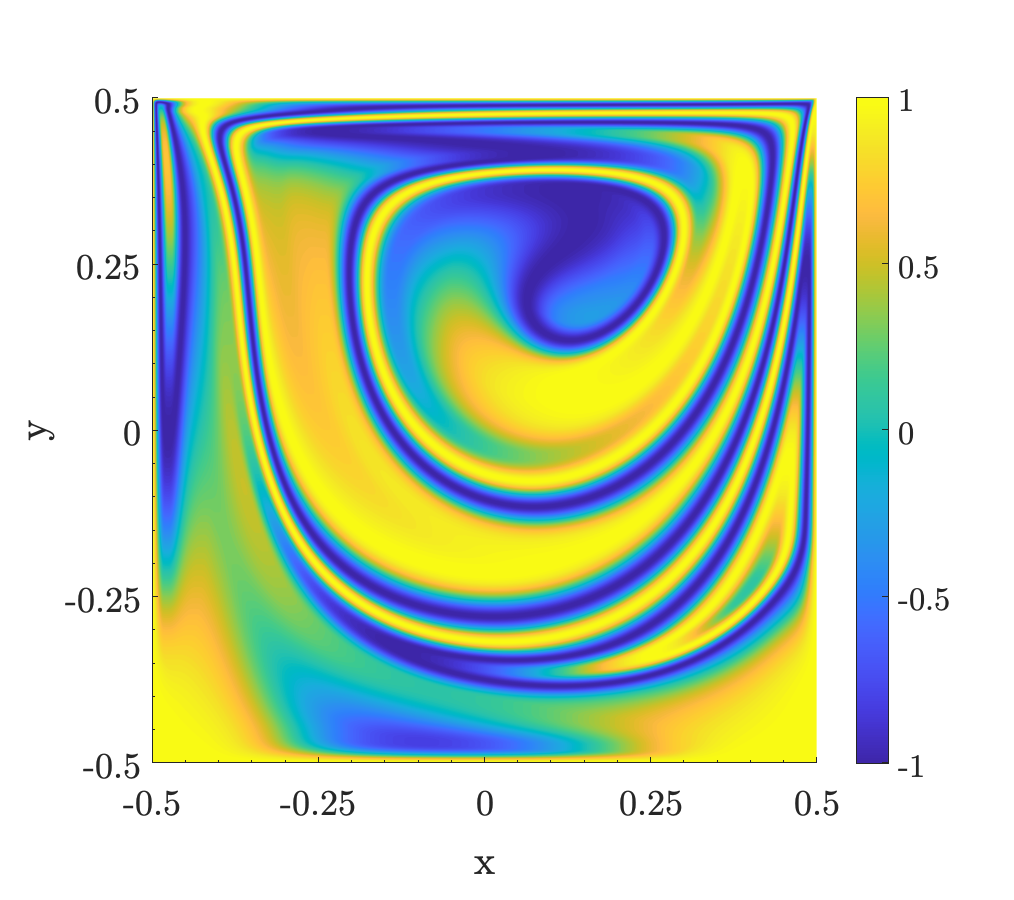}}
	\caption{\footnotesize Filled contours of one component of the distortion field $\A_1$, namely of the $A_{1,12}$ component, for the lid-driven cavity problem at $\reynolds = 100$, at times $t=1.0,  t=2.0, $ (top), $t=4.0,  t=6.0,$ (center) and $t=8.0,  t=10.0$ (bottom).}
	\label{fig:Lid2}
\end{figure}

Figure \ref{fig:Lid1} shows the computational results obtained using the
approach described in section \ref{remarks.dist.field}. This approach separates
the evolution of the two types of information encoded in $\A_a$ and leverages
the capabilities of a semi-analytical solver to efficiently solve the equations
in the stiff relaxation regime. Excellent agreement between the numerical
solution and the Navier-Stokes reference solution of Ghia et al. \cite{Ghia1982}
was obtained. Also for this test, Fig. \ref{fig:Lid2} shows the time evolution
of the $A_{1,12}$ component of the distortion field. It can again be seen that
the distortion field components are excellent candidates for flow
visualisation, revealing in detail the evolution of the flow and keeping track
of the rotations that the fluid element undergoes over time.
\subsection{Elastic vibrations of a beryllium plate}\label{berylium.plate} From
this section, we begin to test the multiphase property of the governing
equations and the ability of the model to describe the solid and fluid branches
of continuum mechanics in a single PDE system. In the following test problem, we simulate the vibrations of an edely elastic beryllium plate subjected to an
initial velocity perturbation. The setup follows \cite{ShashkovCellCentered,
Mairesolid, burton2015reduction, BoscheriDumbser2016a, Hyper-Hypo2018}, but with
the notable modification. Namely, instead of considering a solid body in a vacuum, we consider a solid body in a gas, and thus we need to define two separate density fields via their respective volume fractions.

Compared to the Lagrangian setup, used in the previously mentioned works, the computational domain considered here is larger, as in \cite{frontiers}, and is assumed to be  $\Omega=[-4.0 ; 4.0]\times[-2.0 ; 2.0]$ and the initial conditions for the first phase (the solid) are 
\begin{align}
\begin{split}
&\alpha_1(x,y)  = \begin{cases}	1- 2 \epsilon  &\mathrm{if}  \ \  \xx \in \Omega_1 , \\
						\epsilon\ \ & \mathrm{if}  \ \  \xx \notin \Omega_1,
			\end{cases} \quad 
	    \vv_1(x,y) = \begin{cases}	(0, \vv_{1,2})  &\mathrm{if}  \ \  \xx \in \Omega_1 , \\
						(0, 0) \ \ & \mathrm{if}  \ \  \xx \notin \Omega_1,
			\end{cases}\\ &\rho_1 = 1.845 ,  \quad p_1 = 10^{-4}, \quad \A_1 = \Id, 
\end{split}
\end{align} 
while the second phase (the gas) is initialised as follows
\begin{align}
\begin{split}
&\alpha_2(x,y)  = \begin{cases}	\epsilon  &\mathrm{if}  \ \  \xx \in \Omega_1 , \\
						1- 2 \epsilon\ \ & \mathrm{if}  \ \  \xx \notin \Omega_1,
			\end{cases}  \\
	   & \rho_2 = 10^{-3} , \quad \vv_2 = \vec{0} \quad p_2 = 10^{-4}, \quad \A_2 = \Id, 
\end{split}
\end{align}
where $\Omega_1 = [-3.0 ; 3.0]\times[-0.5 ; 0.5] $ is the subdomain that defines the initial geometry of the beryllium bar, and the initial vertical velocity component $ \vv_{1,2}$, according to Boscheri et al. \cite{BoscheriDumbser2016a}, is given as
\begin{equation}
\vv_{1,2}(x) = C_1 \omega  \big(  C_2 (\mathrm{sinh}(C_3(x+3)) +  \mathrm{sin}(C_3(x+3))  ) 
	 - C_4( \mathrm{cosh}(C_3(x+3)) +  \mathrm{cos}(C_3(x+3))   )  \big) 
\end{equation}
with $C_3 = 0.7883401241$, $C_2 = 0.2359739922$,  $C_1 = 0.004336850425$, $C_4 = 57.64552048$ and $C_2 = 56.53585154$. The third phase has $\alpha_3 = \epsilon$. The other parameters and physical quantities that define the properties of the beryllium material and that are required to use the stiffened-gas EOS \eqref{eqn.stiffened} are chosen as $\gamma_1 = 1.4$, $\Cv_1 = 1000$, $\Cs_1 = 0.905$, $\rhoo_a = 1.845$, $\Co_1 = 1.287$ and $\po_1= p_1$. Additinally, to have an ideal elastic material we set $\tau^\mathrm{e}_1 = 10^{14}$. For the gas phase surrounding the solid phase, the EOS of ideal gases is used and the physical parameters are $\gamma_2 = 1.2$, $\Cv_2 = 1000$, $\Cs_2 = 1.0$ and $\nu_2 = 10^{-4}$.

The simulation is carried out up to the final time $t_f =53.2$ which corresponds approximately to two complete periods of vibration and the computational domain is discretized with an uniform Cartesian mesh composed of $1024\times512$ control volumes. In contrast to Lagrangian schemes, it is not necessary to impose boundary conditions on the surface of the solid, as the solid-gas boundary condition is directly taken into account within the governing PDE system. Hence, in our simulation, periodic boundaries are set everywhere. 

\begin{figure}[!htbp]
	\renewcommand{\figurename}{\footnotesize{Fig.}}
	\centering
	{\includegraphics[width=.45\textwidth]{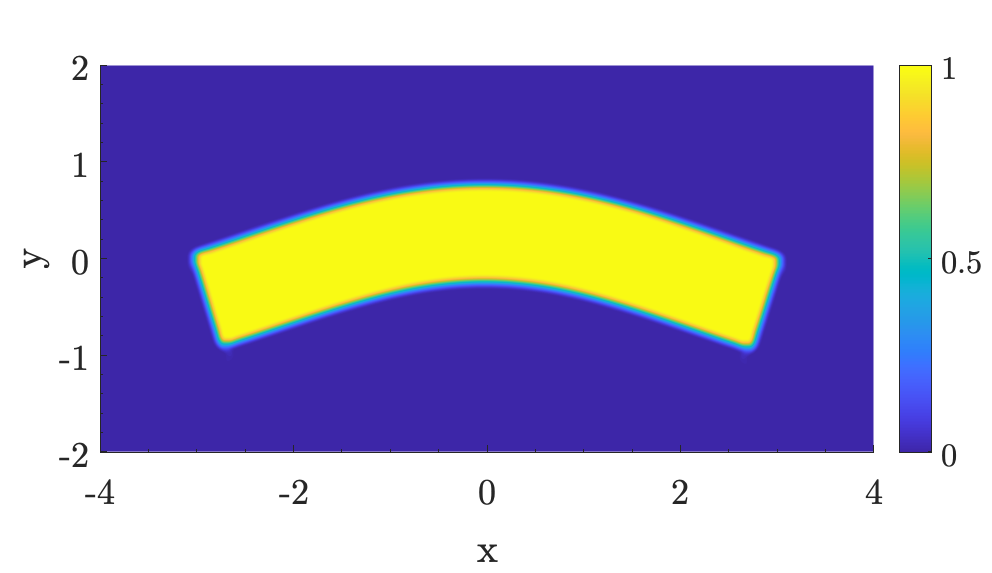}}
	{\includegraphics[width=.45\textwidth]{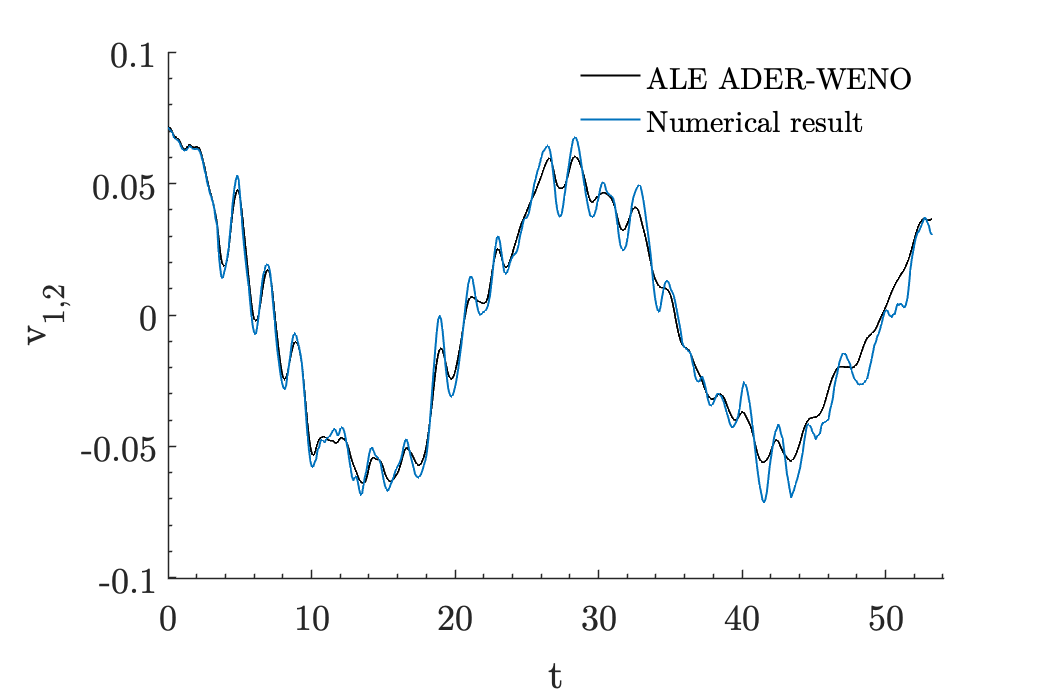}} 	
	\caption{\footnotesize Filled contour map of the volume fraction function $\alpha_1$ for the first phase, which represents the geometry of the beryllium bar at time $t = 8$ (left). The time evolution of the vertical velocity component $v_{1,2}(0, 0, t)$ at $\vec{\bar{x}} = (0, 0)$, i.e. in the barycenter of the bar (right).}
	\label{fig:Bp1}
\end{figure}

\begin{figure}[!htbp]
	\renewcommand{\figurename}{\footnotesize{Fig.}}
	\centering
	{\includegraphics[width=.45\textwidth]{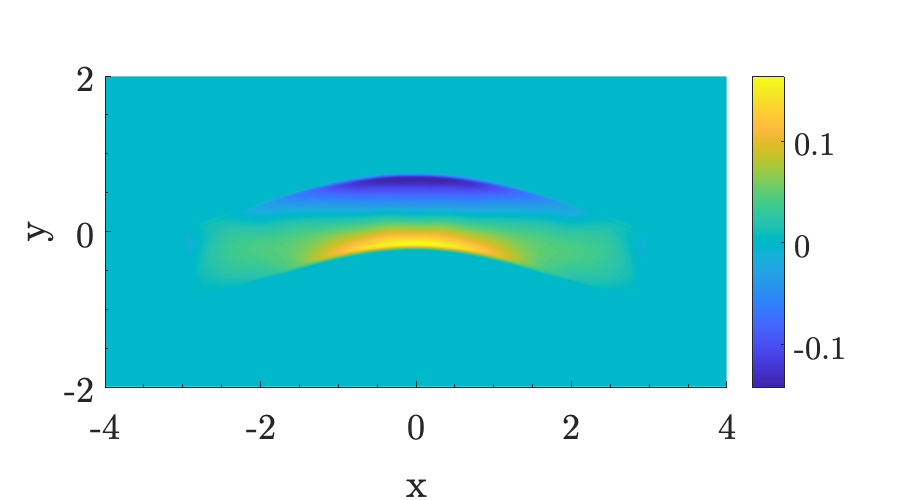}}
	{\includegraphics[width=.45\textwidth]{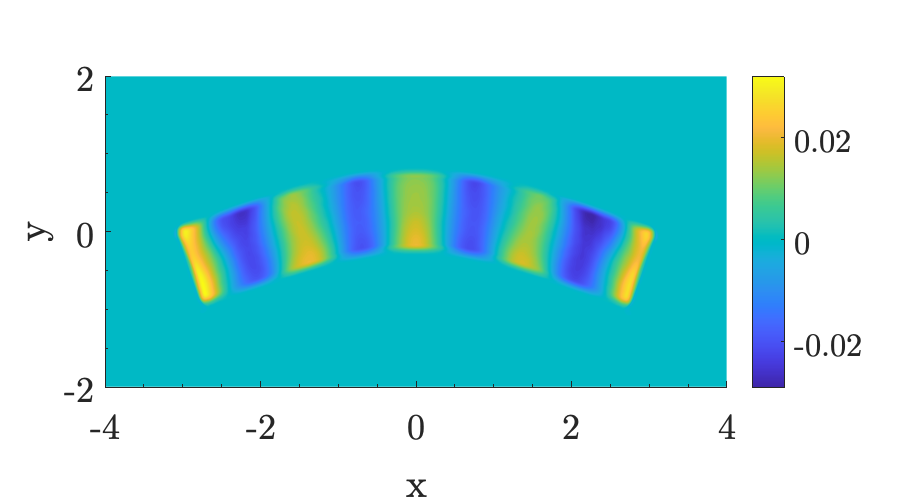}} 	
	{\includegraphics[width=.45\textwidth]{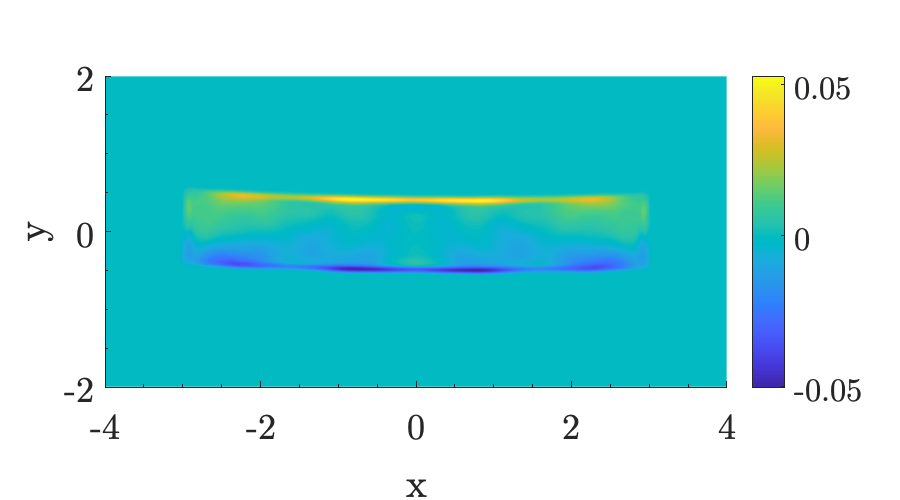}}
	{\includegraphics[width=.45\textwidth]{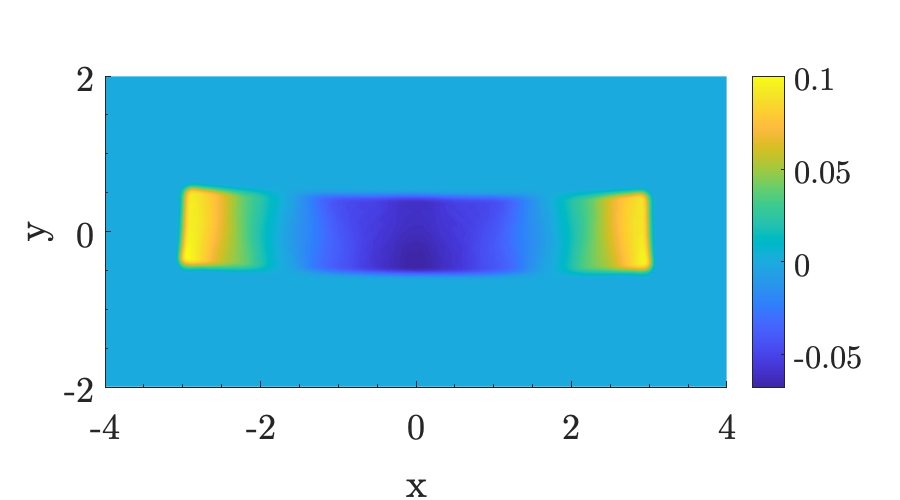}} 
	{\includegraphics[width=.45\textwidth]{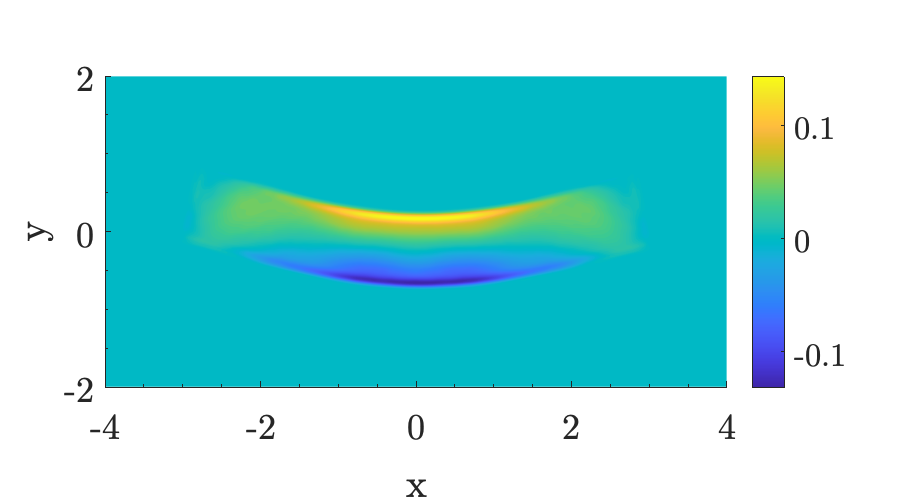}}
	{\includegraphics[width=.45\textwidth]{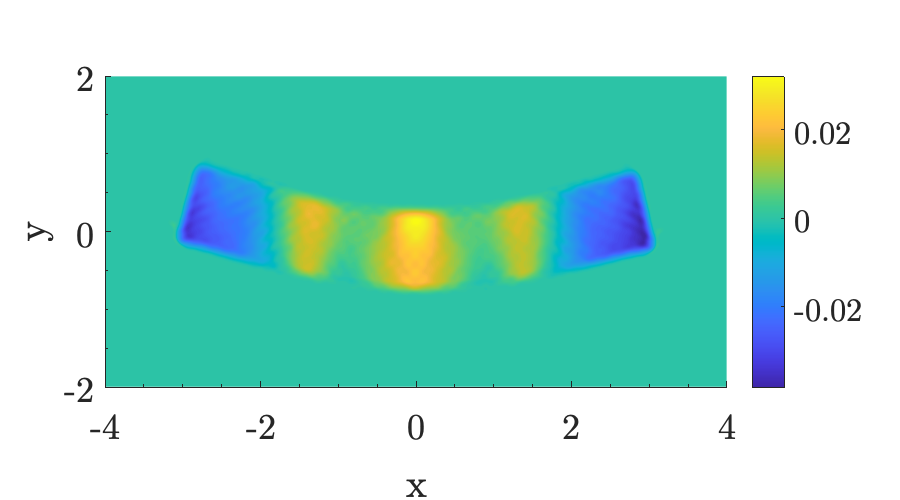}}
	{\includegraphics[width=.45\textwidth]{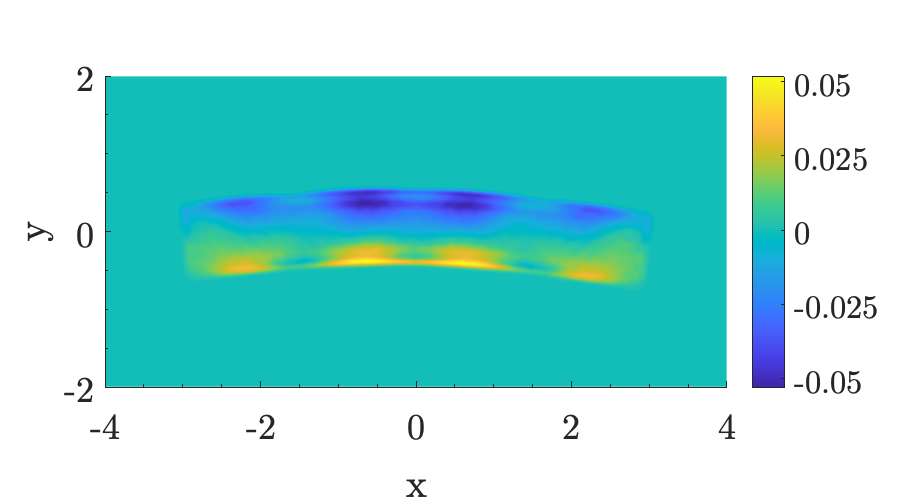}}
	{\includegraphics[width=.45\textwidth]{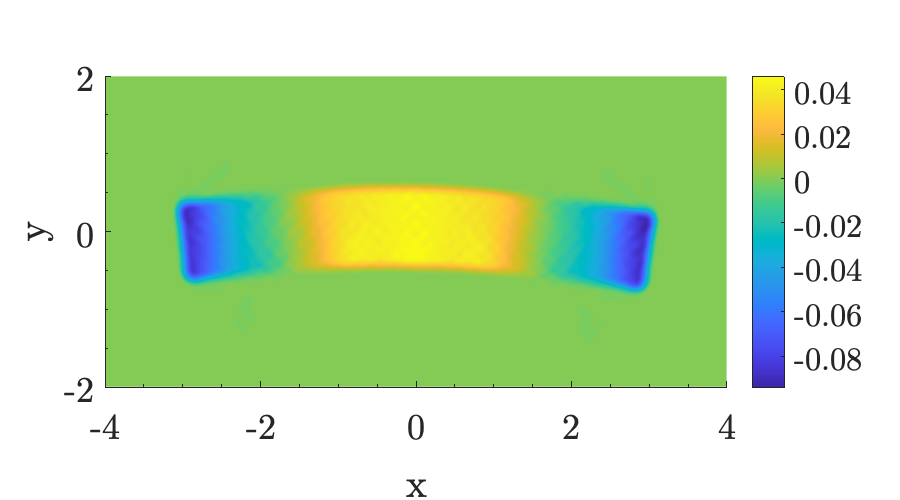}}
	\caption{\footnotesize  Results for the elastic vibrations of a beryllium plate, at times $t = 8$, $t = 15$, $t = 23$ and $t = 30$ (from top to bottom), for the first component of the stress tensor $\sigma_{1,11}$ and the vertical component of the velocity $v_{1,2}$.}
	\label{fig:Bp2}
\end{figure}

In Fig.\,\ref{fig:Bp1}, we represent the contour map of the volume fraction
$\alpha_1$, which represents the geometry of the bar at time $t = 8$ and in the
same figure, we also depict the time evolution of the vertical velocity
component $v_{1,2}(0, 0, t)$ at $\bar{\xx} = (0, 0)$, i.e. in the barycenter
of the bar. For comparison, we also show the results obtained with a
third-order ALE ADER-WENO scheme (black line), with which our numerical
solution (blue line) is in good agreement.
In Fig.\,\ref{fig:Bp2}, the first component of the stress tensor $\sigma_{1,11}$ and the vertical component of the velocity $v_{1,2}$ are shown on the left and right panels respectively for the intermediate times $t = 8$, $t = 15$, $t = 23$ and $t = 30$, covering approximately one bending period.
Despite the fact that in our case we use a diffuse interface approach on a fixed
Cartesian grid, our computational results compare visually well with the
reference solutions available in the literature, e.g.
\cite{ShashkovCellCentered, Mairesolid, burton2015reduction,
BoscheriDumbser2016a, Hyper-Hypo2018}, which all have been simulated with a pure
Lagrangian or arbitrary Lagrangian-Eulerian scheme on moving meshes.

%
%
\subsection{Taylor bar impact}\label{Taylor.bar} In the previous test, we
considered an ideal elastic material, which is the limit case for
$\tau_1\rightarrow \infty$. In the following test, we show how a material can
also exhibit non-linear elastic-plastic behaviour. Here we consider the Taylor
bar impact problem, which is a classical benchmark for an elasto-plastic target
that impacts on a rigid solid wall, e.g. see \cite{ShashkovCellCentered, Mairesolid,
2014JCoPh, BoscheriDumbser2016a} for pure Lagrangian or ALE schemes on moving
meshes and \cite{frontiers} for an Eulerian diffuse interface approach. 

As in the previous test, we define two separate density fields representing a
gas and a solid phase through their respective volume fractions. The
computational domain considered here is larger than in the Lagrangian setup in
order to include the space occupied by the gas phase around the solid. The
computational domain under consideration, as in \cite{frontiers}, is $\Omega=
[-150, +150]\times[0,600]$ and the initial conditions for the solid phase are 
\begin{align}
\begin{split}
&\alpha_1(x,y)  = \begin{cases}	1- 2 \epsilon  &\mathrm{if}  \ \  \xx \in \Omega_1 , \\
						\epsilon\ \ & \mathrm{if}  \ \  \xx \notin \Omega_1,
			\end{cases}   \quad
	    \vv_1(x,y) = \begin{cases}	(0, \vv_{1,2} )  &\mathrm{if}  \ \  \xx \in \Omega_1 , \\
						(0, 0) \ \ & \mathrm{if}  \ \  \xx \notin \Omega_1,
			\end{cases} \\ & \rho_1 = 2.785 , \quad p_1 = 10^{-4}, \quad \A_1 = \Id, 
\end{split}
\end{align}
where $\Omega_1 = [-50,+50] \times[0, 500] $ is the subdomain that defines the
initial geometry of the solid bar, and the initial vertical velocity component
is $ \vv_{1,2} = -0.015$; while the gas phase is initialised as follows
\begin{align}
\begin{split}
&\alpha_2(x,y)  = \begin{cases}	\epsilon  &\mathrm{if}  \ \  \xx \in \Omega_1 , \\
						1- 2 \epsilon\ \ & \mathrm{if}  \ \  \xx \notin \Omega_1,
			\end{cases}  \\
	   &  \rho_2 = 10^{-3} ,\quad  \vv_2 = \vec{0} \quad P_2 = 10^{-4}, \quad \A_2 = \Id.
\end{split}
\end{align}
According to \cite{ShashkovCellCentered, Mairesolid, 2014JCoPh,
BoscheriDumbser2016a}, the solid projectile is assumed to be an aluminium bar,
with following material parameters $\gamma_1 = 1.4$, $\Cv_1 = 1000$, $\Cs_1 =
0.305$, $\rhoo_a = 2.785$, $\Co_1 = 0.533$ and $\po_1= p_1$. To obtain a
non-linear elasto-plastic material behaviour the relaxation time
$\tau^\mathrm{e}_1$ has to be choosen as a non-linear function of an invariant
of the shear stress tensor as follows \cite{BartonRom2010, Hyper-Hypo2018}
\begin{align}\label{rel.time.plast}
\tau^\mathrm{e}_1 = \tau_{\mathrm{\scriptsize{o}}}  \left(\frac{{\sigma}_{\mathrm{\scriptsize{o}}} }{\bar{\sigma}_1}\right)^m,
\end{align}
where $\tau_{\mathrm{\scriptsize{o}}}= 1$ is the scaling constant, ${\sigma}{_\mathrm{\scriptsize{o}}} = 0.003$ is the yield stress of the material under quasi static conditions, the exponent parameter is chosen equal to $m = 20$ (the higher $m$ is the less rate-dependent the effective yield strength is \cite{BartonRom2010,Hyper-Hypo2018}) and the von Mises stress ${\bar{\sigma}_1}$ is given by
\begin{align}\label{eq.Mises}
{\bar{\sigma}_1} = \bigg( \frac{1}{2}\big( (\sigma_{1,11} - \sigma_{1,22})^2 + (\sigma_{1,22} - \sigma_{1,33})^2 + (\sigma_{1,33} - \sigma_{1,11})^2  
+ 6( \sigma_{1,21}^2 + \sigma_{1,31}^2 + \sigma_{1,32}^2 )  \big) \bigg)^{1/2}.
\end{align}
For the surrounding gas phase, the ideal gas EOS with the parameters $\gamma_2 = 1.2$, $\Cv_2 = 1000$, $\Cs_2 = 1.0$ and $\nu_2 = 10^{-4}$ is used.

The simulation is carried out up to the final time $t_f =5000$, and the
computational domain is discretized with a uniform Cartesian mesh composed of
$2048\times1024$ control volumes. In contrast to the Lagrangian schemes, it is
not necessary to impose boundary conditions on the surface of the solid; in our
simulation, periodic boundaries are set in the $x$-direction while reflective
slip wall boundary conditions are set in the $y$-direction.
\begin{figure}[!htbp]
	\renewcommand{\figurename}{\footnotesize{Fig.}}
	\centering
	{\includegraphics[width=.325\textwidth]{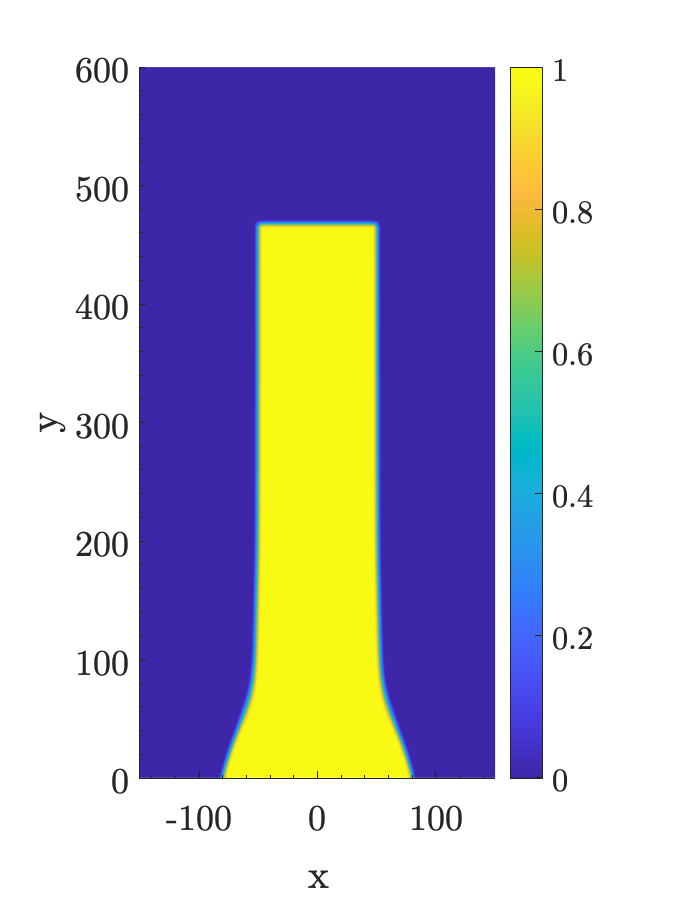}}
	{\includegraphics[width=.325\textwidth]{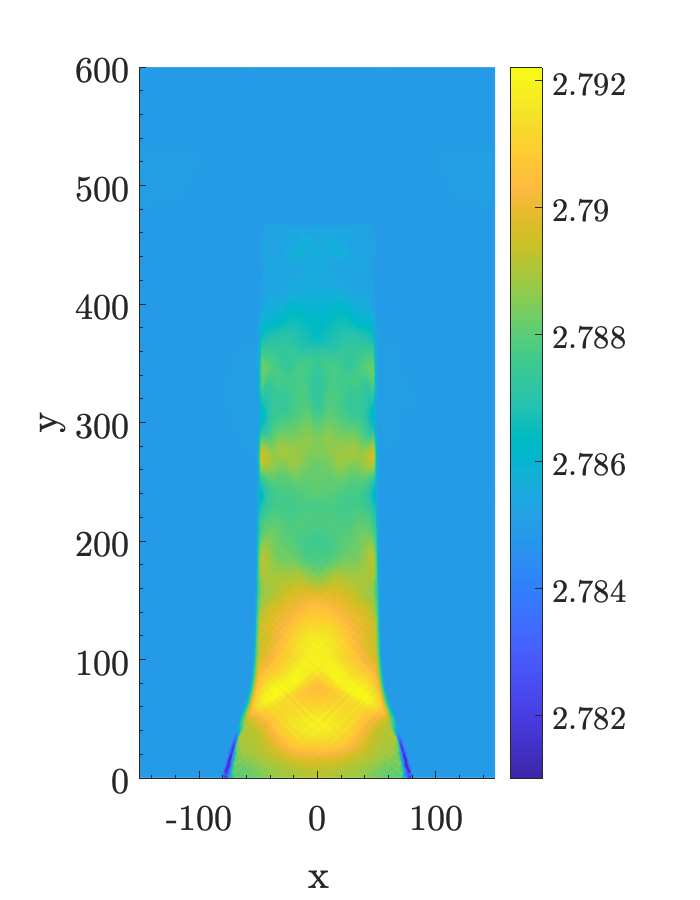}} 	
	{\includegraphics[width=.325\textwidth]{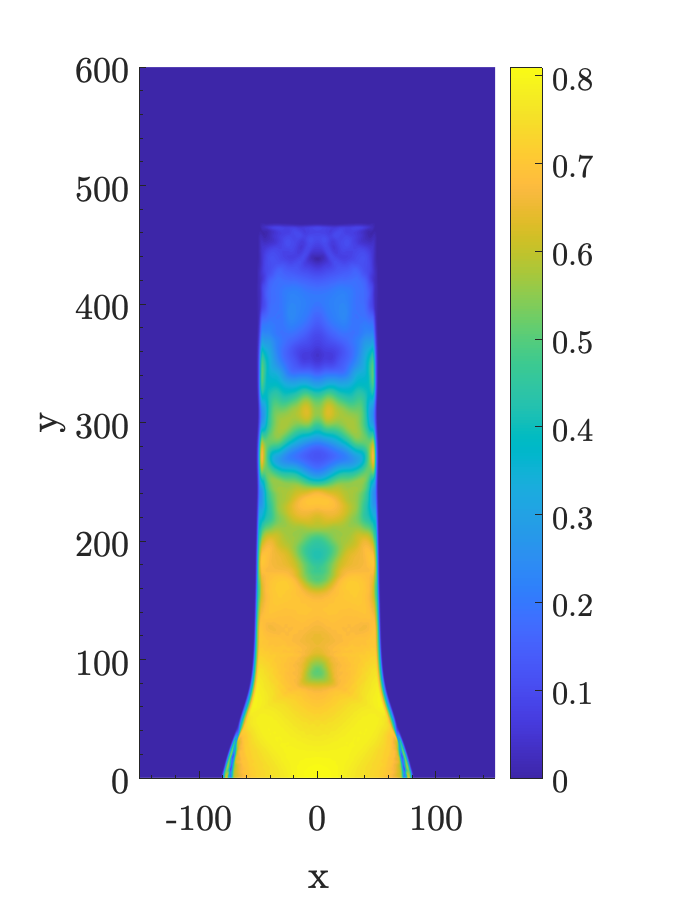}}
	{\includegraphics[width=.325\textwidth]{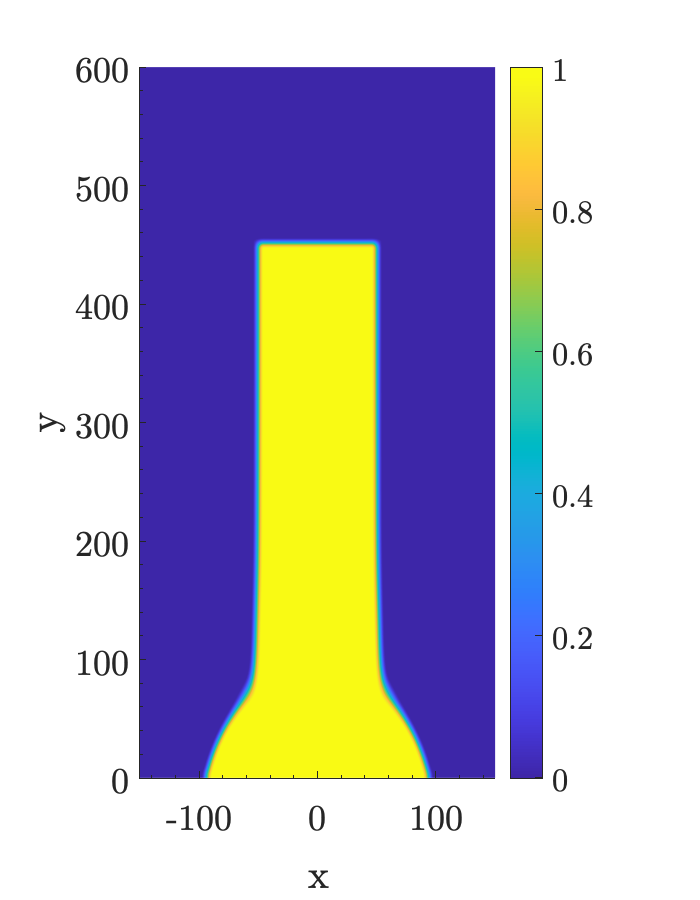}} 
	{\includegraphics[width=.325\textwidth]{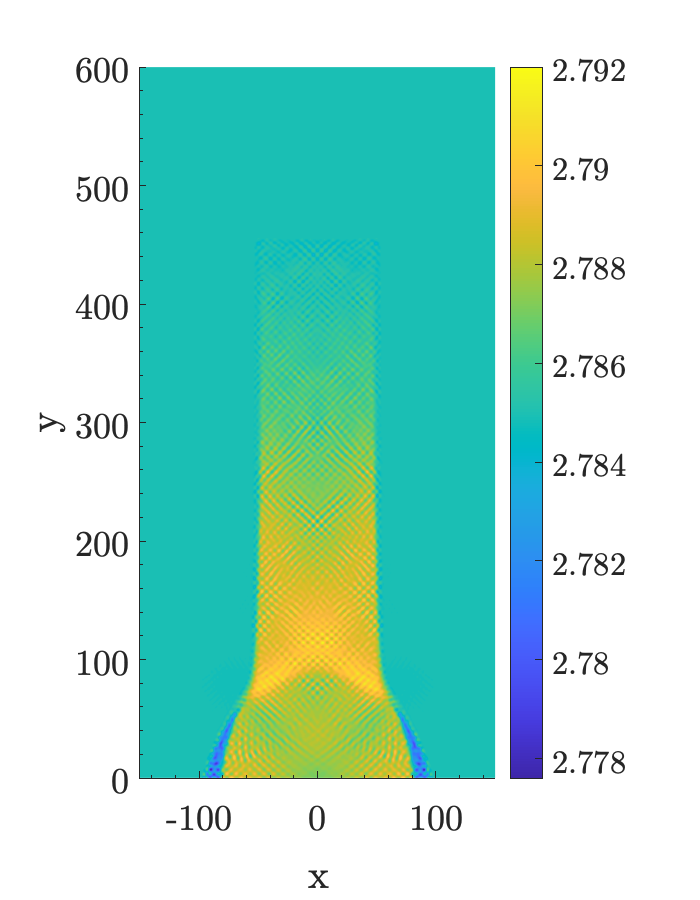}}
	{\includegraphics[width=.325\textwidth]{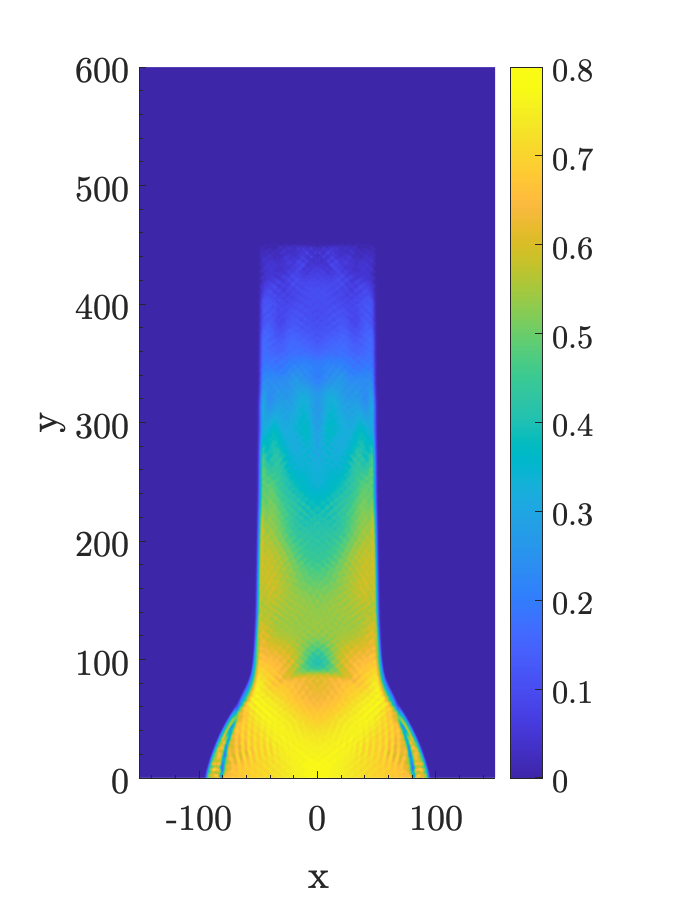}}
	\caption{\footnotesize Results for the non-linear elasto-plastic Taylor bar impact, at times $t = 2500$ and $t = 5000$ (from top to bottom): the volume fraction (left), the density distribution (center) and the plastic rate $\eta_1 = {\bar{\sigma}_1}/ {\sigma}{_\mathrm{\scriptsize{o}}}$ (right).}
	\label{fig:Tb1}
\end{figure}

In Fig.\,\ref{fig:Tb1}, we present the results computed at output times $t =
2500$ and $t = 5000$. The volume fraction (left), the density distribution
(center) and the plastic rate $\eta_1 = {\bar{\sigma}_1}/
{\sigma}{_\mathrm{\scriptsize{o}}}$ (right) are depicted. It can be observed
that the numerical solution is in good agreement with that presented in
\cite{Mairesolid,BoscheriDumbser2016a}, although the models used are
significantly different. Moreover at time $t = 5000$, the final length of the
aluminium bar is $L_f = 455$, which fits the results achieved in
\cite{Mairesolid,BoscheriDumbser2016a} within the $2\%$ error.
\subsection{Multiphase Rayleigh-Taylor instability}\label{Rayleigh.T} The two
previous tests focused mainly on the validation of the unified model for the
elastic and elasto-plastic type behavior, and not much emphasis was placed on
the dynamics of the fluid phase. In this test, however, we will finally put to
the test one of the main features of the model and the numerical scheme, namely
the ability to describe several interacting phases. To this end, we will
simulate a true \textit{two-phase} low-Mach, $M_1 \simeq M_2 \simeq 0.1$, viscous Rayleigh--Taylor
instability. The approach follows \cite{Chiocchetti2023}, but with the notable
modification that instead of initialising a single fluid with a jump in the
density, we define two fluids through the volume fraction and each fluid has a
constant phase density. This makes the problem more challenging because the
quasi-vacuum states of either phase are introduced almost throughout the entire
computational domain, however it introduces much more freedom in defining the
material characteristics of each phase. The computational domain under
consideration, as in \cite{Chiocchetti2023}, is $\Omega=  [0, 1/3]\times[0,1]$
and the initial conditions for the hevier phase (on top) are 
\begin{align}\label{in.cond.RT1}
\begin{split}
&\alpha_1(x,y)  = 	\bar{s}(1- 2 \epsilon) + (1-\bar{s}) \epsilon,     \qquad \rho_1 = 2.0,  \qquad \vv_1 = 	\vec{0}, \\
& p_1 = \bar{s}p_t + ( 1 - \bar{s} )p_b, \qquad \A_1 = \Id, 
\end{split}
\end{align}
and for the lighter fluid (at the bottom) are
\begin{align}\label{in.cond.RT2}
\begin{split}
&\alpha_2(x,y)  = 	1- \alpha_1- \epsilon,     \qquad \rho_2 = 1.0,  \qquad \vv_2 = \vec{0}, \\
& p_2 = \bar{s}p_t + ( 1 - \bar{s} )p_b, \qquad \A_2 = \Id, 
\end{split}
\end{align}
where $\bar{s}$ is a switch function introduced to impose a smooth transition between the two states and to avoid inaccurate representation of the initial condition on the discrete Cartesian grid. This function $\bar{s}$ is defined as
\begin{align}
	\bar{s} = \halb + \halb \mathrm{erf} \left( \frac{y-y_{\up{I}}}{\delta}		\right),
\end{align}
where $y_{\up{I}} = 0.5 + 0.01\mathrm{cos}(6\pi x)$ is the initially perturbed interface between phases and  $\delta= \max (0.004, 6 \Delta x)$ can be seen as the thickness of this interface. The initial top and bottom phase pressures, in \eqref{in.cond.RT1}, \eqref{in.cond.RT2}, are defined as 
\begin{align}
\begin{split}
& p_t  = 1+ \rho_1 (1 - y) \boldsymbol{g},  \\
& p_b = 1 + 0.5  \rho_1  \boldsymbol{g} +  \rho_2  ( 0.5 - y ) \boldsymbol{g},
\end{split}
\end{align}
where the gravity vector is $\boldsymbol{g}= (0, -0.1, 0)^\transpose$.
The other parameters and physical quantities are equal for the two gas phases and using the ideal-gas EOS are set as $\gamma = 1.4$, $\Cv = 1000$, $\Cs = 0.3$ and the dynamic viscosity $\mu = 6\times10^{-5}$, which translates to $\tau^\mathrm{e}_1 =2\times10^{-3} $ and $\tau^\mathrm{e}_2 =4\times10^{-3}$ so that the Reynolds number of the test problem is $\reynolds \simeq 2000$ for the both phases.
\begin{figure}[!htbp]
	\renewcommand{\figurename}{\footnotesize{Fig.}}
	\centering
	{\includegraphics[width=.32\textwidth]{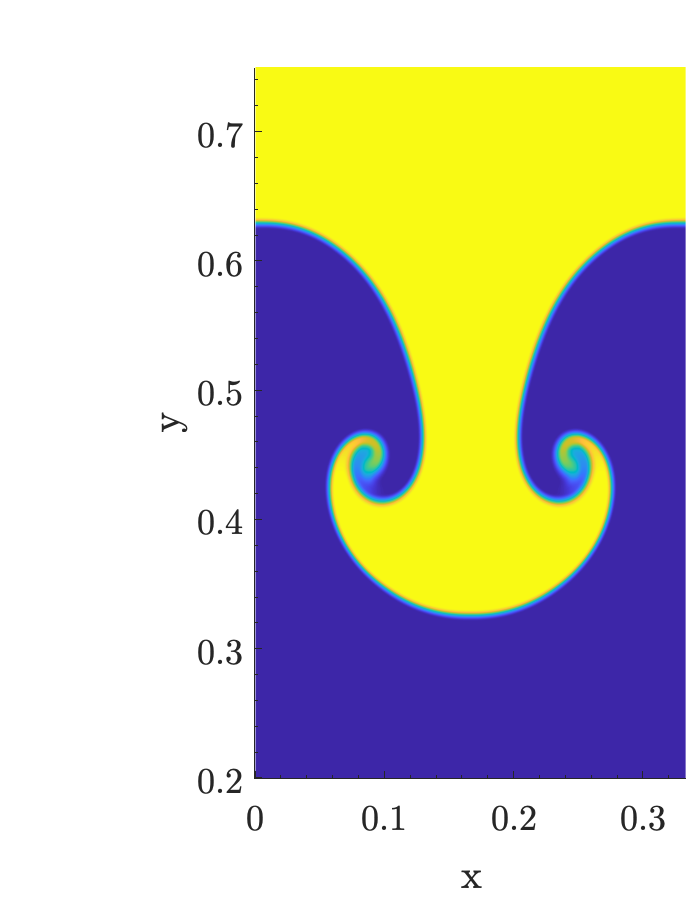}}
	{\includegraphics[width=.32\textwidth]{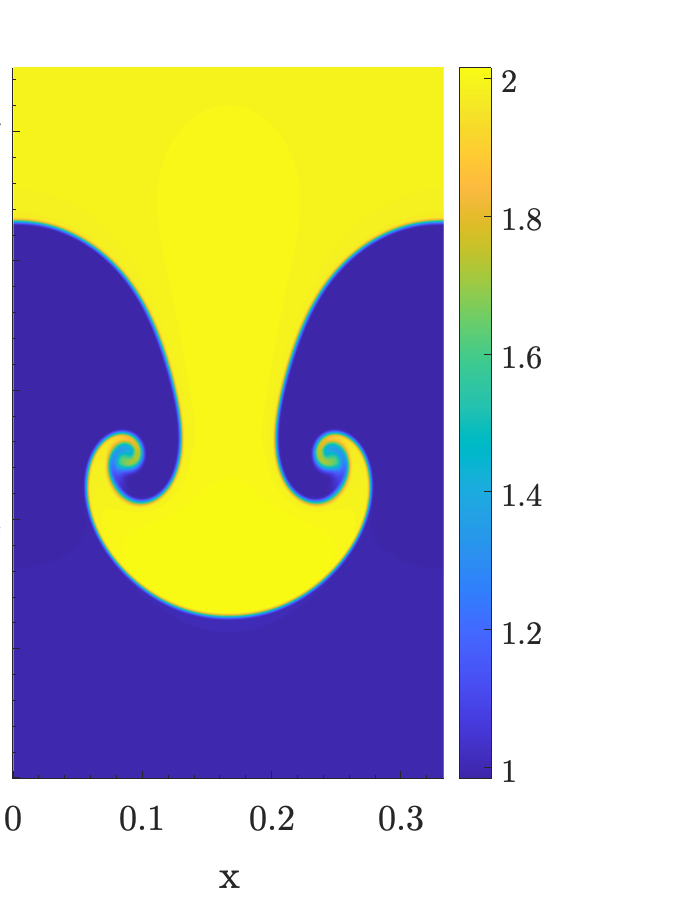}} 	
	\caption{\footnotesize Mixture density $\rho$. Mesh convergence test for the Rayleigh-Taylor instability problem, at times $t = 7$; on the left the result obtained with a mesh consisting of $512 \times 1536 $ is shown, and on the right the result obtained by doubling the mesh resolution.}
	\label{fig:Rt1}
		\vspace{-1em}
\end{figure}

Two simulations are carried out up to the final time $t_f =10$ on two different
uniform Cartesian meshes in order to verify mesh convergence of the solution
algorithm. In these simulations, periodic boundaries are set in $x$-direction
while reflective slip wall boundary conditions are set in $y$-direction.
Fig.\,\ref{fig:Rt1} shows on the left the result for the mixture density
$\rho = \alpha_1\rho_1+ \alpha_2\rho_2+ \alpha_3\rho_3$ obtained with a
mesh consisting of $512 \times 1536 $, while on the right the result obtained by
doubling the mesh resolution, both at time $t = 7$. It is possible to see that
mesh convergence has already been achieved with the coarsest mesh, since the
macroscopic structure of the flow does not depend on mesh size.

\begin{figure}[!htbp]
	\renewcommand{\figurename}{\footnotesize{Fig.}}
	\centering
	{\includegraphics[width=.32\textwidth]{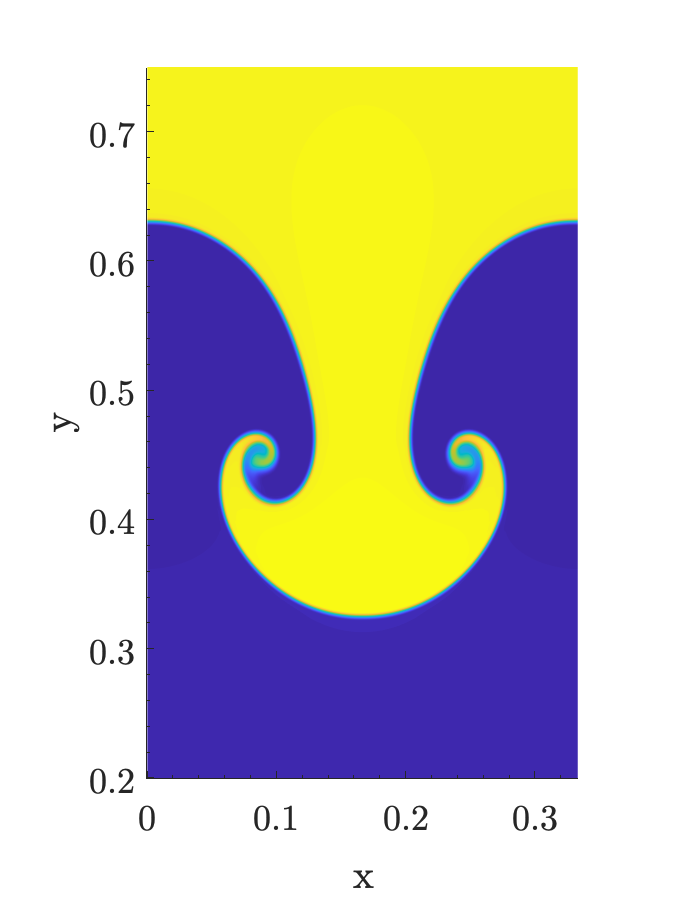}}
	{\includegraphics[width=.32\textwidth]{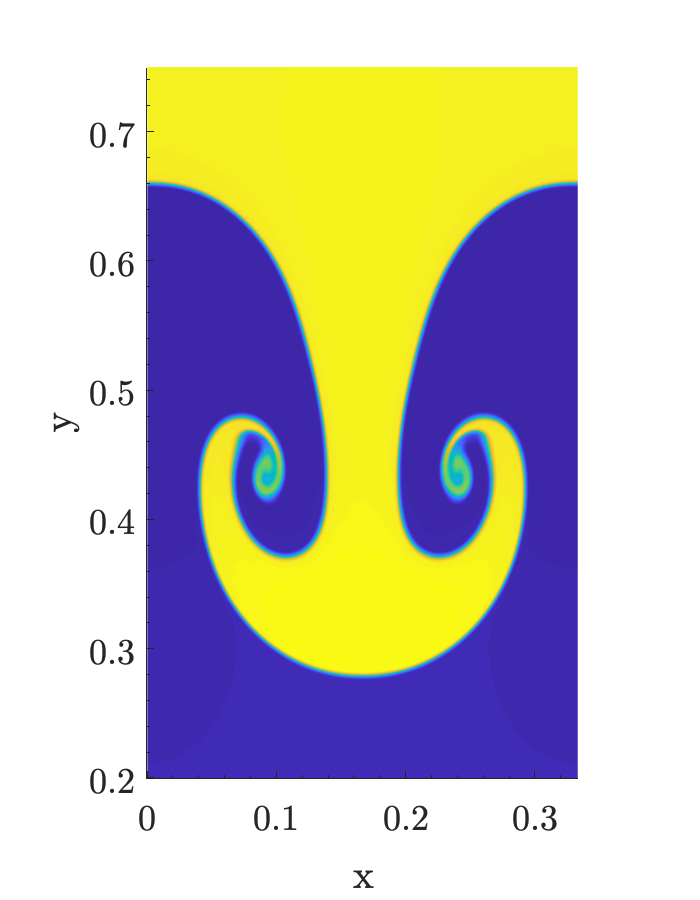}} 	
	{\includegraphics[width=.32\textwidth]{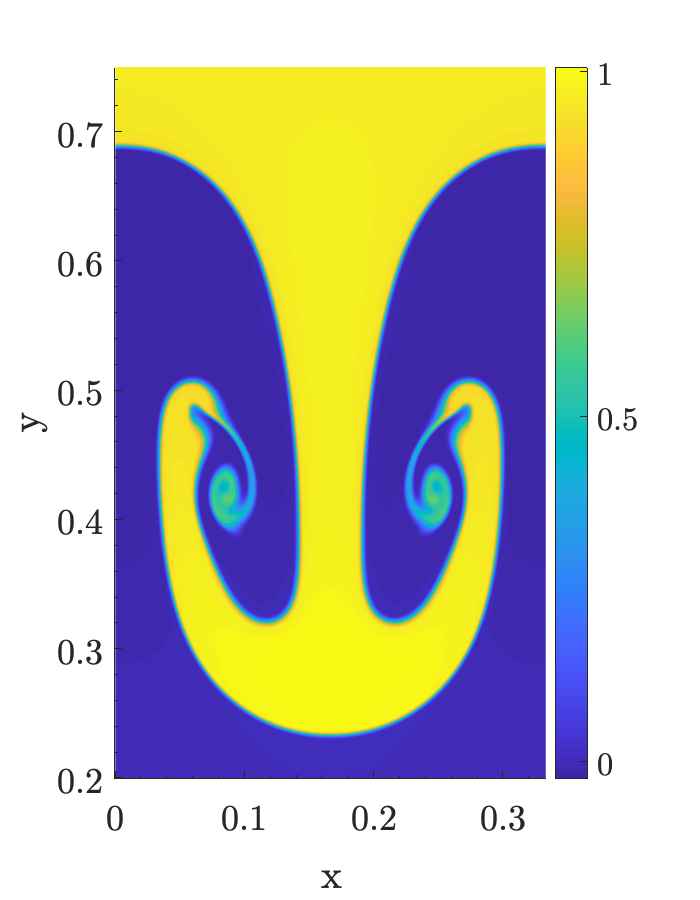}}
	{\includegraphics[width=.32\textwidth]{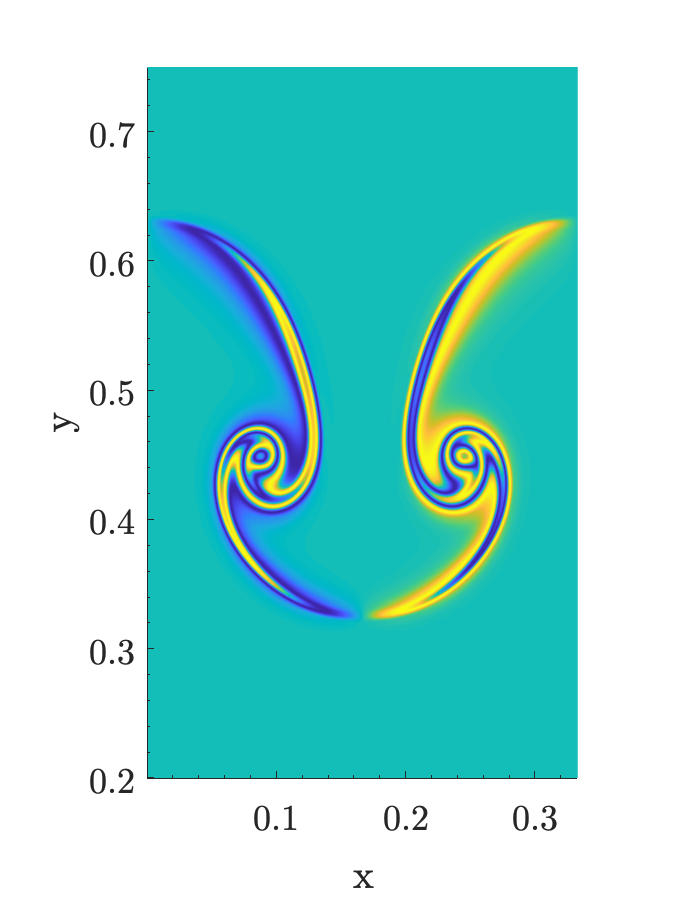}} 
	{\includegraphics[width=.32\textwidth]{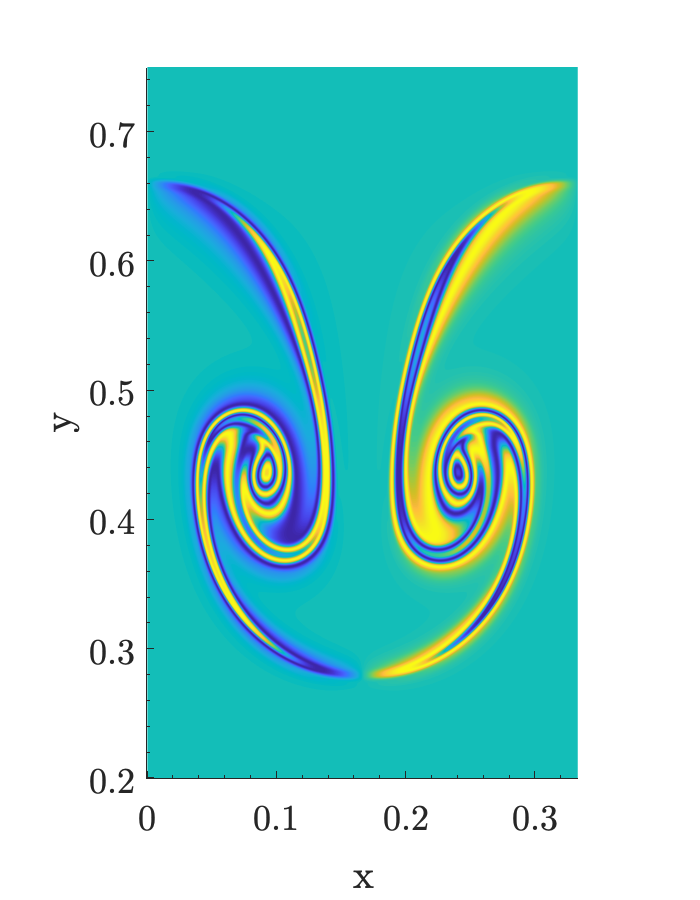}}
	{\includegraphics[width=.32\textwidth]{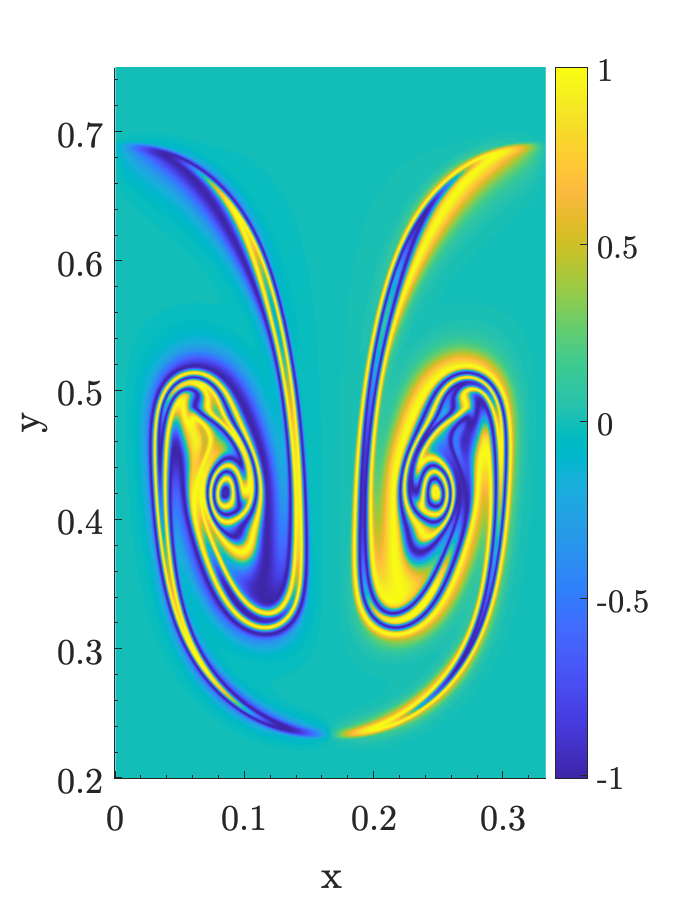}} 
	\caption{\footnotesize Results for the multiphase Rayleigh-Taylor instability problem, at times $t = 6$, $t=7$ and $t = 8$ (from left to right);  the volume fraction $\alpha_1$ (top) and the $A_{1,12}$  component of the distortion field for the first phase (bottom) are represented.}
	\label{fig:Rt2}
\end{figure}

Fig.\,\ref{fig:Rt2} shows the time evolution of both the volume fraction
function $\alpha_1$ and the $A_{1,12}$ component of the distortion field, for
the first phase. 
It is interesting to note that it for such a low Mach number as in this test,
the mixture density, depicted in Fig.\,\ref{fig:Rt1} is macroscopically
proportional to the volume fraction structure in Fig.\,\ref{fig:Rt2}. Moreover,
due to the velocity relaxation, both distortion fields encode the same flow
structure, except that they satisfy two different algebraic constraints, so that
each phase mass conservation equation is the consequence of the time evolution
of each phase distortion field, see \ref{det.constr}.

Our computational results, in particular the temporal evolution in
Fig.\,\ref{fig:Rt2} compare visually well with the reference solutions available
in the literature, see \cite{MaireCyl2}, obtained in this case with ALE schemes
on moving meshes.
\subsection{ Three-phase triple point problem}
Finally, in this section, we will test all the capabilities of the model and the
numerical scheme developed, i.e. the ability to describe several, up to three in
this paper, interacting phases. The problem that will be addressed is a typical
test in the ALE community, namely the so-called triple point problem. This test
is a three state, two material, 2D Riemann problem in a rectangular domain that
generates vorticity, which is why it is very popular in the ALE community for
testing the ability of a code to handle complex fluid motion on moving mesh. It was
introduced in \cite{ReALE2010} and was used to compare ALE approaches in the
case of a two-material Riemann problem in
\cite{MaireCyl2,ShashkovMultiMat3,Kucharik2014} or the simplified one-material
case in \cite{boscheri2015direct,ReALE2015,arepo}. In the Eluerian context, this
problem has been addressed with an interface-capturing approach considering
three immiscible compressible fluids in \cite{menshov1}. In the context of this
work, the main aim of this problem is to verify the ability of the code and
model to correctly propagate shock waves over multiphase and multi-material
regions.

Specifically, we follow the setting presented in \cite{ReALE2010}, but with the significant modification that instead of initialising a two-material Riemann problem, we set up an initial problem involving three phases, where two of them have the same material parameters. The computational domain of the triple point problem $\Omega = [0; 7] \times[0; 3]$ is divided into three subdomains filled with three phases describing different perfect gases, thus yielding a three-phase, two-material problem. The initial condition, in our diffuse interface framework, can be easily set by means of jumps in volume fraction as follows. The first phase, with a state of high pressure and high density, is initialized as
\begin{align}
\begin{split}
&\alpha_1(x,y)  = 	\begin{cases}	1- 2 \epsilon  &\mathrm{if}  \ \  \xx \in \Omega_1 , \\
						\epsilon\ \ & \mathrm{if}  \ \  \xx \notin \Omega_1,
			\end{cases}     \qquad  \mathrm{with} \qquad \Omega_1 =  [0; 1] \times [0; 3] \\
			&\rho_1 = 1.0,  \quad \vv_1 = \vec{0}, \quad p_1 = 1.0, \quad \A_1 = \Id, 
\end{split}
\end{align}
the second phase, with a state of low pressure high density, as 
\begin{align}
\begin{split}
&\alpha_2(x,y)  = 	\begin{cases}	1- 2 \epsilon  &\mathrm{if}  \ \  \xx \in \Omega_2 , \\
						\epsilon\ \ & \mathrm{if}  \ \  \xx \notin \Omega_2,
			\end{cases}     \qquad  \mathrm{with} \qquad \Omega_2 =  [1; 7]\times[0; 1.5] \\
			&\rho_2 = 1.0,  \quad \vv_2 = \vec{0}, \quad p_2 = 0.1, \quad \A_2 = \Id, 
\end{split}
\end{align}
and the third, with an initial low pressure and low density state, is initialized as
\begin{align}
\begin{split}
&\alpha_3(x,y)  = 	\begin{cases}	1- 2 \epsilon  &\mathrm{if}  \ \  \xx \in \Omega_3 , \\
						\epsilon\ \ & \mathrm{if}  \ \  \xx \notin \Omega_3,
			\end{cases}     \qquad  \mathrm{with} \qquad \Omega_3 =  [1; 7]\times[1.5;3.0] \\
			&\rho_3 = 0.125,  \quad \vv_3 = \vec{0}, \quad p_3 = 0.1, \quad \A_3 = \Id.
\end{split}
\end{align}
All the phases represent ideal inviscid gases, thus we are in the stiff limit of
the model and the relaxation times are $\tau^\mathrm{e}_1 = \tau^\mathrm{e}_2 =
\tau^\mathrm{e}_3 = 14^{-14}$. Furthermore, according to \cite{ReALE2010}, since
the first phase and the third phase represent the same material with the material parameters $\gamma_1 = \gamma_3 =
1.5$, $\Cv_1 = \Cv_3 = 1$, $\Cs_1 = \Cs_3 = 1.0$, while the second phase has the paramers $\gamma_2 = 1.4$, $\Cv_2 = 1$, $\Cs_2 = 1.0$. 

The simulation is carried out up to the final time $t_f =5$ discretizing the computational domain with a uniform Cartesian mesh composed of $3584\times1536$ control volumes; reflective slip wall boundary conditions are set in all the directions.

\begin{figure}[!b]
	\renewcommand{\figurename}{\footnotesize{Fig.}}
	\centering
	{\includegraphics[width=.64\textwidth]{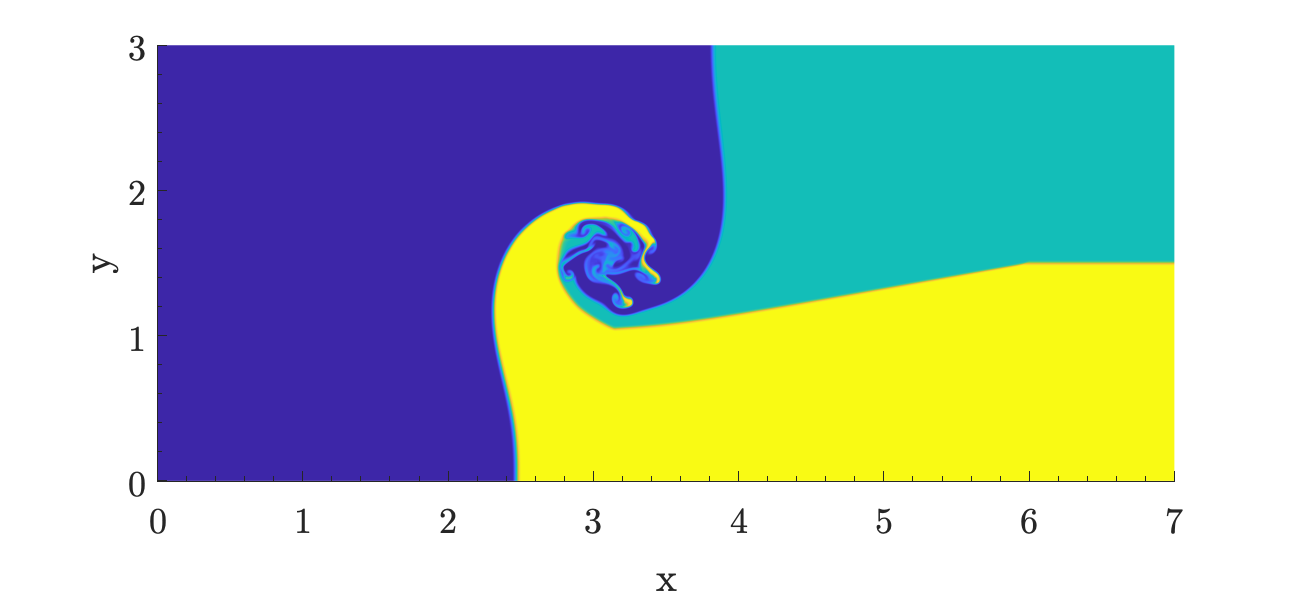}}
	{\includegraphics[width=.35\textwidth]{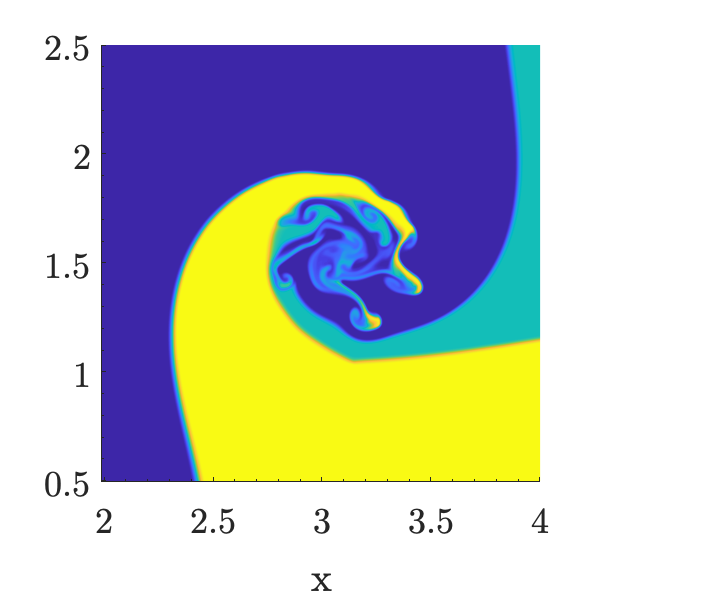}} 	
	{\includegraphics[width=.64\textwidth]{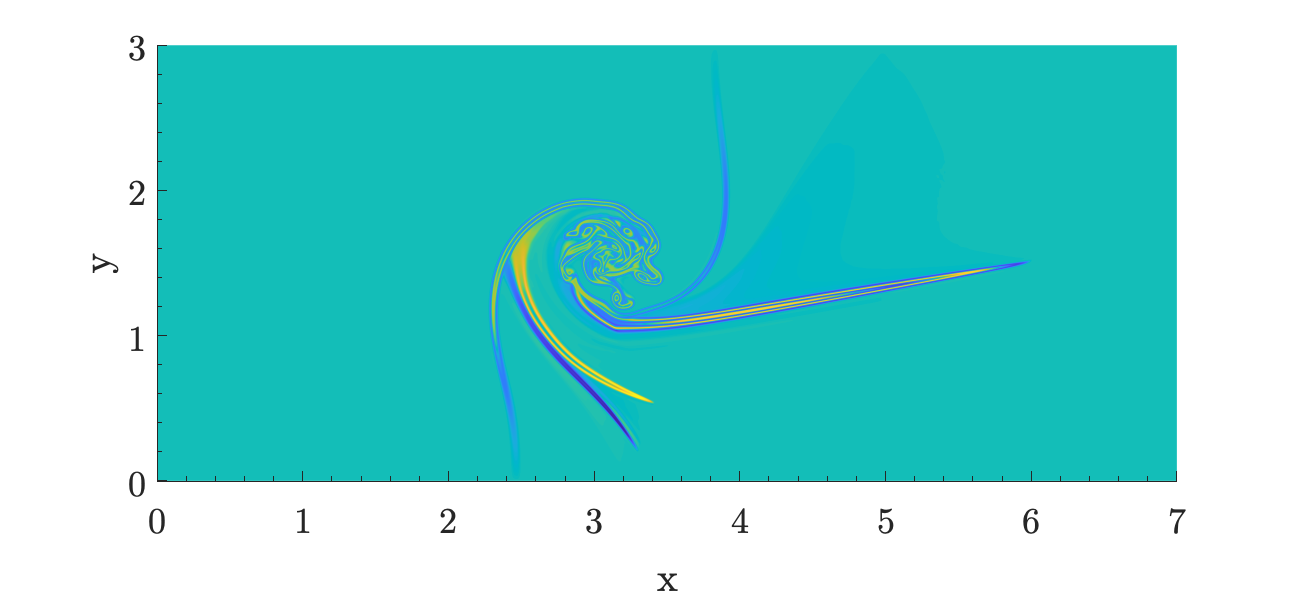}}
	{\includegraphics[width=.35\textwidth]{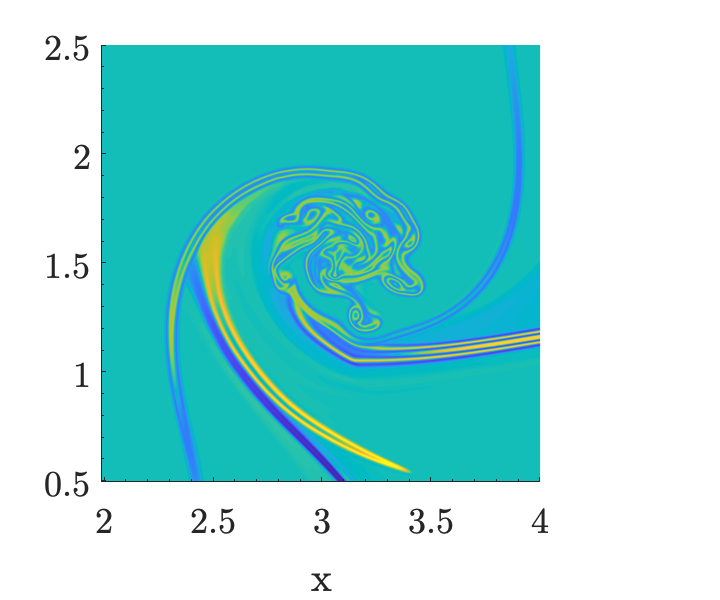}} 
	{\includegraphics[width=.64\textwidth]{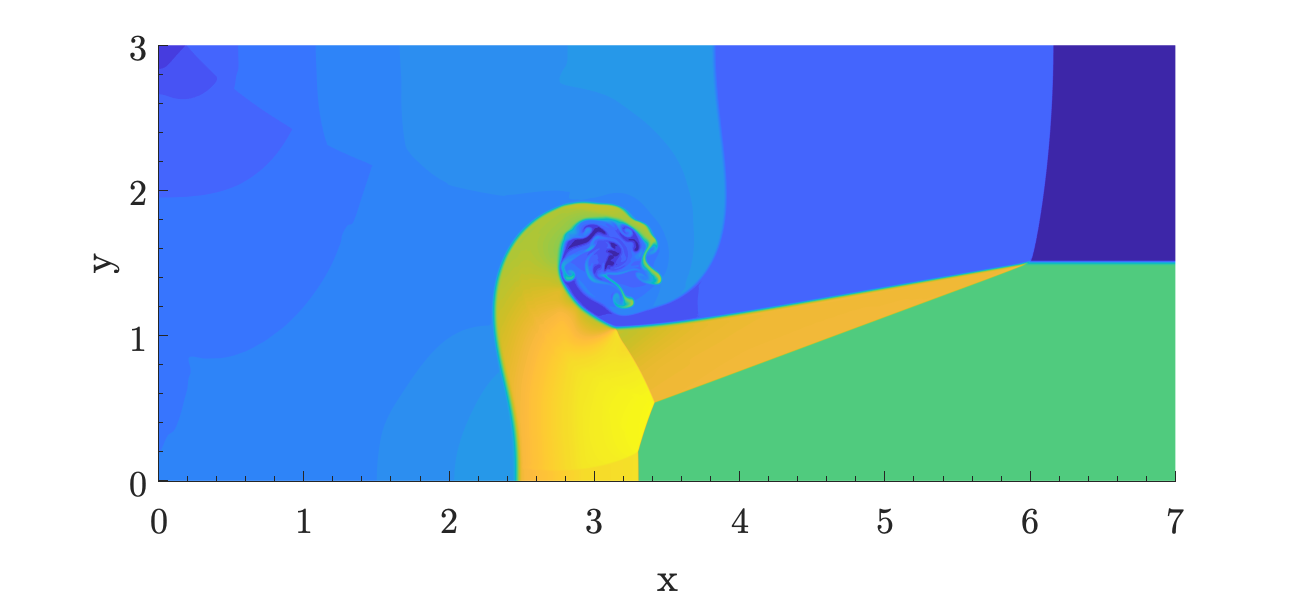}}
	{\includegraphics[width=.35\textwidth]{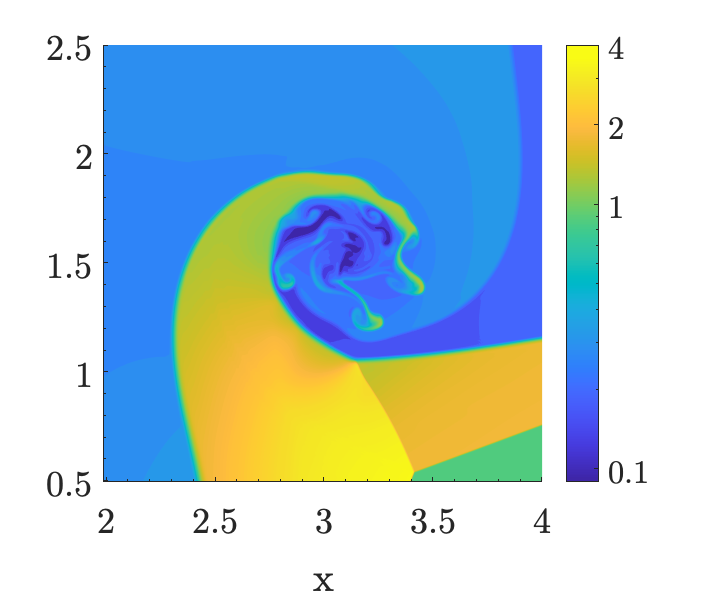}} 
	\caption{\footnotesize Results for the multiphase and multi-material triple point problem at time $t = 3$: the contour plots of the volume fractions (first phase in blue, second phase in yellow and the third one in blue-green) (top), the $A_{1,12}$ component of the distortion field for the first phase (center) and the mixture density (bottom). A zoom in the region $[2.0, 4.0] \times[0.5, 2.5]$, which illustrates the formation of vorticity resulting from the initial contact discontinuity, is shown (right).}
	\label{fig:3P3}
\end{figure}

\begin{figure}[!b]
	\renewcommand{\figurename}{\footnotesize{Fig.}}
	\centering
	{\includegraphics[width=.64\textwidth]{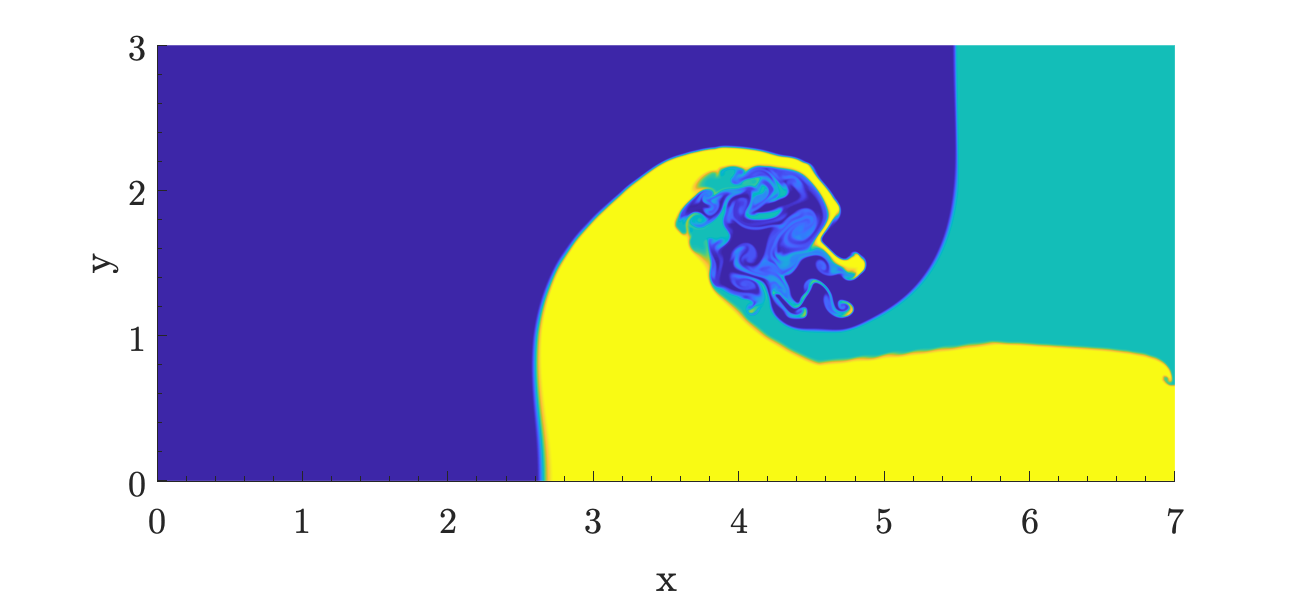}}
	{\includegraphics[width=.35\textwidth]{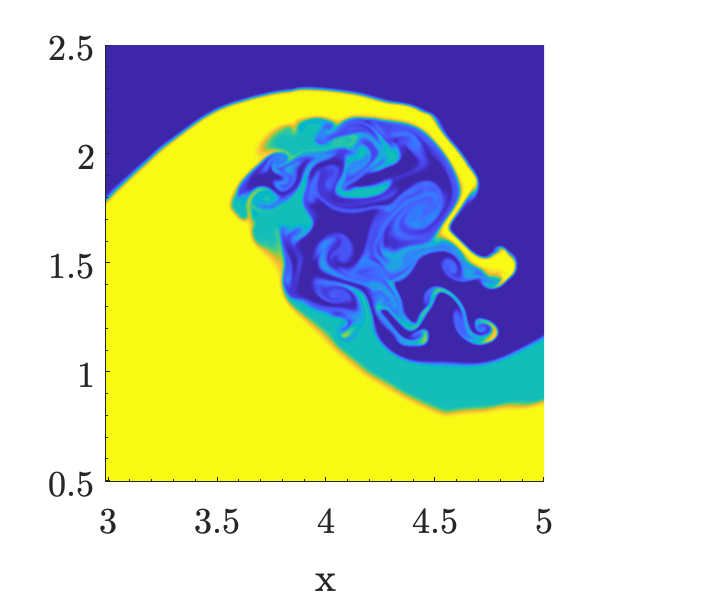}} 	
	{\includegraphics[width=.64\textwidth]{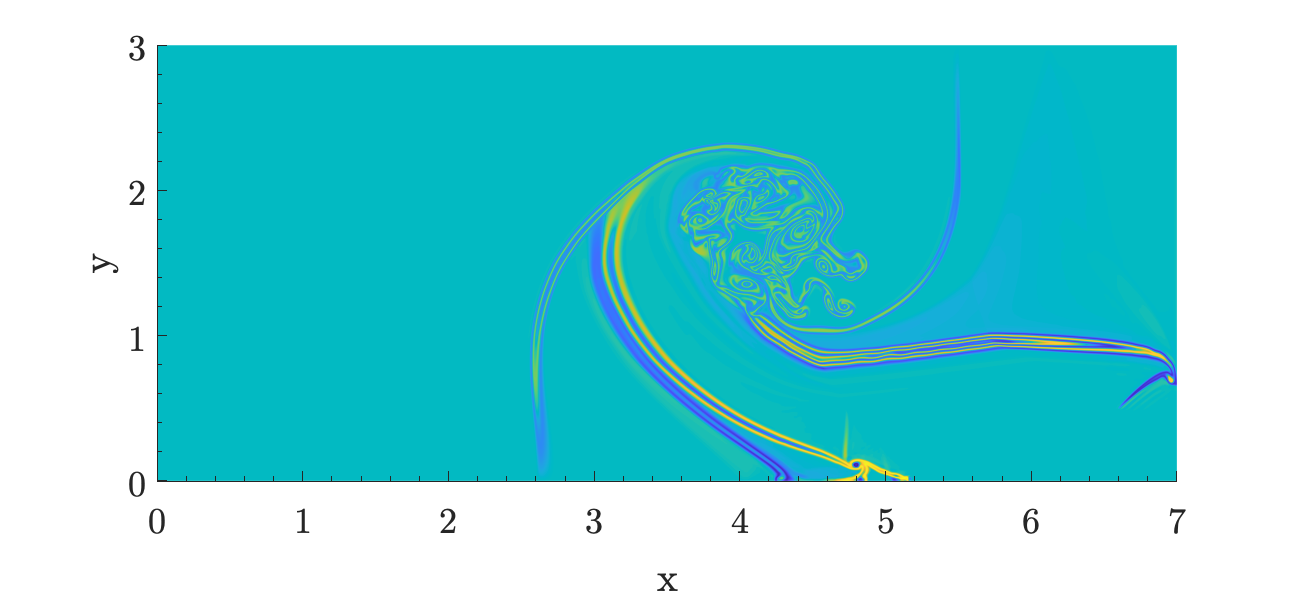}}
	{\includegraphics[width=.35\textwidth]{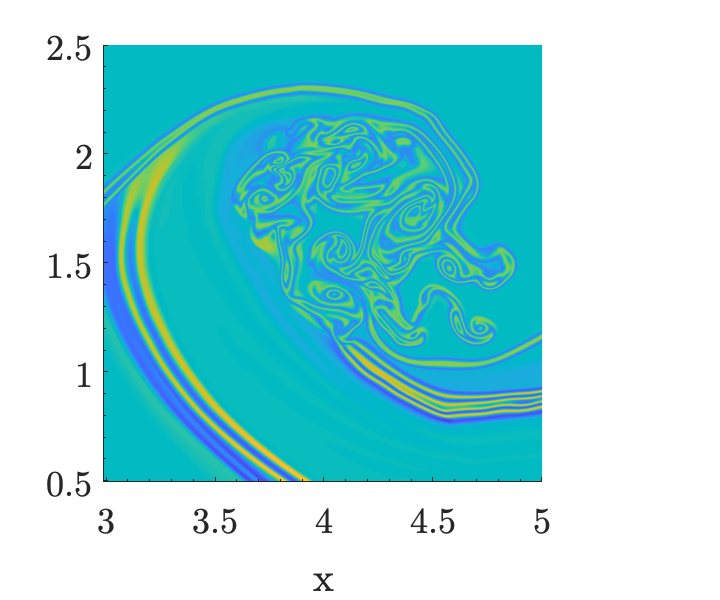}} 
	{\includegraphics[width=.64\textwidth]{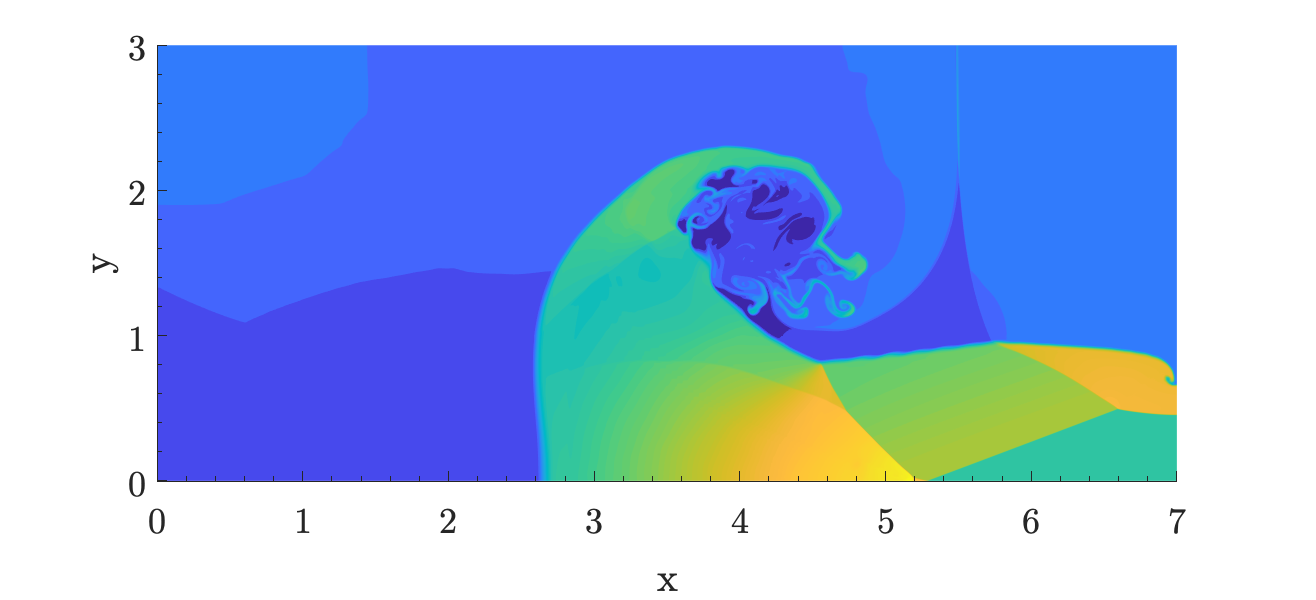}}
	{\includegraphics[width=.35\textwidth]{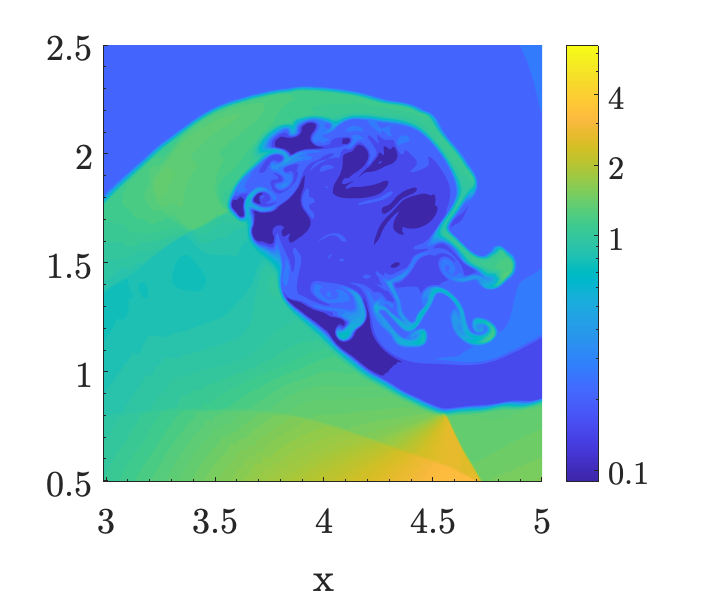}} 
	\caption{\footnotesize Results for the multiphase and multi-material triple point problem at time $t = 5$: the contour levels of the volume fractions (first phase in blue, second phase in yellow and the third one in blue-green) (top), the $A_{1,12}$ component of the distortion field for the first phase (center) and the mixture density (bottom). A zoom in the region $[3.0, 5.0] \times [0.5, 2.5]$, which illustrates the formation of vorticity resulting from the initial contact discontinuity, is shown (right).}
	\label{fig:3P5}
\end{figure}

In Fig.\,\ref{fig:3P3} and \ref{fig:3P5}, we present the results obtained for
the multiphase and multi-material triple point problem at time $t = 3$ and $t =
5$, respectively. In particular, we illustrate the evolution of the different
volume fractions (top), where the first phase is shown in blue, the second in
yellow and the third in blue-green. These contour levels clearly show that the
vortex shape is well resolved and a zoom is shown on the right to better
visualise the vorticity formation resulting from the initial contact
discontinuity. From these results, it is clear how suitable a modelling and
numerical description via an Eulerian diffuse interface approach is for
describing these mixing areas characterised by high vorticity. 

From the representation of the density field in Fig.\,\ref{fig:3P3} and
\ref{fig:3P5}, the entire dynamics of the problem can be well understood. The
fluid flow after the initial discontinuity has broken is characterised by a
rarefaction wave pointing to the left and two shock waves pointing to the right,
separated by a horizontal contact discontinuity. Moreover, these two shock waves
has different velocities, as the densities of the materials are different, and
this leads to the formation of a strong vortex. Our computational results
compare visually well with reference solutions available in the literature
\cite{ReALE2010, MaireCyl2, arepo}, proving the ability of the code and model to
correctly propagate shock waves over multiphase and multi-material regions,
despite the results being obtained on a simple fixed Cartesian grid.

Additionally, in Fig.\,\ref{fig:3P3} and \ref{fig:3P5}, we present the evolution
of the $A_{1,12}$ component of the first phase distortion field. However, as
emphasised in the previous test case, the distortion fields of all phases encode
the same flow structure, except that they satisfy three different algebraic
constraints \eqref{constrain.A}, thus the representation in Fig.\,\ref{fig:3P3}
and \ref{fig:3P5} is indicative of the flow structure of all phases. We note
once again an extremely useful ability of the distortion field (encoded in its
rotational component $\boldsymbol{R}_1$) to demonstrate the details of the flow
structures hidne otherwise. Thus, thanks to this ability of the
distortion field, it is possible to identify a strong shear zone along a contact
discontinuity inside the second phase (the yellow one) that otherwise would not
have been visible using the other state variables.
\subsection{Water entry of a symmetric wedge}
\begin{figure}[!h]
	\renewcommand{\figurename}{\footnotesize{Fig.}}
	\centering
	{\includegraphics[width=0.35\textwidth]{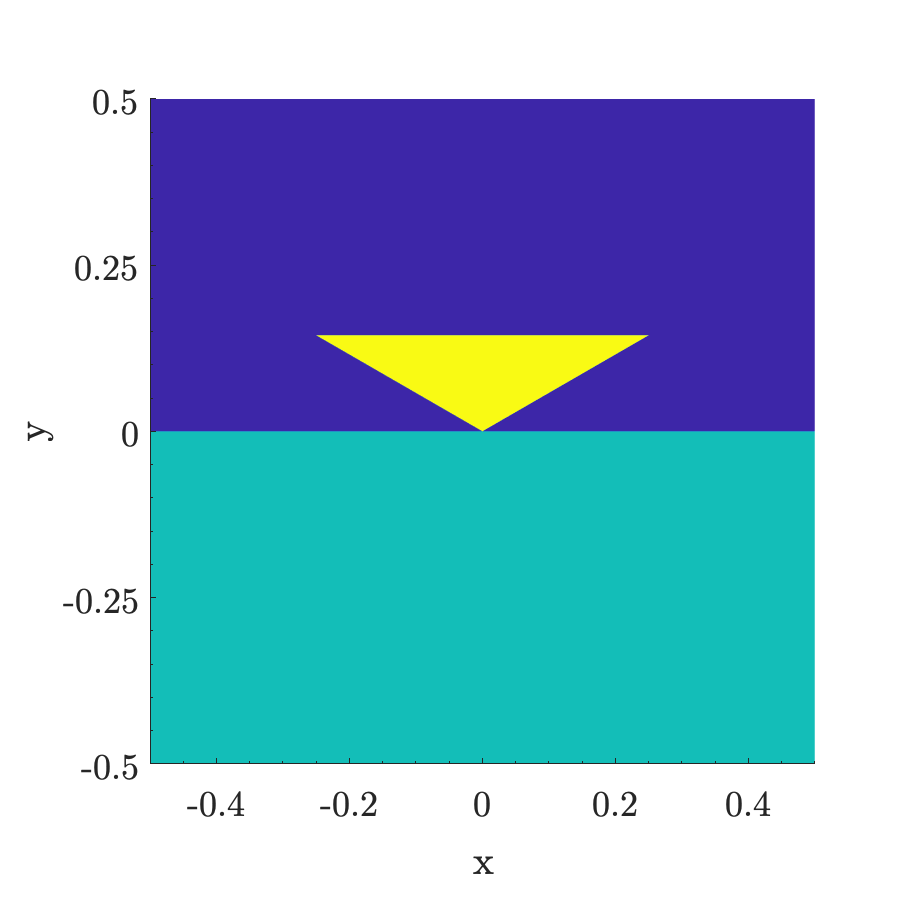}}
	{\includegraphics[width=0.35\textwidth]{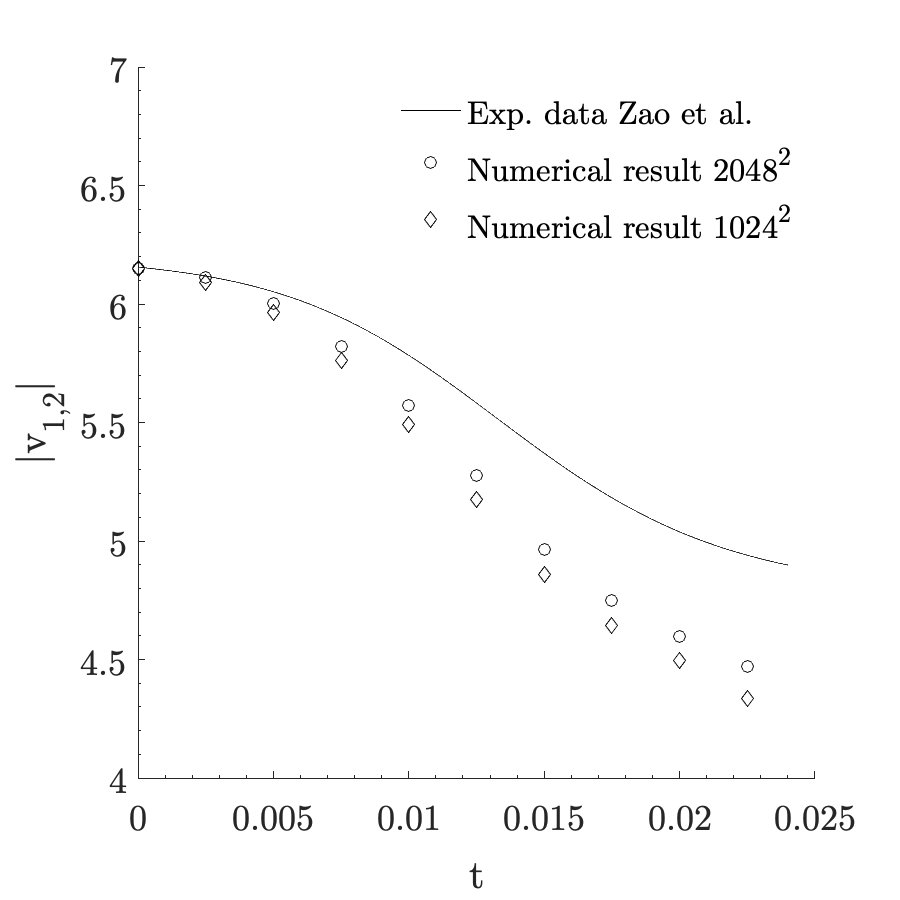}}
	\caption{\footnotesize Filled contour map of the different volume fractions of the constituents to represent the initial test condition (left). Comparison with the reference vertical wedge fall velocity experimentally recorded by Zhao et al. \cite{Wedge1} of the average components of the vertical velocity of the solid wedge $v_{1,2}$ evaluated with the explicit FV scheme for the complete compressible multiphase fluid and solid mechanics on two different uniform Cartesian meshes (right).}
	\label{fig:wedge_comp}
\end{figure}

In all the previous test cases all three phases were formally considered, but
some of them had identical material parameters. In this sense, the test
addressed in this section is more general, and all three phases will be
considered having very different material parameters. Namely, we consider one
solid, one liquid and one gaseous phase. The aim of this numerical test is to
reproduce, as fully as possible, the experiment conducted by Zhao et al.
\cite{Wedge1}, which consists of the impact of a symmetrical wedge with a free
surface. This test has already been considered in the validation of the
semi-implicit numerical scheme presented by the authors in \cite{SIVOF2}. Unlike
the previous work, using the explicit FV scheme developed in this paper for the
compressible multiphase fluid and solid mechanics, we have a possibility to
reproduce the experimental setup \cite{Wedge1} in its completeness. Namely, by
imposing only the initial velocity of the solid wedge, we can now evaluate the
complete dynamics of the three phases resulting from their interaction.
Therefore, the vertical velocity of the wedge recorded experimentally by Zhao et
al. \cite{Wedge1} and reported in Fig.\,\ref{fig:wedge_comp}, which was
\textit{prescribed} in the previous work \cite{SIVOF2} and in the numerical test
by Oger et al. \cite{Wedge2}, now becomes an excellent indicator to judge about
the validity of the multiphase simulation presented here.

In this numerical experiment, we follow the geometric setup used in our previous
work \cite{SIVOF2} but with the significant modification of defining an initial
problem involving three phases. The computational domain $\Omega=[ -0.5; 0.5]
\times [-0.5; 0.5]$ is divided into three subdomains occupied with three phases
describing an ideal elastic solid, a viscous liquid phase and a viscous gaseous
phase. The initial condition for each phases are set by means of jumps in volume
fraction. The first phase, the solid one, is defined according to the the
geometry of the experimental section, which consider the wedge with a dead-rise
angle of 30$^\circ$ illustrated in Fig.\,\ref{fig:wedge_comp} (left);
specifically the initial conditions for this phase read
\begin{align}
&\alpha_1(x,y)  = 	\begin{cases}	1- 2 \epsilon  &\mathrm{if}  \ \   |x| \leq {0.25} \  \land \ y \geq |x \tan(\frac{\pi}{6})| \ \land  \ y \leq 0.25 \tan(\frac{\pi}{6})  ,  \\ \nonumber
						\epsilon\ \ & \mathrm{if}  \ \  \mathrm{otherwise},
			\end{cases}     \qquad   \\
			&\rho_1 = 7{\times} 10^3,  \quad \vv_1 = (0, -6.15), \quad p_1 = 10^{3}, \quad \A_1 = \Id,
\end{align}
where the density is evaluated to obtain the total weight of the measuring section used in the experiment by Zhao et al. \cite{Wedge1}, which corresponds to $255.5 \ kg$; thus the density is calculated as the weight of the instrumental tools divided by the effective area described by the wedge in this numerical setup.
The initial condition for the second phase, defining the viscous liquid phase, are
\begin{align}
\begin{split}
&\alpha_2(x,y)  = 	\begin{cases}	1- 2 \epsilon  &\mathrm{if}  \ \  \xx \in \Omega_2 , \\
						\epsilon\ \ & \mathrm{if}  \ \  \xx \notin \Omega_2,
			\end{cases}     \quad  \mathrm{with} \quad \Omega_2 =  [-0.5; 0.5]{\times}[-0.5;0.0] \\
			&\rho_2 = 10^3,  \quad \vv_2 = \vec{0}, \quad p_2 = 10^{3}, \quad \A_2 = \Id, 
\end{split}
\end{align}
while the third one, defining the gaseous phase, is initialized as
\begin{align}
\alpha_3  = 1 -\alpha_1 - \alpha_2, \quad
			\rho_3 = 1,  \quad \vv_3 = \vec{0}, \quad p_3 = 10^{3}, \quad \A_3 = \Id.
\end{align}
For the solid and liquid phase, the stiffened gas EOS is used; the other material parameters are $\gamma_1 = \gamma_2 = 1.4$, $\Cv_1 = \Cv_2 = 1$, $\Cs_1 = 120$, $\Cs_2 = 100$, $\Co_1 = \Co_2 = 120$ and $\po_1 = \po_2 = p_1$. For the viscous gas, the ideal gas EOS is used with the following parameters $\gamma_3 = 1.4$, $\Cv_3 = 1$, $\Cs_3 = 60.0$ and $\nu_3 = 10^{-1}$. To consider an ideal elastic material, the relaxation time for the first phase is chosen to be $\tau^\mathrm{e}_1 = 10^{14}$, while a kinematic viscosity $\nu_1= 10^{-6}$ is adopted for the viscous liquid phase. 

Two simulations are carried out up to the final time $t_f =0.025$ on two different uniform Cartesian meshes composed of $1024\times1024$ and $2048\times2048$ control volumes, in order to verify qualitatively the mesh convergence of the solution. In these simulations, reflective slip wall boundary conditions are set in all the directions. 

Fig.\,\ref{fig:wedge_comp} (left) shows the distribution of the different volume
fractions in accordance with the initial conditions describing the geometric and
experimental setup of the test water entry of a symmetric wedge. The solid phase
is shown in yellow, the liquid phase in blue-green and the gas phase in blue.
Moreover, in Fig.\,\ref{fig:wedge_comp} (right) we present a comparison that
verifies the validity of the results obtained. In this comparison, the vertical
velocity of the wedge experimentally recorded by Zhao et al. \cite{Wedge1} is
compared with the purely elastic solid body velocity computed in this test by
evaluating an averaged vertical velocity using the volume fraction, in
accordance with the following definition 
\begin{align}
		|\vv_{1,2}| = \frac{ | 	\sum^{N_1 N_2}_{ij} \alpha_{1,ij} v_{1,2}	|}{\sum^{N_1 N_2}_{ij} \alpha_{ij}},
\end{align}
where $N_1$ and $N_2$ are the discrete elements in the first and in the second directions, respectively. It is possible to observe how qualitatively the dynamics of the impact is well represented, in particular the deceleration over time follows the correct trend, i.e. deceleration increases in modulus until about half the simulation time and then tends to decrease. 
Furthermore, both the results obtained with a $1024 \times1024$ mesh and that
obtained by doubling the mesh resolution are represented, and it can be
seen that the numerical solution is getting closer to that recorded experimentally by Zhao
et al. \cite{Wedge1} as mesh is getting more finer. The main reason for the discrepancy from the experimental
result has to be found in the strongly low Mach nature of the test. Indeed, this
impact, in which the solid must maintain a particularly rigid behaviour,
represents a complex test for an explicit numerical scheme. Moover, it
should be noted that this is the first time this test has been solved by considering the interaction of three phases through a monolithic mathematical
model for compressible multiphase fluid and solid mechanics.

\begin{figure}[!htbp]
	\renewcommand{\figurename}{\footnotesize{Fig.}}
	\centering
	{\includegraphics[width=.35\textwidth]{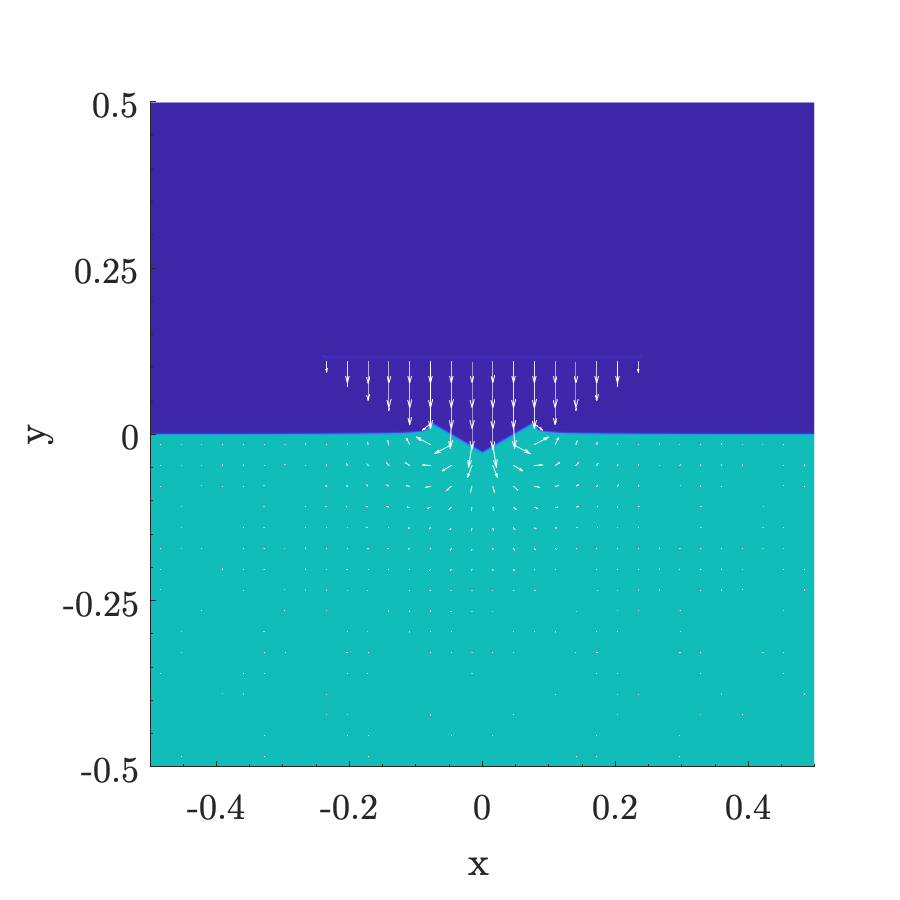}}
	{\includegraphics[width=.35\textwidth]{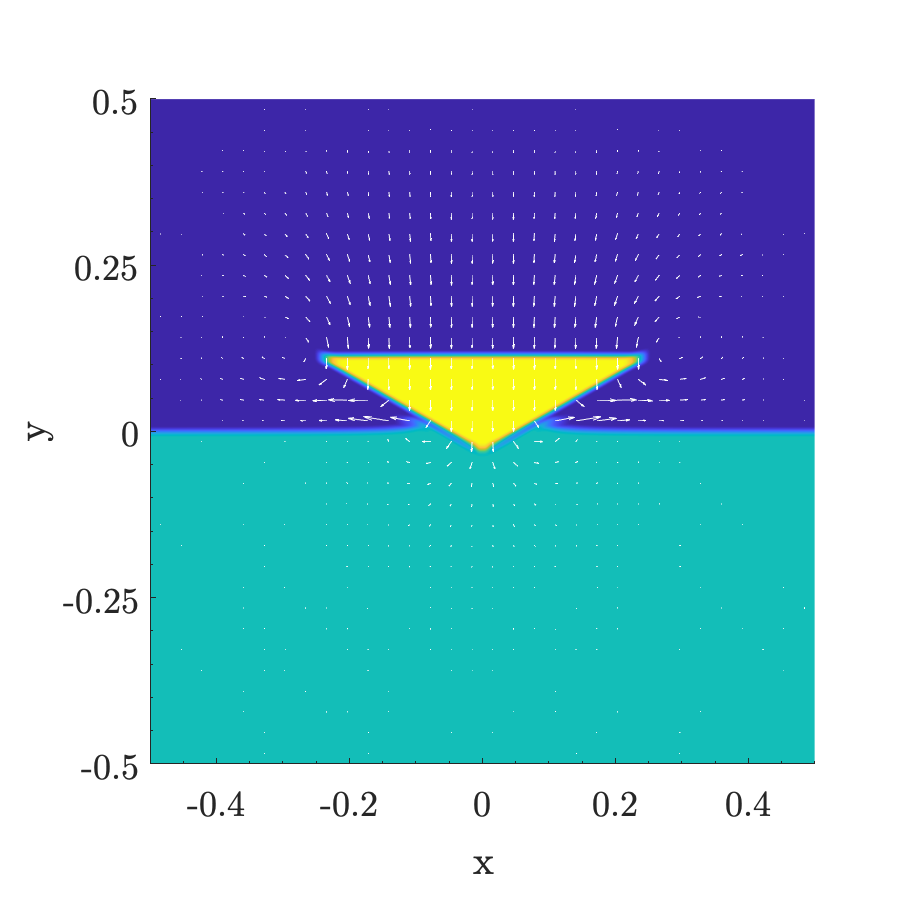}} 	
	{\includegraphics[width=.35\textwidth]{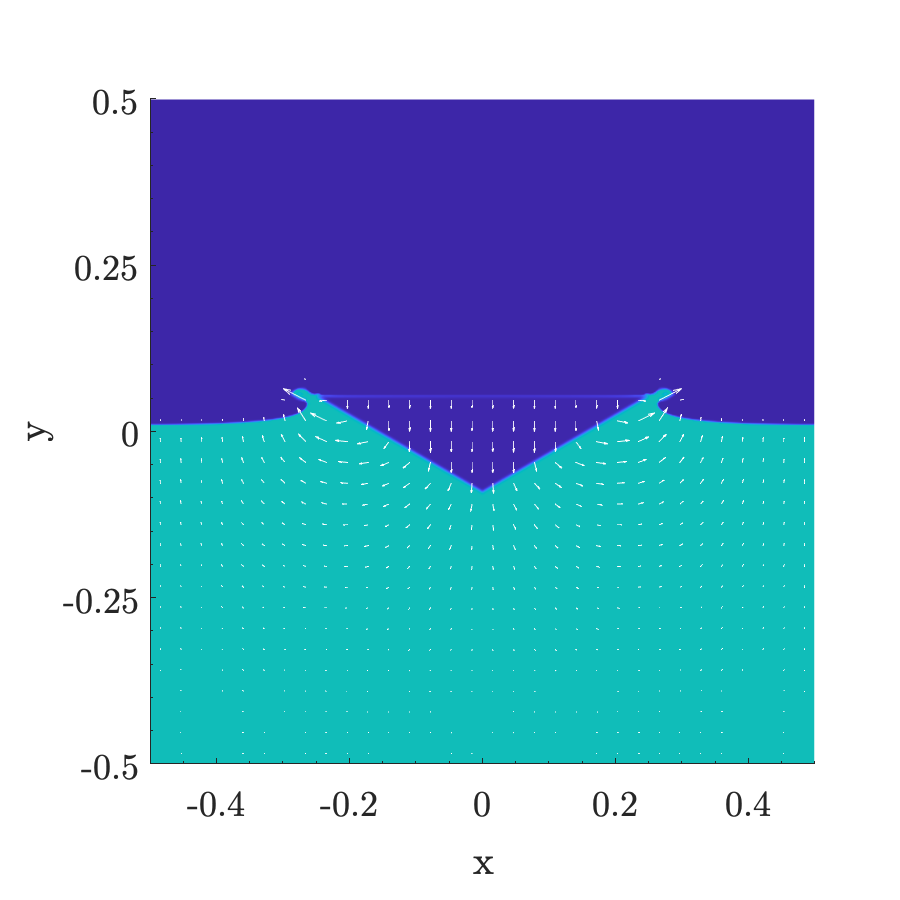}}
	{\includegraphics[width=.35\textwidth]{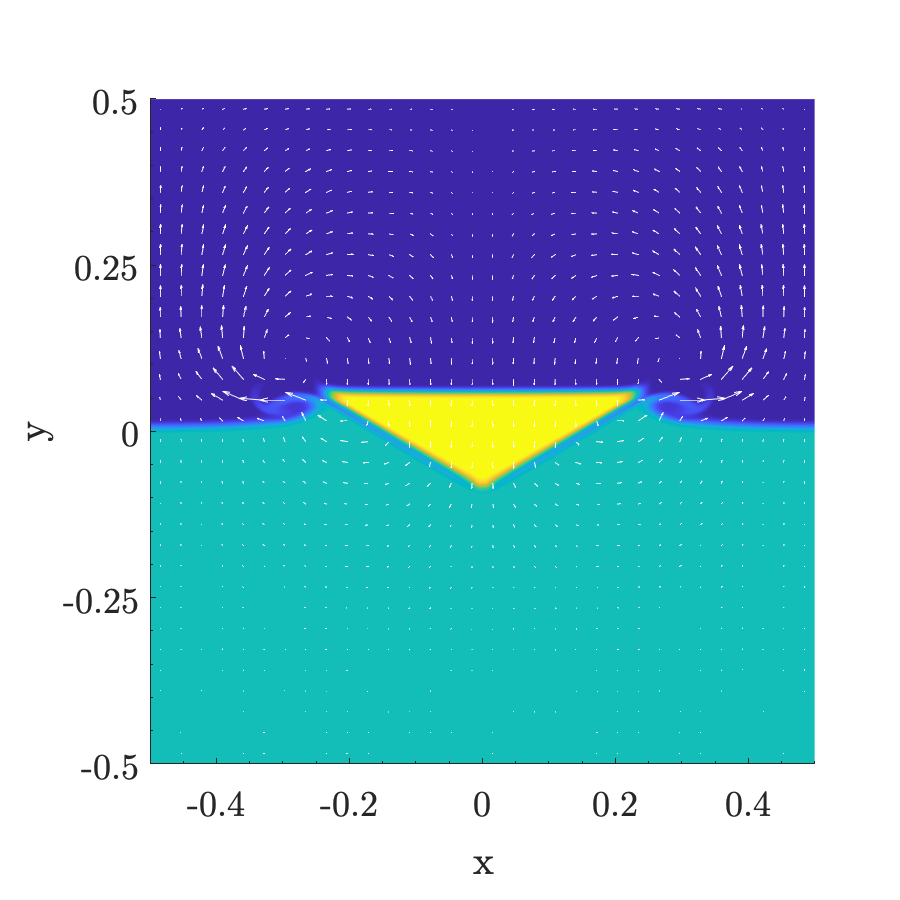}} 
	{\includegraphics[width=.35\textwidth]{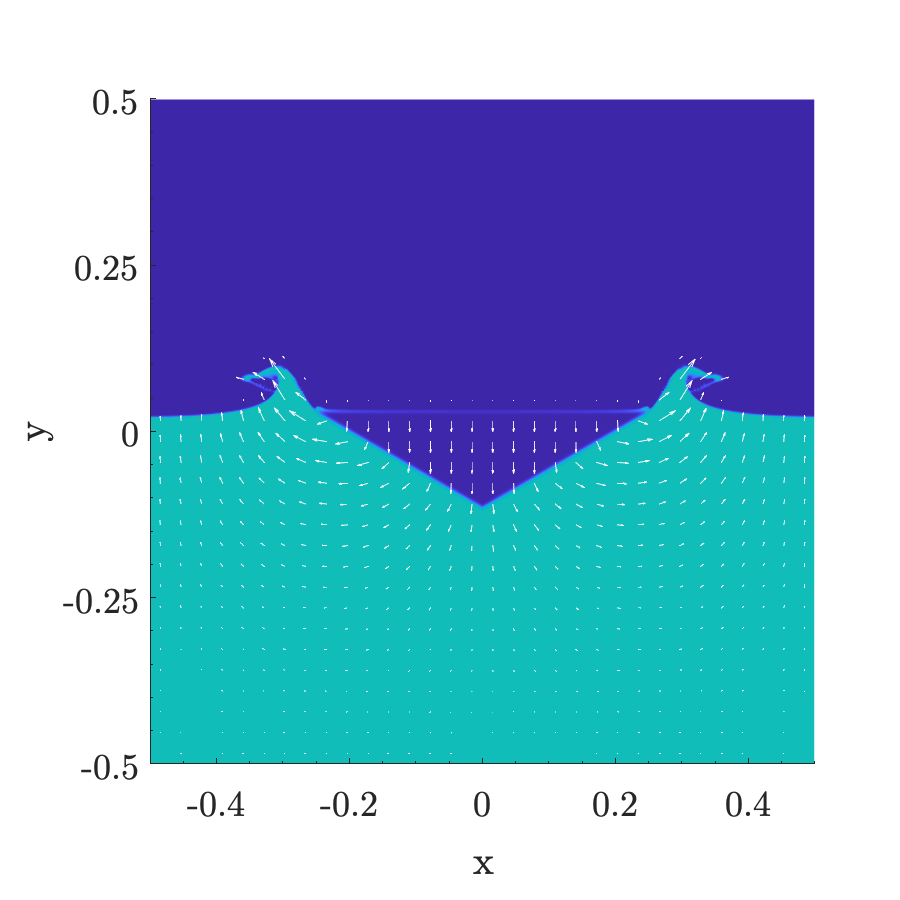}}
	{\includegraphics[width=.35\textwidth]{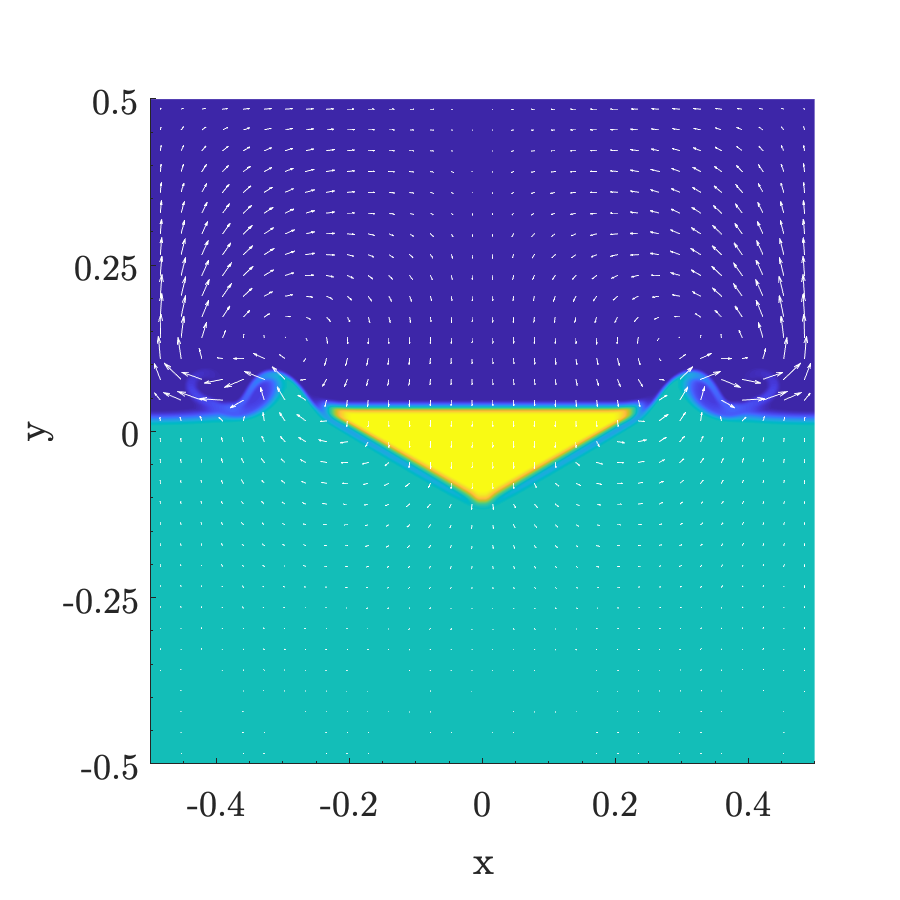}} 
	\caption{\footnotesize Filled contour map of the different volume fractions of the constituents and velocity field of the mixture, the solid phase is shown in yellow, the liquid phase in blue-green and the gas phase in blue. Results obtained with the explicit FV scheme for compressible three-phase fluid and solid mechanics (right), and the two-phase results obtained with a semi-implicit numerical method from \cite{SIVOF2} (left). 
	From the top to the bottom, three different time instants are shown: $t = 0.005$, $t = 0.015$ and $t = 0.020$.}
	\label{fig:wedge_evol}
\end{figure}

Fig.\,\ref{fig:wedge_evol} shows the temporal evolution of volume fractions
obtained with the explicit FV scheme for compressible multiphase fluid and solid
mechanics (right column). For the sake of comparison, the results previously
obtained in \cite{SIVOF2} with a semi-implicit numerical method are shown along
side with the new results, see the left column. From the top to the bottom, the
results for three different instants are shown: $t = 0.005$, $t = 0.015$ and $t
= 0.020$. 
It can be observed that the phenomenological evolution of the free surface
during the entry of the wedge, obtained with the new numerical method that
solves the entire three-phase dynamics, is in agreement with what was previously
obtained in \cite{SIVOF2} with a semi-implicit solver. The first time instant
shows quite well the formation of two jets escaping along the edges of the
wedge. At time $t = 0.015$, these two jets reach the point of separation, which
corresponds to the end of the edge. It can be seen that the jets leave the edge
almost tangentially at this initial phase of flow separation. Then, at $t =
0.020$, the jets tend to move more vertically as well as breaking. This shape
is qualitatively similar to the experimental illustration in the article by Zhao
et al. \cite{Wedge1}.

Furthermore, Fig.\,\ref{fig:wedge_evol} shows the velocity fields obtained with
both numerical methods. One can see that the gas phase was not considered in
\cite{SIVOF2} presented in the left column, while the dynamics of all three
phases is taken into account in this paper. The interraction of the liquid jets
and the gas phase might in particular be responsible for the slight differences
in jets shape between the two simulations.
\subsection{ Multiphase and multi-material solid impact}
This is the last test presented in this work, the aim of which is to show from a
qualitative point of view the wide applicability of the model and numerical
scheme presented. As in the previous test, all the capabilities of the developed
model and numerical scheme are tested. Namely, an impact of true multi-material
\textit{three-phases} will be simulated, which exhibits elastic and also
elasto-plastic behaviour in a inviscid fluid environment. 

In this test, we roughly follow the approach presented in \cite{DEBRAUER2016121,
Iollo2017}, which has been modified by not paying particular attention to
physical material characteristics, as the interest is to qualitatively verify
the approach ability to solve such a complex test. Thus, while maintaining the
geometry similar to the tests in the literature, the parameters and physical
quantities that define the properties of the aluminium bar in \ref{Taylor.bar}
are used. Therefore, this test case presents an impact of an aluminium ball into
an aluminium plate embedded in the surrounding inviscid fluid modeld with the
perfect gas EOS. The computational domain $\Omega = [-0.5; 0.5] \times[-0.5;
0.5]$ is divided into three subdomains filled with three phases describing two
solids with the same properties and a perfect gas. As for the previous test, the
initial conditions are set by means of jumps in volume fraction; the first
phase, that defines the aluminium ball, is initialized as
\begin{align}
\begin{split}
&\alpha_1(x,y)  = 	\begin{cases}	1- 2 \epsilon  &\mathrm{if}  \ \  ((-0.125-x)^2 + (0.5-y)^2)^{0.5} \leq0.075  , \\
						\epsilon\ \ & \mathrm{if}  \ \  \mathrm{otherwise},
			\end{cases}     \qquad   \\
			&\rho_1 = 2.785,  \quad \vv_1 = (0.002, 0), \quad p_1 = 10^{-6}, \quad \A_1 = \Id, 
\end{split}
\end{align}
the second phase, defining the aluminium plate, as
\begin{align}
\begin{split}
&\alpha_2(x,y)  = 	\begin{cases}	1- 2 \epsilon  &\mathrm{if}  \ \  \xx \in \Omega_2 , \\
						\epsilon\ \ & \mathrm{if}  \ \  \xx \notin \Omega_2,
			\end{cases}     \quad  \mathrm{with} \quad \Omega_2 =  [0.05; 0.225]\times[-0.35;0.35] \\
			&\rho_2 = 2.785,  \quad \vv_2 = \vec{0}, \quad p_2 = 10^{-6}, \quad \A_2 = \Id, 
\end{split}
\end{align}
and the third, defining the surrounding perfect gas, is initialised as
\begin{align}
\alpha_3  = 1 -\alpha_1 - \alpha_2, \quad
			\rho_3 = 10^{-3},  \quad \vv_3 = \vec{0}, \quad p_3 = 10^{-6}, \quad \A_3 = \Id.
\end{align}
Since the first and second phases represent the same material, i.e. aluminium,
the other parameters and physical quantities that define the properties of such
a material using the stiffened gas EOS are $\gamma_1 = \gamma_2 = 1.4$, $\Cv_1 =
\Cv_2 = 1000$, $\Cs_1 = \Cs_2 = 0.305$, $\Co_1 = \Co_2 = 0.533$ and $\po_1 =
\po_2 = p_1$. For the inviscid gas phase surrounding the solid phases, the ideal
gas EOS is used and the physical parameters are $\gamma_3 = 1.2$, $\Cv_3 =
1000$, $\Cs_3 = 0.0$ and $\tau^\mathrm{e}_3 = 10^{-14}$. In a first simulation,
the solid materials are assumed to have purely elastic behaviour, so the
relaxation time is assumed to be $\tau^\mathrm{e}_1 = \tau^\mathrm{e}_2 =
10^{14}$ for both phases. Subsequently, to obtain a non-linear elasto-plastic
behaviour of the material, the relaxation time $\tau^\mathrm{e}_1,
\tau^\mathrm{e}_2$ is chosen as a non-linear function of an invariant of the
shear stress tensor as done in the previous test, see \eqref{rel.time.plast}. In
this case, however, the yield stress of the material is set to a lower number,
i.e. ${\sigma}{_\mathrm{\scriptsize{o}}} = 2.5{\times}10^{-4}$, for the sake of
making the plastic deformations more visible. 

Two simulations are carried out up to the final time $t_f =200$ discretizing the computational domain with a uniform Cartesian mesh composed of $2048\times2048$ control volumes; periodic boundary conditions are set in all the directions.
\begin{figure}[!ht]
	\renewcommand{\figurename}{\footnotesize{Fig.}}
	\centering
	{\includegraphics[width=.325\textwidth]{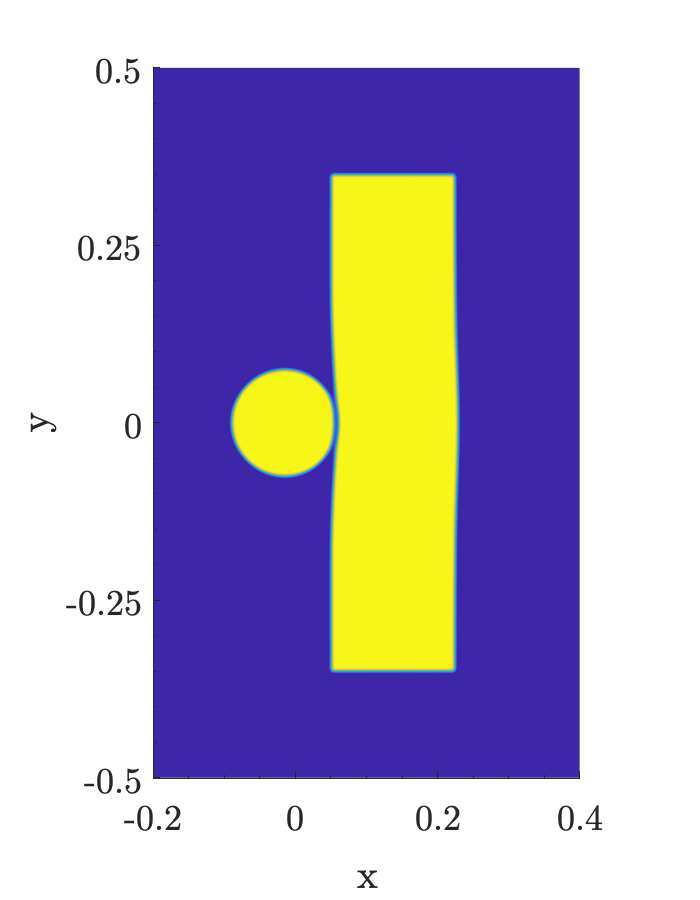}}
	{\includegraphics[width=.325\textwidth]{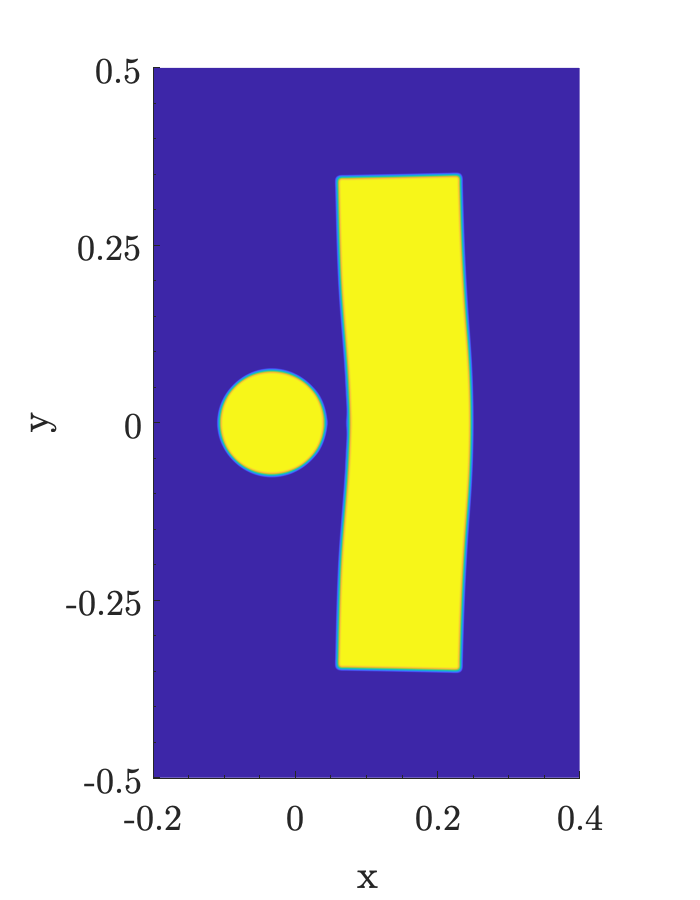}} 	
	{\includegraphics[width=.325\textwidth]{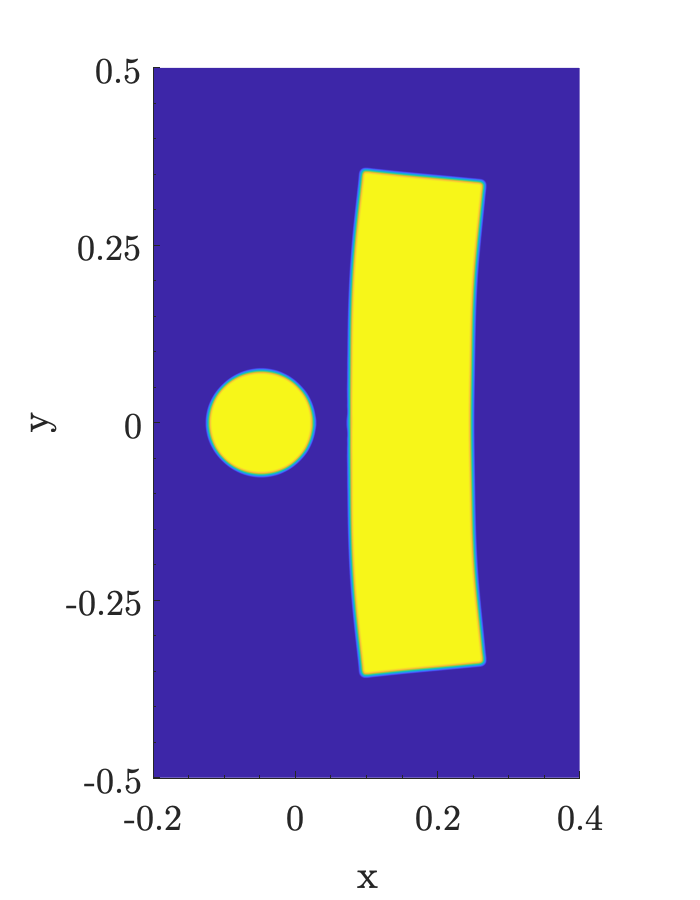}}
	\caption{\footnotesize Results for the multiphase and multi-material elastic solid impact. The contour levels of the volume fractions of the first and the second phases at times $t = 60$, $t=100$, and $t = 140$ (from left to right).}
	\label{fig:Solid_impact1}
\end{figure}
\begin{figure}[!ht]
	\renewcommand{\figurename}{\footnotesize{Fig.}}
	\centering
	{\includegraphics[width=.325\textwidth]{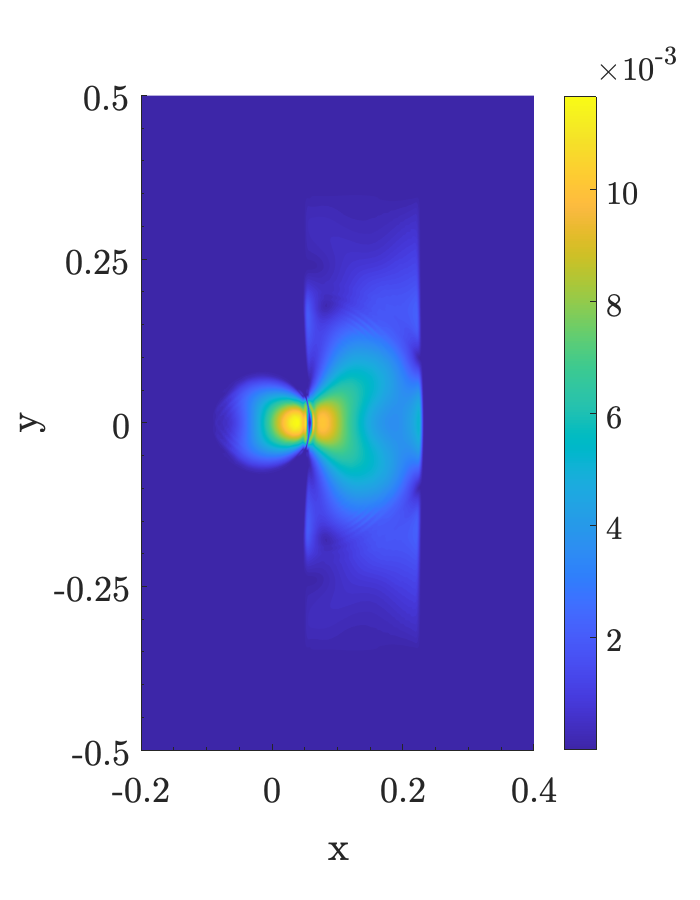}}
	{\includegraphics[width=.325\textwidth]{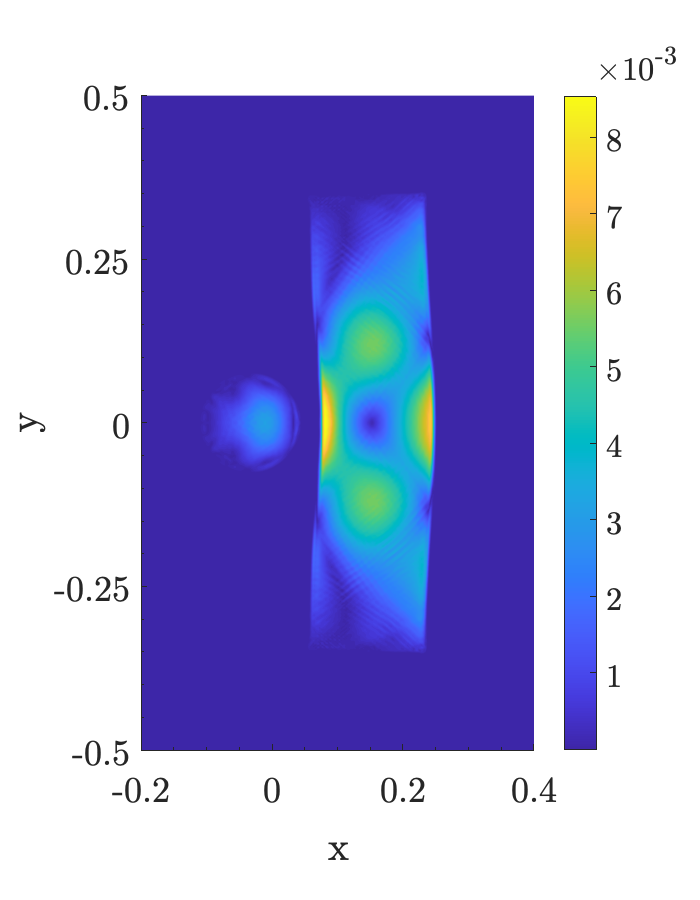}} 	
	{\includegraphics[width=.325\textwidth]{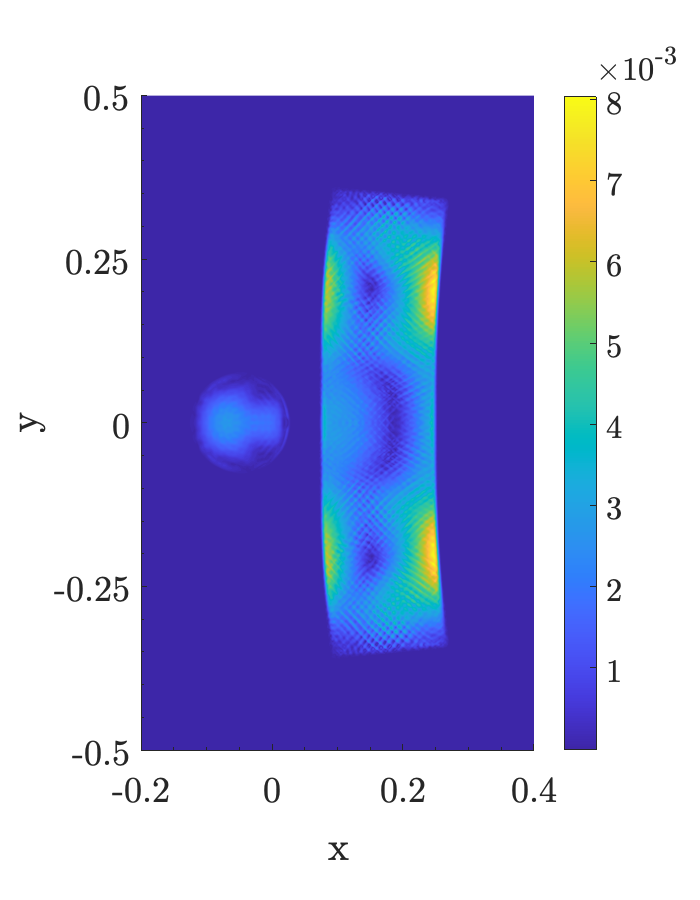}}
	\caption{\footnotesize Results for the multiphase and multi-material elastic solid impact. The contour levels of the von Mises stress of the first ${\bar{\sigma}_1}$ and the second ${\bar{\sigma}_2}$ at times $t = 60$, $t=100$ and $t = 140$ (from left to right).}
	\label{fig:Solid_impact3}
\end{figure}
\begin{figure}[!htbp]
	\renewcommand{\figurename}{\footnotesize{Fig.}}
	\centering
	{\includegraphics[width=.32\textwidth]{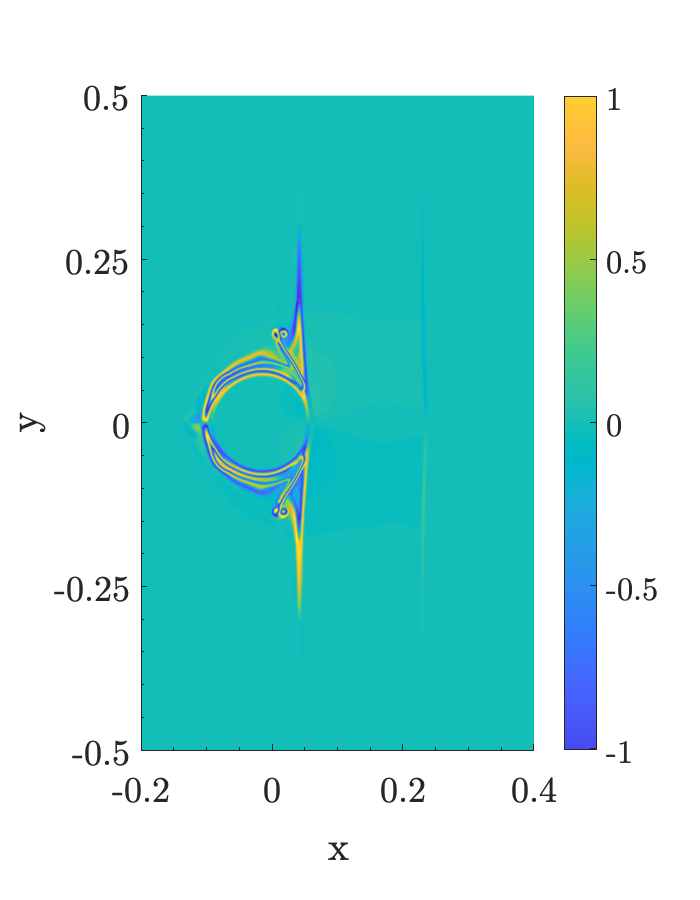}}
	{\includegraphics[width=.32\textwidth]{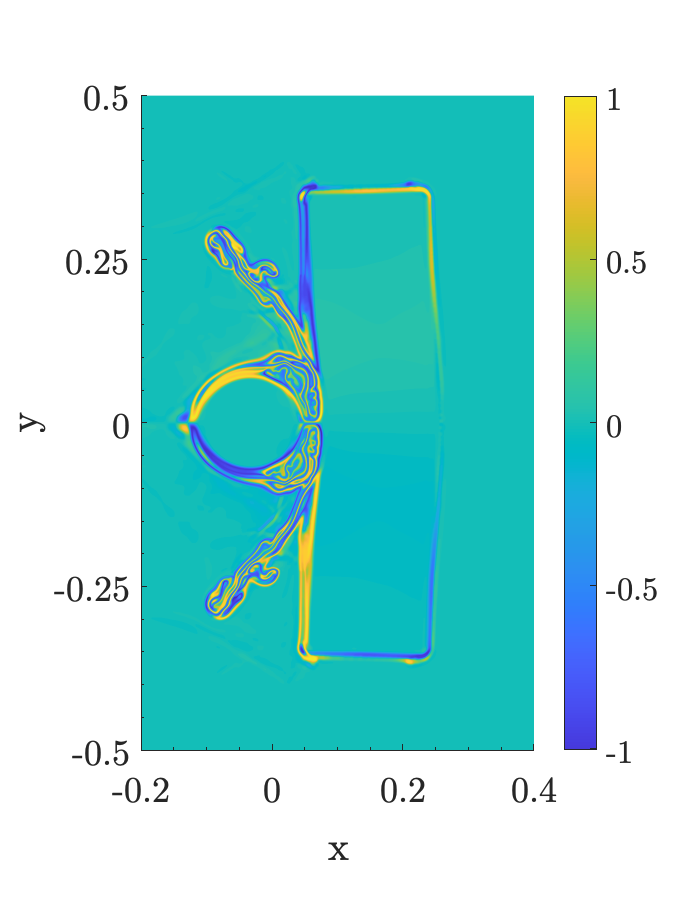}} 	
	{\includegraphics[width=.32\textwidth]{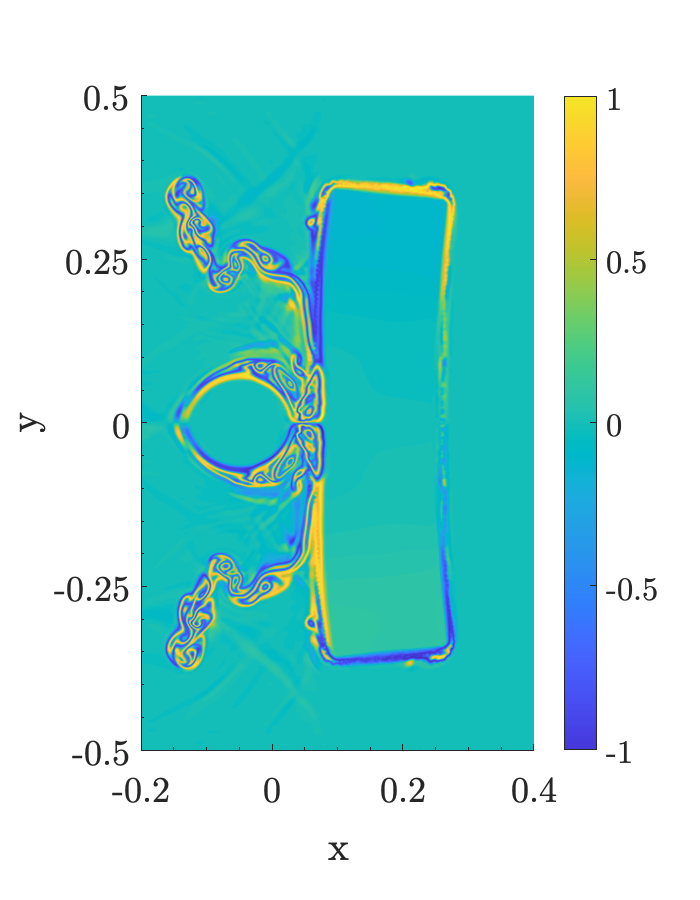}}
	\caption{\footnotesize Results for the multiphase and multi-material elastic solid impact. The contour levels of the  $A_{12}$ component of the mixture  distortion field $\A$ are presented at times $t = 60$, $t=100$, and $t = 140$ (from left to right).}
	\label{fig:Solid_impact5}
\end{figure}

\begin{figure}[!htbp]
	\renewcommand{\figurename}{\footnotesize{Fig.}}
	\centering
	{\includegraphics[width=.325\textwidth]{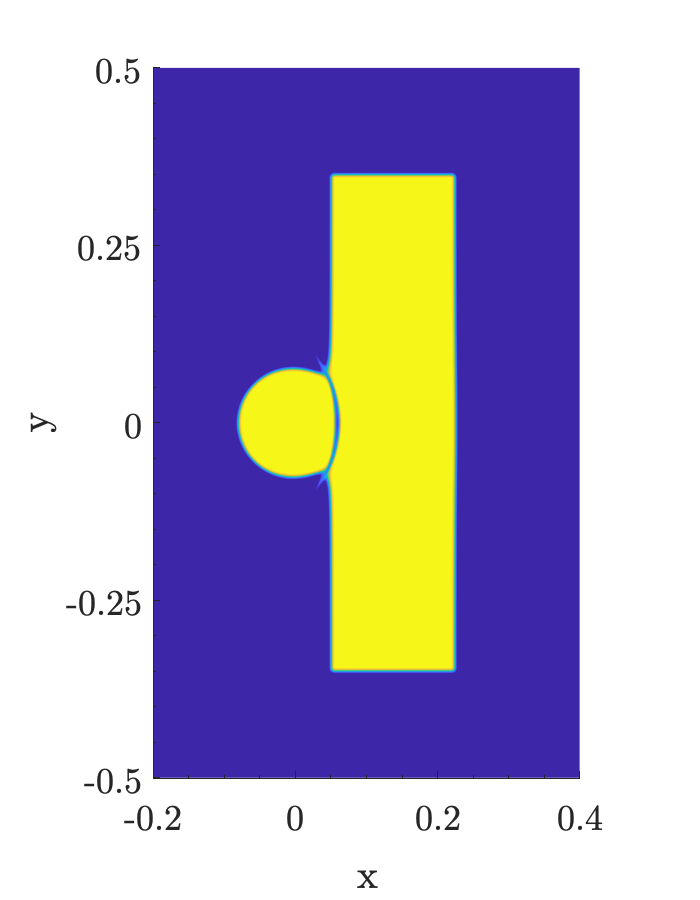}}
	{\includegraphics[width=.325\textwidth]{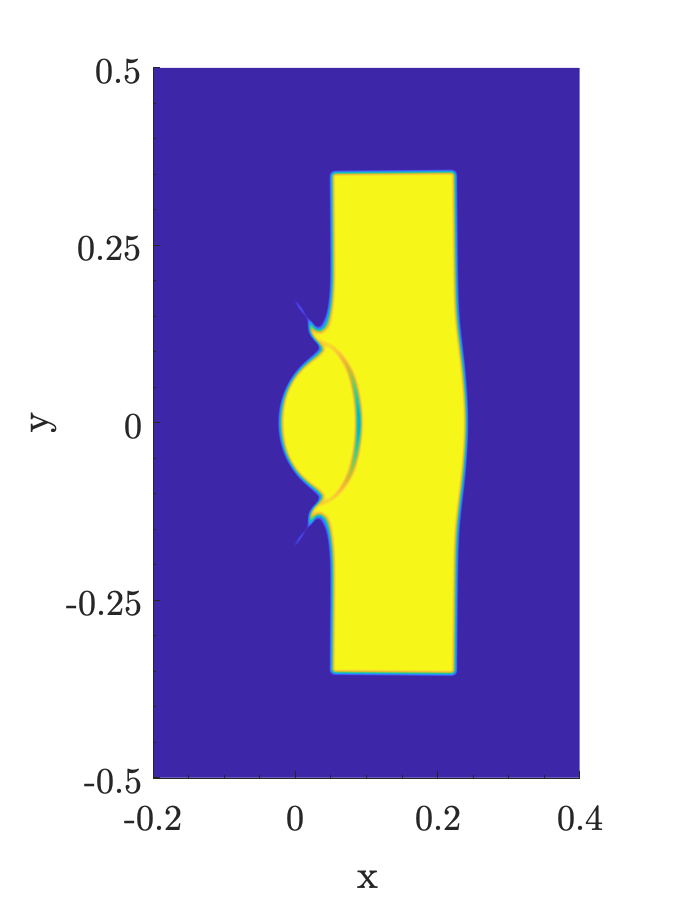}} 	
	{\includegraphics[width=.325\textwidth]{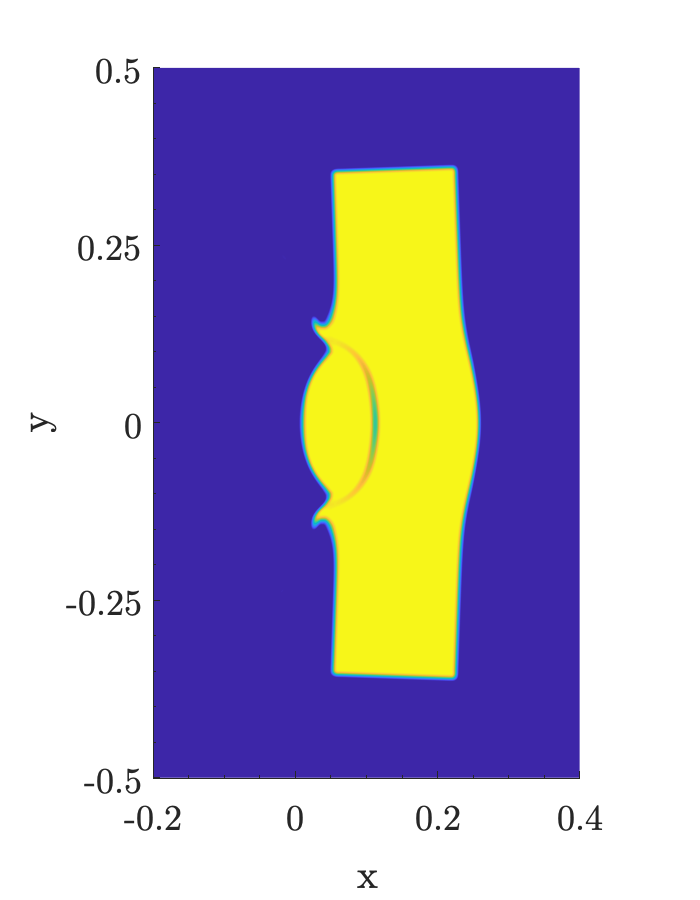}}
	\caption{\footnotesize Results for the multiphase and multi-material elasto-plastic solid impact. The contour levels of the volume fractions of the first and the second phases at times $t = 60$, $t=100$, and $t = 140$ (from left to right).}
	\label{fig:Solid_impact2}
\end{figure}
\begin{figure}[!htbp]
	\renewcommand{\figurename}{\footnotesize{Fig.}}
	\centering
	{\includegraphics[width=.325\textwidth]{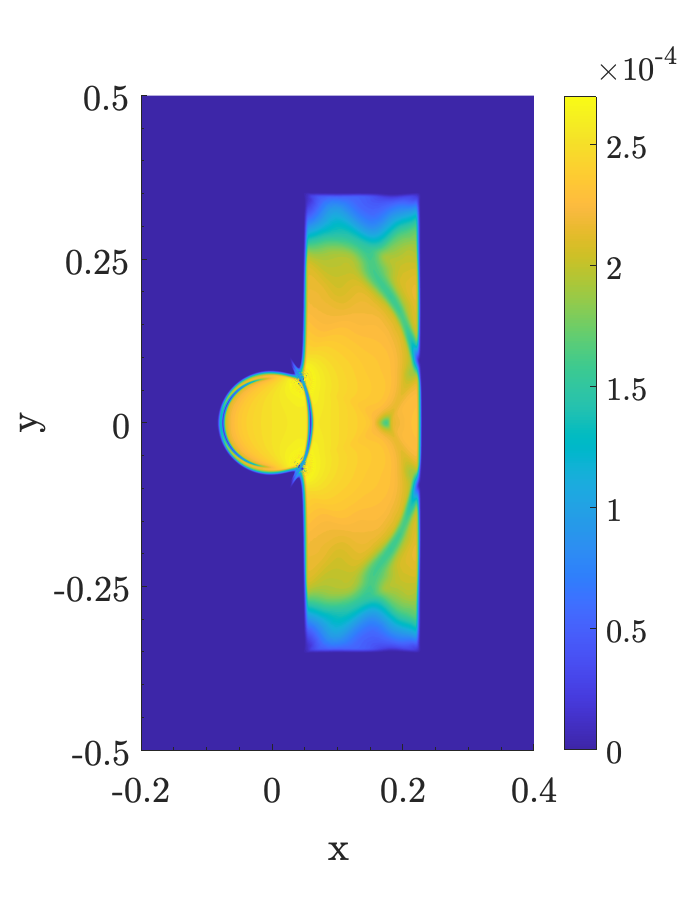}}
	{\includegraphics[width=.325\textwidth]{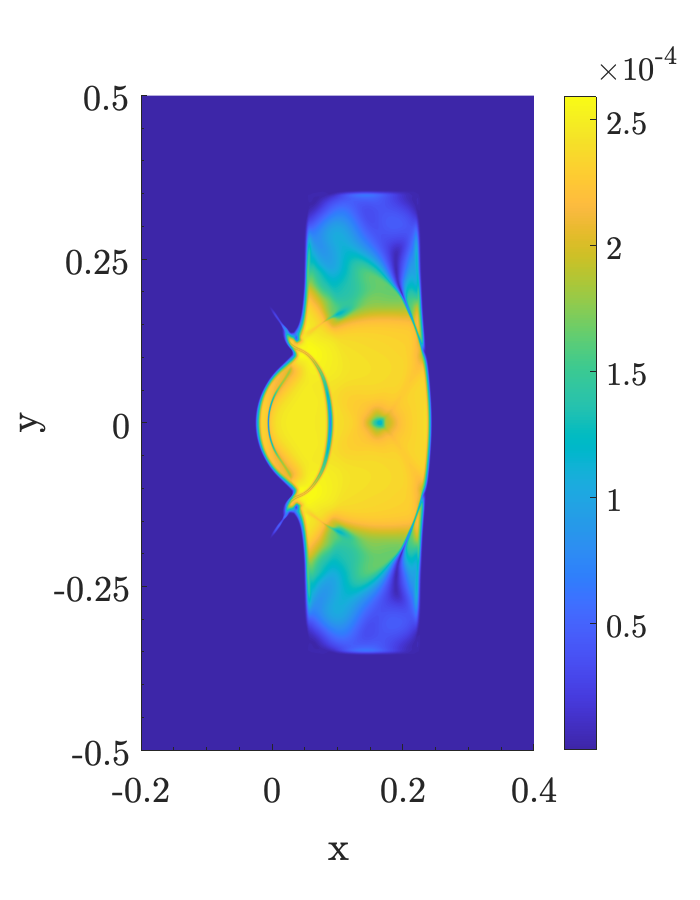}} 	
	{\includegraphics[width=.325\textwidth]{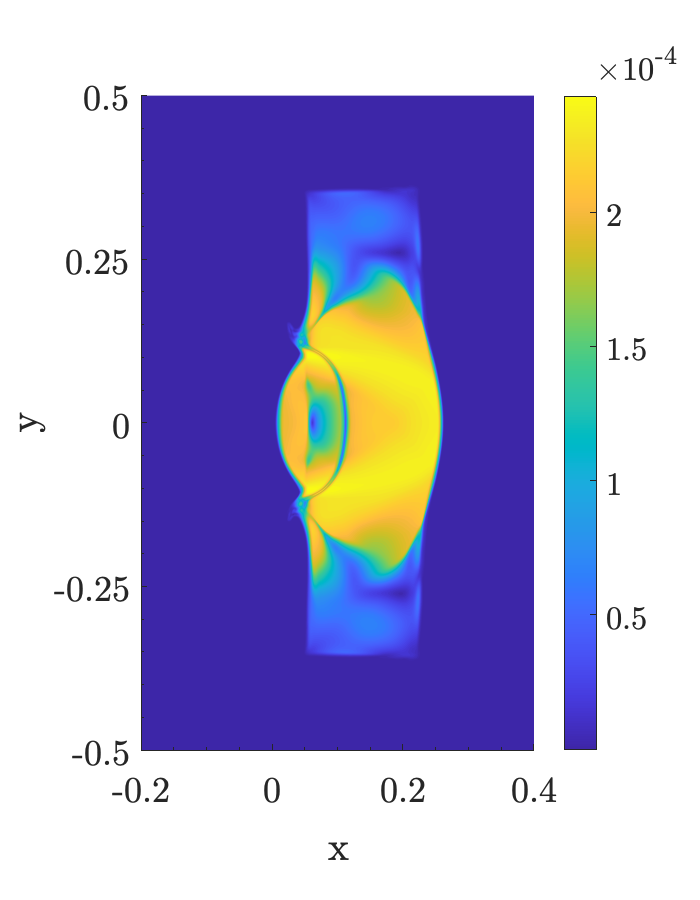}}
	\caption{\footnotesize Results for the multiphase and multi-material elasto-plastic solid impact. The contour levels of the von Mises stress of the first ${\bar{\sigma}_1}$ and the second ${\bar{\sigma}_2}$ at times $t = 60$, $t=100$, and $t = 140$ (from left to right).}
	\label{fig:Solid_impact4}
\end{figure}

Fig.\,\ref{fig:Solid_impact1} shows the volume fraction of the first and second
phase at times $t=60$, $t= 100$ and $t = 140$ from left to right, respectively.
The first time instant represents the moment of impact with the plate. As one
can see from the subsequent instants the behavior of an elastic collision is
qualitatively well represented by the numerical solution. It should be
emphasized that, in a diffuse interface approach, and if both the solid objects
are represented by the same volume fraction function, it is not obvious that the
two solids would bounce instead of sticking to each other. The results for
multibody problems in which the solids are carrying their own volume fractions
is, therefore, of considerable interest. Additionally,
Fig.\,\ref{fig:Solid_impact5} shows the $A_{12}$ component of the mixture
distortion field, obtained through the following relation
\begin{equation}
A_{12} = \alpha_1 A_{1,12}+ \alpha_2 A_{2,12}+ \alpha_3 A_{3,12}.
\end{equation}
This allows the dynamics of the gas phase to be clearly shown as well, making it
evident that the dynamics of all three phases have been resolved through a
distortion field for each phase. It is possible to see the two fluid jets with
non-trivial vorticity being generated at the moment of impact. 

Similarly,
Fig.\,\ref{fig:Solid_impact2} shows a collision of the same solid objects but
withe plasticity effect taken into account. The obtained results visually
compare well with the one in \cite{Iollo2017}. 
To better understand how the different definitions of material properties in
these two tests affect the behavior of solids, it is useful to observe the von
Mises stress of the first phase ${\bar{\sigma}_1}$ and the second phase
${\bar{\sigma}_2}$, evaluated as in \eqref{eq.Mises}, and presented in
Fig.\,\ref{fig:Solid_impact3} and \ref{fig:Solid_impact4}. It can be seen that
the stress in an ideal elastic material propagates through the body by means of
waves, which are reflected over time from the body boundaries. On the contrary,
in the case of an elasto-plastic material, it can be observed that the stress
reaches a lower magnitude than in the ideal elastic case, due to the stress
relaxation process in the inelastic deformations. Furthermore, it is evident
that over time, the highest stress values are localised in the area undergoing
plastic deformations, while the regions far from the impact, in this case, are
less stressed. 

\begin{figure}[!htbp]
	\renewcommand{\figurename}{\footnotesize{Fig.}}
	\centering
	{\includegraphics[width=.32\textwidth]{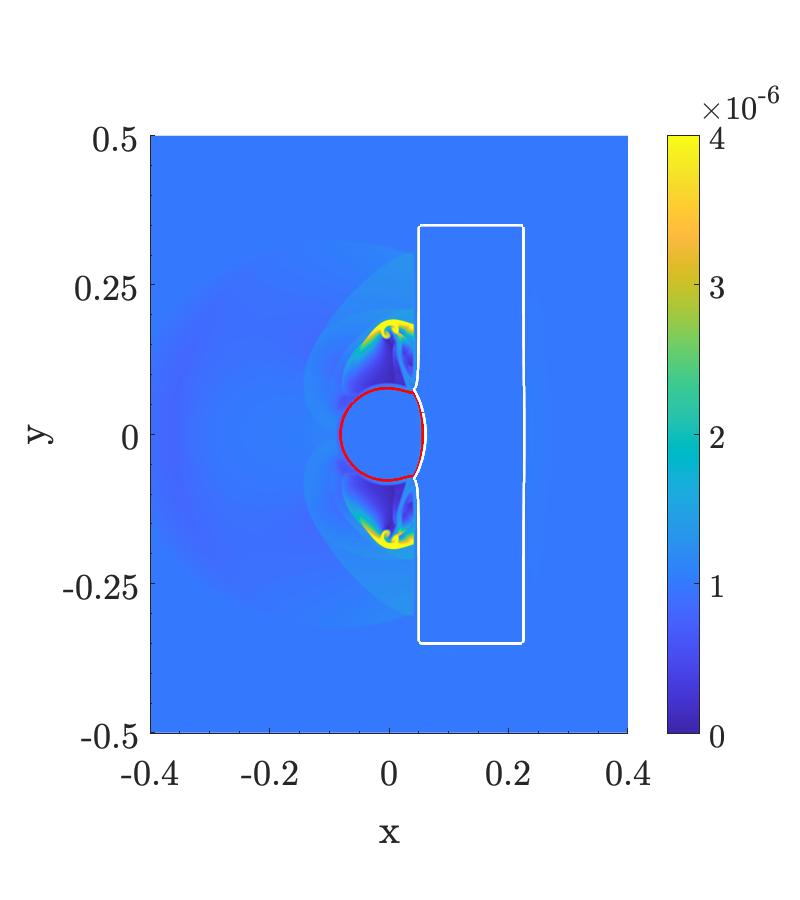}}
	{\includegraphics[width=.32\textwidth]{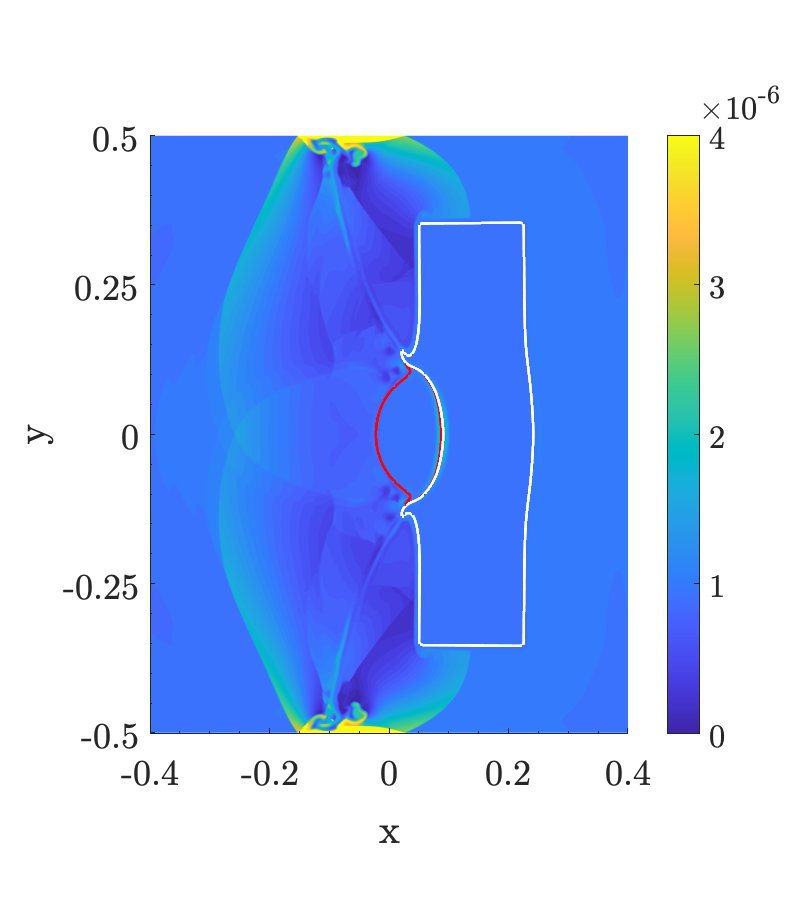}} 	
	{\includegraphics[width=.32\textwidth]{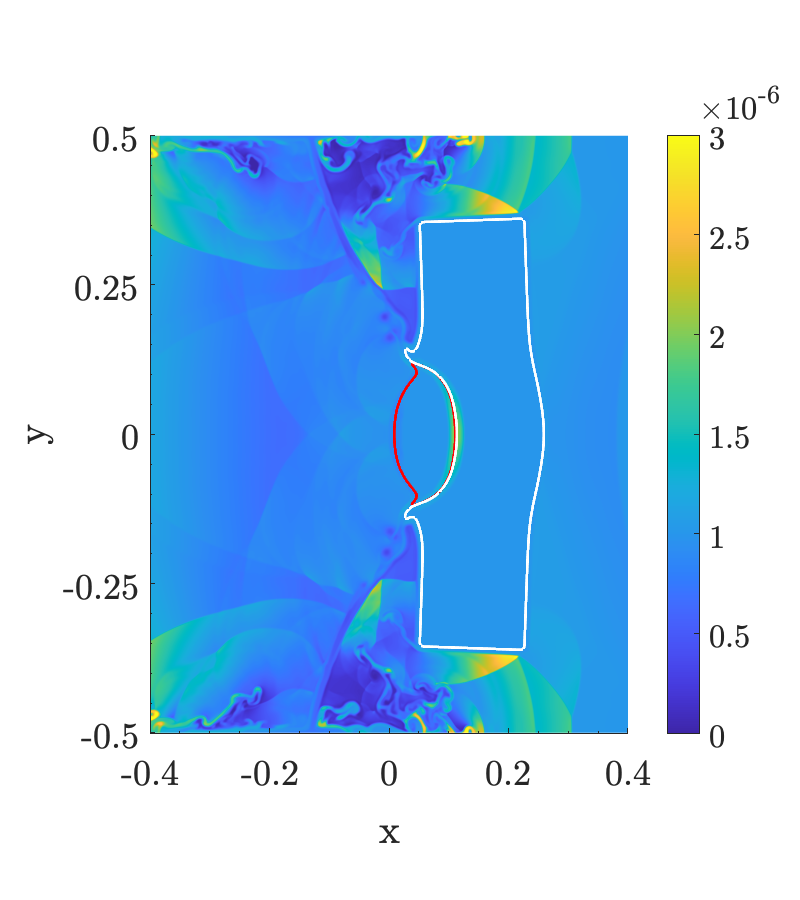}}
	\caption{\footnotesize Results for the multiphase and multi-material elastic-plastic solid impact. The contour levels of the  the gas pressure $p_3$ and the contours of the solids presented at times $t = 60$, $t=100$, and $t = 140$ (from left to right).}
	\label{fig:Solid_impact6}
\end{figure}

Finally, to emphasize the multimaterial character of the test,
Fig.\,\ref{fig:Solid_impact6} shows the gas pressure $p_3$ and the contours of
the solid objects. One can see quite complicated flow structures consisting of
multiple shock waves interracting with the boundaries of the solid bodies and
with each other. 

\section{Conclusion} \label{sec.conclusion} 

The goal of the presented research was to investigate the applicability of the
continuous mixture theory to the numerical simulation of multiphase and
multi-material flows, and in particular to the fluid structure interaction
problems. Within various formulations for the mixture theory, we opted for the
SHTC formulation of continuum mechanics and multiphase flows. The finite volume
MUSCL-Hancock method was selected due to its robustness and easy adaptability to
complex PDE systems like the one presented in this paper.

The multiphase SHTC model \cite{Romenski2007,Romenski2016,Poroelast2020} was
presented in a general non-equilibrium formulation with different phase
pressure, velocities, temperatures, etc., and for arbitrary number of phases.
The phase can be of different nature, viscous gas or liquid, or elastoplastic
solid, and the material interfaces are treated as diffuse interfaces represented
by a rapid but smooth change of the corresponding volume fraction function
across the interface. This presented level of unification of various materials
has been achieved by introducing the unified model of continuum fluid and solid
mechanics \cite{HPR2016,DUMBSER2016824} which also belongs to the SHTC class of
equations.

The SHTC formulation was put into a form resembling the Baer-Nunziato (BN)
model, which is a popular model within the compressible multiphase community.
The motivation for such a reformulation of the SHTC equations was twofold.
First, we wanted to demonstrate how the BN model can be extended beyond its
two-fluid formulation. The second reason is due to our experience with the
numerical discretiztion of the BN model \cite{kemm2020simple} in the past.
However, we find out that in a true multiphysics (e.g. including heat
conduction, elasticity, phase transition, etc.) and multiphase (more than 2
phase) setting, a BN formulation might be not an optimal choice for building a
computational code due to the complexity of the system of equations.
Nevertheless, it was a simplified version of the BN form of the SHTC equations
that was solved numerically in this work.

Despite being a standard finite volume method, the MUSCL-Hancock scheme required some adaptation to the complicated PDE system presented in this paper. To deal with the presence of non-conservative product terms a path-conservative variant of the scheme is employed. Moreover, in order to cope with the stiff character of the relaxation source terms, the MUSCL-Hancock scheme is implemented in an operator splitting manner with the aid of a specially designed implicit discretization of the sources. Particular attention is also given to the reconstruction phase through a specific reconstruction procedure in the primitive variable space that also performs positivity preserving limiting of certain quantities.

This research will be continued in the future in the following directions. We
plan to solve the original SHTC multiphase equations with the help of the new
class of Hyperbolic Thermodynamically Compatible (HTC) schemes
\cite{Busto2021b,ABGRALL2023127629,Thomann2023}, which will allow to respect
most of the structure of the continuous SHTC equations at the discrete level,
and to deal with non-isentropic flows due to the natural ability of the HTC
schemes to deal with hyperbolic models with multiple entropy inequalities, e.g
see \cite{Thomann2023,Busto2021b}.

The multiphase SHTC model presented in Sec.\,\ref{sec.SHTC.model} can be further
extended to include phenomena such as {phase transition} and {surface tension}
\cite{SHTC_surf2023,SHTCSurfaceTension}, likewise to include the SHTC coupling
with electrodynamics \cite{DUMBSER2017298}, e.g. for multifluid plasma models.
Also, a proper SHTC multi-distortion formulation for multiphase flows is
missing. In this paper, the multi-distortion formulation was achieved in a
\textit{ad hoc} manner, and further research is required.

\section{Acknowledgments}

This work was financially supported by the Italian Ministry of Education, University 
and Research (MIUR) in the framework of the PRIN 2022 project \textit{High order structure-preserving semi-implicit schemes for hyperbolic equations}, the PRIN 2022 project No. 2022N9BM3N \textit{Efficient numerical schemes and optimal control methods for time-dependent partial differential equations} and via the  Departments of Excellence  Initiative 2018--2027 attributed to DICAM of the University of Trento (grant L. 232/2016).
\noindent
M.D. was also co-funded by the European Union NextGenerationEU (PNRR, Spoke 7 CN
HPC).  Views and opinions expressed are however those of the author(s) only and
do not necessarily reflect those of the European Union or the European Research
Council. Neither the European Union nor the granting authority can be held
responsible for them. 
\noindent
E.R. was supported by the state contract of the Sobolev Institute of Mathematics
(project no. FWNF-2022-0008)
\noindent
M.D. and I.P. are members of the INdAM GNCS group in Italy.  



\appendix

\section{Eigenvalue estimates}\label{sec.eigenvalues}
From the very nature of the hyperbolic equations, it is important to understand
the \textit{characteristic} structure of the PDE system \eqref{eqn.BN.NumMeth}
under consideration. However, even if simplified with respect to the full SHTC
model in \eqref{eqn.SHTC.system3}, this system remains too complex for
analytical calculation of all eigenvalues due to coupling of convective,
acoustic, and shear parts, as well as due to coupling between the phases.
Therefore, since our FV method requires the knowledge of the maximum sound
speeds, we shall use some estimates for the eigenvalues discussed below. 

As has been suggested in \cite{Chiocchetti2023} in teh context of -two-fluid model coupled with the GPR model, a practical and effective
choice for estimating the spectral radius of the Jacobian matrix for moderate Mach numbers is
$\lambda^{\scriptscriptstyle{M}}_k$ in each direction $x_k$ as
\begin{equation}
\label{eqn.lambda.max}
	\lambda^{\scriptscriptstyle{M}}_k = \max_a\left( \mathopen| \vv_a \cdot \hat{\boldsymbol{n}}_k  + \lambda  \mathclose|,  \mathopen| \vv_a \cdot \hat{\boldsymbol{n}}_k  - \lambda  \mathclose| \right), \qquad \mathrm{with} \qquad \lambda = \sqrt{\lambda_p^2 + \lambda_s^2}
\end{equation}
where 
\begin{equation}
\lambda_p = C_a = p_a \gamma_a/\rho_a
\end{equation}
accounts for the pressure waves of the multiphase model, while
\begin{equation}\label{estimate.shear}
\lambda_s = \sqrt{4 \ \Cs_{a}^2 /3}
\end{equation}
 is a linearised estimate for contribution by the shear modes. Through our
 numerical experiments, we found out that \eqref{estimate.shear} is quite safe,
 leading to only occasional slight overestimates.

\section{Variational formulation}\label{sec.variation}

Here, we provide a variational formulation of the SHTC multifluid heat conducting equations in the
Lagrangian coordinates, which extends the variational scheme from
\cite{SHTC-GENERIC-CMAT} for two-fluid systems to the multifluid case. With the
variational principle we can only obtain the Eulerian homogeneous multifluid
system \eqref{eqn.SHTC.system2} that represents the reversible part of the
overall time evolution. The irreversible part is represented by the relaxation
source terms, which are added afterwards in accordance with the second law of
thermodynamics and Onsager's principle \cite{SHTC-GENERIC-CMAT}. Additionally,
the following formulation is limited to the case of single distortion field
because it remains unclear to us how to extend it to multiple distortion fields.

We consider a vector potential $x_i(t,X_K)$, $i=1,2,3$, $\mathrm{N}-1$-scalar
potentials $\chi_a(t,X_K)$, and another set of $\mathrm{N}$ scalar potentials
$\phi_a(t,X_K)$. Here, $t$ is the time, and $x_i(t,X_K)$ is the map (also called
motion) between the Lagrangian $X_K$, $K=1,2,3$ and Eulerian position $x_i$ of the mixture element. The
physical meaning of the other potentials might be not clear (similar to the
vector potential in the electrodynamics). Yet, their time derivatives are
assumed to have the following meaning
\begin{align}\label{eqn.pot2}
\begin{split}
	V_i = \frac{\rmd x_i}{\rmd t}, \qquad F_{iK} = \frac{\pd x_i}{\pd 
	X_K},\\
	-\mu_a = \frac{\rmd \chi_a}{\rmd t}, \qquad w_{a,K} = \frac{\pd 
	\chi_a}{\pd X_K},\\
	-\theta_a = \frac{\pd \phi_a}{\pd t}, \qquad J_{a,K} = \frac{\pd 
	\phi_a}{\pd X_K},
\end{split}
\end{align}
where $\mu_a$, $a=1,\ldots,\Nph-1$ stands for the phase chemical potential, while $\theta_a$, $a=1,\ldots,\Nph$ is the phase temperature. Also, $\rmd/\rmd t$ should be understood as the material time derivative.

The action integral can be defined in a general form
\begin{equation}\label{eqn.pot3}
	\mathcal{L} = \int \Lambda(x_i,\chi_a,\phi_a;\rmd_t x_i,\pd_K x_i,\rmd_t \chi_a,\pd_K \chi_a,\rmd_t \phi_a,\pd_K \phi_a)  \rmd t \rmd \XX,
\end{equation} 
with the Lagrangian density $\Lambda$ being a function (unspecified at the
moment) of the potentials $x_i$, $\chi_a$, $\phi_a$, and their time and space
gradients denoted with the symbols $\rmd_t = \rmd/\rmd t$ and $\pd_K = \pd/\pd X_K$. However, for our purposes, the explicit dependence of $\Lambda$ on
$x_i$, $\chi_a$, and $\phi_a$ is not required.

One can immediately write down the Euler-Lagrange equations corresponding to this action integral:
\begin{subequations}\label{eqn.EL0}
	\begin{align}
		\frac{\rmd \Lambda_{V_i} }{\rmd t}+ \frac{\pd \Lambda_{F_{iK}} 
			}{\pd X_K} = 0,\label{eqn.pot5}\\
		-\frac{\rmd \Lambda_{\mu_a} }{\rmd t}+ \frac{\pd \Lambda_{w_{a,K} } 
			}{\pd X_K} = 0,\label{eqn.pot6}\\
		-\frac{\rmd \Lambda_{\theta_a} }{\rmd t}+ \frac{\pd \Lambda_{J_{a,K} } 
			}{\pd X_K} = 0, \label{eqn.pot17}
	\end{align}
\end{subequations}
and the following space-time integrability conditions
\begin{subequations}\label{eqn.IntCond0}
	\begin{align}
		&\dfrac{\rmd F_{iJ}}{\rmd t} - \dfrac{\pd {V_i} }{\pd X_K} = 0, 
		\label{eqn.pot7} \\
		&\dfrac{\rmd {w_{a,K} } }{\rmd t} + \dfrac{\pd  \mu_a }{\pd X_K} = 0, 
		\label{eqn.pot8}\\
		&\dfrac{\rmd {J_{a,K} } }{\rmd t} + \dfrac{\pd  \theta_a }{\pd X_K} = 
		0. 
		\label{eqn.pot81}
	\end{align}
\end{subequations}
and pure spatial integrability conditions
\begin{equation}
	\dfrac{\pd F_{iI}}{\pd X_K} - \dfrac{\pd F_{iK}}{\pd X_I} 
	= 0,
	\qquad  
	\dfrac{\pd w_{a,I}}{\pd X_K} - \frac{\pd w_{a,K} }{\pd X_I} = 0,  
	\qquad
	\dfrac{\pd J_{a,I}}{\pd X_K} - \dfrac{\pd J_{a,K} }{\pd X_I} = 0.
\end{equation}
which are trivial consequences of definitions \eqref{eqn.pot2}.

After introducing new unknowns (mixture momentum $U_i$)
\begin{equation}
	U_i = \Lambda_{V_i}, \quad \crho_a =-\Lambda_{\mu_a}, \quad \eta_a = -\Lambda_{\theta_a},
\end{equation}
a new potential
\begin{equation}
	\mathcal{U}(U_i,F_{iJ},\crho_a,w_{a,K},\eta_a,J_{a,K}) = V_i \Lambda_{V_i}
	- \mu_a \Lambda_{\mu_a} - \theta_a \Lambda_{\theta_a}  - \Lambda,
\end{equation}
and noting that the derivatives of the old potential and new potential are related by
\begin{equation}
	\mathcal{U}_{F_{iI}} = -\Lambda_{F_{iI}}, \quad \mathcal{U}_{w_{a,I}} = \Lambda_{w_{a,I}}, \quad \mathcal{U}_{J_{a,I}} = \Lambda_{J_{a,I}},
\end{equation}
\begin{equation}
	U_i = \mathcal{U}_{V_i}, 
	\quad
	\mu_a = \mathcal{U}_{\crho_a},
	\quad
	\theta_a = \mathcal{U}_{\eta_a}, 
\end{equation}
we can rewrite the Euler-Lagrange equations \eqref{eqn.EL0} and the
integrability conditions \eqref{eqn.IntCond0} in the following form
\begin{subequations}\label{eqn.SHTC.Lag1}
	\begin{align}
		&\frac{\rmd U_i 		}{\rmd t} - \frac{ \pd \mathcal{U}_{ 
		F_{iK} } 	}{\pd X_K} = 0, \label{eqn.pot10} \\
		&\frac{\rmd F_{iK}  	}{\rmd t} - \frac{\pd \mathcal{U}_{ 
		U_i} 		} {\pd X_K} = 0, \label{eqn.pot11} \\
		&\frac{\rmd \crho_a 	}{\rmd t} + \frac{\pd \mathcal{U}_{ 
		w_{a,K} }	}{\pd X_K} = 0, \label{eqn.pot13}\\
		&\frac{\rmd w_{a,K} 	}{\rmd t} + \frac{\pd \mathcal{U}_{ 
		\crho_a } 	}{\pd X_K} = 0,  \label{eqn.pot12} \\
		&\frac{\rmd \eta_a 	}{\rmd t} + \frac{\pd \mathcal{U}_{ 
		J_{a,K} }	}{\pd X_K} = 0, \label{eqn.pot132}\\
		&\frac{\rmd J_{a,K} 	}{\rmd t} + \frac{\pd \mathcal{U}_{ 
		\eta_a } 	}{\pd X_K} = 0, \label{eqn.pot121}
	\end{align}
\end{subequations}

The conversion from the SHTC Lagrangian master system \eqref{eqn.SHTC.Lag1},
which admits the original Godunov structure \cite{God1961}, to the Eulerian
system \eqref{eqn.SHTC.system2} can be carried out by means of the
Lagrange-to-Euler coordinate transformation which also concerns the change of
the state variables, potential $\mathcal{U}$, and is described in detail in
\cite{SHTC-GENERIC-CMAT,GodRom1996b}.

\printbibliography

\end{document}